\documentclass[twoside]{article}
\usepackage[latin1]{inputenc}
\usepackage[T1]{fontenc}
\usepackage{parskip}
\usepackage{amsmath,amsfonts,amssymb,amsxtra}
\usepackage{latexsym}
\usepackage{theorem}
\usepackage{graphicx,psfrag}
\usepackage{verbatim}
\usepackage{hyperref}

\newtheorem{theo}{Theorem}
\newtheorem{prop}{Proposition}
\newtheorem{lemma}{Lemma}
\newtheorem{cor}{Corollary}

{\theoremstyle{plain}                       % Remark
\theorembodyfont{\rmfamily}
\newtheorem{Remark}{Remark}}

{\theoremstyle{plain}                       % Example
\theorembodyfont{\rmfamily}
\newtheorem{Example}{Example}}

{\theoremstyle{plain}                       % Definition
\theorembodyfont{\rmfamily}
\newtheorem{Definition}{Definition}}

{\theoremstyle{plain}                       % Condition
\theorembodyfont{\rmfamily}
}

\def\A{\mathcal{A}}
\def\B{\mathcal{B}}
\def\C{\mathcal{C}}
\def\D{\mathcal{D}}
\def\E{\mathcal{E}}
\def\F{\mathcal{F}}
\def\H{\mathcal{H}}
\def\J{\mathcal{J}}
\def\I{\mathcal{I}}
\def\K{\mathcal{K}}
\def\L{\mathcal{L}}
\def\P{\mathcal{P}}

\def\N{\mathbb{N}}
\def\Z{\mathbb{Z}}
\def\R{\mathbb{R}}
\def\U{\mathcal{U}}
\def\e{\epsilon}
\def\d{\delta}
\def\a{\alpha}
\def\b{\beta}
\def\g{\gamma}
\def\G{\Gamma}
\def\t{\tau}
\def\T{\Theta}
\def\o{\omega}
\def\O{\Omega}

\begin{document}

\title{Lower bounds for the dynamically defined measures}

\author{Ivan Werner
\footnote{ Email: ivan\_werner@mail.ru}} 
\date{February 2, 2022}
\maketitle

\begin{abstract}\noindent
 The dynamically defined measure (DDM) $\Phi$ arising from a finite measure $\phi_0$ on an initial $\sigma$-algebra on a set and an invertible map acting on the latter is considered. Several lower bounds for it are obtained and sufficient conditions for its positivity are deduced under the general assumption that there exists an invariant measure $\Lambda$ such that $\Lambda\ll\phi_0$.  
 
In particular, DDMs arising from the Hellinger integral $\mathcal{J}_\alpha(\Lambda,\phi_0)\geq\mathcal{H}^{\alpha,0}(\Lambda,\phi_0)\geq\mathcal{H}_\alpha(\Lambda,\phi_0)$ are constructed with $\mathcal{H}_{0}\left(\Lambda,\phi_0\right)(Q) = \Phi(Q)$, $\mathcal{H}_{1}\left(\Lambda,\phi_0\right)(Q) = \Lambda(Q)$, and
\[\Phi(Q)^{1-\alpha}\Lambda(Q)^{\alpha}\geq\mathcal{J}_{\alpha}\left(\Lambda,\phi_0\right)(Q)\]
for all measurable $Q$ and $\alpha\in[0,1]$, and further computable lower bounds for them are obtained and analyzed. The function $(0,\gamma]\owns\alpha\longmapsto\mathcal{H}_{\alpha}(\Lambda,\phi_0)$ is computed explicitly for $\gamma\geq 1$ such that $\int(d\Lambda/d\phi_0)^{\gamma-1}d\Lambda<\infty$ in the case of a discrete ergodic decomposition of $\Lambda$, and the other two functions are computed under the additional condition of the equivalence of $\phi_0$ and $\Lambda$. In particular, if $\Lambda$ is ergodic, it is shown that the first function is completely determined by the $\Lambda$-essential supremum (infimum) of $d\Lambda/d\phi_0$ for all $0<\a<1$ ($1<\a\leq\g$), and, if it is continuous at $0$, the above inequalities become equalities. The computation of it enables an explicit computation of some DDMs arising as outer measure approximations with respect to it, which demonstrates that this technique allows to obtain new measures, and that such measures can have phase transitions with respect to the DDM specifying the covering sets. In general, regularity properties of all above functions are studied. In particular, all one-sided derivatives of the two greater functions are obtained, and some lower bounds for them by means of the derivatives are given. Some sufficient conditions for the continuity and a one-sided differentiability of the smallest one are provided.


\end{abstract}

\tableofcontents

\section{Introduction}

This article is concerned with the development of general methods for computation of lower bounds for the dynamically defined measures (DDMs) \cite{Wer3},\cite{Wer10},\cite{Wer12},\cite{Wer15} and thus obtaining conditions for their positivity. The latter became particularly required after the recently discovered error in \cite{Wer3}, see \cite{Wer13}.

A DDM is obtained via an extension of the Carath\'{e}odory construction which provides a method for obtaining a measure from a sequence of refining, but not necessarily consistent contents. The inconsistency requires a separate argument on why the obtained measure is not zero. This, in turn, drives a further extension of the measure-theoretic foundation, which leads to new mathematical tools and phenomena opening new possibilities for measure-theoretic modeling and analysis.

Originally, the dynamically defined outer measure $\Phi$ arising from a finite measure $\phi_0$ on an initial $\sigma$-algebra was proposed in \cite{Wer3} as a way to construct the coding map for a contractive Markov system (CMS) \cite{Wer1} almost everywhere with respect to an outer measure which is also obtained constructively (at least on compact sets; in general, it still requires the axiom of choice, but the obtained measure is unique). This outer measure arose in a natural way from the condition of the contraction on average. 

Later, the author  also could not avoid the routine to define the coding map almost everywhere with respect to a measure which is obtained in the canonical, non-constructive and less descriptive way (via the Krylov-Bogolyubov argument) \cite{Wer11}. However, before the dynamically defined outer measure became redundant, it was shown in \cite{Wer10} and \cite{Wer12} that the restriction of the outer measure on the Borel $\sigma$-algebra is a measure the normalization of which provides a construction for equilibrium states for CMSs (the local energy function of which is given by means of the coding map \cite{Wer6}\cite{Wer11}, which makes it highly irregular, so that no other known method, to the author's knowledge, is capable to provide a construction). 

The normalization is, of course, possible only if the measure is not zero. The discovered error in \cite{Wer3} puts it into serious doubts in  a general case. At the time of writing, it has only been shown in \cite{Wer13} that the measure is not zero if all the maps of the CMS are contractions (which does not go far beyond the case accessible by means of a Gibbs measure), with a little comfort that no openness of the Markov partition is required (which makes the local energy function still only measurable in  general). 

The method which is used in  \cite{Wer13} is based on the proof that the logarithm of the supremum of the density function of an invariant measure with respect to the initial one along the trajectories is integrable, which seemed to be a very strong condition in general. 

Trying to weaken that led to the introduction of a {\it relative entropy measure} in this article (Subsection \ref{rems}). The proof that it is a measure is based just on a few of its properties, which are weaker than that of an outer measure. It requires a notion of an {\it outer measure approximation} and a generalization of the Carath\'{e}odory theorem for it. The extension of the {\it Measure Theory} on such constructions in a general setting, based on sequences of {\it measurement pairs}, which can be called the {\it Dynamical Measure Theory}, was developed in \cite{Wer15}. It enables us to compute and analyze all lower bounds for the DDMs in this paper.

All lower bounds for the DDMs here are obtained in the case when the measurement pairs are generated by an invertible map from an initial $\sigma$-algebra and a measure on it. Moreover, for the computations of the lower bounds, we will always assume that there exists an invariant measure $\Lambda$ which is absolutely continuous with respect to the initial one, $\phi_0$. 

In this paper, we develop the dynamical measure theory slightly further by showing that the set functions arising here as outer measure approximations satisfy the regularity property (Lemma \ref{omr}), which, together with the fact that their restrictions on measurable set are measures, implies that they must be outer measures (Proposition \ref{omp}).
 
It became clear after the development of the dynamical measure theory in \cite{Wer15} that it is logical, from the point of view of the structure of the theory in this article, and advantageous for the purpose of obtaining the best practical lower bounds, to introduce first an intermediate family of {\it DDMs arising from the  Hellinger integral} $\H_{\a}\left(\Lambda,\phi_0\right)$, $\a\geq 0$, with $\H_{0}\left(\Lambda,\phi_0\right)(Q) = \Phi(Q)$, and $\mathcal{H}_{1}\left(\Lambda,\phi_0\right)(Q) = \Lambda(Q)$, which provide lower bounds for $\Phi$ through
\[\Phi(Q)^{1-\a}\Lambda(Q)^{\a}\geq\H_{\a}\left(\Lambda,\phi_0\right)(Q)\]
for all measurable $Q$ and $\a\in[0,1]$, and then to obtain a lower bound for $\H_{\a}\left(\Lambda,\phi_0\right)$ via the relative entropy measure (Theorem \ref{lbre}), the local finiteness of which guaranties the positivity of $\H_{\a}\left(\Lambda,\phi_0\right)$.  

Furthermore, this approach allowed us to obtain a practical sufficient condition for the positivity of $[0,1]\owns\a\longmapsto\H_{\a}\left(\Lambda,\phi_0\right)(Q)$ via the limit $\a\to 1$ (Corollary \ref{lpc} (ii)).

In Subsection \ref{subre}, we also provide some natural upper bounds on the relative entropy measure. In particular, in the case of an ergodic $\Lambda$, we show that the finiteness of the relative entropy measure is equivalent to the essential boundedness of $d\Lambda/d\phi_0$ with respect to $\Lambda$ and to the absolute continuity of $\Lambda$ with respect to $\H_{\a}\left(\Lambda,\phi_0\right)$ for all $\a\in[0,1)$ (Corollary \ref{fcem}).

Another advantage of this approach is the possibility for obtaining criteria for the positivity of $\Phi$ via the dependence of $\H_{\a}\left(\Lambda,\phi_0\right)$ on $\a$. This led to the study of other DDMs, in particular, another DDM arising from the  Hellinger integral as an outer measure approximation
\begin{equation}\label{mlbi}
  \Phi(Q)^{1-\a}\Lambda(Q)^{\a}\geq\H^{\a,0}\left(\Lambda,\phi_0\right)(Q)\geq\H_{\a}\left(\Lambda,\phi_0\right)(Q)
\end{equation}
for all measurable $Q$ and $\a\in[0,1]$. 

Clearly, establishing that the functions $[0,1]\owns\a\longmapsto\H_{\a}\left(\Lambda,\phi_0\right)(Q)$ and $[0,1]\owns\a\longmapsto\H^{\a,0}\left(\Lambda,\phi_0\right)(Q)$ have different properties on $[0,1)$ would immediately imply the positivity of $\Phi$. In the case of the first function, we were only able to show that it is positive all the way to the left if is positive at some point in $(0,1)$ and it is zero all the way to the right in the open interval if it is zero at such a point (Lemma \ref{shm} (iii)), but the second is always either zero everywhere on $(0,1)$ or strictly positive on $[0,1]$ (Lemma \ref{pchm} (iv)), due to a certain property of a logarithmic almost convexity of the function. We were not able to establish the continuity of the first function on $(0,1)$ in general, but it holds true for the second (Lemma \ref{pchm} (vii)). 

The continuity of the first function on $(0,1)$ in general could be obtained only under a condition (Proposition \ref{hcp}), which is, in particular, satisfied if $d\Lambda/d\phi_0$ is $\Lambda$-essentially bounded away from zero. In this case, it is also right and left differentiable (with the left derivative not smaller than the right) (Theorem \ref{hrdt} and Theorem \ref{hldt}), which implies, by the well-known result going back to Beppo Levi, that it is differentiable everywhere except at most countably many points (Corollary \ref{hbdc}). 

If $\Lambda$ has a discrete ergodic decomposition, we obtained the function explicitly on $(0,1)$ (Theorem \ref{aedt} (ii)). In particular, it gives an explicit computation of $\Phi$ if the function is continuous at $0$ (Corollary \ref{hecc} (i)), which is, for example, satisfied if $\phi_0$ and $\Lambda$ are equivalent (Corollary \ref{hecc} (ii)). In this case, the function is smooth (Corollary \ref{hcdc}). In particular, it is completely determined by the $\Lambda$-essential supremum of $d\Lambda/d\phi_0$ if $\Lambda$ is ergodic. If, in addition to the ergodicity, $\phi_0$ and $\Lambda$ are equivalent, then the inequalities in \eqref{mlbi} become equalities (Corollary \ref{pebc}). (In the case of a discrete ergodic decomposition of $\Lambda$ and the finiteness of the integral $\int(d\Lambda/d\phi_0)^\g d\phi_0$ for some $\g>1$, we computed the function explicitly also for all $1<\a\leq\g$. For this parameter values, it is completely determined by the $\Lambda$-essential infimum of $d\Lambda/d\phi_0$ in the ergodic case (Theorem \ref{ec2}).)

Also, we obtained a sufficient condition for the continuity of the functions at $1$ (Proposition \ref{ca1}) (which is slightly stronger than the weakest obtained sufficient condition for the positivity of $[0,1]\owns\a\longmapsto\H_{\a}\left(\Lambda,\phi_0\right)(Q)$).

Due to the Lipschitz continuity of the function $(0,1)\owns\a\longmapsto\H^{\a,0}\left(\Lambda,\phi_0\right)(Q)$ on every closed subinterval, it is already differentiable almost everywhere. This all encourages us to investigate other possibly finer regularity (or irregularity) properties of it. To that end, we obtained some (singed) measures which naturally suggest themselves as candidates for the derivatives of it. We showed that the first one is in fact the right derivative (Theorem \ref{dlpl}), but the left one still turned out to be something else (Theorem \ref{ldt}), but also not smaller than the right. However, again as a consequence of the Beppo Levi Theorem, there exists an at most countable set such that the function is differentiable everywhere except at the points in it, and the restriction of the function on the complement of it is continuously differentiable (Corollary \ref{dc}). Moreover, we showed that the logarithmic almost convexity of the function implies that it is strictly smaller than the weighted geometric average in the inequality \eqref{mlbi} at the points where its left derivative is greater than the right (Proposition \ref{ftp}).

The latter inspired us to introduce another DDM arising from the Hellinger integral as an outer measure approximation, with an additional condition on the covering sets, which satisfies the inequalities
\begin{equation*}
	\Phi(Q)^{1-\a}\Lambda(Q)^{\a}\geq\J_{\a}\left(\Lambda,\phi_0\right)\geq\H^{\a,0}\left(\Lambda,\phi_0\right)(Q)
\end{equation*}
for all measurable $Q$ and $\a\in[0,1]$ (Section \ref{tolb}). We showed that $[0,1]\owns\a\longmapsto\J_\a\left(\Lambda,\phi_0\right)(Q)$ is logarithmically convex (which also leads to a general definition of {\it the logarithmic almost convexity} for $[0,1]\owns\a\longmapsto\H^{\a,0}\left(\Lambda,\phi_0\right)(Q)$), but its one-sided derivatives seem to be also different in general (Definition \ref{Jld} and Definition \ref{Jrd}). In any case, the positive derivatives can be used to obtain lower bounds for the functions (Corollary \ref{HlbX} and Corollary \ref{rdlb}).

In the case of a discrete ergodic decomposition of $\Lambda$ and the equivalence of $\phi_0$ and $\Lambda$, all three functions coincide on $[0,1]$, but only the last two coincide on $(1,\g]$ if they are finite there, and the first one is strictly smaller if $\phi_0$ and  $\Lambda$ are different (Theorem \ref{homa} and Theorem \ref{homa2}).

To sum up, so far, it seems that the main result of the approach taken in this paper is the explicit computation of $\Phi$ in terms of the Radon-Nikodym derivative $d\Lambda/d\phi_0$ via an expression of $[0,1]\owns\a\longmapsto\H_{\a}\left(\Lambda,\phi_0\right)(Q)$ as a sum over ergodic components in the case of a discrete ergodic decomposition of $\Lambda$ and the continuity of the function at $0$ (Corollary \ref{pebc}). As a by-product of it and the computation of $(1, \g]\owns\a\longmapsto\H_{\a}\left(\Lambda,\phi_0\right)(Q)$ (Theorem \ref{ec2} (ii)), we obtain explicit computations of some DDMs arising from the Hellinger integral as outer measure approximations  (Corollary \ref{dhsp}, Corollary \ref{hsdc}, Subsection \ref{homs}, and Subsection \ref{homs2}), which, in particular, demonstrate that this technique leads to new measures and start to paint a picture on how an inductively constructed DDM depends on the previous one. In particular, we obtain an example where this dependence exhibits a phase transition of the measure (Theorem \ref{homa}). Furthermore, as another by-product, we give examples of DDMs obtained from sequences of non-additive contents (Subsection \ref{inc}), which suggests that the dynamical measure theory discovered so far is most likely just the tip of the iceberg.    

As indicated by the names of the introduced auxiliary measures, we will need some preliminaries from the information theory, which are collected in Subsection \ref{itps}.

Concluding the introduction, a few words on the notation. All considerations in this article will take place on a set $X$. We will denote  the collection of all subsets of $X$ by $\P(X)$. As usual, $\N$ and  $\Z$ will denote the set of all natural numbers (without zero) and  the set of all integers respectively. We will use the notation '$f|_{A}$' to denote the restriction of a function $f$ on a set $A$, '$1_{A}$' to denote the indicator function of a set $A$, '$\ll$' to denote the absolute continuity relation for set functions, '$A\Delta B$' to denote the symmetric difference between sets, '$f\vee g$'  ('$f\wedge g$') to denote the maximum (minimum) of $f$ and $g$ with the abbreviations $f^+:=f\vee 0$ and $f^-:=-(f\wedge 0)$, and '$x\to^+y$' ('$x\to^-y$') to abbreviate the convergence $x\to y$ and $x>y$ ($x<y$).

\section{The setup for the dynamically defined measure $\Phi$}

In this section, we define the main object of the study in this article - a particular case of the dynamically defined measure as specified in Section 5 in \cite{Wer15}. 

Let $X$ be a set and $S:X\longrightarrow X$ be an invertible map. Let $\A$ be a $\sigma$-algebra on $X$. Let $\A_0$ be the $\sigma$-algebra generated by $\bigcup_{i=0}^\infty S^{-i}\A$ and  $\B$ be the $\sigma$-algebra generated by $\bigcup_{i=-\infty}^\infty S^{-i}\A$. Define
\[\A_m:=S^{-m}\A_0\ \ \ \mbox{ for all }m\in\Z\setminus\N.\]
It is not difficult to verify that $\A_0\subset\A_{-1}\subset...$, $\B$ is generated by $\bigcup_{m\leq 0}\A_m$, and $S$ is $\B$-$\B$ and $\A_0$-$\A_0$-measurable (see Section 5 in \cite{Wer15}).

Let $\phi_0$ be a finite, positive measure on $\A_0$. For $Q\subset X$, define
\[\C(Q):=\left\{(A_m)_{m\leq 0}|\ A_m\in\mathcal{A}_{m}\
\forall m\leq0\mbox{ and }Q\subset\bigcup\limits_{m\leq 0}A_m \right\}\]
and
\[\Phi(Q):=\inf\limits_{(A_m)_{m\leq 0}\in\C(Q)}\sum\limits_{m\leq 0}\phi_{0}(S^mA_m).\]
Then  $\Phi(S^{i}Q)\leq\Phi(S^{i-1}Q)$ for all $i\leq 0$ (see Sections 4 and 5 in \cite{Wer15}). Define
\[\bar\Phi(Q):=\lim\limits_{i\to-\infty}\Phi\left(S^{i}Q\right).\]
Then, by Theorem 16 (i) (Theorem 4 (i) in the arXiv version) in \cite{Wer15}, $\bar\Phi(Q)=\Phi(Q)$ for all $Q\in\B$, and $\Phi$ is a (obviously $S$-invariant) measure on $\B$, which we call the {\it dynamically defined measure (DDM)} associated with $\phi_0$.

\begin{Example}\label{e1}
	Let $P:=(p_{ij})_{1\leq i,j\leq N}$ be a stochastic $N\times N$-matrix. Let $X:=\{1,...,N\}^{\Z}$ (be the set of all $(...,\sigma_{-1},\sigma_0,\sigma_{1},...)$, $\sigma_i\in\{1,...,N\})$ and $S$ be the left shift map on $X$ (i.e. $(S\sigma)_i=\sigma_{i+1}$ for all $i\in\Z$). Let ${_0[a]}$ denote a cylinder set at time $0$ (i.e. the set of all $(\sigma_i)_{i\in\Z}\in X$ such that $\sigma_0=a$ where $a\in\{1,...,N\}$). Let $\A$ be the $\sigma$-algebra generated by the partition $({_0[a]})_{a\in\{1,...,N\}}$.
	
	 Let $\nu$ be a probability measure on all subsets of $\{1,...,N\}$.  Let $\phi_0$  be the probability measures on $\A_0$ given by
	\[\phi_0\left({_0[i_1,...,i_n]}\right):=\nu\{i_1\}p_{i_1i_2}...p_{i_{n-1}i_n}\]
	for all $_0[i_1,...,i_n]\subset\{1,...,N\}^{\Z}$ and $n\geq0$. One easily sees that $\Phi(X)>0$ if $P$ is irreducible and $\nu\{i\}>0$ for all $i\in\{1,...,N\}$ (see Example 2 in \cite{Wer15}).
\end{Example}

For an example in which the positivity of $\Phi$ is not that obvious, see \cite{Wer13}.

In this note, we will use the measure theory developed in \cite{Wer15} to obtain lower bounds for $\Phi$ in terms of various (signed) measures in the case when there exists  $\phi'_0\ll\phi_0$  such that $\phi'_0\circ S^{-1}= \phi'_0$, which will allow us not only to obtain sufficient conditions for the positivity of $\Phi$ (which is another important role which is going to be salvaged from the erroneous Lemma 2 (ii) in \cite{Wer3}), but also it will give several  necessary  and sufficient conditions for $\Phi'|_\B\ll\Phi|_\B$ in the case when $\phi'_0$ is ergodic. By Proposition 11 (Proposition 1 in the arXiv version) in \cite{Wer15}, $\Phi'|_{\A_m}=\phi'_0\circ S^m$ for all $m\leq 0$.

In the following, we will denote by $\Lambda$ a positive and finite measure on $\A_0$ such that $\Lambda\circ S^{-1} = \Lambda$ and $\Lambda\ll\phi_0$.  
Its unique extension on $\B$, which is, for example, given by Proposition 11 in \cite{Wer15}, and the dynamically defined outer measure (in this case, the usual Lebesgue outer measure) will be also denoted by $\Lambda$, since it is always clear what is meant from the set to which it is applied.  

Let $Z$ be a measurable version of the Radon-Nikodym derivative $d\Lambda/d\phi_0$.

\section{Preliminaries}

As a usual convenience in a measure-theoretic analysis, we will use extended real numbers $\{-\infty\}\cup\R\cup\{\infty\}$ with the continuous extensions of common functions on them. However, contrary to the usual way of defining objects in the measure theory, we will consider the product $0\times\infty$ as NOT defined. Instead, a Lebesgue integral $\int_X fg d\mu$ of the product of measurable functions $f$ and $g$ with values in $[-\infty, +\infty]$ will be understood as $\int_X fg d\mu:=\int_{X\setminus\mathcal{N}} (fg)^+ d\mu-\int_{X\setminus\mathcal{N}} (fg)^- d\mu$ with a measurable $\mu$-zero set $\mathcal{N}$ such that the product $f(x)g(x)$ is well defined for all $x\in X\setminus\mathcal{N}$ and $\int_{X\setminus\mathcal{N}} (fg)^+ d\mu<\infty$ or $\int_{X\setminus\mathcal{N}} (fg)^- d\mu<\infty$. 

Also, we will use the continuous extensions $0\log 0:=0$ and $0\log(0/0):=0\log 0-0\log 0=0$. (As a consequence, $0^0 = 1$, since $y^x:=e^{x\log y}$.)

\subsection{Preliminaries for the derivatives of an exponential function}\label{pexp}

In this article, we are going to study, in particular, some functions obtained as some infimums and supremums of the functions $[0,\g)\owns\a\longmapsto\sum_{m\leq 0}\int_{S^mA_m}Z^\a d\phi_0$. In this context, since  $dZ^\a/d\a = Z^\a\log Z$, we will need the following simple lemmas.

\begin{lemma}\label{hfpl} 
	For every $n\in\N$ and $0\leq\a<1$,
	\[\max\limits_{x\in[0,1]}x|\log x|^n = \left(\frac{n}{e}\right)^n\ \left(\mbox{it is attained at }e^{-n}\right),\]
	\[\max\limits_{x\in[0,\infty)}e^{-(1-\a)x}x^n = \left(\frac{n} {e(1-\a)}\right)^n\ \left(\mbox{it is attained at }\frac{n}{1-\a}\right).\]
\end{lemma}
{\it Proof.}
The proof is straightforward.
\hfill$\Box$

\begin{lemma}\label{edl} 
	Let $0<\a_0<\a$, $n\in\N\cup\{0\}$, $Z\geq 0$ and
	\[D_n^{\a,\a_0}(Z):=\frac{Z^{\a}(\log Z)^{n} - Z^{\a_0}(\log Z)^{n}}{\a-\a_0}\mbox{ with }(\log Z)^0:=1.\]
	(i) If $n$ is even, then
	\[Z^{\a_0}(\log Z)^{n+1}\leq D_n^{\a,\a_0}(Z)\leq Z^{\a}(\log Z)^{n+1}.\]
	(ii) If $n$ is odd, then
	\[0\leq D_n^{\a,\a_0}(Z)\leq 1_{\{Z\leq C\}}\frac{1}{1-(\a-\a_0)\log C}Z^{\a_0}(\log Z)^{n+1}+1_{\{Z > C\}}Z^{\a}(\log Z)^{n+1}\]
	for all $C\geq 1$ such that $(\a-\a_0)\log C<1$, and, for $0<\a_0<\a<1$,
	\begin{eqnarray*}
		&& \max\left\{Z^{\a_0}(\log Z)^{n+1}-(\a-\a_0)\left(\frac{n+2}{\a_0e}\right)^{n+2}1_{\{Z\leq 1\}}\right.,\\
		&&\left.Z^{\a}(\log Z)^{n+1}-(\a-\a_0)\left(\frac{n+2}{(1-\a)e}\right)^{n+2}Z1_{\{Z> 1\}}\right\}\\
		&\leq& D_n^{\a,\a_0}(Z)\leq \min\left\{Z^{\a}(\log Z)^{n+1}+(\a-\a_0)\left(\frac{n+2}{\a_0 e}\right)^{n+2}1_{\{Z\leq 1\}},\right.\\
		&&\left.Z^{\a_0}(\log Z)^{n+1}+(\a-\a_0)\left(\frac{n+2}{(1-\a)e}\right)^{n+2}Z1_{\{Z>1\}} \right\}.\ 
	\end{eqnarray*}
\end{lemma}
{\it Proof.} Obviously, (i) and (ii) are correct if $Z=0$. Suppose $Z>0$. 

$(i)$ Observe that
\[Z^{\a_0}(\log Z)^{n+1}=\frac{1}{\a-\a_0}Z^{\a_0}(\log Z)^{n}\log Z^{\a-\a_0}\leq\frac{1}{\a-\a_0}Z^{\a_0}(\log Z)^{n}\left(Z^{\a-\a_0}-1\right).\]
This implies the first inequality in (i). Also,
\[Z^{\a}(\log Z)^{n+1}=-\frac{1}{\a-\a_0}Z^{\a}(\log Z)^{n}\log Z^{\a_0-\a}\geq-\frac{1}{\a-\a_0}Z^{\a}(\log Z)^{n}\left(Z^{\a_0-\a}-1\right).\]
This implies the second inequality in (i).

$(ii)$ The inequality $0\leq D_n^{\a,\a_0}(Z)$ is obvious. Furthermore, observe that for $0\leq Z\leq 1$,
\begin{eqnarray*}
	Z^\a(\log Z)^{n+1} &=& -\frac{1}{\a-\a_0}Z^\a(\log Z)^{n}\log Z^{-\a+\a_0}\\
	&\leq&-\frac{1}{\a-\a_0}Z^\a(\log Z)^{n}\left(Z^{-\a+\a_0}-1\right)\\
	&=& D_n^{\a,\a_0}(Z).
\end{eqnarray*}
For $Z\geq 1$, as in (i),
\[Z^{\a_0}(\log Z)^{n+1}\leq D_n^{\a,\a_0}(Z).\]
Hence, for every $Z\geq 0$,
\begin{equation}\label{oda}
D_n^{\a,\a_0}(Z)\geq 1_{\{Z\leq 1\}}Z^\a(\log Z)^{n+1}+1_{\{Z> 1\}}Z^{\a_0}(\log Z)^{n+1}.
\end{equation}
Then on one hand, by (i) and Lemma \ref{hfpl}, for $\a_0>0$,
\begin{eqnarray}\label{foldi}
D_n^{\a,\a_0}(Z)&\geq& Z^{\a_0}(\log Z)^{n+1} + 1_{ \{Z\leq 1\}}\left(Z^\a-Z^{\a_0}\right)(\log Z)^{n+1}\nonumber\\
&\geq& Z^{\a_0}(\log Z)^{n+1} + 1_{\{ Z\leq 1\}}Z^{\a_0}(\log Z)^{n+2}(\a-\a_0)\nonumber\\
&\geq& Z^{\a_0}(\log Z)^{n+1} - 1_{\{Z\leq 1\}}\left(\frac{n+2}{\a_0e}\right)^{n+2}(\a-\a_0),
\end{eqnarray}
and on the other hand, by (i) and Lemma \ref{hfpl}, for $\a<1$,
\begin{eqnarray}\label{soldi}
D_n^{\a,\a_0}(Z)&\geq& Z^{\a}(\log Z)^{n+1} - 1_{\{ Z> 1\}}\left(Z^\a-Z^{\a_0}\right)(\log Z)^{n+1}\nonumber\\
&\geq& Z^{\a}(\log Z)^{n+1} - 1_{\{ Z> 1\}}Z^{\a}(\log Z)^{n+2}(\a-\a_0)\nonumber\\
&=& Z^{\a}(\log Z)^{n+1} - 1_{\{ Z> 1\}}Ze^{-(1-\a)\log Z}(\log Z)^{n+2}(\a-\a_0)\nonumber\\
&\geq& Z^{\a}(\log Z)^{n+1} - 1_{\{ Z> 1\}}Z\left(\frac{n+2}{(1-\a)e}\right)^{n+2}(\a-\a_0).
\end{eqnarray}
Thus, \eqref{foldi} and \eqref{soldi} imply the first inequality of the second part in (ii).

Let $C\geq 1$ such that $(\a-\a_0)\log C<1$. If $Z\leq C$, then, by (i),
\begin{eqnarray*}
	Z^{\a_0}(\log Z)^{n+1}&=&\frac{1}{\a-\a_0}Z^{\a_0}(\log Z)^{n-1}\left(\log\frac{Z}{C}\right)\log Z^{\a-\a_0}\\
	&&+Z^{\a_0}\log C(\log Z)^{n}\\
	&\geq&\frac{1}{\a-\a_0}Z^{\a_0}(\log Z)^{n-1}\left(\log\frac{Z}{C}\right)\left(Z^{\a-\a_0}-1\right)\\
	&&+Z^{\a_0}\log C(\log Z)^{n}\\
	&=& D_{n}^{\a,\a_0}(Z)+\log C(\log Z)^{n-1}\left(Z^{\a_0}\log Z-\frac{Z^\a-Z^{\a_0}}{\a-\a_0}\right)\\
	&\geq& D_{n}^{\a,\a_0}(Z)+\log C(\log Z)^{n-1}\left(Z^{\a_0}\log Z-Z^{\a}\log Z\right)\\
	&=& D_{n}^{\a,\a_0}(Z)\left(1-(\a-\a_0)\log C \right).
\end{eqnarray*}
If $Z\geq C$, then, as in (i),
\[Z^\a(\log Z)^{n+1}\geq D_n^{\a,\a_0}(Z).\]
Hence, it follows the second inequality of the first part in (ii). 

Then, as above, by (i) and Lemma \ref{hfpl}, on one hand, for $\a<1$,
\begin{equation}\label{fordi}
D_n^{\a,\a_0}(Z)\leq  Z^{\a_0}(\log Z)^{n+1} +(\a-\a_0)1_{\{Z> 1\}}Z\left(\frac{n+2}{(1-\a)e}\right)^{n+2},
\end{equation}
and on the other hand, for $\a_0>0$,
\begin{equation}\label{sordi}
D_n^{\a,\a_0}(Z)\leq Z^{\a}(\log Z)^{n+1} +(\a-\a_0)1_{\{Z\leq 1\}}\left(\frac{n+2}{\a_0e}\right)^{n+2}.
\end{equation}
Thus, \eqref{fordi} and \eqref{sordi} imply the second inequality in (ii).
\hfill$\Box$

\subsection{Information-theoretic preliminaries}\label{itps}

In this article, we will also make use of some generalizations and  derivations of some relations between measures which were developed in the information theory. We collect the required preliminary material in this subsection. 

Let $(X, \A, \Lambda)$ be a finite measure space, i.e. $\A$ is a $\sigma$-algebra, and $\Lambda$ is a positive and finite measure on it.

Let $\phi$ be another positive and finite measure on $\A$ such that $\Lambda\ll \phi$. Let $f$ be a measurable version of the Radon-Nikodym derivative $d\Lambda/d\phi$. (Note that $\Lambda\{f=0\} = 0$.)

\begin{Definition}
	Let $A\in\A$. Define
	\[K\left(\Lambda|\phi\right)(A):=\int\limits_{A}  \log  fd\Lambda,\ \ \ \mbox{ and }\ \ \ K\left(\Lambda|\phi\right):=K\left(\Lambda|\phi\right)(X).\]
	The latter is called the {\it Kullback-Leibler divergence} of $\Lambda$ with respect to $\phi$. For $\alpha\geq 0$, define
	\[H_\alpha(\Lambda,\phi)(A):=\int\limits_{A} f^\alpha d\phi,\ \ \ \mbox{ and }\ \ \ H_\alpha(\Lambda,\phi):=H_\alpha(\Lambda,\phi)(X).\]  
	The latter is called the {\it Hellinger integral}.
\end{Definition}

For some properties of them, e.g. see \cite{B}. Since $x\log x\geq x-1$ for all $x\geq 0$, $K\left(\Lambda|\phi\right)(A)\geq\Lambda(A)-\phi(A)$. In particular, $K\left(\Lambda|\phi\right)(A)\geq 0$ if $\Lambda(A)\geq\phi(A)$.
Obviously, by the concavity of $x\longmapsto x^\alpha$, $0\leq H_\alpha(\Lambda,\phi)(A)\leq\phi(A)^{1-\alpha}\Lambda(A)^\alpha$ for all $0\leq \alpha\leq 1$. 

In this article, we are going, in particular, to extend the following relation of the measures to that of the corresponding DDMs which allows to obtain a lower bound for the DDM of the main concern. 
\begin{lemma}\label{HKr} 
	Let $A\in\A$ such that $\Lambda(A)>0$. Then
	\[K\left(\Lambda|\phi\right)(A)\geq-\frac{\Lambda(A)}{\alpha}\log\frac{H_{1-\alpha}(\Lambda,\phi)(A)}{\Lambda(A)}\ \ \ \mbox{ for all }0<\alpha\leq 1,\mbox{ and }\]
	\[K\left(\Lambda|\phi\right)(A)=-\lim\limits_{\alpha\to 0}\frac{\Lambda(A)}{\alpha}\log\frac{H_{1-\alpha}(\Lambda,\phi)(A)}{\Lambda(A)}.\]
\end{lemma}
{\it Proof.}
First, observe that, by the convexity of $x\longmapsto e^{-x}$,
\[H_{1-\alpha}(\Lambda,\phi)(A) \geq \int\limits_{A} e^{-\alpha\log f}d\Lambda\geq \Lambda(A)e^{-\frac{\alpha}{\Lambda(A)}\int\limits_{A}\log fd\Lambda}=\Lambda(A)e^{-\frac{\alpha}{\Lambda(A)} K\left(\Lambda|\phi\right)(A)}\]
for all $0<\alpha\leq 1$. This implies the first part of the assertion.

Now, one easily checks that $1/\alpha(x- x^{1-\alpha})\uparrow x\log x$ as $\alpha\to 0$ for all $x\geq 0$, and that the approximating 
functions are equibounded from below. Hence, by the Monotone Convergence Theorem,
\begin{eqnarray*}
	&& -\lim\limits_{\alpha\to 0}\frac{\Lambda(A)}{\alpha}\log\frac{H_{1-\alpha}(\Lambda,\phi)(A)}{\Lambda(A)}\geq\lim\limits_{\alpha\to 0}\frac{1}{\alpha}\left(\Lambda(A)-H_{1-\alpha}(\Lambda,\phi)(A)\right)\\
	&=&\lim\limits_{\alpha\to 0}\int\limits_{A}\frac{1}{\alpha}(f-f^{1-\alpha})d\phi=\int\limits_{A} f\log fd\phi.
\end{eqnarray*}
\hfill$\Box$

\begin{Definition}
	Let $A\in\A$ such that $\Lambda(A)>0$. Let $\Lambda_A$ and ${\phi}_A$ denote the measures on $\A$ given by
	\[\Lambda_A(B):=\frac{\Lambda(B\cap A)}{\Lambda(A)}\ \ \ \mbox{ and }\ \ \ \phi_A(B):=\frac{\phi(B\cap A)}{\phi(A)}\mbox{ for all }B\in\A.\]
	Set $K\left(\Lambda_A|\phi_A\right) := 0$ if $\Lambda(A)=0$.     
\end{Definition}

\begin{lemma}\label{ml}   Let $A\in\A$. Then\\
	(i) 
	\begin{equation*}
    	\Lambda(A)\log\frac{\Lambda(A)}{\phi(A)}+\Lambda(A)K\left(\Lambda_A|\phi_A\right)=K\left(\Lambda|\phi\right)(A),
	\end{equation*}
	(ii) 
	\[H_\alpha\left(\Lambda_A,{\phi}_A\right) = \frac{H_\alpha(\Lambda,\phi)(A)}{\Lambda(A)^{\alpha}\phi(A)^{1-\alpha}}\ \ \ \mbox{ for all }0\leq\alpha\leq 1\mbox{ if }\Lambda(A)>0,\mbox{ and}\]
	(iii) 
	\[\Lambda(A)\log\frac{\Lambda(A)}{\phi(A)}-\Lambda(A)\frac{1}{\alpha}\log  H_{1-\alpha}(\Lambda_A,{\phi}_A)\leq K\left(\Lambda|\phi\right)(A)\]
	for all $0<\alpha\leq 1$ if $\Lambda(A)>0$, and in the limit, as $\alpha\to 0$, holds true the equality.
	
	(iv) For every $\b,\g\in[0,1]$ such that $\b>0$ if $\g>0$ and $0<\alpha\leq 1$,
	\begin{eqnarray*}
	  &&\int\limits_Af^\g d\phi\log\frac{\int_Af^\g d\phi}{\int_Af^\b d\phi}  - \frac{\int_Af^\g d\phi}{\a}\log\frac{\int_Af^{(1-\a)\g+\a\b} d\phi}{\left(\int_Af^\g d\phi\right)^{1-\a}\left(\int_Af^\b d\phi\right)^{\a}}\\
	  &\leq&(\g-\b)\int\limits_Af^\g\log f d\phi,
	\end{eqnarray*}
	and in the limit, as $\alpha\to 0$, holds true the equality.
\end{lemma}
{\it Proof.} $(i)$ Clearly, we can assume that $\Lambda(A)>0$.  Let $f_A$ be a measurable version of the Radon-Nikodym derivative $d\Lambda_A/d\phi_A$. A straightforward computation, using the uniqueness of the Radon-Nikodym derivative, shows that
\begin{eqnarray}\label{ldf}
f_A=\frac{\phi(A)}{\Lambda(A)}f\ \ \ \ {\phi}_A\mbox{-a.e.}
\end{eqnarray}
Therefore, 
\begin{eqnarray*}
	\int f_A\log f_A d{\phi}_A&=&\frac{1}{\Lambda(A)}\int\limits_Af\left(\log\frac{\phi(A)}{\Lambda(A)}+\log f\right)d\phi\\
	&=&\log\frac{\phi(A)}{\Lambda(A)} +\frac{1}{\Lambda(A)}\int\limits_Af\log fd\phi.
\end{eqnarray*}
The multiplication by $\Lambda(A)$ implies (i).

$(ii)$ The assertion follows immediately from \eqref{ldf}. 

$(iii)$  The assertion follows from (i) and Lemma \ref{HKr}.

$(iv)$ Clearly, we only need to proof the case $\b>0$.  Define $\phi'(A):=\int_Af^\b d\phi$ and $\Lambda'(A):=\int_Af^\g d\phi$ for all $A\in\A$. Then, one easily sees that $\Lambda'\ll\phi'$, $\phi'\{f=0\}=0$ and
\[\frac{d\Lambda'}{d\phi'}=f^{\g-\b}\ \phi'\mbox{-a.e.}\]
Thus, the assertion follows from (ii) and (iii) applied to $\phi'$ and $\Lambda'$.
\hfill$\Box$ 

\begin{Remark}\label{pcr} 
	Obviously, by Lemma \ref{ml} (i) or (iii),
	\begin{equation*}
	\Lambda(A)\log\frac{\Lambda(A)}{\phi(A)}\leq\int\limits_A\log f d\Lambda.
	\end{equation*}
	Furthermore, recall that the  sum $\sum_{m}\Lambda(A_m)\log(\Lambda(A_m)/\phi(A_m))$ converges\\ monotonously to $\int\limits\log f d\Lambda$ with a converging refinement of the partitions $(A_m)$ if  $\Lambda$ and $\phi$ are probability measures (e.g. see Theorem 4.1 in  \cite{Sl}). Hence, in the stationary information theory,  the second term in Lemma \ref{ml} (i) makes no contribution in the limit.  The contribution of that term in the limit in the dynamical generalization of it, which we develop in this article, is unknown. However, despite the fact that, by Lemma \ref{HKr}, the term can be well approximated in terms of the density function (which makes it easier to estimate), the author was not able to make any use of it so far. 
\end{Remark}

\section{Lower bounds for $\Phi$ via the DDMs arising from the Hellinger integral $\H_{\a}\left(\Lambda,\phi_0\right)$ and $\H^{\a,\b}\left(\Lambda,\phi_0\right)$}\label{olbs}

First, we are going to obtain some inequalities which can be used for inferring a residual relation between $\Lambda$ and $\Phi$ from $\Lambda\ll\phi_0$ (or $K(\Lambda|\phi_0)<\infty$) which gives a lower bound for $\Phi$. 

 The following lemma lists a hierarchy of methods which can be used for a deduction of the positivity of $\Phi$.

\begin{lemma}\label{lbl}
  Let $Q\in \P(X)$  and $(A_m)_{m\leq 0}\in\C(Q)$. 
  
(i) For every $m\leq 0$ such that $\Lambda(A_m)>0$,
\begin{eqnarray*}
    &&\Lambda(A_m)^{1-\alpha}\phi_0  (S^mA_m)^\alpha\geq\int\limits_{S^mA_m}Z^{1-\alpha} d\phi_0\geq\Lambda(A_m)e^{-\frac{\alpha}{\Lambda(A_m)}\int\limits_{S^mA_m}\log Z d\Lambda}
\end{eqnarray*}
for all $0<\alpha\leq 1$, with the definition $\log(0):=-\infty$.

(ii) If $\sum_{m\leq 0}\Lambda(A_m)>0$, then, for every $0<\alpha< 1$,
\begin{eqnarray*}
	&&\sum\limits_{m\leq 0}\Lambda(A_m)^{1-\alpha}\phi_0(S^mA_m)^\alpha\leq\left(\sum\limits_{m\leq 0}\Lambda(A_m)\right)^{1-\alpha}\left(\sum\limits_{m\leq 0}\phi_0(S^mA_m)\right)^\alpha.
\end{eqnarray*}
\end{lemma}
    {\it Proof.}
    $(i)$  The first inequality of (i) is obvious, by the concavity of $x\longmapsto x^{1-\alpha}$ or the H\"{o}lder inequality, 
and the second follows immediately by the convexity of $x\longmapsto e^{-x}$, as in the proof of Lemma \ref{HKr}.

$(ii)$  By the assumption, there exists $m\leq 0$ such that $\Lambda(S^mA_m)>0$, and therefore,  $\phi_0(S^mA_m)>0$. Thus, the product on the right hand side is well defined, and one sees the assertion by the H\"{o}lder inequality, or simply the concavity of $x\longmapsto x^{1-\alpha}$.
\hfill$\Box$

Clearly, by Lemma \ref{ml}, one can obtain a stronger inequality than that of Lemma \ref{lbl} (i), since
\[\Lambda(A_m)^{1-\alpha}\phi_0  (S^mA_m)^\alpha=\Lambda(A_m)e^{-\frac{\alpha}{\Lambda(A_m)}
	\Lambda(A_m)\log\frac{\Lambda(A_m)}{\phi_0\left(S^mA_m\right)}}\]
for all $A_m\in\A_m$ with $\Lambda(A_m)>0$ and $m\leq 0$. However, because of Remark \ref{pcr}, we proceed guided by Lemma \ref{lbl} and start with the following object for the computation of lower bounds for $\Phi$, which leads to the best practical estimates which we could obtain so far.

\begin{Definition}
Let $\a\geq 0$ and $Q\in\P(X)$.  Define
\[\H_{\a}\left(\Lambda,\phi_0\right)(Q):=\inf\limits_{(A_m)_{m\leq 0}\in\C(Q)}\sum\limits_{m\leq 0}\int\limits_{S^mA_m}Z^\a d\phi_0.\]
\end{Definition}
Obviously, $\H_{0}\left(\Lambda,\phi_0\right)(Q)=\Phi(Q)$, and $\H_{1}\left(\Lambda,\phi_0\right)(Q)=\Lambda(Q)$ by Proposition 11 (i) (Proposition 1 (i) in the arXiv version) in \cite{Wer15}. For $0\leq\a\leq 1$, the following provides an approach to computations of lower bounds for $\Phi$ on $\B$.

\begin{lemma}\label{hmlb}
(i) For $0\leq\a\leq 1$,
\begin{equation*}
\H_{\a}\left(\Lambda,\phi_0\right)(Q)\leq\Phi(Q)^{1-\a}\Lambda(Q)^{\a}\ \ \ \mbox{ for all }Q\in\P(X).
\end{equation*}

(ii) $\H_{\a}\left(\Lambda,\phi_0\right)$ is a finite, $S$-invariant measure on $\B$ for all $\a\in[0,\infty)$ such that $\int Z^\a d\phi_0<\infty$, in particular for all $\a\in[0,1]$.

(iii) $\H_{\a}\left(\Lambda,\phi_0\right)\ll\Phi$ for all $\a\in[0,1)$, and $\H_{\a}\left(\Lambda,\phi_0\right)\ll\Lambda$ for all $\a\in(0,1]$.
\end{lemma}
{\it Proof.}
$(i)$ Clearly, we only need to consider the case $0<\a<1$. Let $Q\in\P(X)$, $\e>0$  and $(A_m)_{m\leq 0}\in\C(Q)$ such that
\[\sum\limits_{m\leq 0}\phi_0\left(S^mA_m\right)<\Phi(Q)+\e.\]
If $\sum_{m\leq 0}\Lambda(A_m)=0$, then $\Lambda(Q)=0$. Hence, the right hand side of the inequality is, obviously, zero, and the left hand side is also zero, since $\int_{S^mA_m}Z^\a d\phi_0\leq \int_{S^mA_m}Z^{\a-1} d\Lambda$ for all $m\leq 0$. Let $\sum_{m\leq 0}\Lambda(A_m)>0$.
Then, by Lemma \ref{lbl},
 \[\left(\Phi(Q)+\e\right)^{1-\a}\left(\sum\limits_{m\leq 0}\Lambda(A_m)\right)^{\a}\geq\sum\limits_{m\leq 0}\int\limits_{S^mA_m}Z^{\a} d\phi_0\geq \H_{\a}\left(\Lambda,\phi_0\right)(Q). \]
Hence, by the $S$-invariance of $\Lambda$, the assertion follows by Proposition 12 (i) (Proposition 2 (i) in the arXiv version) in \cite{Wer15} if $Q\in\B$, but, as we are going to show next, the same argument applies also for general $Q$, see Lemma \ref{icp} (ii). 

$(ii)$ It follows by Theorem 16 (ii) (Theorem 4 (ii) in the arXiv version) in \cite{Wer15} and (i), as $\H_{\a}\left(\Lambda,\phi_0\right)(X)$ is finite also for $\a>1$ such that $\int Z^{\a}d\phi_0<\infty$, since $(...,\emptyset,\emptyset,X)\in\C(X)$.

$(iii)$ It follows by (i).
\hfill$\Box$ 

\begin{Remark}
	Note that $\H_{\a}\left(\Lambda,\phi_0\right)\ll\Lambda$ also for $\a=0$ if $\phi_0\ll\Lambda$, by Lemma 19 (Lemma 10 in the arXiv version) in \cite{Wer15}. However, in this paper, we only assume $\Lambda\ll\phi_0$ if some additional conditions are not stated explicitly, and it does not imply $\Lambda\ll\Phi$ in general, as we will see later (Subsection \ref{satss}).
\end{Remark}

\subsection{Inductive constructions}

It turns out that one can obtain greater DDMs arising from the Hellinger integral via the inductive construction from Subsection 4.1 in \cite{Wer15}. They generalize $\H_{\a}\left(\Lambda,\phi_0\right)$ and also provide lower bounds for $\Phi$, but the main purpose for their introduction so far is their usefulness for obtaining criteria for the positivity of $\Phi$ via their dependence on the parameter.

In order to cover the needs also of other inductive constructions of measures in this paper, we will define such a construction first for general measures to the depth of 3 and then specify the case arising from the Hellinger integral.
\begin{Definition}\label{icfw}
	Let $\xi$ and $\psi$ be non-negative measures on $\A_0$. Let $Q\in\P(X)$ and $\e>0$. Let $\Xi(Q)$ be defined as $\Phi(Q)$ with $\xi$ in place of $\phi_0$.  Suppose $\Xi(Q)<\infty$. Define 
	\[\C^{\xi}_{\e}(Q):=\left\{(A_m)_{m\leq 0}\in\C(Q)|\ \sum\limits_{m\leq 0}\xi(S^mA_m)<\Xi(Q)+\e\right\}.\] 
	Further, we will use the following abbreviations for our special cases, 
	\[\C^{0}_{\e}(Q):=\C^{\phi_0}_{\e}(Q),\ \C^1_{\e}(Q):=\C^{\Lambda}_{\e}(Q),\mbox{ and,}\] 
	in general, for every $\a\geq 0$ such that $\H_{\a}(\Lambda|\phi_0)(Q)<\infty$, 
	\[\C^{\a}_{\e}(Q):=\C^{h_{\a}}_{\e}(Q)\mbox{ where }h_\a(A):=\int_A Z^\a d\phi_0 \mbox{ for all } A\in\A_0.\]  Define
	\[\Psi^\xi_{\e}(Q):=\inf\limits_{(A_m)_{m\leq 0}\in\C^\xi_\e(Q)}\sum\limits_{m\leq 0}\psi\left(S^mA_m\right)\mbox{ and}\]
	\[\Psi^\xi(Q):=\lim\limits_{\e\to 0}\Psi^\xi_{\e}(Q).\]
	If $\psi=h_\a$ for some $\a\geq 0$, we will use the special notation
	\[\H^{\a,\xi}_{\e}\left(\Lambda,\phi_0\right)(Q):=\Psi^\xi_{\e}(Q)\mbox{ and } \H^{\a,\xi}\left(\Lambda,\phi_0\right)(Q):=\Psi^\xi(Q), \mbox{ and}\]
	\[\H^{\a,\g}_{\e}\left(\Lambda,\phi_0\right)(Q):=\H^{\a,\xi}_{\e}\left(\Lambda,\phi_0\right)(Q)\mbox{ and } \H^{\a,\g}\left(\Lambda,\phi_0\right)(Q):=\H^{\a,\xi}\left(\Lambda,\phi_0\right)(Q)\mbox{ if}\]
	$\xi=h_\g$ for $\g\geq 0$ such that $\H_{\g}(\Lambda|\phi_0)(Q)<\infty$.
	Suppose $\Psi^\xi(Q)<\infty$. Define
	\[\C^{\psi, \xi}_{\e}(Q):=\left\{(A_m)_{m\leq 0}\in\C^\xi_\e(Q)|\ \sum\limits_{m\leq 0}\psi(S^mA_m)<\Psi^\xi(Q)+\e\right\}.\]
	If $\psi=h_\a$ and $\H^{\a,\xi}\left(\Lambda,\phi_0\right)(Q)<\infty$ for some $\a\geq 0$, we will use the abbreviations
	\[\C^{\a, \xi}_{\e}(Q):=\C^{\psi, \xi}_{\e}(Q),\mbox{ and }\C^{\a, \g}_{\e}(Q):=\C^{\a, \xi}_{\e}(Q)\mbox{ if }\xi=h_\g\]
	for $\g\geq 0$ such that $\H_{\g}(\Lambda|\phi_0)(Q)<\infty$. For a non-negative measure $\o$ on $\A_0$, define
	\[\O^{\psi,\xi}_{\e}(Q):=\inf\limits_{(A_m)_{m\leq 0}\in\C^{\psi,\xi}_\e(Q)}\sum\limits_{m\leq 0}\o\left(S^mA_m\right)\mbox{ and}\]
	\[\O^{\psi,\xi}(Q):=\lim\limits_{\e\to 0}\O^{\psi,\xi}_{\e}(Q).\]
	Analogously, in the special case where $\xi=h_\a$, $\psi=h_\b$, and $\o=h_\g$, we will use the notation 
	\[\H^{\g,\b,\a}_{\e}\left(\Lambda,\phi_0\right)(Q):=\O^{\psi,\xi}_{\e}(Q)\mbox{ and }\H^{\g,\b,\a}\left(\Lambda,\phi_0\right)(Q):=\O^{\psi,\xi}(Q).\]
\end{Definition}

The inductive continuation of the construction is obvious. By Theorem 16 (Theorem 4 in the arXiv version) in \cite{Wer15}, each of the set functions $\Xi$, $\Psi^\xi$, and $\O^{\psi,\xi}$ in the above definition is an invariant measures on $\B$ if it is finite (and each of the previous ones in the inductive construction is also a finite measure on $\B$). 

The property (i) in the next lemma is called {\it the regularity} of the outer measure, see \cite{B} for more on it. We show now that this property extends also on set functions obtained as outer measure approximations \cite{Wer15}, as constructed in Definition \ref{icfw}.
\begin{lemma}\label{omr} 
	Let $\xi$, $\psi$, and $\o$ be non-negative measures on $\A_0$ such that $\Xi(X)<\infty$. Let $\Xi$, $\Psi^\xi$, and $\O^{\psi,\xi}$ be as in Definition \ref{icfw}. Let $Q\in\P(X)$. Then there exists $B\in\B$ such that $Q\subset B$,
	
	(i)  $\Xi(Q)=\Xi(B)$,
	
	(ii) $\Psi^\xi(Q)=\Psi^\xi(B)$, and
	
	(iii) if $\Psi^\xi(X)<\infty$, then $\O^{\psi,\xi}(Q)=\O^{\psi,\xi}(B)$.
\end{lemma}
{\it Proof.} 
$(i)$ For each $n\in\N$, choose $(A_m^n)_{m\leq 0}\in\C^\xi_{1/n}(Q)$ and define $B:=\bigcap_{n\in\N}\bigcup_{m\leq 0}A^n_m$. Then, obviously, $B\in\B$ and $Q\subset B$. Since, for each $n\in\N$, $(A_m^n)_{m\leq 0}\in\C(B)$, 
\[\Xi(B)\leq\sum\limits_{m\leq 0}\xi\left(S^mA_m^n\right)<\Xi(Q)+\frac{1}{n}.\] 
Hence, it follows (i), since, by Lemma 2 (Lemma 1 in the arXiv version) in \cite{Wer15}, $\Xi$ is an outer measure on $X$.

$(ii)$ By the above, for every $\e>0$,
\[\C^\xi_\e(B)\subset\C^\xi_\e(Q),\]
and therefore,
\begin{equation*}
   \Psi^\xi(Q)\leq\Psi^\xi(B).
\end{equation*}
Clearly, the equality of (ii) holds if $\Psi^\xi(Q)=\infty$.

Now, suppose $\Psi^\xi(Q)<\infty$. Then, for each $n\in\N$, there exists $(B_m^n)_{m\leq 0}\in\C^{\psi,\xi}_{1/n}(Q)$. Define $\tilde B:=\bigcap_{n\in\N}\bigcup_{m\leq 0}B^n_m$. Then $\tilde B\in\B$ and $Q\subset\tilde B$. Since, for each $n\in\N$, $(B_m^n)_{m\leq 0}\in\C(\tilde B)$, and $\Xi$ is an outer measure on $X$, 
\[\Xi\left(\tilde B\right)\leq\sum\limits_{m\leq 0}\xi\left(S^mB_m^n\right)<\Xi(Q)+\frac{1}{n}\leq\Xi\left(\tilde B\right)+\frac{1}{n}.\]
Hence,
\[\Xi(Q)=\Xi\left(\tilde B\right),\mbox{ and }(B_m^n)_{m\leq 0}\in\C^\xi_{1/n}(\tilde B) \mbox{ for all } n\in\N.\]
Hence, for each $n\in\N$,
\[\Psi^\xi_{1/n}\left(\tilde B\right)\leq\sum\limits_{m\leq 0}\psi\left(S^mB^n_m\right)<\Psi^\xi(Q)+1/n,\]
and therefore,
\[\Psi^\xi\left(\tilde B\right)\leq\Psi^\xi(Q).\]
On the other hand, since $\C^\xi_\e(\tilde B)\subset\C^\xi_\e(Q)$ for all $\e>0$, it follows that $\Psi^\xi\left(\tilde B\right)\geq\Psi^\xi\left(Q\right)$. Thus
\[\Psi^\xi\left(\tilde B\right)=\Psi^\xi\left(Q\right),\]
and we obtain $B:=\tilde B$ with $Q\subset B$ and the properties of (i) and (ii).

$(iii)$ Suppose $\Psi^\xi(X)<\infty$. By (i) and (ii), $\C^{\psi, \xi}_\e(B)\subset\C^{\psi, \xi}_\e(Q)$ for all $\e>0$. Hence,
\[\O^{\psi,\xi}(Q)\leq\O^{\psi,\xi}\left(B\right).\]
Thus the equality of (iii) is true if $\O^{\psi,\xi}(Q)=\infty$.

Finally, assuming $\O^{\psi,\xi}(Q)<\infty$, taking $(B_m^n)_{m\leq 0}\in\C^{\psi,\xi}_{1/n}(Q)$ such that
\[\sum\limits_{m\leq 0}\o\left(S^mB ^n_m\right)<\O^{\psi,\xi}(Q)+1/n\] for all $n\in\N$, and proceeding the same way as in (ii) gives $B\in\B$ with $Q\subset B$ which satisfies (i), (ii), and (iii).
\hfill$\Box$ 

It is obvious from the proof of Lemma \ref{omr} that it can be also proved to any depth of the inductive construction if needed.

One of the consequences of the regularity of the obtained set functions and the fact that they are measures on $\B$ when they are finite is that they must be outer measures on $X$ in this case.  
\begin{prop}\label{omp} 
	Let $\xi$, $\psi$, and $\o$ be non-negative measures on $\A_0$ such that $\Xi(X)<\infty$ and $\Psi^\xi(X)<\infty$. Then 
	
	(i) $\Psi^\xi$ is an outer measure on $X$, and 
	
	(ii) $\O^{\psi,\xi}$ is also an outer measure on $X$ if it is finite.
\end{prop}
{\it Proof.}
  We will give a simultaneous proof (i) and (ii) (it will be obvious in it that no assumption on $\O^{\psi,\xi}$ is required for the proof of (i)). Suppose $\O^{\psi,\xi}(X)<\infty$.
  
  First, we show the monotonicity of the set functions. Let $Q_1\subset Q_2\subset X$. By Lemma \ref{omr}, there exist $B_1, B_2\in\B$ such that $Q_1\subset B_1$, $Q_2\subset B_2$, $\Xi(Q_i)=\Xi(B_i)$, $\Psi^\xi(Q_i)=\Psi^\xi(B_i)$, and $\O^{\psi,\xi}(Q_i)=\O^{\psi,\xi}(B_i)$ for all $i\in\{1,2\}$. Define $B:=B_1\cap B_2$. Then $Q_1=Q_1\cap Q_2\subset B$, and therefore, $C(B)\subset C(Q_1)$. Let $\e>0$ and $(A_m)_{m\leq 0}\in\C_\e^\xi(B)$. Then
  \[\sum\limits_{m\leq 0}\xi\left(S^mA_m\right)<\Xi(B)+\e\leq\Xi(B_1)+\e=\Xi(Q_1)+\e.\]
  Hence, $(A_m)_{m\leq 0}\in\C_\e^\xi(Q_1)$, i.e. 
 \begin{equation}\label{ommp}
    \C_\e^\xi(B)\subset\C_\e^\xi(Q_1),
 \end{equation} and therefore, since $\Psi^\xi$ is a measure on $\B$,
  \[\Psi^\xi(Q_1)\leq\Psi^\xi(B)\leq\Psi^\xi(B_2)=\Psi^\xi(Q_2).\]
  This proves the monotonicity of $\Psi^\xi$. Let $(B_m)_{m\leq 0}\in\C_\e^{\psi,\xi}(B)$. Then
  \[\sum\limits_{m\leq 0}\psi\left(S^mB_m\right)<\Psi^\xi(B)+\e\leq\Psi^\xi(B_1)+\e=\Psi^\xi(Q_1)+\e,\]
  which together with \eqref{ommp} shows that $\C_\e^{\psi,\xi}(B)\subset\C_\e^{\psi,\xi}(Q_1)$. Hence, since $\O^{\psi,\xi}$ is a measure on $\B$,
  \[\O^{\psi,\xi}(Q_1)\leq\O^{\psi,\xi}(B)\leq\O^{\psi,\xi}(B_2)=\O^{\psi,\xi}(Q_2),\]
  which proves the monotonicity of $\O^{\psi,\xi}$.
  
  Finally, we turn to the countable subadditivity of the set functions. Let $(Q_n)_{n\in\N}\subset\P(X)$. By Lemma \ref{omr}, there exist $(B_n)_{n\in\N}\subset\B$ such that $Q_n\subset B_n$, $\Psi^\xi(Q_n)=\Psi^\xi(B_n)$, and $\O^{\psi,\xi}(Q_n)=\O^{\psi,\xi}(B_n)$ for all $n\in\N$.  Then, by the monotonicity of $\Psi^\xi$ and the fact that it is a measure on $\B$,
  \[\Psi^\xi\left(\bigcup\limits_{n\in\N}Q_n\right)\leq\Psi^\xi\left(\bigcup\limits_{n\in\N}B_n\right)
  \leq\sum\limits_{n\in\N}\Psi^\xi\left(B_n\right)=\sum\limits_{n\in\N}\Psi^\xi\left(Q_n\right),\]
  which proves the countable subadditivity of $\Psi^\xi$. Obviously, the same argument works also for $\O^{\psi,\xi}$.
\hfill$\Box$

The following lemma extends Proposition 12 (Proposition 2 in the arXiv version) and Proposition 13 (Proposition 3 in the arXiv version) in \cite{Wer15}. In particular, it shows that the outer measure $\Psi^\xi$ (obtained in Proposition \ref{omp} (i)) coincides with the outer measure $\Psi$ if $\psi$ or $\xi$ is invariant.
\begin{lemma}\label{icp} 
	Let $\xi$ and $\psi$ be non-negative measures on $\A_0$ such that $\Xi(X)<\infty$ and $\Psi(X)<\infty$. Let $Q\in\P(X)$. Then 
		
	(i) $\Xi^\Lambda(Q)=\Xi(Q)$,
	
	(ii) $\Lambda^\xi_\e(Q)=\Lambda(Q)$ for all $\e>0$,
	
	(iii) $\Psi^{\xi,\Lambda}(Q)=\Psi^\xi(Q)$ if $\Psi^\xi(X)<\infty$, and 
	
	(iv)  $\Lambda^{\psi,\xi}_\e(Q)=\Lambda(Q)$ for all $\e>0$ if $\Psi^\xi(X)<\infty$.
\end{lemma}
{\it Proof.}
$(i)$ By Proposition 13 (Proposition 3 in the arXiv version) in \cite{Wer15}, $\Xi^\Lambda(B)=\Xi(B)$ for all $B\in\B$. By Lemma \ref{omr} (i), there exists $B'\in\B$ such that $Q\subset B'$ and $\Xi(Q) =\Xi(B')$.  By Proposition \ref{omp} (i), $\Xi^\Lambda$ is an outer measure on $X$, and therefore,
\[\Xi^\Lambda(Q)\leq\Xi^\Lambda(B')=\Xi(B')=\Xi(Q).\]
Thus
\[\Xi^\Lambda(Q)\leq\Xi(Q),\]
and the converse inequality is obvious.

$(ii)$ Let $\e>0$. By (i), there exists $(A_m)_{m\leq 0}\in\C^\Lambda_\e(Q)$ such that
\[\sum\limits_{m\leq 0}\xi\left(S^mA_m\right)<\Xi^\Lambda(Q)+\e=\Xi(Q)+\e.\]
Hence, $(A_m)_{m\leq 0}\in\C^\xi_\e(Q)$, and therefore,
\[\Lambda ^\xi_\e(Q)\leq\sum\limits_{m\leq 0}\Lambda\left(A_m\right)<\Lambda(Q)+\e,\]
which implies that
\[\Lambda^{\xi}(Q)\leq\Lambda(Q).\]
Also, obviously, $\Lambda(Q)\leq\Lambda ^\xi_\e(Q)\leq\Lambda^{\xi}(Q)$. It completes the proof of (ii).

Now, let $\Psi^\xi(X)<\infty$.

$(iii)$  By Proposition 13 (Proposition 3 in the arXiv version) in \cite{Wer15}, $\Psi^{\xi,\Lambda}(B')=\Psi^\xi(B')$ for all $B'\in\B$. By Lemma \ref{omr} (ii), there exists $B\in\B$ such that $Q\subset B$ and $\Psi^{\xi}(Q)=\Psi^\xi(B)$. Hence, by (i) and Proposition \ref{omp} (ii),
\[\Psi^{\xi}(Q)\leq\Psi^{\xi,\Lambda}(Q)\leq\Psi^{\xi,\Lambda}(B)=\Psi^{\xi}(B)=\Psi^\xi(Q).\]

$(iv)$ By (iii), there exists $(B_m)_{m\leq 0}\in\C^{\xi,\Lambda}_\e(Q)$ such that
\[\sum\limits_{m\leq 0}\psi\left(S^mB_m\right)<\Psi^{\xi,\Lambda}(Q)+\e=\Psi^{\xi}(Q)+\e.\]
Also, by (i),
\[\sum\limits_{m\leq 0}\xi\left(S^mB_m\right)<\Xi^\Lambda(Q)+\e=\Xi(Q)+\e.\]
Hence, $(B_m)_{m\leq 0}\in\C^{\psi,\xi}_\e(Q)$, and therefore,
\[\Lambda ^{\psi,\xi}_\e(Q)\leq\sum\limits_{m\leq 0}\Lambda\left(B_m\right)<\Lambda(Q)+\e,\]
which implies that
\[\Lambda^{\psi,\xi}(Q)\leq\Lambda(Q).\]
Since, obviously, $\Lambda(Q)\leq\Lambda^{\psi,\xi}_\e(Q)\leq\Lambda^{\psi,\xi}(Q)$, it follows (iv).
\hfill$\Box$

Now, we show that, for $0\leq\a_0\leq\a< 1$, $\H^{\a,\a_0}\left(\Lambda,\phi_0\right)$ also provides a lower bound for $\Phi$. The result also enables us to shed some light on the dependence of $\H_{\a}\left(\Lambda,\phi_0\right)(Q)$ on $\a$.

\begin{lemma}\label{shm} Let $Q\in\P(X)$.
	
	(i)  \[\H^{\a,\a_0}\left(\Lambda,\phi_0\right)(Q)\leq\H_{\a_0}\left(\Lambda,\phi_0\right)(Q)^{\frac{1-\a}{1-\a_0}}\Lambda(Q)^\frac{\a-\a_0}{1-\a_0}\  \mbox{ for all } 0\leq\a_0\leq\a<1.\] 
	
	(ii) For every $0\leq\a_0\leq\a\leq1$, $\H^{\a,\a_0}\left(\Lambda,\phi_0\right)$ is a finite, $S$-invariant measure on $\B$.
	
	(iii)  If $\H_{\a}(\Lambda,\phi_0)(Q)>0$ for some $\a\in(0,1)$, then $\H_{\a_0}(\Lambda,\phi_0)(Q)>0$ for all $\a_0\in[0,\a]\cup\{1\}$.
	If $\H_{\a_0}(\Lambda,\phi_0)(Q)=0$ for some $\a_0\in[0,1)$, then $\H_{\a}(\Lambda,\phi_0)(Q)=0$ for all $\a\in[\a_0, 1)$. 
\end{lemma}
{\it Proof.}
$(i)$ Clearly, we can assume that $\a_0<\a$.  Let $\e>0$ and $(A_m)_{m\leq 0}\in\C^{\a_0,1}_{\e}(Q)$. Let us first consider the case $\sum_{m\leq 0}\Lambda(A_m)>0$. By the convexity of $x\longmapsto x^{(1-\a_0)/(1-\a)}$ and Lemma \ref{icp} (i),
\begin{eqnarray*}
\H_{\a_0}\left(\Lambda,\phi_0\right)(Q)+\e&>&\sum\limits_{m\leq 0}\int\limits_{S^mA_m}Z^{\a_0} d\phi_0\\
&\geq&\sum\limits_{m\leq 0}\int\limits_{S^mA_m}\left(Z^{\a-1}\right)^\frac{1-\a_0}{1-\a} d\Lambda\\
&\geq&\left(\sum\limits_{m\leq 0}\Lambda(A_m)\right)^{1-\frac{1-\a_0}{1-\a}}\left(\sum\limits_{m\leq 0}\int\limits_{S^mA_m}Z^\a d\phi_0\right)^{\frac{1-\a_0}{1-\a}}\\
&\geq&\left(\Lambda(Q)+\e\right)^{1-\frac{1-\a_0}{1-\a}}\H^{\a,\a_0}_\e\left(\Lambda,\phi_0\right)(Q)^{\frac{1-\a_0}{1-\a}}.
\end{eqnarray*}
If $\sum_{m\leq 0}\Lambda(A_m)=0$, then $\sum_{m\leq 0}\int_{S^mA_m}Z^\a d\phi_0=0$, and the inequality 
\[\H_{\a_0}\left(\Lambda,\phi_0\right)(Q)+\e>\left(\Lambda(Q)+\e\right)^{1-\frac{1-\a_0}{1-\a}}\left(\sum\limits_{m\leq 0}\int\limits_{S^mA_m}Z^\a d\phi_0\right)^{\frac{1-\a_0}{1-\a}}\]
is obviously correct also. It implies (i).

$(ii)$ It follows immediately from (i), Lemma \ref{icp} (ii), and Theorem 16 (ii) (Theorem 4 (ii) in the arXiv version) in \cite{Wer15}.

$(iii)$ It follows immediately from (i) and Lemma \ref{hmlb} (iii).
\hfill$\Box$

The previous lemma does not say anything about $\H^{\a,\a_0}\left(\Lambda,\phi_0\right)$ if $\a<\a_0$. However, as we have already seen in Lemma \ref{icp} and also it will be shown later (Subsection \ref{homs}) that, in some cases, the hypothesis of the following lemma is satisfied.
\begin{lemma}\label{hsl} 
		Let $\xi$ and $\psi$ be non-negative measures on $\A_0$. Let $Q\in\P(X)$ such that $\Xi(Q)<\infty$ and $\Psi(Q)<\infty$. Then 
		\[\Xi^\psi(Q)=\Xi(Q)\implies\Psi^\xi(Q)=\Psi(Q).\]
\end{lemma}
{\it Proof.} Let $\Xi^\psi(Q)=\Xi(Q)$. Let $\e>0$. Then $\C^{\xi,\psi}_\e(Q)\subset\C^{\xi}_\e(Q)$. Hence, for every $(A_m)_{m\leq 0}\in\C^{\xi,\psi}_\e(Q)$,
 \[\Psi^\xi_\e(Q)\leq\Psi^{\xi,\psi}_\e(Q)\leq\sum\limits_{m\leq 0}\psi(S^mA_m)<\Psi(Q)+\e.\]
 Thus
 \[\Psi^\xi(Q)\leq\Psi(Q).\]
The converse inequality is obvious.
\hfill$\Box$

\subsection{A lower bound for $\H_{\a}\left(\Lambda,\phi_0\right)$ via a relative entropy measure}\label{rems}

For the purpose of obtaining a lower bound for $\H_{\a}\left(\Lambda,\phi_0\right)(Q)$, first observe that, by Lemma \ref{edl} (i), for every $0\leq\a<\g$ and $(A_m)_{m\leq 0}\in\C(Q)$ such that $\sum_{m\leq 0}\int_{S^mA_m}Z^\a d\phi_0<\infty$,
\begin{eqnarray}\label{hrds}
&&(\g-\a)\sum\limits_{m\leq 0, \int_{S^mA_m}Z^\g\log Z d\phi_0<0}\int\limits_{S^mA_m}Z^\g\log Z d\phi_0\\
&\geq&\sum\limits_{m\leq 0, \int_{S^mA_m}Z^\g\log Z d\phi_0<0}\int\limits_{S^mA_m}\left(Z^\g-Z^\a\right) d\phi_0\geq-\sum\limits_{m\leq 0}\int\limits_{S^mA_m}Z^\a d\phi_0>-\infty.\nonumber
\end{eqnarray}
Therefore, the sum in the following expression is well defined.
\begin{Definition}
	For $0\leq\a<\g\leq 1$, $Q\in\P(X)$ and $\e>0$, define
	\[\D^{\g,\a}_\e(Q):=\inf\limits_{(A_m)_{m\leq 0}\in\C^\a_{\e}(Q)}\sum\limits_{m\leq 0}\int\limits_{S^mA_m}Z^\g\log Z d\phi_0\]
	and
	\[\D^{\g,\a}(Q):=\lim\limits_{\e\to 0}\D^{\g,\a}_\e(Q).\]
	The same way as in the proof of  Lemma 5 (Lemma 3 in the arXiv version) in \cite{Wer15}, on sees that
	\begin{equation*}
	\D^{\g,\a}_\e(Q)\leq\D^{\g,\a}_\e(S^{-1}Q)\ \ \ \mbox{ for all }Q\in\P(X)\mbox{ and }\e>0.
	\end{equation*}
	Therefore, we can define 
	\[\bar\D^{\g,\a}_\e(Q):=\lim\limits_{n\to\infty}\D^{\g,\a}_\e(S^{-n}Q)\ \ \ \mbox{ for all }Q\in\P(X)\mbox{ and }\e>0, \mbox{ and }\]
	\[\bar\D^{\g,\a}(Q):=\lim\limits_{\e\to 0}\bar\D^{\g,\a}_\e(Q)\ \ \ \mbox{ for all }Q\in\P(X).\]
	One easily sees that
	\[\bar\D^{\g,\a}(Q)=\lim\limits_{n\to\infty}\D^{\g,\a}(S^{-n}Q)\ \ \ \mbox{ for all }Q\in\P(X).\]
	
	Let $\dot\C^\a_\e(Q)$ denote the set of all $(A_m)_{m\leq 0}\in\C^\a_{\e}(Q)$ such that $A_{m}$'s are pairwise disjoint. By Lemma 3 (Lemma 2 in the arXiv version) in \cite{Wer15}, $\dot\C^\a_{\e}(Q)$ is not empty. Define $\dot \D^{\g,\a}_\e(Q)$ the same way as $\D^{\g,\a}_\e(Q)$ with the infimum taken over $\dot\C^\a_\e(Q)$ and $\dot\D^{\g,\a}(Q)$ analogously.
	
	For the important case $\g=1$, we will use the special notation 
	\[\K_{\a,\e}\left(\Lambda,\phi_0\right):=\D^{1,\a}_\e\mbox{ and }\K_{\a}\left(\Lambda,\phi_0\right):=\D^{1,\a}.\]
\end{Definition}

\begin{Definition}
For every $0\leq\a<\g\leq 1$ and $A\in\A_0$, define
\[\kappa^{\g,\a}(A):=\int\limits_A\left(Z^\g\log Z+\frac{1}{\g-\a}Z^\a\right)d\phi_0,\]
and let $\K^{\g,\a}_\e$, $\K^{\g,\a}$ and $\bar\K^{\g,\a}$ be defined the same way as $\D^{\g,\a}_\e$, $\D^{\g,\a}$ and $\bar\D^{\g,\a}$ with $\int_{A}Z^\g\log Zd\phi_0$ replaced by $\kappa^{\g,\a}(A)$.
\end{Definition}

The obtained set functions have the following properties.

\begin{lemma}\label{eml} Let $0\leq\a<\g\leq 1$. Then the following holds true.
	
	(i) For every $Q\in\P(X)$,
	\[\frac{1}{\g-\a}\H_{\g}\left(\Lambda,\phi_0\right)(Q)\leq\K^{\g,\a}(Q),\]
	and for every $Q\in\B$  and $\a<\g<1$,
	\[\K^{\g,\a}(Q)\leq\frac{1}{1-\g}\left(\Lambda(Q)-\H_{\g}\left(\Lambda,\phi_0\right)(Q)\right)+\frac{1}{\g-\a}\H_{\a}\left(\Lambda,\phi_0\right)(Q).\]
	
	(ii) \[\D^{\g,\a}(Q)=  \K^{\g,\a}(Q) -\frac{1}{\g-\a}\H_{\a}\left(\Lambda,\phi_0\right)(Q)\ \ \ \mbox{ for all }Q\in\B.\]
	
	(iii) \[\D^{\g,\a}(Q)= \dot \D^{\g,\a}(Q)\ \ \ \mbox{ for all }Q\in\B.\]
	
	(iv) $\bar\D^{\g,\a}$ is a $S$-invariant, signed measure on $\B$.
	
	(v) $\bar\D^{\g,\a}(Q)=\D^{\g,\a}(Q)$ for all $Q\in\B$ if $\g<1$, and 
	\[\K_\a\left(\Lambda|\phi_0\right)(Q)=\bar\K_\a\left(\Lambda|\phi_0\right)(Q) \mbox{  for all }Q\in\B\mbox{ if } \K_\a\left(\Lambda|\phi_0\right)(X)<\infty.\]
	
	(vi) $\K_\a\left(\Lambda|\phi_0\right)(X) = K\left(\Lambda|\phi_0\right)$ if $\phi_0\circ S^{-1} = \phi_0$.
\end{lemma}
{\it Proof.}
$(i)$ Let $Q\in\P(X)$, $\e>0$ and $(A_m)_{m\leq 0}\in\C^\a_{\e}(Q)$. Then,  by Lemma \ref{edl} (i),
\begin{eqnarray*}
  &&\frac{1}{\g-\a}\H_{\g}\left(\Lambda,\phi_0\right)(Q)\\
  &\leq&\frac{1}{\g-\a}\sum\limits_{m\leq 0}\int\limits_{S^mA_m}Z^\g d\phi_0\leq\sum\limits_{m\leq 0}\int\limits_{S^mA_m}\left(Z^\g\log Z+\frac{1}{\g-\a}Z^\a\right)d\phi_0\\
  &=&\sum\limits_{m\leq 0}\kappa^{\g,\a}\left(S^mA_m\right).
\end{eqnarray*}
Thus, the first inequality of (i) follows.

Now, let $Q\in\B$ and $\g<1$. By Lemma \ref{icp} (ii), we can choose $(B_m)_{m\leq 0}\in\C^\a_{\e}(Q)$ such that $\sum_{m\leq 0}\Lambda(B_m)<\Lambda(Q)+\e$. Then, by Lemma \ref{edl} (i), 
\begin{eqnarray*}
&&\K^{\g,\a}_\e(Q)\\
&\leq&\frac{1}{1-\g}\left(\sum\limits_{m\leq 0}\Lambda(B_m)-\sum\limits_{m\leq 0}\int\limits_{S^mB_m}Z^\g d\phi_0\right) +\frac{1}{\g-\a}\sum\limits_{m\leq 0}\int\limits_{S^mB_m}Z^\a d\phi_0\\
&\leq&\frac{1}{1-\g}\left(\Lambda(Q)+\e-\H_{\g}\left(\Lambda,\phi_0\right)(Q)\right) +\frac{1}{\g-\a}\left(\H_{\a}\left(\Lambda,\phi_0\right)(Q)+\e\right).
\end{eqnarray*}
Thus, the second inequality of (i) follows.
 
$(ii)$  It follows immediately by Lemma 10 (i) (Lemma 6 (i) in the arXiv version) in \cite{Wer15}.

$(iii)$ It follows immediately by (ii) and Lemma 10 (ii) ( Lemma 6 (ii) in the arXiv version) in \cite{Wer15}.

$(iv)$ By (ii),
\[\bar\D^{\g,\a}(Q) =  \bar\K^{\g,\a}(Q) -\frac{1}{\g-\a}\H_{\a}\left(\Lambda,\phi_0\right)(Q)\ \ \ \mbox{ for all }Q\in\B.\]
Thus, (iv) follows by Theorem 7 (Theorem 3 in the arXiv version) in \cite{Wer15}.

$(v)$ The assertion follows immediately by (i), (ii) and Theorem 16 (Theorem 4 in the arXiv version) in \cite{Wer15}.

$(vi)$ Observe that, by the hypothesis, $Z\circ S^{-1} = Z$ $\phi_0$-a.e. Therefore, for every $\e>0$ and $(A_m)_{m\leq 0}\in\dot\C^\a_{\e}(X)$,
\[\sum\limits_{m\leq 0}\int\limits_{S^mA_m}\log Zd\Lambda = \sum\limits_{m\leq 0}\int\limits_{A_m}\log Zd\Lambda=\int\log Zd\Lambda.\]
Thus, the assertion follows by (iii).
\hfill$\Box$

\begin{Remark}\label{cuf}
	Note that  $\K_{\a}(\Lambda|\phi_0)(Q)$ can be infinite. However, by Lemma 17 (Lemma 9 in the arXiv version) in \cite{Wer15}, for every $\e>0$, $\K_{\a,\e}(\Lambda|\phi_0)(Q)$ is finite for a broad class of topological dynamical systems if $K(\Lambda|\phi_0)$ is finite and $Q$ is compact.
\end{Remark}

The following theorem gives some lower bounds for $\H_{\a}\left(\Lambda,\phi_0\right)$ by capturing some residual of the relation from Lemma \ref{HKr}.

\begin{theo}\label{lbre}
	Let $Q\in\B$ and $0\leq\a<\g\leq 1$.
	
	(i) Let $\e>0$ such that $\H^{\g,\a}_\e\left(\Lambda,\phi_0\right)(Q)>0$. Then
	\begin{equation*}
	\H_{\a}\left(\Lambda,\phi_0\right)(Q)\geq\H^{\g,\a}_\e\left(\Lambda,\phi_0\right)(Q)\min\left\{e^{-\frac{\g-\a}{\H^{\g,\a}_\e\left(\Lambda,\phi_0\right)(Q)}\D^{\g,\a}_\e(Q)},e\right\}-\e,
	\end{equation*}
	and
		\[\Phi(Q)\geq\Lambda(Q)\min\left\{e^{-\frac{1}{\Lambda(Q)}\K_{\a}(\Lambda|\phi_0)(Q)},e^\frac{1}{1-\a}\right\}\mbox{ if }\Lambda(Q)>0.\]
	
	(ii) 
	\begin{equation*}
	   \H_{\a}\left(\Lambda,\phi_0\right)(Q)\geq\Lambda(Q)e^{-\frac{1-\a}{\Lambda(Q)}\K_{\a, \e}(\Lambda|\phi_0)(Q)}-\e 	
	\end{equation*}
	for all $0<\e<\Lambda(Q)\left(e-(\H_{\a_0}\left(\Lambda,\phi_0\right)(Q)/\Lambda(Q))^{(1-\a)/(1-\a_0)}\right)$ and $0\leq\a_0\leq\a$, \mbox{and}
	\[\H_{\a}\left(\Lambda,\phi_0\right)(Q)\geq\Lambda(Q)e^{-\frac{1-\a}{\Lambda(Q)}\K_{\a}(\Lambda|\phi_0)(Q)}\]
	if $\H_{\a_0}\left(\Lambda,\phi_0\right)(Q)<\Lambda(Q)e^{(1-\a_0)/(1-\a)}$ for some $0\leq\a_0\leq\a$.
\end{theo}
{\it Proof.}
$(i)$ Clearly, we can assume that $\K_{\a,\e}(\Lambda|\phi_0)(Q)<\infty$ and $\H_{\a}\left(\Lambda,\phi_0\right)(Q)<\H^{\g,\a}_\e\left(\Lambda,\phi_0\right)e-\e$.  

Suppose $\D^{\g,\a}_\e(Q) = 0$. Let $\tau>0$. Then there exists $(B_m)_{m\leq 0}\in\C^{\a}_{\e}(Q)$ such that
\[\sum\limits_{m\leq 0}\int\limits_{S^mB_m}Z^\g\log Z d\phi_0<\tau.\]
Therefore, by Lemma \ref{edl} (i),
\begin{eqnarray}\label{zclb}
\H_{\a}\left(\Lambda,\phi_0\right)(Q)+\e&>&\sum\limits_{m\leq 0}\int\limits_{S^mB_m}Z^{\a} d\phi_0\nonumber\\
&\geq& \sum\limits_{m\leq 0}\int\limits_{S^mB_m}(Z^\g -(\g-\a) Z^\g\log Z)d\phi_0\\
&\geq&\H^{\g,\a}_\e\left(\Lambda,\phi_0\right)(Q)-(\g-\a)\tau.\nonumber
\end{eqnarray}
Hence,
\[\H_{\a}\left(\Lambda,\phi_0\right)(Q)+\e\geq\H^{\g,\a}_\e\left(\Lambda,\phi_0\right)(Q). \]
This proves the assertion in the case $\D^{\g,\a}_\e(Q) = 0$.

Now, suppose $\D^{\g,\a}_\e(Q) \neq 0$. Let $\tau_0>0$  be such that $\D^{\g,\a}_\e(Q) +\tau$ has the same sign as $\D^{\g,\a}_\e(Q)$ for all $0<\tau<\tau_0$. Let $0<\tau<\tau_0$ and $(A_m)_{m\leq 0}\in\C^{\a}_{\e}(Q)$ such that
\begin{equation*}
\D^{\g,\a}_\e(Q) +\tau > \sum\limits_{m\leq 0}\int\limits_{S^mA_m}Z^\g\log Z d\phi_0.
\end{equation*}
Then, as in \eqref{zclb}, one sees that $\sum_{m\leq 0}\int_{S^mA_m}Z^\g d\phi_0<\infty$. Therefore, by Lemma \ref{ml} (iv) and the convexity of $x\longmapsto e^{-x}$,
\begin{eqnarray*}
	\H_{\a}\left(\Lambda,\phi_0\right)(Q)+\e&>&\sum\limits_{m\leq 0}\int\limits_{S^mA_m}Z^{\a} d\phi_0\\
	&\geq&\sum\limits_{m\leq 0}\int\limits_{S^mA_m}Z^{\g} d\phi_0e^{-\frac{\g-\a}{\sum\limits_{m\leq 0}\int\limits_{S^mA_m}Z^{\g} d\phi_0}\sum\limits_{m\leq 0} \int\limits_{S^mA_m}Z^{\g}\log Z d\phi_0}\\
	&\geq&\sum\limits_{m\leq 0}\int\limits_{S^mA_m}Z^{\g} d\phi_0e^{-\frac{\g-\a}{\sum\limits_{m\leq 0}\int\limits_{S^mA_m}Z^{\g} d\phi_0}\left(\D^{\g,\a}_\e(Q)+\tau\right)}.
\end{eqnarray*}
That is
\begin{equation}\label{ebWa}
\frac{1}{\sum\limits_{m\leq 0}\int\limits_{S^mA_m}Z^{\g} d\phi_0}e^{\frac{\g-\a}{\sum\limits_{m\leq 0}\int\limits_{S^mA_m}Z^{\g} d\phi_0}\left(\D^{\g,\a}_\e(Q)+\tau\right)}>\frac{1}{\H_{\a}\left(\Lambda,\phi_0\right)(Q)+\e}.
\end{equation}
Observe that by the assumption on $\H_{\a}\left(\Lambda,\phi_0\right)(Q)$, this implies that
\[\frac{(\g-\a)\left(\D^{\g,\a}_\e(Q) +\tau\right)}{\sum\limits_{m\leq 0}\int\limits_{S^mA_m}Z^{\g} d\phi_0}>\log\frac{\H^{\g,\a}_\e\left(\Lambda,\phi_0\right)(Q)}{\H_{\a}\left(\Lambda,\phi_0\right)(Q)+\e}>-1.\] 
Hence, since the principal branch of Lambert's $W$ function is monotonously increasing, \eqref{ebWa} implies (regardless of the sign of $\D^{\g,\a}_\e(Q)+\tau$) that
\[\H_{\g}\left(\Lambda,\phi_0\right)(Q)\leq\sum\limits_{m\leq 0}\int\limits_{S^mA_m}Z^{\g} d\phi_0<\frac{(\g-\a)\left(\D^{\g,\a}_\e(Q)+\tau\right)}{W\left(\frac{(\g-\a)\left(\D^{\g,\a}_\e(Q)+\tau\right)}{\H_{\a}\left(\Lambda,\phi_0\right)(Q)+\e}\right)}.\]
Therefore, since, by the definition of $W$, $x/W(x)=e^{W(x)}$ for all $x\in[-1/e,\infty)\setminus\{0\}$,
\begin{eqnarray*}
  \log\frac{\H^{\g,\a}_\e\left(\Lambda,\phi_0\right)(Q)}{\H_{\a}\left(\Lambda,\phi_0\right)(Q)+\e}&<&\log\frac{\frac{(\g-\a)\left(\D^{\g,\a}_\e(Q)+\tau\right)}{\H_{\a}\left(\Lambda,\phi_0\right)(Q)+\e}}{W\left(\frac{(\g-\a)\left(\D^{\g,\a}_\e(Q)+\tau\right)}{\H_{\a}\left(\Lambda,\phi_0\right)(Q)+\e}\right)}\\
  &=&W\left(\frac{(\g-\a)\left(\D^{\g,\a}_\e(Q)+\tau\right)}{\H_{\a}\left(\Lambda,\phi_0\right)(Q)+\e}\right).
\end{eqnarray*}
Hence, applying the inverse of $W$ implies that
\[\H^{\g,\a}_\e\left(\Lambda,\phi_0\right)(Q)\log\frac{\H^{\g,\a}_\e\left(\Lambda,\phi_0\right)(Q)}{\H_{\a}\left(\Lambda,\phi_0\right)(Q)+\e}<(\g-\a)\left(\D^{\g,\a}_\e(Q)+\tau\right).\]
Thus, letting $\tau\to0$ proves the first inequality of (i). The second follows immediately from the first, in the case $\g=1$, by Lemma \ref{icp} (ii) and Lemma \ref{hmlb} (i), after letting $\e\to 0$.

$(ii)$ The condition on $\e$ implies that 
\[\H_{\a_0}\left(\Lambda,\phi_0\right)(Q)^{\frac{1-\a}{1-\a_0}}\Lambda(Q)^\frac{\a-\a_0}{1-\a_0}<\Lambda(Q) e-\e.\]  Hence, by Lemma \ref{shm} (i), $\H_{\a}\left(\Lambda,\phi_0\right)(Q)<\Lambda(Q) e-\e$, and therefore, the first inequality of (ii) follows from that of (i) in the case $\g=1$.

The second inequality of (ii) follows from the first, after letting $\e\to 0$.
\hfill$\Box$

 The following corollary can be used for obtaining criteria for the positivity of $\Phi$. 

\begin{cor}\label{lpc}
   Let $Q\in\B$ such that $\Lambda(Q)>0$. 
   
   (i) Suppose there exist $0<\e<e\Lambda(Q)$ and $\g\in[0,1)$ such that 
   \[\K_{\g,\e}(\Lambda|\phi_0)(Q)<\frac{\Lambda(Q)}{1-\g}\log\frac{\Lambda(Q)}{\e}.\]
   Then $\H_{\a}\left(\Lambda,\phi_0\right)(Q)>0$ for all $\a\in[0,\g]$.
   
   (ii) Suppose there exists a function $\tau:(0,1]\longrightarrow [0,\infty)$ which is continuous at $1$ such that $\tau(1)=0$, $\tau(\a)>0$ for all $\a\in(0,1)$ and
   \[\liminf\limits_{\a\to^-1}(1-\a)\K_{\a,\tau(\a)}(\Lambda|\phi_0)(Q)<\infty.\]
   Then $\H_{\a}\left(\Lambda,\phi_0\right)(Q)>0$ for all $\a\in[0,1]$.
\end{cor}
{\it Proof.}
  $(i)$ By the hypothesis, 
  \[\Lambda(Q)e^{-\frac{1-\g}{\Lambda(Q)}\K_{\g,\e}(\Lambda|\phi_0)(Q)}>\e.\]  
  Thus, the assertion follows by Theorem \ref{lbre} (i) and Lemma \ref{shm} (iii).

  $(ii)$ For all $\a\in(0,1)$ large enough,
  \[\tau(\a)<\Lambda(Q)\left(e-\left(\frac{\Phi(Q)}{\Lambda(Q)}\right)^{1-\a}\right).\]
  Therefore, by Theorem \ref{lbre} (ii), $\limsup_{\a\to^-1}\H_{\a}\left(\Lambda,\phi_0\right)(Q)>0$. Hence, by Lemma \ref{shm} (iii), $\H_{\a}\left(\Lambda,\phi_0\right)(Q)>0$ for all $\a\in[0,1)$.
\hfill$\Box$

\subsection{Upper bounds for the relative entropy measure}\label{subre}

Clearly, choosing a good and easy computable upper bound for $\K_\a(\Lambda|\phi_0)(Q)$ most likely depends on the particular application. 
However, there are some natural general upper bounds, which might suggest a direction in a particular case via some weakening or generalization.

\subsubsection{Restricting the set of covers via an invariant measure}\label{rims}

A natural way to obtain an upper bound on $\K_\a(\Lambda|\phi_0)(Q)$ is of course by a further restriction of the set of covers of $Q$ over which the infimum is taken. 

Since the main approach of this paper is a reduction of the proof of the positivity of $\Phi$ to the fact of the existence of $\Lambda$, via an estimation of an integral expression of $Z$, it suggests itself a further restriction of the set of covers via additional conditions in terms of $\Lambda$.

Recall, that, by Lemma \ref{icp} (i),
\begin{equation}\label{hmvi}
\H^{\a,1}\left(\Lambda,\phi_0\right)(Q) = \H_{\a}\left(\Lambda,\phi_0\right)(Q)
\end{equation}
for all $Q\in\P(X)$ and $\a\geq 0$ such that $\H_{\a}\left(\Lambda,\phi_0\right)(X)<\infty$, which suggests the following definition, via the inductive construction from Definition \ref{icfw}.

\begin{Definition}
	Let $0\leq\a< 1$, $Q\in\P(X)$ and $\e>0$. Define
	\[\K_{\a,\Lambda, \e}\left(\Lambda,\phi_0\right)(Q):=\inf\limits_{(A_m)_{m\leq 0}\in\C^{\a,1}_{\e}(Q)}\sum\limits_{m\leq 0}\int\limits_{S^mA_m}Z\log Z d\phi_0.\]
	Also, define $\K_{\a,\Lambda}\left(\Lambda,\phi_0\right)(Q)$, $\bar\K_{\a, \Lambda}\left(\Lambda,\phi_0\right)(Q)$, $\K_{\a, \Lambda}(Q)$ and $\bar\K_{\a, \Lambda}(Q)$ analogously to $\K_{\a}\left(\Lambda,\phi_0\right)(Q)$, $\bar\K_{\a}\left(\Lambda,\phi_0\right)(Q)$, $\K_{\a}(Q)$ and $\bar\K_{\a}(Q)$.
\end{Definition}

Then, since, by \eqref{hmvi}, $\C^{\a,1}_{\e}(Q)\subset\C^\a_{\e}(Q)$,
\[\K_{\a,\e}\left(\Lambda,\phi_0\right)(Q)\leq\K_{\a,\Lambda,\e}\left(\Lambda,\phi_0\right)(Q)\]
for all $Q\in\B$ and $\e>0$. 

However, by Lemma \ref{icp} (iii), such an additional condition on the covers does not change $\K_{\a}$ if it is finite. The next lemma deduces that for $\K_{\a}\left(\Lambda,\phi_0\right)$.

\begin{lemma}\label{urep} Let $\a\in[0,1)$ and $Q\in\B$.
	
(i)	\[\K_{\a,\Lambda}\left(\Lambda|\phi_0\right)(Q)=  \K_{\a,\Lambda}(Q) -\frac{1}{1-\a}\H_{\a}\left(\Lambda,\phi_0\right)(Q).\]

(ii) $\bar\K_{\a,\Lambda}\left(\Lambda|\phi_0\right)$ is a $S$-invariant, signed measure on $\B$.	
	
(iii) If $\K_{\a}\left(\Lambda,\phi_0\right)(X)<\infty$, then 
\[\K_{\a, \Lambda}(\Lambda|\phi_0)(Q)=\K_{\a}\left(\Lambda,\phi_0\right)(Q).\]
\end{lemma}
{\it Proof.} $(i)$ Let $\e>0$ and $(A_m)_{m\leq 0}\in\C^{\a,1}_\e(Q)$. Then 
\begin{eqnarray*}
&&\K_{\a,\Lambda,\e}(\Lambda|\phi_0)(Q) +\frac{1}{1-\a}\H_{\a}\left(\Lambda,\phi_0\right)(Q)\\
&\leq&\sum\limits_{m\leq 0}\int\limits_{S^mA_m}Z\log Z d\phi_0 + \frac{1}{1-\a}\sum\limits_{m\leq 0}\int\limits_{S^mA_m} Z^\a d\phi_0\\
&=&\sum\limits_{m\leq 0}\kappa_{\a}\left(S^mA_m\right)\\
&\leq&\sum\limits_{m\leq 0}\int\limits_{S^mA_m}Z\log Z d\phi_0+\frac{1}{1-\a}\left(\H_{\a}\left(\Lambda,\phi_0\right)(Q)+\e\right).
\end{eqnarray*}
Thus, taking the infimum and letting $\e\to0$ implies (i).

$(ii)$ The proof of (ii) is the same as that of Lemma \ref{eml} (iv).

$(iii)$ By Lemma \ref{eml} (ii), the assumption implies that $\K_{\a}(X)<\infty$. Hence, by Lemma \ref{icp} (iii), $\K_{\a,\Lambda}(Q)=\K_{\a}(Q)$. Thus, (iii) follows by (i) and Lemma \ref{eml} (ii).
\hfill$\Box$

The additional condition on the covers allows us to obtain a slightly more elegant version of Theorem \ref{lbre}, which is also much easier to prove.  (By Lemma 17 (Lemma 9 in the arXiv version) in \cite{Wer15}, for every $\e>0$, $\K_{\a,\Lambda,\e}(\Lambda|\phi_0)(Q)$ is also finite for a broad class of topological dynamical systems if $K(\Lambda|\phi_0)$ is finite and $Q$ is compact.)

For $0\leq\a< 1$, $\e>0$ and $Q\in\B$, define $\lambda_{\a,\e}(Q):=\Lambda(Q)$ if $\K_{\a,\Lambda,\e}(\Lambda|\phi_0)(Q)>0$ and $\lambda_{\a,\e}(Q):=\Lambda(Q)+\e$ otherwise. Obviously, $\Lambda(Q)\leq\lambda_{\a,\e}(Q)\leq\Lambda(Q)+\e$.

\begin{theo}\label{imt}
	Let $Q\in\B$ such that $\Lambda(Q)>0$ and $0\leq\a< 1$.
	
	(i)
	\begin{equation*}
	\H_{\a}\left(\Lambda,\phi_0\right)(Q)\geq\Lambda(Q)e^{-\frac{1-\a}{\lambda_{\a,\e}(Q)}\K_{\a,\Lambda,\e}(\Lambda|\phi_0)(Q)}-\e \ \ \ \mbox{ for all }\e>0,\mbox{ and}
	\end{equation*} 
	\[\Phi(Q)\geq\Lambda(Q)e^{-\frac{1}{\Lambda(Q)}\K_{\a,\Lambda}(\Lambda|\phi_0)(Q)}.\]
	
	(ii) If $\K_{\a}(\Lambda|\phi_0)(X)<\infty$ and $\B$ is generated by a sequence of finite partitions, then
	\[\Phi(X)\geq e^{K\left(\Lambda|\hat\Phi\right) -\K_\a(\Lambda|\phi_0)(X)}\mbox{ where }\hat\Phi := \frac{\Phi}{\Phi(X)}\]
	(hence, $K(\Lambda|\hat\Phi)\leq \K_\a(\Lambda|\phi_0)(X)$ if $\phi_0$ is a probability measure).
\end{theo}
{\it Proof.}
$(i)$ Let $\e>0$. Clearly, we can assume that $\K_{\a,\Lambda,\e}(\Lambda|\phi_0)(Q)<\infty$.   Let $\tau>0$ such that $\K_{\a,\Lambda,\e}(\Lambda|\phi_0)(Q)+\tau$ has the same sign as $\K_{\a,\Lambda,\e}(\Lambda|\phi_0)(Q)$ (we assign to zero '+'). Let $(A_m)_{m\leq 0}\in\C^{\a,1}_{\e}(Q)$ such that
\begin{equation*}
\K_{\a,\Lambda,\e}(\Lambda|\phi_0)(Q)+\tau>\sum\limits_{m\leq 0}\int\limits_{S^mA_m}Z\log Zd\phi_0.
\end{equation*}
 Then, as in the proof of Theorem \ref{lbre} (i), by \eqref{hmvi},
\begin{eqnarray*}
	\H_{\a}\left(\Lambda,\phi_0\right)(Q)+\e&>&\sum\limits_{m\leq 0}\int\limits_{S^mA_m}Z^{\a} d\phi_0\\
	&\geq&\sum\limits_{m\leq 0}\Lambda(A_m)e^{-\frac{1-\a}{\sum_{m\leq 0}\Lambda(A_m)}\left(\K_{\a,\Lambda, \e}(\Lambda|\phi_0)(Q)+\tau\right)}\\
	&\geq&\Lambda(Q)e^{-\frac{1-\a}{\lambda_{\a,\e}(Q)}\left(\K_{\a,\Lambda, \e}(\Lambda|\phi_0)(Q)+\tau\right)}.
\end{eqnarray*}
Thus, letting $\tau\to0$ implies the first inequality of (i).

The second inequality of (i) follows from the first by Lemma \ref{hmlb} (i) after letting $\e\to 0$.

$(ii)$ By second inequality of (i), Lemma \ref{urep} (iii) and Lemma \ref{eml} (v), 
\[\sum\limits^n_{k=1}\Lambda\left(Q_k\right)\log\frac{\Lambda(Q_k)}{\hat\Phi(Q_k)} -\log\Phi(X)\leq \K_\a(\Lambda|\phi_0)(X)\]
for every  $\B$-measurable partition  $(Q_k)_{1\leq k\leq n}$ of $X$.
Using the well-know fact that the sum in the inequality converges to $K(\Lambda|\hat\Phi)$ if one has a sequence of partitions which is increasing with respect to the refinement and generates the $\sigma$-algebra (e.g. Theorem 4.1 in  \cite{Sl}), it follows that
\[K\left(\Lambda|\hat\Phi\right) - \K_\a(\Lambda|\phi_0)(X) \leq\log\Phi(X),\]
which proves (ii).
\hfill$\Box$

\subsubsection{Taking supremum along trajectories}\label{satss}

Obviously, the finiteness of $K(\Lambda|\phi_0)$ implies only that $\Lambda\{Z>n\}\to 0$ as $n\to\infty$. The next corollary shows that the latter does not imply in general that $\Lambda\ll\Phi$. Therefore, by Theorem \ref{imt}, $K(\Lambda|\phi_0)$  is not an upper bound for $\K_\a(\Lambda|\phi_0)(X)$ in general.

A straightforward way to obtain an upper bound on $\K_\a(\Lambda|\phi_0)(X)$, which appears also to be quite practical (see \cite{Wer13}, where it was introduced and used), is the following.
\begin{Definition}
	Define 
	\[ Z^*:=\sup\limits_{m\leq 0}Z\circ S^m\ \ \ \mbox{ and}\]
	\[K^*(\Lambda|\phi_0):=\int  \log  Z^*d\Lambda.\]
\end{Definition}

Since $\int  \log^-  Z^*d\Lambda\leq\int\log^-Zd\Lambda=\int Z \log^-  Zd\phi_0<\infty$, $\int  \log  Z^*d\Lambda$ is well defined.  Obviously, $K(\Lambda|\phi_0)\leq K^*(\Lambda|\phi_0)$, and $K(\Lambda|\phi_0)= K^*(\Lambda|\phi_0)$ if $\phi_0\circ S^{-1}=\phi_0$. 

\begin{lemma}\label{reup} 
	\[\K_\a(\Lambda|\phi_0)(X)\leq K^*(\Lambda|\phi_0)\ \ \ \mbox{ for all }0\leq\a< 1.\]
\end{lemma}
{\it Proof.}
Let $0\leq\a<1$ and $\e>0$.  Let $(B_m)_{m\leq 0}\in\dot\C^\a_{\e}(X)$. Then, by Lemma 10 (ii) (Lemma 6 (ii) in the arXiv version) in \cite{Wer15},
\begin{eqnarray*}
	\K_{\a,\e}(\Lambda|\phi_0)(X)&\leq&\inf\limits_{(A_m)_{m\leq 0}\in\dot\C^\a_{\e}(X)}\sum\limits_{m\leq 0}\int\limits_{S^mB_m}Z\log Zd\phi_0\\
	&\leq&\sum\limits_{m\leq 0}\int\limits_{B_m}\log Z\circ S^md\Lambda\\
	&\leq&\int\log Z^*d\Lambda.
\end{eqnarray*}
Thus, the assertion follows.
\hfill$\Box$

Though $K^*(\Lambda|\phi_0)$ appears to be a very rough upper bound for $\K_\a(\Lambda|\phi_0)(X)$, the next corollary shows that it is quite adequate in some important cases.

\begin{cor}\label{fcem}
	Suppose $\Lambda$ is an ergodic probability measure. Let $0\leq\a< 1$. Then the following are equivalent.
	
	(i) $\Lambda\ll \H_{\a}\left(\Lambda,\phi_0\right)$ on $\B$.
	
	(ii) $Z$ is essentially bounded with respect to $\Lambda$.
	
	(iii) $K^*(\Lambda|\phi_0)<\infty$.
	
	(iv) $\K_{\a}(\Lambda|\phi_0)(X)<\infty$.
\end{cor}
{\it Proof.}
$(i)\Rightarrow (ii)$: Suppose (ii) is not true. Then $\Lambda\{Z>n\}>0$ for all $n\in\N$.  For each $n\in\N$ and $m\in\Z\setminus\N$, define $B^n_m:=S^{-m}\{Z>n\}$. By the ergodicity of $\Lambda$, $\Lambda\left(\bigcup_{m\leq 0}B^n_m\right) = 1$ for all $n\in\N$. Set $B:=\bigcap_{n\in\N}\bigcup_{m\leq 0}B^n_m$. Then
\begin{equation}\label{oifm}
\Lambda(B)=1.
\end{equation}
Set $A^n_0:=B^n_0$ and $A^n_m:=B^n_m\setminus(B^n_{m+1}\cup...\cup B^n_0)$ for all $m\leq -1$ and $n\in\N$. Then,  for each $n\in\N$, $A^n_m$'s are pairwise disjoint, each $A^n_m\in\A_m$ and $\bigcup_{m\leq 0}A^n_m=\bigcup_{m\leq 0}B^n_m$. Therefore,
\begin{eqnarray}
	1&=&\Lambda\left(\bigcup\limits_{m\leq 0}A^n_m\right)=\sum\limits_{m\leq 0}\Lambda\left(S^mA^n_m\right)=\sum\limits_{m\leq 0}\int\limits_{S^mA^n_m}Zd\phi_0\nonumber\\
	&\geq& n^{1-\a}\sum\limits_{m\leq 0}\int\limits_{S^mA^n_m}Z^\a d\phi_0
	\geq n^{1-\a}\H_{\a}\left(\Lambda,\phi_0\right)\left(B\right)\label{ecae}
\end{eqnarray}
for all $n\in\N$. Hence, $\H_{\a}\left(\Lambda,\phi_0\right)(B)=0$, which together with \eqref{oifm} contradicts to (i).

$(ii)\Rightarrow (iii)$ is obvious.

$(iii)\Rightarrow (iv)$ by Lemma \ref{reup}.

$(iv)\Rightarrow (i)$ follows by Theorem \ref{imt} (i), Lemma \ref{eml} (ii) and the fact that $\bar\K_{\a}$ is a measure on $\B$.
\hfill$\Box$

The following corollary covers, in particular, Example \ref{e1}.
\begin{cor}\label{1hme}
	 Suppose $\Lambda$ is an ergodic Borel probability measure such that $\phi_0\ll\Lambda$ (in addition to $\Lambda\ll\phi_0$). Then the following are equivalent.
	
	(i) There exists $0\leq\a<1$ such that $\K_\a(\Lambda|\phi_0)(X)<\infty$.
	
	(ii) For every $0\leq\g\leq 1$,  $\H_{\g}\left(\Lambda,\phi_0\right)(X)>0$ and\\ $\H_{\g}\left(\Lambda,\phi_0\right)(Q)/\H_{\g}\left(\Lambda,\phi_0\right)(X) = \Lambda(Q)$ for all $Q\in\B$.
	
	(iii) There exists $0\leq\a<1$ such that  $\H_{\a}\left(\Lambda,\phi_0\right)(X)>0$ and\\ $\H_{\a}\left(\Lambda,\phi_0\right)(Q)/\H_{\a}\left(\Lambda,\phi_0\right)(X) = \Lambda(Q)$ for all $Q\in\B$.
\end{cor}
{\it Proof.}
$(i)\Rightarrow (ii)$: Let $0\leq\g<1$. By Corollary \ref{fcem}, $\K_\g(\Lambda|\phi_0)(X)<\infty$. Hence, by Theorem \ref{imt} (i) and Lemma \ref{eml} (ii), $\H_{\g}\left(\Lambda,\phi_0\right)(X)>0$. By Lemma 19 (Lemma 10 in the arXiv version) in \cite{Wer15}, $\H_{\g}\left(\Lambda,\phi_0\right)\ll\Lambda$. Since $\H_{\g}\left(\Lambda,\phi_0\right)/\H_{\g}\left(\Lambda,\phi_0\right)(X)$ is a $S$-invariant probability measure on $\B$, the ergodicity of $\Lambda$ implies that $\H_{\g}\left(\Lambda,\phi_0\right)/\H_{\g}\left(\Lambda,\phi_0\right)(X)=\Lambda$ on $\B$. 

$(ii)\Rightarrow (iii)$ is obvious.

$(iii)\Rightarrow (i)$: It follows by (i) $\Rightarrow$ (iv) of Corollary \ref{fcem}.
\hfill$\Box$

\subsection{The regularity of $\a\longmapsto\H_{\a}\left(\Lambda,\phi_0\right)$ and $\a\longmapsto\H^{\a,0}\left(\Lambda,\phi_0\right)$ }

Now, we turn our attention to the regularity of the dependence of $\H_{\a}\left(\Lambda,\phi_0\right)$ and $\H^{\a,0}\left(\Lambda,\phi_0\right)$ on $\a$, which is another way to obtain conditions for their positivity.

\subsubsection{An almost convexity of $\a\longmapsto\H_{\a}\left(\Lambda,\phi_0\right)$ and $\a\longmapsto\H^{\a,0}\left(\Lambda,\phi_0\right)$ }\label{acs}

A natural approach to obtain some regularity properties of the functions $\a\longmapsto\H_{\a}\left(\Lambda,\phi_0\right)$ and $\a\longmapsto\H^{\a,0}\left(\Lambda,\phi_0\right)$ is to try to deduce them from the convexity of $\a\longmapsto Z^\a$. 

This requires another DDM arising from the Hellinger integral via the inductive construction from Definition \ref{icfw}, which also generalizes $\H_{\a}\left(\Lambda,\phi_0\right)$ and provides lower bounds for $\Phi$.

Recall that, by Lemma \ref{shm} (i), for every $Q\in\P(X)$, $\H^{\g,0}\left(\Lambda,\phi_0\right)(Q)<\infty$ for all $0\leq\g\leq 1$, and, by Definition \ref{icfw}, for every $\a,\g\geq 0$, the set function $\H^{\a,\g,0}\left(\Lambda,\phi_0\right)(Q)$ is well defined if $\H^{\g,0}\left(\Lambda,\phi_0\right)(Q)<\infty$.

Obviously, $\H_{\a}\left(\Lambda,\phi_0\right)(Q)\leq\H^{\a,0}\left(\Lambda,\phi_0\right)(Q)\leq\H^{\a,\g,0}\left(\Lambda,\phi_0\right)(Q)$ for all $\a\geq 0$. Also, one easily sees that $\H^{0,\g,0}\left(\Lambda,\phi_0\right)(Q)=\Phi(Q)$, $\H^{\g,\g,0}\left(\Lambda,\phi_0\right)(Q)=\H^{\g,0}\left(\Lambda,\phi_0\right)(Q)$ and, by Lemma \ref{icp} (iv), $\H^{1,\g,0}\left(\Lambda,\phi_0\right)(Q)=\Lambda(Q)$.

The new set functions allow us to formulate the following properties of $\H^{\a,0}(\Lambda,\phi_0)$.

\begin{lemma}\label{pchm} 
	Let $Q\in\P(X)$.	Let $\tilde\H_{\a}(\Lambda,\phi_0)$ and $\tilde\H^{\b,\a}(\Lambda,\phi_0)$ denote either $\H^{\a,0}(\Lambda,\phi_0)$ and $\H^{\b,\a,0}(\Lambda,\phi_0)$ or $\H_{\a}(\Lambda,\phi_0)$ and $\H^{\b,\a}(\Lambda,\phi_0)$.
		
	(i) Let $0\leq\b\leq\a_0<\a\leq\g$. If $\tilde\H_{\a}(\Lambda,\phi_0)(Q)<\infty$, then
	\[\tilde\H^{\a_0,\a}_\e(\Lambda,\phi_0)(Q)\leq{\tilde\H^{\b,\a}_\e\left(\Lambda,\phi_0\right)(Q)}^{1-\frac{\a_0-\b}{\a-\b}}\left(\tilde\H_{\a}(\Lambda,\phi_0)(Q)+\e\right)^{\frac{\a_0-\b}{\a-\b}}\]
	for all $\e>0$. If $\tilde\H_{\a_0}(\Lambda,\phi_0)(Q)<\infty$, then
	\[\tilde\H^{\a,\a_0}_\e(\Lambda,\phi_0)(Q)\leq\left(\tilde\H_{\a_0}(\Lambda,\phi_0)(Q)+\e\right)^{1-\frac{\a-\a_0}{\g-\a_0}}{\tilde\H^{\g,\a_0}_\e\left(\Lambda,\phi_0\right)(Q)}^{\frac{\a-\a_0}{\g-\a_0}}\]
	for all $\e>0$. In particular, for every $0\leq\a_0<\a\leq 1$,
	\begin{equation}\label{gpdl}
	\H^{\a_0,\a,0}(\Lambda,\phi_0)(Q)\leq{\Phi(Q)}^{1-\frac{\a_0}{\a}}{\H^{\a,0}(\Lambda,\phi_0)(Q)}^{\frac{\a_0}{\a}}\mbox{, and}
	\end{equation}
	\begin{equation}\label{hpdl}
	\H^{\a,\a_0,0}(\Lambda,\phi_0)(Q)\leq{\H^{\a_0,0}(\Lambda,\phi_0)(Q)}^{1-\frac{\a-\a_0}{1-\a_0}}{\Lambda(Q)}^{\frac{\a-\a_0}{1-\a_0}}.
	\end{equation}
	
	(ii) $\H^{\a,\b,0}\left(\Lambda,\phi_0\right)(Q)\leq\Phi(Q)^{1-\a}\Lambda(Q)^\a$ for all $\a,\b\in[0,1]$. 
	
	(iii) For every $\a,\b\in[0,1]$, $\H^{\a,\b,0}\left(\Lambda,\phi_0\right)$ is a finite, $S$-invariant measure on $\B$.
	
	(iv) Suppose there exists $0<\tau<1$ such that $\H^{\tau,0}(\Lambda,\phi_0)(Q)>0$. Then $\H^{\a,0}(\Lambda,\phi_0)(Q)>0$ for all $\a\in[0,1]$.
	
	(v) Let $0\leq\b<\a_0<\a<\g\leq 1$. Suppose $\H^{\b,\a}\left(\Lambda,\phi_0\right)(Q)<\infty$. Then
	\begin{eqnarray*}
		&&\max\left\{\frac{\tilde\H_{\a_0}\left(\Lambda,\phi_0\right)(Q)}{\a_0-\b}\log\frac{\tilde\H_{\a_0}\left(\Lambda,\phi_0\right)(Q)}{\tilde\H^{\b,\a}\left(\Lambda,\phi_0\right)(Q)},\phantom{\left(\frac{\tilde\H_{\a}(\Lambda,\phi_0)(Q)}{\tilde\H^{\b,\a}\left(\Lambda,\phi_0\right)(Q)}\right)^\frac{\a_0-\b}{\a-\b}}\right.\\  
		&&\ \ \ \ \ \ \ \left.\frac{\tilde\H^{\b,\a}\left(\Lambda,\phi_0\right)(Q)^{1-\frac{\a_0-\b}{\a-\b}}\tilde\H_{\a}(\Lambda,\phi_0)(Q)^\frac{\a_0-\b}{\a-\b}}{\a-\b}\log\frac{\tilde\H_{\a}(\Lambda,\phi_0)(Q)}{\tilde\H^{\b,\a}\left(\Lambda,\phi_0\right)(Q)}\right\}\\
		&\leq&\frac{\tilde\H_{\a}\left(\Lambda,\phi_0\right)(Q)-\tilde\H_{\a_0}\left(\Lambda,\phi_0\right)(Q)}{\a-\a_0}\\
		&\leq&\min\left\{\frac{\tilde\H_{\a}\left(\Lambda,\phi_0\right)(Q)}{\g-\a_0}\log\frac{\tilde\H^{\g,\a_0}\left(\Lambda,\phi_0\right)(Q)}{\tilde\H_{\a_0}\left(\Lambda,\phi_0\right)(Q)},\phantom{\frac{\tilde\H_{\a_0}\left(\Lambda,\phi_0\right)(Q)}{\g-\a_0}\left(\frac{\tilde\H^{\g,\a_0}(\Lambda,\phi_0)(Q)}{\tilde\H_{\a_0}\left(\Lambda,\phi_0\right)(Q)}\right)^{\frac{\a-\a_0}{\g-\a_0}}\log\frac{\tilde\H^{\g,\a_0}(\Lambda,\phi_0)(Q)}{\tilde\H_{\a_0}\left(\Lambda,\phi_0\right)(Q)}}\right.\\
		&&\ \ \ \ \ \ \ \frac{\tilde\H_{\a}\left(\Lambda,\phi_0\right)(Q)}{\g-\a}\log\frac{\tilde\H^{\g,\a_0}\left(\Lambda,\phi_0\right)(Q)}{\tilde\H_{\a}\left(\Lambda,\phi_0\right)(Q)},\\
		&&\ \ \ \ \ \ \ \left.\frac{\tilde\H_{\a_0}\left(\Lambda,\phi_0\right)(Q)^{1-\frac{\a-\a_0}{\g-\a_0}}\tilde\H^{\g,\a_0}(\Lambda,\phi_0)(Q)^\frac{\a-\a_0}{\g-\a_0}}{\g-\a_0}\log\frac{\tilde\H^{\g,\a_0}(\Lambda,\phi_0)(Q)}{\tilde\H_{\a_0}\left(\Lambda,\phi_0\right)(Q)}\right\}.
	\end{eqnarray*}
	
	(vi) Let $0\leq\b<\a_0<\a<\g\leq 1$.  Then
	\begin{eqnarray*}
		&&\max\left\{\frac{\tilde\H_{\a_0}\left(\Lambda,\phi_0\right)(Q)-\tilde\H^{\b,\a}\left(\Lambda,\phi_0\right)(Q)}{\a_0-\b},\right.\\
		&&\ \ \ \ \ \ \ \left.\frac{\tilde\H_{\a}\left(\Lambda,\phi_0\right)(Q)-\tilde\H^{\b,\a}\left(\Lambda,\phi_0\right)(Q)}{\a-\b}\right\}\\
		&\leq&\frac{\tilde\H_{\a}\left(\Lambda,\phi_0\right)(Q)-\tilde\H_{\a_0}\left(\Lambda,\phi_0\right)(Q)}{\a-\a_0}\\
		&\leq&\min\left\{\frac{\tilde\H_{\a}\left(\Lambda,\phi_0\right)(Q)}{\tilde\H_{\a_0}\left(\Lambda,\phi_0\right)(Q)}\frac{\tilde\H^{\g,\a_0}\left(\Lambda,\phi_0\right)(Q)-\tilde\H_{\a_0}\left(\Lambda,\phi_0\right)(Q)}{\g-\a_0}\right.\\
		&&\ \ \ \ \ \ \ \mbox{ if }\tilde\H_{\a_0}\left(\Lambda,\phi_0\right)(Q)>0,\\
		&&\ \ \ \ \ \ \ \frac{\tilde\H^{\g,\a_0}\left(\Lambda,\phi_0\right)(Q)-\tilde\H_{\a}\left(\Lambda,\phi_0\right)(Q)}{\g-\a},\\
		&&\ \ \ \ \ \ \ \left.\frac{\tilde\H^{\g,\a_0}\left(\Lambda,\phi_0\right)(Q)-\tilde\H_{\a_0}\left(\Lambda,\phi_0\right)(Q)}{\g-\a_0}\right\}.
	\end{eqnarray*}
	
	(vii) Let $0\leq\a_0<\a\leq 1$. Then
	\begin{eqnarray*}
	   &&(\a-\a_0)\frac{\tilde\H_{\a}\left(\Lambda,\phi_0\right)(Q)-\tilde\H^{0,\a}\left(\Lambda,\phi_0\right)(Q)}{\a}\leq\tilde\H_{\a}\left(\Lambda,\phi_0\right)(Q)-\tilde\H_{\a_0}\left(\Lambda,\phi_0\right)(Q)\\
	   &&\leq(\a-\a_0)\frac{\Lambda(Q)-\tilde\H_{\a_0}\left(\Lambda,\phi_0\right)(Q)}{1-\a_0}
	\end{eqnarray*}
\end{lemma}
{\it Proof.}
 We will prove the statements involving $\tilde\H$ for $\H_{\a}(\Lambda,\phi_0)$ and $\H^{\b,\a}(\Lambda,\phi_0)$. The proofs of those with $\H^{\a,0}(\Lambda,\phi_0)$ and $\H^{\b,\a,0}(\Lambda,\phi_0)$ are analogous.

$(i)$ Suppose $\H_{\a}(\Lambda,\phi_0)(Q)<\infty$. Clearly, for a proof of the first inequality, we can assume that $\a_0>0$. Let us abbreviate 
\[\tau:=\frac{\a_0-\b}{\a-\b}.\]
Obviously, $0\leq\tau< 1$.   Let $\e>0$ and  $(A_m)_{m\leq 0}\in\C^\a_\e(Q)$. Let us first consider the case $\sum_{m\leq 0}\int_{S^mA_m}Z^{\a} d\phi_0>0$. Restricting $Z$ on the set $\{Z>0\}$ if necessary, the concavity of $[0,\infty)\owns x\longmapsto x^{1-\tau}$ implies that
\begin{eqnarray}
\sum\limits_{m\leq 0}\int\limits_{S^mA_m}Z^{\a_0} d\phi_0&=&\sum\limits_{m\leq 0}\int\limits_{S^mA_m}\left(Z^{\b-\a}\right)^{1-\tau} Z^{\a}d\phi_0\nonumber\\
&\leq&\sum\limits_{m\leq 0}\left(\int\limits_{S^mA_m}Z^{\b} d\phi_0\right)^{1-\tau}\left(\int\limits_{S^mA_m}Z^{\a} d\phi_0\right)^\tau\nonumber\\
&\leq&\left(\sum\limits_{m\leq 0}\int\limits_{S^mA_m}Z^{\b} d\phi_0\right)^{1-\tau}\left(\sum\limits_{m\leq 0}\int\limits_{S^mA_m}Z^{\a} d\phi_0\right)^\tau.\ \label{hce}
\end{eqnarray}
Hence,
\[\H^{\a_0,\a}_\e(\Lambda,\phi_0)(Q)\leq\left(\sum\limits_{m\leq 0}\int\limits_{S^mA_m}Z^{\b} d\phi_0\right)^{1-\tau}\left(\H_{\a}(\Lambda,\phi_0)(Q)+\e\right)^\tau.\]
If $\sum_{m\leq 0}\int_{S^mA_m}Z^{\a} d\phi_0=0$, then $\sum_{m\leq 0}\int_{S^mA_m}Z^{\a_0} d\phi_0=0$, and the last inequality is obviously correct also. It implies the first inequality of (i).

Now, for a proof of the second inequality, suppose $\H_{\a_0}(\Lambda,\phi_0)(Q)<\infty$. Let $(B_m)_{m\leq 0}\in\C^{\a_0}_\e(Q)$. Set
\[\eta:=\frac{\a-\a_0}{\g-\a_0}.\]
Then $0<\eta\leq 1$. Suppose $\sum_{m\leq 0}\int_{S^mA_m}Z^{\a_0} d\phi_0>0$. Then, by the concavity of $[0,\infty)\owns x\longmapsto x^\eta$, 
\begin{eqnarray*}
	\sum\limits_{m\leq 0}\int\limits_{S^mB_m}Z^{\a} d\phi_0&=&\sum\limits_{m\leq 0}\int\limits_{S^mB_m}\left(Z^{\g-\a_0}\right)^\eta Z^{\a_0}d\phi_0\\
	&\leq&\left(\sum\limits_{m\leq 0}\int\limits_{S^mB_m}Z^{\a_0} d\phi_0\right)^{1-\eta}\left(\sum\limits_{m\leq 0}\int\limits_{S^mB_m}Z^{\g} d\phi_0\right)^\eta\\
	&\leq&\left(\H_{\a_0}(\Lambda,\phi_0)(Q)+\e\right)^{1-\eta}\left(\sum\limits_{m\leq 0}\int\limits_{S^mB_m}Z^{\g} d\phi_0\right)^\eta.
\end{eqnarray*}
If $\sum_{m\leq 0}\int_{S^mA_m}Z^{\a_0} d\phi_0=0$, then $\sum_{m\leq 0}\int_{S^mA_m}Z^{\a} d\phi_0=0$, and we still have
\[\H^{\a,\a_0}_\e(\Lambda,\phi_0)(Q)\leq\left(\H_{\a_0}(\Lambda,\phi_0)(Q)+\e\right)^{1-\eta}\left(\sum\limits_{m\leq 0}\int\limits_{S^mB_m}Z^{\g} d\phi_0\right)^\eta.
\]
It implies the second inequality of (i).

The third and fourth inequalities follow from the first and the second (by Lemma \ref{icp} (iv)) respectively.

$(ii)$ It follows from \eqref{gpdl}, \eqref{hpdl}, Lemma \ref{shm} (i), and Lemma \ref{hmlb} (i).

$(iii)$ It follows immediately by (ii) and Theorem 16 (ii) (Theorem 4 (ii) in the arXiv version) in \cite{Wer15}.

$(iv)$ Let $\a\in[0,1]$. If $\tau<\a\leq 1$, then $\H^{\a,0}(\Lambda,\phi_0)(Q)>0$ by \eqref{gpdl}. If $0\leq\a<\tau$, then $\H^{\a,0}(\Lambda,\phi_0)(Q)>0$ by \eqref{hpdl}.

$(v)$ By (iv) and Lemma \ref{shm} (iii), the first inequality in (v) is obviously correct if $\tilde\H_{\b}(\Lambda,\phi_0)(Q)=0$ (as then $\tilde\H_{\a_0}(\Lambda,\phi_0)(Q)=0$ and $\tilde\H_{\a}(\Lambda,\phi_0)(Q)=0$). Suppose $\H_{\b}(\Lambda,\phi_0)(Q)>0$. By Lemma \ref{edl} (i), $Z^a\leq Z -(1-a)Z^a\log Z$, $0\leq a\leq 1$, which is equivalent to $Y^{1/a}\geq Y+(1/a-1)Y\log Y$. Applying the former to the first inequality of (i) (after letting $\e\to 0$) implies
\begin{eqnarray*}
	&&\H_{\a_0}(\Lambda,\phi_0)(Q)\leq\H_{\a}(\Lambda,\phi_0)(Q) - \left(1-\frac{\a_0-\b}{\a-\b}\right)\H^{\b,\a}\left(\Lambda,\phi_0\right)(Q)\\
	&&\ \ \ \ \ \ \ \ \ \ \ \ \ \ \ \ \ \ \ \ \ \ \ \ \ \ \ \   \times\left(\frac{\H_{\a}(\Lambda,\phi_0)(Q)}{\H^{\b,\a}\left(\Lambda,\phi_0\right)(Q)}\right)^\frac{\a_0-\b}{\a-\b}\log\frac{\H_{\a}(\Lambda,\phi_0)(Q)}{\H^{\b,\a}\left(\Lambda,\phi_0\right)(Q)}.
\end{eqnarray*}
Applying the latter to (the equivalent)
\[\H_{\a}(\Lambda,\phi_0)(Q)\geq\H^{\b,\a}\left(\Lambda,\phi_0\right)(Q)\left(\frac{\H_{\a_0}(\Lambda,\phi_0)(Q)}{\H^{\b,\a}\left(\Lambda,\phi_0\right)(Q)}\right)^\frac{\a-\b}{\a_0-\b}\]
implies that
\begin{eqnarray*}
	&&\H_{\a}(\Lambda,\phi_0)(Q)\\	&\geq&\H_{\a_0}\left(\Lambda,\phi_0\right)(Q)+\left(\frac{\a-\b}{\a_0-\b}-1\right)\H_{\a_0}\left(\Lambda,\phi_0\right)(Q)\log\frac{\H_{\a_0}\left(\Lambda,\phi_0\right)(Q)}{\H^{\b,\a}\left(\Lambda,\phi_0\right)(Q)}.
\end{eqnarray*}
This proves the first inequality in (v).

The second inequality in (v) is obviously correct (by (iv) and Lemma \ref{shm} (iii)) if $\tilde\H_{\a_0}(\Lambda,\phi_0)(Q)=0$. Suppose $\tilde\H_{\a_0}(\Lambda,\phi_0)(Q)>0$. By the second inequality of (i),
\begin{equation}\label{lchm}
\frac{\H_{\a}(\Lambda,\phi_0)(Q)}{\H_{\a_0}(\Lambda,\phi_0)(Q)}\leq\left(\frac{\H^{\g,\a_0}(\Lambda,\phi_0)(Q)}{\H_{\a_0}(\Lambda,\phi_0)(Q)}\right)^\frac{\a-\a_0}{\g-\a_0}.
\end{equation}
Hence,
\begin{eqnarray*}
	\frac{\H_{\a}(\Lambda,\phi_0)(Q)-\H_{\a_0}(\Lambda,\phi_0)(Q)}{\H_{\a}(\Lambda,\phi_0)(Q)}&\leq&\log\frac{\H_{\a}(\Lambda,\phi_0)(Q)}{\H_{\a_0}(\Lambda,\phi_0)(Q)}\\
	&\leq&\frac{\a-\a_0}{\g-\a_0}\log\frac{\H^{\g,\a_0}(\Lambda,\phi_0)(Q)}{\H_{\a_0}(\Lambda,\phi_0)(Q)},
\end{eqnarray*}
which implies the first part of the second inequality of (v). 

Inequality \eqref{lchm} implies also that
\[\left(\frac{\H_{\a}(\Lambda,\phi_0)(Q)}{\H_{\a_0}(\Lambda,\phi_0)(Q)}\right)^\frac{\g-\a}{\g-\a_0}\leq\left(\frac{\H^{\g,\a_0}(\Lambda,\phi_0)(Q)}{\H_{\a}(\Lambda,\phi_0)(Q)}\right)^\frac{\a-\a_0}{\g-\a_0}.\]
The linearization of the left side of the logarithmic version of it, as above, gives the second part of the second inequality of (v).

Applying $Z^a\leq 1+aZ^a\log Z$ to \eqref{lchm} gives the third part of the second inequality of (v).

$(vi)$ Clearly, for the proof of the first inequality of (vi), we can assume that $\H^{\b,\a}\left(\Lambda,\phi_0\right)(Q)<\infty$. Then the first part of it follows immediately from that of (v), since $x\log x\geq x-1$ for all $x\geq 0$. 

Let $(A_m)_{m\leq 0}\in\C^{\a}_\e(Q)$. Then, by $(Z^{\a_0}-Z^\b)/(\a_0-\b)\leq(Z^{\a}-Z^{\a_0})/(\a-\a_0)$ (which follows from the convexity of $x\longmapsto Z^x$ for $x>0$),
\begin{eqnarray*}
	&&\frac{\H_{\a}(\Lambda,\phi_0)(Q)-\sum_{m\leq 0}\int_{S^mA_m}Z^{\b} d\phi_0}{\a-\b}\\
	&\leq&\frac{\sum_{m\leq 0}\int_{S^mA_m}Z^{\a} d\phi_0-\sum_{m\leq 0}\int_{S^mA_m}Z^{\b} d\phi_0}{\a-\b}\\
	&\leq&\frac{\sum_{m\leq 0}\int_{S^mA_m}Z^{\a} d\phi_0-\sum_{m\leq 0}\int_{S^mA_m}Z^{\a_0} d\phi_0}{\a-\a_0}\\
	&\leq&\frac{\H_{\a}(\Lambda,\phi_0)(Q)+\e-\H_{\a_0}(\Lambda,\phi_0)(Q)}{\a-\a_0},
\end{eqnarray*}
which implies the second part of the first inequality of (vi).

For the proof of the second inequality of (vi), by (iv) and Lemma \ref{shm} (iii), we can assume that $\tilde\H_{\a_0}(\Lambda,\phi_0)(Q)>0$. Then the first and the second parts of it follow immediately from those of (v), as $\log x\leq x-1$ for all $x>0$. 

The third part of the second inequality of (vi) follow from the the convexity of $x\longmapsto Z^x$ similarly to the proof of the second part of the first inequality. 

$(vii)$ The assertion for $0<\a_0<\a< 1$ follows immediately from (vi), by setting $\b=0$ and $\g=1$. 

Let $\a_0=0$. Since $\Phi(Q)\leq\H^{0,\a}\left(\Lambda,\phi_0\right)(Q)$, it follow the fist inequality of (vii). By (i), 
\begin{equation}\label{hmla}
  \H_{\a}\left(\Lambda,\phi_0\right)(Q)\leq\Phi(Q)^{1-\a}\Lambda(Q)^\a\leq(1-\a)\Phi(Q)+\a\Lambda(Q),
\end{equation}
which implies the second inequality of (vii). 

Finally, let $\a=1$. Then the first inequality follows form \eqref{hmla} written for $\a_0$ and Lemma \ref{icp} (i), and the second is an equality.
\hfill$\Box$

As a by-product of Lemma \ref{pchm} (i), we obtain also the following method for computation of lower bounds for $\H_{\a_0}(\Lambda,\phi_0)$.

\begin{prop}\label{clbh} Let $Q\in\P(X)$.
	
	(i) Let $0\leq\b\leq 1<\a$. Suppose $\H_{\a}(\Lambda,\phi_0)(Q)<\infty$ and $\H^{\b,\a}\left(\Lambda,\phi_0\right)(Q)<\infty$. Then
	\[\Lambda(Q)\leq{\H^{\b,\a}\left(\Lambda,\phi_0\right)(Q)}^{1-\frac{1-\b}{\a-\b}}{\H_{\a}(\Lambda,\phi_0)(Q)}^{\frac{1-\b}{\a-\b}}.\]
	
	(ii) Let $0\leq\a_0<1\leq\g$. Suppose $\H^{\g,\a_0}\left(\Lambda,\phi_0\right)(Q)<\infty$. Then
	\[\Lambda(Q)\leq{\H_{\a_0}(\Lambda,\phi_0)(Q)}^{1-\frac{1-\a_0}{\g-\a_0}}{\H^{\g,\a_0}\left(\Lambda,\phi_0\right)(Q)}^{\frac{1-\a_0}{\g-\a_0}}.\]
\end{prop}
{\it Proof.}
    $(i)$ It follows from the first inequality of Lemma \ref{pchm} (i) and Lemma \ref{icp} (ii), by setting $\a_0=1$.
    
    $(ii)$ It follows from the second inequality of Lemma \ref{pchm} (i) and Lemma \ref{icp} (ii), by setting $\a=1$.
\hfill$\Box$

\subsubsection{The continuity of $(0,1)\owns\a\longmapsto\H_{\a}\left(\Lambda,\phi_0\right)$}\label{chm}

	 Obviously, Lemma \ref{pchm} would also imply some continuity properties of $\a\longmapsto\H_{\a}\left(\Lambda,\phi_0\right)(Q)$ if we knew that $\H^{\b,\a}\left(\Lambda,\phi_0\right)(Q)<\infty$ for some $0\leq\b<\a< 1$. This can happen. For example, suppose the exists $c>0$ such that $Z\geq c$ $\Lambda$-a.e. (as in Example \ref{e1}). Let $\e>0$ and $(A_m)_{m\leq 0}\in\C^\a_\e(Q)$. Then 
	 \begin{eqnarray*}
	 	&&\H_{\a}(\Lambda,\phi_0)(Q)+\e>\sum\limits_{m\leq 0}\int\limits_{S^mA_m}Z^\a d\phi_0=\sum\limits_{m\leq 0}\int\limits_{S^mA_m}Z^{\a-\b}Z^{\b-1} d\Lambda\\
	 	&\geq& c^{\a-\b}\sum\limits_{m\leq 0}\int\limits_{S^mA_m}Z^\b d\phi_0\geq
	 	c^{\a-\b}\H^{\b,\a}_\e\left(\Lambda,\phi_0\right)(Q).
	 \end{eqnarray*}
	 Hence,
	 \[\H^{\b,\a}\left(\Lambda,\phi_0\right)(Q)\leq\frac{\H_{\a}(\Lambda,\phi_0)(Q)}{c^{\a-\b}}\]
	 for all $0\leq\b<\a< 1$. Therefore, in this case, by Lemma \ref{pchm} (vii), the function $(0,1)\owns\a\longmapsto\H_{\a}\left(\Lambda,\phi_0\right)(Q)$ is continuous (see Theorem \ref{homa} for more cases). This already clarifies the behavior of the function in the case of Example \ref{e1}. 
	 
	 Now, we are going to investigate conditions for the continuity of the function more closely. 
	 
	 First, observe that, for every $0<\b<\a\leq 1$ and $(A_m)_{m\leq 0}\in\C(Q)$ such that $\sum_{m\leq 0}\int_{S^mA_m}Z^\a d\phi_0<\infty$, by Lemma \ref{edl} (i),
	 \begin{eqnarray}\label{hdub}
	 	&&(\a-\b)\sum\limits_{m\leq 0, \int_{S^mA_m}Z^\b\log Z d\phi_0\geq 0}\int\limits_{S^mA_m}Z^\b\log Z d\phi_0\\
	 	&\leq&\sum\limits_{m\leq 0, \int_{S^mA_m}Z^\b\log Z d\phi_0\geq 0}\int\limits_{S^mA_m}\left(Z^\a-Z^\b\right) d\phi_0\leq\sum\limits_{m\leq 0}\int\limits_{S^mA_m}Z^\a d\phi_0<\infty.\nonumber
	 \end{eqnarray}
	 Hence, the sum $\sum_{m\leq 0}\int_{S^mA_m}Z^\b\log Z d\phi_0$ is well defined for all $(A_m)_{m\leq 0}\in\C^\a_{\e}(Q)$ and $\e>0$. Therefore, we can make the following definition. 
	 \begin{Definition}\label{hccd}
	 	Let $0<\b<\a\leq 1$. For $Q\in\P(X)$ and $\e>0$, define 
	 	\[\E^{\b,\a}_{\e}(Q):=\sup\limits_{(A_m)_{m\leq 0}\in\C^\a_{\e}(Q)}\sum\limits_{m\leq 0}\int\limits_{S^mA_m}Z^\b\log Zd\phi_0 \ \ \ \mbox{ and}\]
	 	\[\E^{\b,\a}(Q):=\lim\limits_{\e\to 0}\E^{\b,\a}_{\e}(Q),\]
	 	as, obviously, $\E^{\b,\a}_{\e}(Q)\geq\E^{\b,\a}_{\d}(Q)$ for all $0<\d\leq\e$.
	 \end{Definition}
	 
	 Obviously, by \eqref{hdub}, $\E^{\b,\a}(Q)<\infty$ for all $0<\b<\a\leq 1$ and $Q\in\P(Q)$. (Also, since, by \eqref{hdub},
	 $-\sum_{m\leq 0}\int_{S^mA_m}Z^\b\log Z d\phi_0+1/(\a-\b)\sum_{m\leq 0}\int_{S^mA_m}Z^\a d\phi_0\geq 0$ for all $(A_m)_{m\leq 0}\in\C^\a_{\e}(Q)$ and $\e>0$, one can show, similarly to Lemma \ref{eml} (iv), that $\lim_{i\to\infty}\E^{\b,\a}(S^{-i}.)$ is a signed measure on $\B$, but we will not need it.)
	 
	 The following lemma lists some criteria for the finiteness of  $\H^{\b,\a}\left(\Lambda,\phi_0\right)(Q)$ for $0<\b<\a<1$ via the finiteness from below of $\E^{\b,\a}(Q)$.
	 
	 It seems to be natural to make the following definitions, in order to obtain computable criteria for the latter.
	 
	 \begin{Definition}
	 	For $0\leq\a\leq 1$, $Q\in\P(X)$ and $\e>0$, define 
	 	\[\L_{\a,\e}\left(\Lambda|\phi_0\right)(Q):=\inf\limits_{(A_m)_{m\leq 0}\in\C^{\a}_{\e}(Q)}\sum\limits_{m\leq 0}\int\limits_{S^mA_m\cap\{Z<1\}}\frac{1}{Z}\log\frac{1}{Z}d\Lambda,\]
	 	\[\L_{\a}(\Lambda|\phi_0)(Q):=\lim\limits_{\e\to 0}\L_{\a,\e}(\Lambda|\phi_0)(Q),\]
	 	\[\U_{\a,\e}\left(\Lambda|\phi_0\right)(Q):=\sup\limits_{(A_m)_{m\leq 0}\in\C^{\a}_{\e}(Q)}\sum\limits_{m\leq 0}\int\limits_{S^mA_m\cap\{Z<1\}}\frac{1}{Z}\log\frac{1}{Z}d\Lambda\ \ \ \mbox{ and}\]
	 	\[\U_{\a}(\Lambda|\phi_0)(Q):=\lim\limits_{\e\to 0}\U_{\a,\e}(\Lambda|\phi_0)(Q).\]
	 \end{Definition}
	 
	 Obviously, $\L_{\a}\left(\Lambda|\phi_0\right)(Q)\leq\U_{\a}\left(\Lambda|\phi_0\right)(Q)$.

	 \begin{lemma}\label{fhcl} Let  $Q\in\B$.
	 	
	 	(i) For $0<\b<\a\leq 1$,
	 	\[\H^{\b,\a}\left(\Lambda,\phi_0\right)(Q)\leq\H_{\a}\left(\Lambda,\phi_0\right)(Q)-(\a-\b)\E^{\b,\a}(Q).\]
	 	
	 	(ii) $\E^{\b,\a}(Q)\leq\E^{\g,\a}(Q)$ for all $0<\b\leq\g<\a\leq 1$.
	 	
	 	(iii) For $0<\b<\a\leq 1$ and $\e>0$,
	 	\[\E^{\b,\a}(Q)\geq-\left(\frac{\H_{\a}\left(\Lambda,\phi_0\right)(Q)+\e}{1-\a}\right)^\b {\L_{\a,\e}\left(\Lambda|\phi_0\right)(Q)}^{1-\b}.\] 	 
	 	
	 	(iv) If there exists $0<c<1$ such that $Z\geq c$ $\Lambda$-a.e., then, for $0\leq\a\leq 1$,
	 	\[\U_{\a}\left(\Lambda|\phi_0\right)(Q)\leq \frac{\H_{\a}\left(\Lambda,\phi_0\right)(Q)}{c^\a}\log\frac{1}{c}.\]	
	 \end{lemma}
	 {\it Proof.} Let $0\leq\a\leq 1$, $\e>0$ and $(A_m)_{m\leq 0}\in\C^{\a}_{\e}(Q)$.
	 
	 $(i)$ Clearly, we can assume $\E^{\b,\a}(Q)>-\infty$. By Lemma \ref{edl} (i),
	 \begin{eqnarray*}
	   \H_{\a}\left(\Lambda,\phi_0\right)(Q)+\e&>&\sum\limits_{m\leq 0}\int\limits_{S^mA_m}Z^\a d\phi_0\\
	   &\geq&\sum\limits_{m\leq 0}\int\limits_{S^mA_m}Z^\b d\phi_0+(\a-\b)\sum\limits_{m\leq 0}\int\limits_{S^mA_m}Z^\b\log Zd\phi_0\\
	   &\geq&\H^{\b,\a}_\e\left(\Lambda,\phi_0\right)(Q)+(\a-\b)\sum\limits_{m\leq 0}\int\limits_{S^mA_m}Z^\b\log Zd\phi_0,
	 \end{eqnarray*}
	 which implies (i).
	 
	 $(ii)$ It follows immediately from Lemma \ref{edl} (ii).
	 
	 $(iii)$ First, observe that, by Lemma \ref{edl} (i), 
	 \begin{eqnarray*}
	 	&&\sum\limits_{m\leq 0}\int\limits_{S^mA_m\cap\{Z<1\}}\log\frac{1}{Z}d\Lambda\\
	 	&=&-\sum\limits_{m\leq 0}\int\limits_{S^mA_m\cap\{Z<1\}}Z\log Zd\phi_0\\
	 	&\leq&-\frac{1}{1-\a}\left(\sum\limits_{m\leq 0}\int\limits_{S^mA_m\cap\{Z<1\}}Zd\phi_0-\sum\limits_{m\leq 0}\int\limits_{S^mA_m\cap\{Z<1\}}Z^\a d\phi_0\right)\\
	 	&\leq&\frac{1}{1-\a}\left(\H_{\a}\left(\Lambda,\phi_0\right)(Q)+\e\right).
	 \end{eqnarray*}
	 Therefore, by the concavity of $x\longmapsto x^{1-\b}$,
	 \begin{eqnarray}\label{hdlb}
	    &&\E^{\b,\a}_{\e}(Q)\nonumber\\
	    &\geq&\sum\limits_{m\leq 0}\int\limits_{S^mA_m}Z^\b\log Zd\phi_0\nonumber\\
	    &\geq&-\sum\limits_{m\leq 0}\int\limits_{S^mA_m\cap\{Z<1\}}\left(\frac{1}{Z}\right)^{1-\b}\log\frac{1}{Z}d\Lambda\nonumber\\
	    &\geq&-\left(\sum\limits_{m\leq 0}\int\limits_{S^mA_m\cap\{Z<1\}}\log\frac{1}{Z}d\Lambda\right)^\b\left(\sum\limits_{m\leq 0}\int\limits_{S^mA_m\cap\{Z<1\}}\frac{1}{Z}\log\frac{1}{Z}d\Lambda\right)^{1-\b}\nonumber\\
	    &\geq&-\left(\frac{\H_{\a}\left(\Lambda,\phi_0\right)(Q)+\e}{1-\a}\right)^\b\left(\sum\limits_{m\leq 0}\int\limits_{S^mA_m\cap\{Z<1\}}\frac{1}{Z}\log\frac{1}{Z}d\Lambda\right)^{1-\b},
	 \end{eqnarray}
	 which implies (iii).
	 
	 $(iv)$ Observe that
	 \begin{eqnarray*}
	    \sum\limits_{m\leq 0}\int\limits_{S^mA_m\cap\{Z<1\}}\frac{1}{Z}\log\frac{1}{Z}d\Lambda&=&\sum\limits_{m\leq 0}\int\limits_{S^mA_m\cap\{c\leq Z<1\}}\frac{1}{Z}\log\frac{1}{Z}d\Lambda\\
	    &\leq&\log\frac{1}{c}\sum\limits_{m\leq 0}\int\limits_{S^mA_m\cap\{c\leq Z<1\}}\frac{1}{Z^\a}Z^\a d\phi_0\\
	    &\leq&\frac{1}{c^\a}\log\frac{1}{c}\left(\H_{\a}\left(\Lambda,\phi_0\right)(Q)+\e\right),
	 \end{eqnarray*}
	 which implies the assertion.
	 \hfill$\Box$
	 
	 Now, we are able to shed some light on the continuity of the function $(0,1)\owns\a\longmapsto\H_{\a}\left(\Lambda,\phi_0\right)(Q)$ by means of $\E^{\b,\a}(Q)$.
	 \begin{prop}\label{hcp} Let $0<\b<\a\leq 1$ and $Q\in\B$.
	 	
	 (i)	\[(\a-\b)\E^{\b,\a}(Q)\leq\H_{\a}\left(\Lambda,\phi_0\right)(Q)-\H_{\b}\left(\Lambda,\phi_0\right)(Q)\leq(\a-\b)\frac{\Lambda(Q)-\H_{\b}\left(\Lambda,\phi_0\right)(Q)}{1-\b}.\]
	 (ii) Let $0<\a<1$. If there exists $0<\b<\a$ such that $\E^{\b,\a}(Q)>-\infty$, then $(0,1)\owns x\longmapsto\H_{x}\left(\Lambda,\phi_0\right)(Q)$ is continuous at $\a$ from the left.
	 \end{prop}
	 {\it Proof.} $(i)$ It follows by Lemma \ref{pchm} (vii) and Lemma \ref{fhcl} (i). 
	 
	 $(ii)$ It follows immediately from (i), since $\E^{\b,\a}(Q)\leq\E^{\g,\a}(Q)$ for all $\b\leq\g$.
	 \hfill$\Box$

	 However, there are no problems with the continuity of $\a\longmapsto\H^{\a,0}\left(\Lambda,\phi_0\right)$ (compare also Lemma \ref{shm} (iii) and Lemma \ref{pchm} (iv)) (Lemma \ref{pchm} (vii) shows that the function $(0,1)\owns\a\longmapsto\H^{\a,0}\left(\Lambda,\phi_0\right)(Q)$ is continuous for all $Q\in\B$). This suggests that the functions are different in general. In such a case, it follows immediately that $\H^{\a,0}\left(\Lambda,\phi_0\right)(Q)>0 $ for all $\a\in[0,1]$.

\subsubsection{The continuity of $[0,1]\owns\a\longmapsto\H_{\a}\left(\Lambda,\phi_0\right)$ and $[0,1]\owns\a\longmapsto\H^{\a,0}\left(\Lambda,\phi_0\right)$ at $0$ and $1$ }\label{chm1}

Obviously, the continuity of the function $[0,1]\owns\a\longmapsto\H_{\a}\left(\Lambda,\phi_0\right)(Q)$ at $1$ implies, by Lemma \ref{shm} (iii), that it is strictly positive if $\Lambda(Q)>0$. The same argument can be also applied to the function $[0,1]\owns\a\longmapsto\H^{\a,0}\left(\Lambda,\phi_0\right)$, by Lemma \ref{pchm} (iv). 

Now, we give a sufficient condition for the continuity at $1$ for the functions which follows from Lemma \ref{pchm} (vii) and Theorem \ref{lbre}.  In particular, it immediately clarifies the behavior of the functions at the point $1$ in an essentially bounded case, as e.g. in Example \ref{e1}.

\begin{prop}\label{ca1}
	Let $Q\in\B$. Suppose $\K_{\a,\tau}(\Lambda|\phi_0)(Q)$ is finite for all $\tau>0$ and there exists a function $\tau:(0,1]\longrightarrow [0,\infty)$ which is continuous at $1$ such that $\tau(1)=0$, $\tau(\a)>0$ for all $\a<1$, and
	\[\lim\limits_{\a\to^-1}(1-\a)\K_{\a,\tau(\a)}(\Lambda|\phi_0)(Q)=0.\]
	Then the functions $[0,1]\owns\alpha\longmapsto \H_{\a}\left(\Lambda,\phi_0\right)(Q)$ and $[0,1]\owns\alpha\longmapsto \H^{\a,0}\left(\Lambda,\phi_0\right)(Q)$ are continuous at $1$.
\end{prop}
{\it Proof.} 
By Lemma \ref{pchm} (vii), for $\a\in(0,1)$,
\begin{equation*}
-(1-\a)\Phi(Q)\leq\Lambda(Q)-\H^{\a,0}\left(\Lambda,\phi_0\right)(Q)\leq\Lambda(Q)-\H_{\a}\left(\Lambda,\phi_0\right)(Q).
\end{equation*}

If $\Lambda(Q)=0$, then the continuity holds true by Lemma \ref{shm} (iii). Otherwise, for $\a\in(0,1)$ large enough,
\[\tau(\a)<\Lambda(Q)\left(e-\left(\frac{\Phi(Q)}{\Lambda(Q)}\right)^{1-\a}\right).\]
Therefore, by Theorem \ref{lbre} (ii) and the inequality $e^x\geq x+1$,
\[\Lambda(Q)-\H_{\a}\left(\Lambda,\phi_0\right)(Q)\leq(1-\a)\K_{\a,\tau(\a)}(\Lambda|\phi_0)(Q)+\tau(\a)\]
for all such $\a\in[0,1)$. Thus, the assertion follows.
\hfill$\Box$

Now, we turn to the continuity at zero. In particular, it will enable us to compute $\Phi$ in some cases (see Corollary \ref{pebc}). 
The following lemma shows also that $\Phi(X)>0$ if one of the functions is discontinuous at $0$. 

Note that, in our case, the relation $\Lambda\ll\phi_0$ is equivalent to $Z>0$ $\phi_0$-a.e.
\begin{prop}\label{ca0}
	Let $Q\in\B$. The functions $[0,1]\owns\alpha\longmapsto \H_{\a}\left(\Lambda,\phi_0\right)(Q)$ and $[0,1]\owns\alpha\longmapsto \H^{\a,0}\left(\Lambda,\phi_0\right)(Q)$ are continuous at $0$ if at least one of the following is true:
	
	(a) $\Phi(Q)=0$, or 
	
	(b) $Z>0$ $\phi_0$-a.e., and there exists a function $\e:[0,1)\longrightarrow [0,\infty)$ which is continuous at $0$ such that $\e(0)=0$, $\e(\a)>0$ for all $\a>0$, and 
	     \[\lim_{\a\to 0}\a\L_{\a,\e(\a)}(\Lambda|\phi_0)(Q)=0.\] 
	
\end{prop}
{\it Proof.} 
 (a)  Let $0<\a<1$. By Lemma \ref{shm} (i),
  \[\H^{\a,0}\left(\Lambda,\phi_0\right)(Q)\leq(1-\a)\Phi(Q)+\a\Lambda(Q).\]
  Hence,
  \begin{equation}\label{zcfb}
     \H_{\a}\left(\Lambda,\phi_0\right)(Q)-\Phi(Q)\leq\H^{\a,0}\left(\Lambda,\phi_0\right)(Q)-\Phi(Q)\leq\a\Lambda(Q).
  \end{equation}
  This implies, in particular the continuity of the functions at $0$ if $\Phi(Q)=0$.
  
  (b) Now, suppose $\Phi(Q)>0$. Clearly, we can assume that $\L_{\a,\e(\a)}(\Lambda|\phi_0)(Q)<\infty$ for sufficiently small $\a$. Observe that the integral $\int_A\log Zd\phi_0$ is well-defined for all $A\in\A_0$, as
  $\int_{A\cap\{Z\geq 	1\}}\log Zd\phi_0=-\int_{A\cap\{Z\geq 1\}}1/Z\log(1/Z)d\Lambda<\Lambda(A)/e$. Let $(A_m)_{m\leq 0}\in\C^{\a}_{\e(\a)}(Q)$ such that 
  \[\sum_{m\leq 0}\int_{\{Z<1\}\cap S^mA_m}1/Z\log(1/Z)d\Lambda<\L_{\a,\e(\a)}(\Lambda|\phi_0)(Q)+\e(\a).\] 
  Then, by proceeding via finite sums and then taking the limit, 
  \begin{eqnarray*}
     &&\H_{\a}\left(\Lambda,\phi_0\right)(Q)+\e(\a) \\
     &>&\sum\limits_{m\leq 0}\int\limits_{S^mA_m}e^{\a\log Z}d\phi_0\\
     &\geq &\sum\limits_{m\leq 0}\phi_0\left(S^mA_m\right)e^{\frac{\a}{\sum_{m\leq 0}\phi_0(S^mA_m)}\sum\limits_{m\leq 0}\int\limits_{\{Z<1\}\cap S^mA_m}\log Zd\phi_0}\\
     &\geq &\Phi(Q)e^{-\frac{\a}{\Phi(Q)}\left(\L_{\a,\e(\a)}(\Lambda|\phi_0)(Q)+\e(\a)\right)}.
  \end{eqnarray*}
  Hence, using $e^x\geq x+1$, it follows that
  \begin{eqnarray*}
     \H^{\a,0}\left(\Lambda,\phi_0\right)(Q)-\Phi(Q)&\geq&\H_{\a}\left(\Lambda,\phi_0\right)(Q)-\Phi(Q)\\
     &\geq&-\a\L_{\a,\e(\a)}(\Lambda|\phi_0)(Q)-(1+\a)\e(\a),
  \end{eqnarray*}
  which, together with \eqref{zcfb}, implies the assertion.
\hfill$\Box$

\subsubsection{The right differentiability of $(0,1)\owns\a\longmapsto\H_{\a}\left(\Lambda,\phi_0\right)$}

Clearly, the function cannot be zero everywhere if it is not differentiable at some point.

In this subsection, we will give a sufficient condition for the right differentiability of $(0,1)\owns\a\longmapsto\H_{\a}\left(\Lambda,\phi_0\right)(Q)$ for all $Q\in\B$.

By \eqref{hrds}, we can make the following definition.
\begin{Definition} Let $0<\a\leq 1$, $0\leq\b\leq 1$, $Q\in\P(X)$ and $\e>0$. Define
    \[\D^{\a,\b}_{1,\e}(Q):=\inf\limits_{(A_m)_{m\leq 0}\in\C^{\b,1}_{\e}(Q)}\sum\limits_{m\leq 0}\int\limits_{S^mA_m}Z^\a\log Zd\phi_0\ \ \ \mbox{ and}\]
	\[\D^{\a,\b}_{1}(Q):=\lim\limits_{\e\to 0}\D^{\a,\b}_{1,\e}(Q).\]
\end{Definition}

Obviously, $\D^{1,\b}_{1}(Q)=\K_{\b,\Lambda}\left(\Lambda,\phi_0\right)(Q)$. The following lemma indicates that $\D^{\a,\a}_{1}(Q)$ might be a derivative of the function if it is greater than minus infinity.
\begin{lemma}\label{fdl} (i) Let $0<\a_0<\a\leq 1$ and $Q\in\B$. Let $\e_0,\e>0$.  Then
	\begin{eqnarray*}
		(\a-\a_0)\D^{\a_0,\a}_{1,\e_0}(Q)-\e_0&<&\H_{\a}\left(\Lambda,\phi_0\right)(Q)-\H_{\a_0}\left(\Lambda,\phi_0\right)(Q)\\
		&<&(\a-\a_0)\D^{\a,\a_0}_{1,\e}(Q)+\e.
	\end{eqnarray*} 
	(ii) Let $0<\b<1$, $0\leq\a< 1$ and $Q\in\B$. Then
	\[-\left(\frac{\H_{\a}\left(\Lambda,\phi_0\right)(Q)+\e}{1-\a}\right)^\b {\U_{\a}\left(\Lambda|\phi_0\right)(Q)}^{1-\b}\leq\D^{\b,\a}_{1,\e}(Q)\leq\frac{\Lambda(Q)}{e(1-\b)}\]
	for all $\e>0$.
\end{lemma}
{\it Proof.}
$(i)$ Let $(A_m)_{m\leq 0}\in\C^{\a,1}_{\e_0}(Q)$. Then, by Lemma \ref{edl} (i) and \eqref{hmvi},
\begin{eqnarray*}
	(\a-\a_0)\D^{\a_0,\a}_{1,\e_0}(Q)&\leq&(\a-\a_0)\sum\limits_{m\leq 0}\int\limits_{S^mA_m}Z^{\a_0}\log Zd\phi_0\\
	&\leq&\sum\limits_{m\leq 0}\int\limits_{S^mA_m}Z^\a d\phi_0-\sum\limits_{m\leq 0}\int\limits_{S^mA_m}Z^{\a_0} d\phi_0\\
	&<&\H_{\a}\left(\Lambda,\phi_0\right)(Q)+\e_0-\H_{\a_0}\left(\Lambda,\phi_0\right)(Q).
\end{eqnarray*} 
This implies the first inequality of (i).

Let $(B_m)_{m\leq 0}\in\C^{\a_0,1}_{\e}(Q)$. Then, by Lemma \ref{edl} (i) and \eqref{hmvi},
\begin{eqnarray*}
	\H_{\a}\left(\Lambda,\phi_0\right)(Q)-\H_{\a_0}\left(\Lambda,\phi_0\right)(Q)-\e_0&<&\sum\limits_{m\leq 0}\int\limits_{S^mB_m}Z^\a d\phi_0-\sum\limits_{m\leq 0}\int\limits_{S^mB_m}Z^{\a_0} d\phi_0\\
	&\leq&(\a-\a_0)\sum\limits_{m\leq 0}\int\limits_{S^mB_m}Z^\a\log Zd\phi_0.
\end{eqnarray*} 
This implies the second inequality (i).

$(ii)$ Let $\e>0$ and $(C_m)_{m\leq 0}\in\C^{\a,1}_{\e}(Q)$. Then,  as in \eqref{hdlb} (the restrictions for $\a$ and $\b$ in \eqref{hdlb} were determined only by Definition \ref{hccd}), by \eqref{hmvi},
\begin{eqnarray*}
&&\sum\limits_{m\leq 0}\int\limits_{S^mC_m}Z^\b\log Zd\phi_0\\
&\geq&-\left(\frac{\H_{\a}\left(\Lambda,\phi_0\right)(Q)+\e}{1-\a}\right)^\b\left(\sum\limits_{m\leq 0}\int\limits_{S^mC_m\cap\{Z<1\}}\frac{1}{Z}\log\frac{1}{Z}d\Lambda\right)^{1-\b}\\
&\geq&-\left(\frac{\H_{\a}\left(\Lambda,\phi_0\right)(Q)+\e}{1-\a}\right)^\b{\U_{\a}\left(\Lambda|\phi_0\right)(Q)}^{1-\b},
\end{eqnarray*}
which implies the first inequality of (ii). The second follows by Lemma \ref{hfpl}, since $\D^{\b,\a}_{1,\e}(Q)\leq\D^{\b,\a}_{1}(Q)$.
\hfill$\Box$

The following lemma gives a condition for the continuity of $\D^{\a,\b}_{1}(Q)$ with respect to the first parameter.

\begin{lemma}\label{hrdc}
	Let $0<\a_0\leq\a\leq 1$, $0\leq\b\leq 1$ and $Q\in\B$. Suppose there exists $0<\d<\a_0$ such that $\D^{\a_0-\d,\b}_{1}(Q)>-\infty$. Then
	\begin{eqnarray*}
		0&\leq&\D^{\a,\b}_{1}(Q)-\D^{\a_0,\b}_{1}(Q)\\
		&\leq&(\a-\a_0)\left[-\frac{1}{\d}\D^{\a_0-\d,\b}_{1}(Q)+\left(\frac{1}{\d e(1-\a_0+\d)}+\left(\frac{	2}{e(1-\a)}\right)^2\right)\Lambda(Q)\right].
	\end{eqnarray*}
\end{lemma}
{\it Proof.}
By Lemma \ref{edl} (ii), $\D^{\a_0-\d,\b}_{1}(Q)\leq\D^{\a_0,\b}_{1}(Q)\leq\D^{\a,\b}_{1}(Q)$, which implies the first inequality.

Let $\e>0$ and $(A_m)_{m\leq 0}\in\C^{\b,1}_{\e}(Q)$. Then, by Lemma \ref{edl} (ii) and Lemma \ref{hfpl},
\begin{eqnarray*}
	&&\frac{1}{\a-\a_0}\left(\sum\limits_{m\leq 0}\int\limits_{S^mA_m}Z^\a\log Zd\phi_0-\sum\limits_{m\leq 0}\int\limits_{S^mA_m}Z^{\a_0}\log Zd\phi_0\right)\\
	&\leq&\sum\limits_{m\leq 0}\int\limits_{S^mA_m\cap\{Z<1\}}Z^{\a_0}\left(\log Z\right)^2d\phi_0+\left(\frac{2}{e(1-\a)}\right)^2\left(\Lambda(Q)+\e\right).
\end{eqnarray*}
Now, observe that, by Lemma \ref{hfpl},
\begin{eqnarray*}
	&&\sum\limits_{m\leq 0}\int\limits_{S^mA_m\cap\{Z<1\}}Z^{\a_0}\left(\log Z\right)^2d\phi_0\\
	&=&\frac{1}{\d}\sum\limits_{m\leq 0}\int\limits_{S^mA_m\cap\{Z<1\}}Z^{\a_0}\left(\log\frac{1}{Z}\right)\left(\log\frac{1}{Z^\d}\right) d\phi_0\\
	&\leq&-\frac{1}{\d}\sum\limits_{m\leq 0}\int\limits_{S^mA_m\cap\{Z<1\}}Z^{\a_0-\d}\log Zd\phi_0+\frac{1}{\d}\sum\limits_{m\leq 0}\int\limits_{S^mA_m\cap\{Z<1\}}Z^{\a_0}\log Zd\phi_0\\
	&\leq&-\frac{1}{\d}\sum\limits_{m\leq 0}\int\limits_{S^mA_m}Z^{\a_0-\d}\log Zd\phi_0\\
	&&+\frac{1}{\d}\sum\limits_{m\leq 0}\int\limits_{S^mA_m\cap\{Z\geq 1\}}e^{-(1-\a_0+\d)\log Z}\log Zd\Lambda\\
	&\leq&-\frac{1}{\d}\D^{\a_0-\d,\b}_{1,\e}(Q)+\frac{1}{\d e(1-\a_0+\d)}\left(\Lambda(Q)+\e\right).
\end{eqnarray*}
Therefore,
\begin{eqnarray}\label{dfpc}
	&&\frac{1}{\a-\a_0}\left(\D^{\a,\b}_{1,\e}(Q)-\sum\limits_{m\leq 0}\int\limits_{S^mA_m}Z^{\a_0}\log Zd\phi_0\right)\nonumber\\
	&\leq&\frac{1}{\a-\a_0}\left(\sum\limits_{m\leq 0}\int\limits_{S^mA_m}Z^\a\log Zd\phi_0-\sum\limits_{m\leq 0}\int\limits_{S^mA_m}Z^{\a_0}\log Zd\phi_0\right)\\
	&\leq&-\frac{1}{\d}\D^{\a_0-\d,\b}_{1,\e}(Q)+\frac{1}{\d e(1-\a_0+\d)}\left(\Lambda(Q)+\e\right)+\left(\frac{2}{e(1-\a)}\right)^2\left(\Lambda(Q)+\e\right).\nonumber
\end{eqnarray}
Thus, the second inequality follows.
\hfill$\Box$

Now, we are able to give a sufficient condition for the right differentiability of the function, which, by Lemma \ref{fdl} (ii) and Lemma \ref{fhcl} (iv), is satisfied in the case of Example \ref{e1}.
\begin{theo} \label{hrdt}
	Let $Q\in\B$ and $0<\a_0<1$. Suppose there exists $\d>0$ such that $\D^{\a_0-\d,\a_0}_1(Q)>-\infty$ and $\lim_{\e\downarrow 0}\e\D^{\a_0,\a_0+\e}_{1,\e}(Q)=0$. Then the function $(0,1)\owns x\longmapsto\H_{x}\left(\Lambda,\phi_0\right)(Q)$ is right differentiable at $\a_0$, and 
	\[\left.\frac{d_+}{d_+x}\H_x(\Lambda,\phi_0)(Q)\right|_{x=\a_0}=\D^{\a_0,\a_0}_{1}(Q)=\lim\limits_{x\to^+\a_0}\D^{\a_0,x}_{1}(Q)\]
	where $d_+/d_+x$ denotes the right derivative.
\end{theo}
{\it Proof.} 
    Let $\e>0$ such that $\D^{\a_0,\a}_{1,\a-\a_0}(Q)>-\infty$ for all $0<\a-\a_0\leq\e$. Let $\a:=\a_0+\e$ and $(A_m)_{m\leq 0}\in\C^{\a,1}_{\e}(Q)$. Then, by Lemma \ref{fdl} (i),
    \begin{eqnarray*}
       &&\H_{\a_0}\left(\Lambda,\phi_0\right)(Q) +(\a-\a_0)\D^{\a,\a_0}_{1}(Q)+\e\geq\H_{\a}\left(\Lambda,\phi_0\right)(Q)+\e\\
       &>&\sum\limits_{m\leq 0}\int\limits_{S^mA_m}Z^\a d\phi_0\geq\sum\limits_{m\leq 0}\int\limits_{S^mA_m}Z^{\a_0} d\phi_0+(\a-\a_0)\sum\limits_{m\leq 0}\int\limits_{S^mA_m}Z^{\a_0}\log Z d\phi_0\\
       &\geq&\sum\limits_{m\leq 0}\int\limits_{S^mA_m}Z^{\a_0} d\phi_0+(\a-\a_0)\D^{\a_0,\a}_{1,\e}(Q).
    \end{eqnarray*}
    Hence, since, by Lemma \ref{fdl} (i), $(\a-\a_0)(\D^{\a,\a_0}_{1}(Q)-\D^{\a_0,\a}_{1,\e}(Q))+\e\geq 0$, $(A_m)_{m\leq 0}\in\C^{\a_0,1}_{(\a-\a_0)(\D^{\a,\a_0}_{1}(Q)-\D^{\a_0,\a}_{1,\e}(Q))+2\e}(Q)$.
    That is 
    \[\C^{\a,1}_{\e}(Q)\subset\C^{\a_0,1}_{(\a-\a_0)\left(\D^{\a,\a_0}_{1}(Q)-\D^{\a_0,\a}_{1,\e}(Q)\right)+2\e}(Q).\]
    Therefore, for every $0<\b\leq 1$,
    \[\D^{\b,\a}_{1}(Q)\geq\D^{\b,\a}_{1,\e}(Q)\geq\D^{\b,\a_0}_{1,(\a-\a_0)\left(\D^{\a,\a_0}_{1}(Q)-\D^{\a_0,\a}_{1,\e}(Q)\right)+2\e}(Q).\]
    Hence, by Lemma \ref{fdl} (i) and Lemma \ref{hrdc},
    \begin{eqnarray*}
       &&\D^{\a_0,\a_0}_{1,(\a-\a_0)\left(\D^{\a,\a_0}_{1}(Q)-\D^{\a_0,\a}_{1,\a-\a_0}(Q)\right)+2(\a-\a_0)}(Q)\leq\D^{\a_0,\a}_{1}(Q)\\
       &\leq&\frac{\H_{\a}\left(\Lambda,\phi_0\right)(Q)-\H_{\a_0}\left(\Lambda,\phi_0\right)(Q)}{\a-\a_0}\leq\D^{\a,\a_0}_{1}(Q)
       \leq\D^{\a_0,\a_0}_{1}(Q)\\
       &&+(\a-\a_0)\left[-\frac{1}{\d}\D^{\a_0-\d,\a_0}_{1}(Q)+\left(\frac{1}{\d e(1-\a_0+\d)}+\left(\frac{	2}{e(1-\a)}\right)^2\right)\Lambda(Q)\right].
    \end{eqnarray*}
    Thus, the hypothesis implies the assertion.
\hfill$\Box$

\subsubsection{The left differentiability of $(0,1)\owns\a\longmapsto\H_{\a}\left(\Lambda,\phi_0\right)$}

Now, we give a sufficient condition for the left differentiability of $(0,1)\owns\a\longmapsto\H_{\a}\left(\Lambda,\phi_0\right)(Q)$ for every measurable $Q$.

 \begin{Definition}
 	Let $0<\a\leq 1$, $0\leq\b\leq 1$, $Q\in\P(X)$ and $\e>0$. Define 
 	\[\E^{\a,\b}_{1,\e}(Q):=\sup\limits_{(A_m)_{m\leq 0}\in\C^{\b,1}_{\e}(Q)}\sum\limits_{m\leq 0}\int\limits_{S^mA_m}Z^\a\log Zd\phi_0 \ \ \ \mbox{ and}\]
 	\[\E^{\a,\b}_1(Q):=\lim\limits_{\e\to 0}\E^{\a,\b}_{1,\e}(Q).\]
 \end{Definition}
 
 Clearly, $\E^{\a,\b}_{1}(Q)\leq\E^{\a,\b}(Q)$ for all $Q\in\B$. However, there still might be a problem with its finiteness from below for $\a\leq\b$.

Similarly to $\D^{\a,\b}_1(Q)$, the set function has the following continuity property with respect to the first parameter.

\begin{lemma}\label{hldc}
	Let $0<\a_0\leq\a\leq 1$, $0\leq\b\leq 1$ and $Q\in\B$.  Suppose there exists $0<\d<\a_0$ such that $\D^{\a_0-\d,\b}_{1}(Q)>-\infty$. Then
	\begin{eqnarray*}
	  0&\leq&\E^{\a,\b}_1(Q)-\E^{\a_0,\b}_1(Q)\\
	  &\leq&(\a-\a_0)\left[-\frac{1}{\d}\D^{\a_0-\d,\b}_1(Q)+\left(\frac{1}{\d e(1-\a_0+\d)}+\left(\frac{	2}{e(1-\a)}\right)^2\right)\Lambda(Q)\right].
	\end{eqnarray*}
\end{lemma}
{\it Proof.}
By the hypothesis and Lemma \ref{edl} (ii), $-\infty<\E^{\a_0,\b}_{1}(Q)\leq\E^{\a,\b}_{1}(Q)$, which implies the first inequality.

Let $\e>0$ and $(A_m)_{m\leq 0}\in\C^{\b,1}_{\e}(Q)$. Then, by \eqref{dfpc},
\begin{eqnarray*}
	&&\frac{1}{\a-\a_0}\left(\sum\limits_{m\leq 0}\int\limits_{S^mA_m}Z^\a\log Zd\phi_0-\E^{\a_0,\b}_{1,\e}(Q)\right)\\
	&\leq&-\frac{1}{\d}\D^{\a_0-\d,\b}_{1,\e}(Q)+\frac{1}{\d e(1-\a_0+\d)}\left(\Lambda(Q)+\e\right)+\left(\frac{2}{e(1-\a)}\right)^2\left(\Lambda(Q)+\e\right),
\end{eqnarray*}
which implies the second inequality.
\hfill$\Box$

Also, similarly to Lemma \ref{fdl} (i), we have the following.
\begin{lemma} \label{hldl}
  (i)	Let $0<\a_0<\a\leq 1$ and $Q\in\B$. Then
	\[(\a-\a_0)\E^{\a_0,\a}_1(Q)\leq\H_{\a}\left(\Lambda,\phi_0\right)(Q)-\H_{\a_0}\left(\Lambda,\phi_0\right)(Q)\leq(\a-\a_0)\E^{\a,\a_0}_1(Q).\]
	
  (ii)	Let $0<\a< 1$, $0\leq\b\leq 1$ and $Q\in\B$. Then
  \[\E^{\a,\b}_{1,\e}(Q)\leq\frac{\Lambda(Q)+\e}{e(1-\a)}\ \ \ \mbox{ for all }\e>0.\]
\end{lemma}
{\it Proof.} $(i)$ Let $\e>0$ and $(A_m)_{m\leq 0}\in\C^{\a,1}_{\e}(Q)$. Then, by Lemma \ref{edl} (i) and \eqref{hmvi},
\begin{eqnarray*}
  &&(\a-\a_0)\sum\limits_{m\leq 0}\int\limits_{S^mA_m}Z^{\a_0}\log Zd\phi_0\leq\sum\limits_{m\leq 0}\int\limits_{S^mA_m}Z^{\a}d\phi_0-\sum\limits_{m\leq 0}\int\limits_{S^mA_m}Z^{\a_0}d\phi_0\\
  &<&\H_{\a}\left(\Lambda,\phi_0\right)(Q)+\e-\H_{\a_0}\left(\Lambda,\phi_0\right)(Q),
\end{eqnarray*}
which implies the first inequality of (i).

 Now, let $(B_m)_{m\leq 0}\in\C^{\a_0,1}_{\e}(Q)$. Then, by Lemma \ref{edl} (i) and \eqref{hmvi},
 \begin{eqnarray*}
 	(\a-\a_0)\E^{\a,\a_0}_{1,\e}(Q)&\geq&(\a-\a_0)\sum\limits_{m\leq 0}\int\limits_{S^mB_m}Z^{\a}\log Zd\phi_0\\
 	&\geq&\sum\limits_{m\leq 0}\int\limits_{S^mB_m}Z^{\a}d\phi_0-\sum\limits_{m\leq 0}\int\limits_{S^mB_m}Z^{\a_0}d\phi_0\\
 	&>&\H_{\a}\left(\Lambda,\phi_0\right)(Q)-\H_{\a_0}\left(\Lambda,\phi_0\right)(Q)-\e,
 \end{eqnarray*}
 which implies the second inequality of (i).
 
 $(ii)$ It follows immediately by Lemma \ref{hfpl}.
\hfill$\Box$

\begin{theo} \label{hldt}
  Let $Q\in\B$ and $0<\a<1$. Suppose there exists $0<\a_0<\a$ such that $\D^{\a_0,\a}_1(Q)>-\infty$. Then the function $(0,1)\owns x\longmapsto\H_{x}\left(\Lambda,\phi_0\right)(Q)$ is left differentiable at $\a$, and 
  \[\left.\frac{d_-}{d_-x}\H_x(\Lambda,\phi_0)(Q)\right|_{x=\a}=\E^{\a,\a}_{1}(Q)=\lim\limits_{x\to^-\a}\E^{\a,x}_{1}(Q)\]
   where $d_-/d_-x$ denotes the left derivative.
\end{theo}
{\it Proof.} Let $\a_0<x<\a$ and $\d>0$ such that $\a_0<x-\d$. Then, by the hypothesis and Lemma \ref{edl} (ii), $\E^{x,\a}_{1}(Q)\geq\E^{x-\d,\a}_{1}(Q)\geq\D^{x-\d,\a}_{1}(Q)\geq\D^{\a_0,\a}_{1}(Q)>-\infty$. Let $\e>0$ and $(A_m)_{m\leq 0}\in\C^{x,1}_{\e}(Q)$. Then, by Lemma \ref{hldl}, \eqref{hmvi}, Lemma \ref{edl} (i) and Lemma \ref{hfpl},
\begin{eqnarray*}
  &&\H_{\a}\left(\Lambda,\phi_0\right)(Q)-(\a-x)\E^{x,\a}_{1}(Q)+\e\geq\H_{x}\left(\Lambda,\phi_0\right)(Q)+\e\\
  &>&\sum\limits_{m\leq 0}\int\limits_{S^mA_m}Z^{x}d\phi_0\geq\sum\limits_{m\leq 0}\int\limits_{S^mA_m}Z^{\a}d\phi_0-(\a-x)\sum\limits_{m\leq 0}\int\limits_{S^mA_m}Z^{\a}\log Zd\phi_0\\
  &\geq&\sum\limits_{m\leq 0}\int\limits_{S^mA_m}Z^{\a}d\phi_0-(\a-x)\E^{\a,x}_{1,\e}(Q).
\end{eqnarray*}
Hence, since, by Lemma \ref{hldl}, $\E^{x,\a}_1(Q)\leq\E^{\a,x}_1(Q)\leq\E^{\a,x}_{1,\e}(Q)\leq(\Lambda(Q)+\e)/(e(1-\a))$,
\[(A_m)_{m\leq 0}\in\C^{\a,1}_{(\a-x)\left((\Lambda(Q)+\e)/(e(1-\a))-\E^{x,\a}_{1}(Q)\right)+\e}(Q).\]
That is
\[\C^{x,1}_{\e}(Q)\subset\C^{\a,1}_{(\a-x)\left((\Lambda(Q)+\e)/(e(1-\a))-\E^{x,\a}_{1}(Q)\right)+\e}(Q).\]
Hence, for every $0<\b\leq 1$,
\[\E^{\b,x}_1(Q)\leq\E^{\b,x}_{1,\e}(Q)\leq\E^{\b,\a}_{1,(\a-x)\left((\Lambda(Q)+\e)/(e(1-\a))-\E^{x,\a}_{1}(Q)\right)+\e}(Q).\]
Therefore, by Lemma \ref{hldc} and Lemma \ref{hldl},
\begin{eqnarray*}
	&&\E^{\a,\a}_1(Q)\\
	&&-(\a-x)\left[-\frac{1}{\d}\D^{x-\d,\a}_1(Q)+\left(\frac{1}{\d e(1-x+\d)}+\left(\frac{2}{e(1-\a)}\right)^2\right)\Lambda(Q)\right]\\
	&\leq&\E^{x,\a}_1(Q)\\
	&\leq&\frac{\H_{\a}\left(\Lambda,\phi_0\right)(Q)-\H_{x}\left(\Lambda,\phi_0\right)(Q)}{\a-x}\\
	&\leq&\E^{\a,x}_1(Q)\\
    &\leq&\E^{\a,\a}_{1,(\a-x)\left((\Lambda(Q)+\e)/(e(1-\a))-\E^{x,\a}_{1}(Q)\right)+\e}(Q).
\end{eqnarray*}
Thus, setting $\e=\a-x$ and letting $x\to \a$ implies the assertion.
\hfill$\Box$

\begin{Remark}
  Observe that the assertion of Theorem \ref{hldt} remains true also for $\a=1$ if also there exists $C<\infty$ such that $\E^{1,x}_{1,\e}(Q)\leq C$ for all $x<1$ sufficiently close to $1$ and all sufficiently small $\e>0$, as in the case of Example \ref{e1}.
\end{Remark}

\subsubsection{The differentiability of $(0,1)\owns\a\longmapsto\H_{\a}\left(\Lambda,\phi_0\right)$}\label{hds}

In this subsection, we shed some light on the differentiability of the function if $Z$ is $\Lambda$-essentially bounded away from zero.

\begin{cor}\label{hbdc}
	Let $Q\in\B$. Suppose $Z$ is $\Lambda$-essentially bounded away from zero. Then the function $(0,1)\owns x\longmapsto\H_{x}\left(\Lambda,\phi_0\right)(Q)$ is left and right differentiable, and 
	\[\left.\frac{d}{dx}\H_x(\Lambda,\phi_0)(Q)\right|_{x=\a}=\D^{\a,\a}_{1}(Q)=\E^{\a,\a}_{1}(Q)\]	
	for all except at most countably many $\a\in(0,1)$. 
\end{cor}
{\it Proof.} By Lemma \ref{fhcl} (iv) and Lemma \ref{fdl} (ii), the hypotheses of Theorem \ref{hrdt} and Theorem \ref{hldt} are satisfied. Therefore, the function is right and left differentiable. Thus, the assertion follows by the well-known Beppo Levi Theorem (e.g. see \cite{To}, p. 143).
\hfill$\Box$

\subsubsection{Candidates for the derivatives of $(0,1)\owns\a\longmapsto\H^{\a,0}\left(\Lambda,\phi_0\right)$}\label{cdhm}

By Lemma \ref{pchm} (ii) and (vii), the function $(0,1)\owns\a\longmapsto\H^{\a,0}\left(\Lambda,\phi_0\right)$ appears to have better continuity properties. We are going now to investigate its differentiability properties. (Clearly, the function cannot be zero everywhere if it has some irregularity at some $\a\in(0,1)$.)

We will use the inductive construction from Subsection 4.1.2 in \cite{Wer15}, to obtain some measures on $\B$ as natural candidates for the derivatives of the function.

\begin{Definition}\label{hmdd}
	Let $0\leq\a\leq 1$, $Q\in\P(X)$, $\e>0$. Define $\C^\a_{0,\e}(Q):=\C^0_{\e}(Q)$ and $\Psi^{\a}_{0}(Q):= \H^{\a,0}\left(\Lambda,\phi_0\right)(Q)$. For $n\in\N$ and $0<\a\leq 1$, define recursively (with $(-\infty)^0:=1$) (it will be shown in the next lemma that each of the following set functions is finite)
	\[\C^\a_{n,\e}(Q):=\left\{(A_m)_{m\leq 0}\in\C^\a_{n-1,\e}(Q)|\ \bar\Psi^\a_{n-1}(Q)>\sum\limits_{m\leq 0}\int\limits_{S^mA_m}Z^\a\left(\log Z\right)^{n-1}d\phi_0-\e\right\},\]
	\[\Psi^\a_{n,\e}(Q):=\inf\limits_{(A_m)_{m\leq 0}\in\C^\a_{n,\e}(Q)}\sum\limits_{m\leq 0}\int\limits_{S^mA_m}Z^\a\left(\log Z\right)^{n}d\phi_0,\]
	\[\bar\Psi^\a_{n,\e}(Q):=\lim\limits_{i\to\infty}\Psi^\a_{n,\e}(S^{-i}Q)\ \ \ \mbox{ and}\]
	\[\bar\Psi^\a_{n}(Q):=\lim\limits_{\e\to 0}\bar\Psi^\a_{n,\e}(Q),\]
	since,  as in the proof of Lemma 3 in \cite{Wer3},  $\Psi^\a_{n,\e}(Q)\leq\Psi^\a_{n,\e}(S^{-1}Q)$ and, obviously, $\Psi^\a_{n,\e}(Q)\leq\Psi^\a_{n,\d}(Q)$ for all $0<\d\leq\e$.
	
	Let $n\in\N$. Let $0\leq\a_0\leq 1$ if $n=1$ and $0<\a_0\leq 1$ otherwise.  Define
	\[\Psi^{\a,\a_0}_{n,\e}(Q):=\inf\limits_{(A_m)_{m\leq 0}\in\C^{\a_0}_{n,\e}(Q)}\sum\limits_{m\leq 0}\int\limits_{S^mA_m}Z^\a\left(\log Z\right)^{n}d\phi_0,\]
	\[\Psi^{\a,\a_0}_{n}(Q):=\lim\limits_{\e\to 0}\Psi^{\a,\a_0}_{n,\e}(Q),\]
	\[\bar\Psi^{\a,\a_0}_{n,\e}(Q):=\lim\limits_{i\to\infty}\Psi^{\a,\a_0}_{n,\e}(S^{-i}Q)\ \ \ \mbox{ and}\]
	\[\bar\Psi^{\a,\a_0}_{n}(Q):=\lim\limits_{\e\to 0}\bar\Psi^{\a,\a_0}_{n,\e}(Q).\]
	
	Let $\dot\C^\a_{n,\e}(Q)$ denote the set of all $(A_m)_{m\leq 0}\in\C^\a_{n,\e}(Q)$ such that $A_m$'s are pairwise disjoint. By Lemma 10 (ii) (Lemma 6 (ii) in the arXiv version) in \cite{Wer15}, $\dot\C^\a_{n,\e}(Q)$ is not empty. Define
	\[\dot\Psi^{\a,\a_0}_{n,\e}(Q):=\inf\limits_{(A_m)_{m\leq 0}\in\dot\C^{\a_0}_{n,\e}(Q)}\sum\limits_{m\leq 0}\int\limits_{S^mA_m}Z^\a\left(\log Z\right)^{n}d\phi_0 \ \ \ \mbox{ and }\]
	$\bar{\dot\Psi}^{\a,\a_0}_{n}(Q)$ the same way as $\bar\Psi^{\a,\a_0}_{n}(Q)$.
\end{Definition}
By Lemma 10 (ii) in \cite{Wer15}, $\bar{\dot\Psi}^{\a,\a_0}_{n}(Q) = \bar\Psi^{\a,\a_0}_{n}(Q)$.

 The set functions $\Psi^{\a,\a_0}_{n}(Q)$, $Q\in\B$, have the following properties.

Let us abbreviate
\[\G_{n}^{\a_0,\a}(Q):=\left(\frac{n}{\a_0 e}\right)^n\Phi(Q)+\left(\frac{n}{(1-\a)e}\right)^n\Lambda(Q)\]
for all $Q\in\B$, $\a_0\in(0,1]$, $\a\in[0,1)$ and $n\in\N$.

\begin{lemma}\label{sdpl}
	Let $n\in\N$, $Q\in\B$ and $\a\in(0,1)$. Let $0\leq\a_0\leq 1$ if $n=1$ and $0<\a_0\leq 1$ otherwise. Then the following holds true.
	
	(i) If $n$ is odd, then
	\[-\left(\frac{n}{\a e}\right)^n\Phi(Q)\leq\Psi^{\a,\a_0}_{n}(Q)\leq\left(\frac{n}{(1-\a)e}\right)^n\Lambda(Q)\ \ \ \mbox{ and}\]
	\[\frac{1}{\a}\left(\H^{\a,0}\left(\Lambda,\phi_0\right)(Q)-\Phi(Q)\right)\leq\Psi^{\a,\a_0}_{1}(Q)\leq\frac{1}{1-\a}\left(\Lambda(Q)-\H^{\a,0}\left(\Lambda,\phi_0\right)(Q)\right).\]
	
	(ii) If $n$ is even, then
	\[0\leq\Psi^{\a,\a_0}_{n}(Q)\leq\G_{n}^{\a,\a}(Q).\]
	
	(iii) \[\Psi^{\a,\a_0}_{n}(Q)=\bar\Psi^{\a,\a_0}_{n}(Q)\ \ \ \mbox{ for all }Q\in\B,\ \ \ \mbox{ and}\]
	$\Psi^{\a,\a_0}_{n}$ is a $S$-invariant (signed) measure on $\B$.
\end{lemma}
{\it Proof.} The proof completes Definition \ref{hmdd} by induction.

$(i)$ Let $\e>0$ and $(A_m)_{m\leq 0}\in\C^{\a_0}_{n,\e}(Q)$. Since, by Lemma \ref{hfpl},
\begin{eqnarray*}
	&&-\left(\frac{n}{\a e}\right)^n(\Phi(Q)+\e)\leq-\left(\frac{n}{\a e}\right)^n\sum\limits_{m\leq 0}\phi_0\left(S^mA_m\right)\leq\sum\limits_{m\leq 0}\int\limits_{S^mA_m}Z^\a(\log Z)^nd\phi_0\\
	&&\leq\sum\limits_{m\leq 0}\int\limits_{S^mA_m\cap\{Z > 1\}}e^{-(1-\a)\log Z}(\log Z)^nd\Lambda
	\leq\left(\frac{n}{(1-\a) e}\right)^n\sum\limits_{m\leq 0}\Lambda\left(A_m\right),
\end{eqnarray*}
the first assertion in (i) follows by Proposition 12 (Proposition 2 in the arXiv version) in \cite{Wer15}. The second and the third assertions in (i) follow by the inequalities $1/\a(Z^\a - 1)\leq Z^\a\log Z\leq1/(1-\a)(Z-Z^\a)$. 

$(ii)$ The first inequality in (ii) is obvious.

By Lemma \ref{hfpl}, 
\begin{eqnarray*}
	\Psi^{\a,\a_0}_{n,\e}(Q)&\leq&\sum\limits_{m\leq 0}\int\limits_{S^mA_m}Z^\a(\log Z)^nd\phi_0\\
	&=&\sum\limits_{m\leq 0}\int\limits_{S^mA_m\cap\{Z\leq 1\}}Z^\a(\log Z)^nd\phi_0\\
	&& + \sum\limits_{m\leq 0}\int\limits_{S^mA_m\cap\{Z > 1\}}e^{-(1-\a)\log Z}(\log Z)^nd\Lambda\\
	&\leq&\left(\frac{n}{\a e}\right)^n(\Phi(Q)+\e) + \left(\frac{n}{(1-\a) e}\right)^n\sum\limits_{m\leq 0}\Lambda\left(A_m\right).
\end{eqnarray*}
Hence, by Proposition 12 in \cite{Wer15},
\[\Psi^{\a,\a_0}_{n,\e}(Q)\leq\left(\frac{n}{\a e}\right)^n(\Phi(Q)+\e) + \left(\frac{n}{(1-\a) e}\right)^n\Lambda\left(Q\right).\]
Thus, the second inequality in (ii) follows.

$(iii)$ Let $A\in\A_0$ and $n\in\N\cap\{0\}$. Define (with $(-\infty)^0:=1$)
\begin{equation*}
c_{\a_0,n}:=\left\{\begin{array}{cc} \left(\frac{n}{\a_0 e}\right)^{n}& \mbox{if }n\mbox{ is odd,}\\
0&\mbox{ otherwise}
\end{array}\right.
\end{equation*}
and
\begin{equation*}
\psi_{\a_0,n}(A):=\left\{\begin{array}{cc} \int\limits_{A}\left(Z^{\a_0}\left(\log Z\right)^{n}+c_{\a_0,n}\right)d\phi_0& \mbox{if }n\mbox{ is odd,}\\
\int\limits_{A}Z^{\a_0}\left(\log Z\right)^{n}d\phi_0  &\mbox{ otherwise.}
\end{array}\right. 
\end{equation*}
Then, by Lemma \ref{hfpl}, $\psi_{\a_0,n}(A)>0$, and
\[\int\limits_{A}Z^{\a_0}\left(\log Z\right)^{n}d\phi_0=\psi_{\a_0,n}(A)-c_{\a_0,n}\phi_0(A)\]
for all $n$. Thus, applying Lemma 10 (i) (Lemma 6 (i) in the arXiv version) in \cite{Wer15} to the families $\psi_{\a_0,0}$,...,$\psi_{\a_0,n}$,$\psi_{\a,n+1}$ and $c_{\a_0,0}$,...,$c_{\a_0,n}$,$c_{\a,n+1}$ implies, by Corollary 8 (ii) (Corollary 1 (ii) in the arXiv version) in \cite{Wer15}, that $\bar\Psi^{\a,\a_0}_{n+1}$ is a (signed) $S$-invariant measure on $\B$. Since, by (i) or (ii) it is finite, it follows by Theorem 16 (ii) (Theorem 4 (ii) in the arXiv version) in \cite{Wer15}, that it is equal to $\Psi^{\a,\a_0}_{n+1}$ on $\B$.
\hfill$\Box$

\subsubsection{The continuity of the candidates for the derivatives of $(0,1)\owns\a\longmapsto\H^{\a,0}\left(\Lambda,\phi_0\right)$}

Now, we show some continuity properties of the obtained measures with respect to the first parameter. 

\begin{lemma}\label{eql}
	Let $n\in\N\cup\{0\}$, $0<\a_0\leq\a<1$, $\g\in[0,1]$ and $Q\in\B$.\\
	(i) In $n$ is even, then
	\begin{eqnarray}\label{sac}
	&&-(\a-\a_0)\left(\frac{n+1}{\a_0e}\right)^{n+1}\Phi(Q)\leq\Psi^{\a,\g}_{n}(Q)-\Psi^{\a_0,\g}_{n}(Q)\nonumber\\
	&\leq&(\a-\a_0)\left(\frac{n+1}{(1-\a)e}\right)^{n+1}\Lambda(Q).
	\end{eqnarray}
		
	(ii) If $n$ is odd, then 
	\begin{eqnarray}\label{odeb}
	0&\leq&\Psi^{\a,\g}_{n,\e}(Q)-\Psi^{\a_0,\g}_{n,\e}(Q)\\
	&\leq& (\a-\a_0)\left(\left(\frac{n+1}{\a_0e}\right)^{n+1}(\Phi(Q)+\e) +\left(\frac{n+1}{(1-\a)e}\right)^{n+1}\Lambda(X)\right)+\e\left(\frac{n}{\a_0 e}\right)^{n}\nonumber
	\end{eqnarray}
	for all $\e>0$, and
	\begin{equation}\label{odcl}
	0\leq\Psi^{\a,\g}_{n}(Q)-\Psi^{\a_0,\g}_{n}(Q)\leq (\a-\a_0)\G_{n+1}^{\a_0,\a}(Q).
	\end{equation}
\end{lemma}
{\it Proof.}
Let  $\a_0<\a$ and $\e>0$. 

$(i)$ Let $(B_m)_{m\leq 0}\in\C^{\g}_{n,\e}(Q)$. Then, by the first inequality of Lemma \ref{edl} (i) and Lemma \ref{hfpl}, 
\begin{eqnarray*}
	&&-(\a-\a_0)\left(\frac{n+1}{\a_0e}\right)^{n+1}(\Phi(Q)+\e)\\
	&\leq&\sum\limits_{m\leq 0}\int\limits_{S^mA_m}Z^{\a}(\log Z)^nd\phi_0-\sum\limits_{m\leq 0}\int\limits_{S^mA_m}Z^{\a_0}(\log Z)^nd\phi_0\\
	&\leq&\sum\limits_{m\leq 0}\int\limits_{S^mA_m}Z^{\a}(\log Z)^nd\phi_0-\Psi^{\a_0,\g}_{n,\e}(Q).
\end{eqnarray*}
Thus, it follow the first inequalities of \eqref{sac}.

Now, let $(A_m)_{m\leq 0}\in\C^{\g}_{n,\e}(Q)$ such that
\[\sum\limits_{m\leq 0}\int\limits_{S^mA_m}Z^{\a_0}(\log Z)^nd\phi_0<\Psi^{\a_0,\g}_{n}(Q)+\e.\]
Then, by the second inequality of Lemma \ref{edl} (i) and Lemma \ref{hfpl},
\begin{eqnarray*}
	&&\Psi^{\a,\g}_{n,\e}(Q)-\Psi^{\a_0,\g}_{n}(Q)-\e\\
	&\leq&\sum\limits_{m\leq 0}\int\limits_{S^mA_m}Z^{\a}(\log Z)^nd\phi_0-\sum\limits_{m\leq 0}\int\limits_{S^mA_m}Z^{\a_0}(\log Z)^nd\phi_0\\
	&\leq&(\a-\a_0)\left(\frac{n+1}{(1-\a)e}\right)^{n+1}\sum\limits_{m\leq 0}\Lambda(A_m).
\end{eqnarray*}
Hence, by Proposition 12 (Proposition 2 in the arXiv version) in \cite{Wer15}, it follows the second inequality of \eqref{sac}.

$(ii)$ Obviously, by Lemma \ref{edl} (ii),
\[0\leq\Psi^{\a,\g}_{n,\e}(Q)-\Psi^{\a_0,\g}_{n,\e}(Q).\]
Let $(B_m)_{m\leq 0}\in\dot\C^{\g}_{n+1,\e}(Q)$. Then, by Lemma \ref{edl} (ii) and Lemma \ref{hfpl},
\begin{eqnarray*}
	&&\dot\Psi^{\a,\g}_{n,\e}(Q)-\sum\limits_{m\leq 0}\int\limits_{S^mB_m}Z^{\a_0}(\log Z)^nd\phi_0\\
	&\leq&\sum\limits_{m\leq 0}\int\limits_{S^mB_m}Z^{\a}(\log Z)^nd\phi_0-\sum\limits_{m\leq 0}\int\limits_{S^mB_m}Z^{\a_0}(\log Z)^nd\phi_0\\
	&\leq&(\a-\a_0)\left(\left(\frac{n+1}{\a_0e}\right)^{n+1}\sum\limits_{m\leq 0}\phi_0\left(S^mA_m\right)+\left(\frac{n+1}{(1-\a)e}\right)^{n+1}\sum\limits_{m\leq 0}\Lambda(A_m)\right)\\
	&\leq&(\a-\a_0)\left(\left(\frac{n+1}{\a_0e}\right)^{n+1}(\Phi(Q)+\e)+\left(\frac{n+1}{(1-\a)e}\right)^{n+1}\Lambda(X)\right).
\end{eqnarray*}
Hence,
\begin{eqnarray*}
	&&\dot\Psi^{\a,\g}_{n,\e}(Q)-\dot\Psi^{\a_0,\g}_{n,\e}(Q)\\
	&\leq&(\a-\a_0)\left(\left(\frac{n+1}{\a_0e}\right)^{n+1}(\Phi(Q)+\e)+\left(\frac{n+1}{(1-\a)e}\right)^{n+1}\Lambda(X)\right).
\end{eqnarray*}
Since $\Psi^{\a,\g}_{n,\e}(Q)\leq\dot\Psi^{\a,\g}_{n,\e}(Q)$ and, by Lemma 10 (ii) (Lemma 6 (ii) in the arXiv version) in \cite{Wer15}, 
\[\dot\Psi^{\a_0,\g}_{n,\e}(Q)\leq\Psi^{\a_0,\g}_{n,\e}(Q)+\e\left(\frac{n}{\a_0 e}\right)^{n},\]
it follows \eqref{odeb}. \eqref{odcl} follows by Lemma \ref{edl} (ii) and Lemma \ref{hfpl}, the same way as in the proof of (i).
\hfill$\Box$

\begin{Remark}
In the case $n=0$, Lemma \ref{eql} (i) gives the following continuity property of $(0,1)\owns\a\longmapsto\H^{\a,0}\left(\Lambda,\phi_0\right)$.
\[-(\a-\a_0)\frac{\Phi(Q)}{\a_0e}\leq\H^{\a,0}\left(\Lambda,\phi_0\right)(Q)-\H^{\a_0,0}\left(\Lambda,\phi_0\right)(Q)\leq(\a-\a_0)\frac{\Lambda(Q)}{(1-\a)e}\]
for all $0<\a_0\leq\a<1$ and $Q\in\B$, which is weaker than that of Lemma \ref{pchm} (vii).
\end{Remark}

\subsubsection{The right derivative of $(0,1)\owns\a\longmapsto\H^{\a,0}\left(\Lambda,\phi_0\right)$}\label{fdhs}

We show now that $\Psi^{\a,\a}_{1}(Q)$ is the right derivative of $(0,1)\owns\a\longmapsto\H^{\a,0}\left(\Lambda,\phi_0\right)(Q)$ for all $Q\in\B$. Also, as a by-product, we obtain another lower bound for $\Phi$ in terms of $\Psi^{\a,\a}_{1}$ and $\H^{\a,0}(\Lambda,\phi_0)$.

\begin{lemma}\label{hmfd} Let $0<\a_0<\a\leq 1$ and $Q\in\B$.
	
	(i)  Let $\e_0,\e>0$. Let $\d_0,\d>0$ such that $\H^{\a_0,0}_{\d_0}\left(\Lambda,\phi_0\right)(Q)>\H^{\a_0,0}\left(\Lambda,\phi_0\right)(Q)-\e_0$ and $\H^{\a,0}_\d\left(\Lambda,\phi_0\right)(Q)>\H^{\a,0}\left(\Lambda,\phi_0\right)(Q)-\e$. Then
	 \begin{eqnarray*}
	 	(\a-\a_0)\Psi^{\a_0,\a}_{1,\d_0}(Q)-\e_0-\d_0&<&\H^{\a,0}\left(\Lambda,\phi_0\right)(Q)-\H^{\a_0,0}\left(\Lambda,\phi_0\right)(Q)\\
	 	&<&(\a-\a_0)\Psi^{\a,\a_0}_{1,\d}(Q)+\e+\d.
	 \end{eqnarray*} 
	 
	(ii) 
	\[\Psi^{\a_0,\a_0}_{1}(Q)\leq\frac{\H^{\a,\a_0,0}\left(\Lambda,\phi_0\right)(Q)-\H^{\a_0,0}\left(\Lambda,\phi_0\right)(Q)}{\a-\a_0}\leq\Psi^{\a,\a_0}_{1}(Q), \]
	\[\Psi^{\a_0,\a}_{1}(Q)\leq\frac{\H^{\a,0}\left(\Lambda,\phi_0\right)(Q)-\H^{\a_0,\a,0}\left(\Lambda,\phi_0\right)(Q)}{\a-\a_0}\leq\Psi^{\a,\a_0}_{1}(Q),\mbox{ and}\]
	\[0\leq\frac{\H^{\a,\a_0,0}\left(\Lambda,\phi_0\right)(Q)-\H^{\a,0}\left(\Lambda,\phi_0\right)(Q)}{\a-\a_0}\leq\Psi^{\a,\a_0}_{1}(Q)-\Psi^{\a_0,\a}_{1}(Q). \]	 	
\end{lemma}
{\it Proof.}
 $(i)$  By Lemma \ref{edl} (i), for any  $(A_m)_{m\leq 0}\in\C(Q)$ with $\sum_{m\leq 0}\phi_0(S^mA_m)<\infty$, 
 \begin{equation}\label{hmde}
 	\sum\limits_{m\leq 0}\int\limits_{S^mA_m}Z^{\a_1} d\phi_0\geq\sum\limits_{m\leq 0}\int\limits_{S^mA_m}Z^{\a_2} d\phi_0+\left(\a_1-\a_2\right)\sum\limits_{m\leq 0}\int\limits_{S^mA_m}Z^{\a_2}\log Z d\phi_0
 \end{equation}
 for all $\a_1\in [0,1]$ and $\a_2\in (0,1]$.
 Hence, putting $\a_1=\a$, $\a_2=\a_0$ and taking $(A_m)_{m\leq 0}\in\C^{\a}_{1,\d_0}(Q)$ implies that
 \begin{eqnarray*}
 	\H^{\a,0}(\Lambda,\phi_0)(Q) +\d_0&>&\H^{\a_0,0}_{\d_0}\left(\Lambda,\phi_0\right)(Q) +(\a-\a_0)\Psi^{\a_0,\a}_{1,\d_0}(Q)\\
 	&>&\H^{\a_0,0}(\Lambda,\phi_0)(Q) -\e_0 +(\a-\a_0)\Psi^{\a_0,\a}_{1,\d_0}(Q),
 \end{eqnarray*}
 which is the first inequality of (i). The same way, putting $\a_1=\a_0$, $\a_2=\a$ and taking infimum over all $(A_m)_{m\leq 0}\in\C^{\a_0}_{1,\d}(Q)$ implies the second inequality.
 
 $(ii)$ Let $(A_m)_{m\leq 0}\in\C^{\a_0}_{1,\d}(Q)$. Substituting $\a_1:=\a_0$ and $\a_2:=\a$ in \eqref{hmde} implies that
 \[\H^{\a_0,0}(\Lambda,\phi_0)(Q) +\d>\H^{\a, \a_0,0}_\d(\Lambda,\phi_0)(Q)-\left(\a-\a_0\right)\sum\limits_{m\leq 0}\int\limits_{S^mB_m}Z^{\a}\log Z d\phi_0.\]
 This gives the second inequality of (ii).
 
 Substituting $\a_1:=\a$ and $\a_2:=\a_0$ in \eqref{hmde} implies that
 \[\sum\limits_{m\leq 0}\int\limits_{S^mB_m}Z^{\a} d\phi_0\geq\H^{\a_0, \a_0,0}_\d(\Lambda,\phi_0)(Q)+\left(\a-\a_0\right)\Psi^{\a_0,\a_0}_{1,\d}(Q).\]
 This gives the first inequality of (ii). 
 
 If $(A_m)_{m\leq 0}\in\C^{\a}_{1,\d}(Q)$, then
 \[\H^{\a,0}(\Lambda,\phi_0)(Q) +\d>\H^{\a_0, \a,0}_\d(\Lambda,\phi_0)(Q)+\left(\a-\a_0\right)\Psi^{\a_0,\a}_{1,\d}(Q).\]
 This implies the third inequality in (ii).
 
 The fourth inequality in (ii) follows from (i), since $\H^{\a_0,0}(\Lambda,\phi_0)(Q)\\\leq\H^{\a_0,\a,0}(\Lambda,\phi_0)(Q)$.
 
 The fifth inequality in (ii) is obvious.
 
 The sixth inequality in (ii) is obvious if $\a=1$ and $\Psi^{1,\a_0}_{1}(Q)=+\infty$. Suppose $\a<1$, or $\Psi^{1,\a_0}_{1}(Q)<+\infty$.
 Let $\eta,\tau>0$. Let $(C_m)_{m\leq 0}\in\C^{\a_0}_{1,\tau}(Q)$ such that
 \[\sum\limits_{m\leq 0}\int\limits_{S^mC_m}Z^{\a}\log Z d\phi_0<\Psi^{\a,\a_0}_{1,\tau}(Q)+\eta.\]
 Then, by (i),
 \begin{eqnarray*}
   &&\H^{\a,0}(\Lambda,\phi_0)(Q)+\tau\\
   &\geq&\H^{\a_0,0}(\Lambda,\phi_0)(Q)+(\a-\a_0)\Psi^{\a_0,\a}_{1}(Q)+\tau\\
   &>&\sum\limits_{m\leq 0}\int\limits_{S^mC_m}Z^{\a_0}d\phi_0+(\a-\a_0)\Psi^{\a_0,\a}_{1}(Q)\\
   &\geq&\sum\limits_{m\leq 0}\int\limits_{S^mC_m}Z^{\a}d\phi_0-(\a-\a_0)\sum\limits_{m\leq 0}\int\limits_{S^mC_m}Z^{\a}\log Zd\phi_0+(\a-\a_0)\Psi^{\a_0,\a}_{1}(Q)\\
   &>&\sum\limits_{m\leq 0}\int\limits_{S^mC_m}Z^{\a}d\phi_0-(\a-\a_0)\left(\Psi^{\a,\a_0}_{1,\tau}(Q)-\Psi^{\a_0,\a}_{1}(Q)+\eta\right).
 \end{eqnarray*}
 Hence,
 \[(C_m)_{m\leq 0}\in\C^{\a}_{1,(\a-\a_0)\left(\Psi^{\a,\a_0}_{1,\tau}(Q)-\Psi^{\a_0,\a}_{1}(Q)+\eta\right)+\tau}(Q).\]
 Therefore,
 \begin{eqnarray*}
   H^{\a,\a_0,0}_\tau(\Lambda,\phi_0)(Q)&\leq&\sum\limits_{m\leq 0}\int\limits_{S^mC_m}Z^{\a}d\phi_0\\&<&\H^{\a,0}(\Lambda,\phi_0)(Q)+(\a-\a_0)\left(\Psi^{\a,\a_0}_{1,\tau}(Q)-\Psi^{\a_0,\a}_{1}(Q)+\eta\right)\\
   &&+\tau.
 \end{eqnarray*}
 Since $\eta,\tau>0$ were arbitrary, this implies the sixth inequality of (ii).
\hfill$\Box$

\begin{prop}\label{lbvd}
     For every $0\leq\b\leq\a_0<\a\leq 1$ and $Q\in\B$,
	\[\H^{\a_0,0}\left(\Lambda,\phi_0\right)(Q)\log\frac{\H^{\a_0,0}\left(\Lambda,\phi_0\right)(Q)}{\H^{\b,\a,0}\left(\Lambda,\phi_0\right)(Q)}\leq(\a_0-\b)\Psi^{\a,\a_0}_{1}(Q).\]
	In particular,
	\[\Phi(Q)\geq \H^{\a_0,0}\left(\Lambda,\phi_0\right)(Q)e^{\frac{-\a_0}{\H^{\a_0,0}\left(\Lambda,\phi_0\right)(Q)}\Psi^{\a,\a_0}_{1}(Q)}\]
	if $\H^{\a_0,0}\left(\Lambda,\phi_0\right)(Q)>0$.
\end{prop}
{\it Proof.} 
    For $\b<\a_0$, the assertion follows by the first inequality of Lemma \ref{pchm} (v) together with the second one of Lemma \ref{hmfd} (i). For $\b=\a_0$, it is obvious, since $\H^{\a_0,0}\left(\Lambda,\phi_0\right)(Q)\leq\H^{\a_0,\a,0}\left(\Lambda,\phi_0\right)(Q)$.
    
    It can be also deduced from Lemma \ref{ml} (iv).
\hfill$\Box$

Now, we are ready to show the right differentiability of $(0,1)\owns\a\longmapsto\H^{\a,0}(\Lambda,\phi_0)(Q)$. In order also to shed some light on the problem for $\Psi^{\a,\a}_{1}(Q)$ being also the left derivative of the function, we need the following definitions.

\begin{Definition}
	Let $Q\in\P(X)$ and $\t>0$. Define
	\begin{eqnarray*}
		\d_\t(\a_1,\a_2)&:=&|\a_1-\a_2|^\frac{\t}{2}\sup\left\{ 0<\d<|\a_1-\a_2|^\frac{\t}{2}:\phantom{\H^{\a_i,0}_{\d}\left(\Lambda,\phi_0\right)}\right.\\ &&\left.\H^{\a_i,0}_{\d}\left(\Lambda,\phi_0\right)(Q)>\H^{\a_i,0}\left(\Lambda,\phi_0\right)(Q)-|\a_1-\a_2|^\t\mbox{ for }i=1,2\vphantom{|\a_1-\a_2|^\frac{\t}{2}}\right\}
	\end{eqnarray*}
	for all $\a_1,\a_2\in[0,1]$. For $0<\a_0\leq\a< 1$, define
	\[\e_\t(\a_0,\a):=(\a-\a_0)\left(\Psi^{\a,\a_0}_{1,\d_\t(\a_0,\a)}(Q)-\Psi^{\a_0,\a}_{1,\d_\t(\a_0,\a)}(Q)\right)+2(\a-\a_0)^\t+3\d_\t(\a_0,\a).\]
\end{Definition}

\begin{theo}\label{dlpl}  
 Let $Q\in\B$. Then the function $(0,1)\owns\a\longmapsto\H^{\a,0}(\Lambda,\phi_0)(Q)$ is right differentiable, and
 \[\left.\frac{d_+}{d_+\a}\H^{\a,0}(\Lambda,\phi_0)(Q)\right|_{\a=\a_0} =  \Psi^{\a_0,\a_0}_{1}(Q)=\lim\limits_{\a\to^+\a_0}\Psi^{\a,\a}_{1}(Q)\]
 for all $0<\a_0<1$ where $d_+/d_+\a$ denotes the right derivative. 
 
 From the left, for every $0<\a<1$ and $\t>0$,
 \[\lim\limits_{\a_0\to^-\a}\Psi^{\a_0,\a_0}_{1,\e_\t(\a_0,\a)}(Q)=\lim_{\a_0\to^-\a}\Psi^{\a,\a_0}_{1,\e_\t(\a_0,\a)}(Q)=\Psi^{\a,\a}_{1}(Q),\mbox{ and}\]
  \begin{eqnarray}\label{ldi}
  \Psi^{\a,\a}_{1}(Q)&\leq&\liminf\limits_{\a_0\to^-\a}\frac{\H^{\a_0,0}\left(\Lambda,\phi_0\right)(Q)-\H^{\a,0}\left(\Lambda,\phi_0\right)(Q)}{\a_0-\a}\leq\liminf\limits_{\a_0\to^-\a}\Psi^{\a,\a_0}_{1,\d_\t(\a_0,\a)}(Q)\nonumber\\
  &=&\liminf\limits_{\a_0\to^-\a}\Psi^{\a_0,\a_0}_{1,\d_\t(\a_0,\a)}(Q)
  \end{eqnarray}
  for all $\t>1$.
\end{theo}
{\it Proof.} 
Let $0<\g_0<\g< 1$ and $\t>0$. Observe that $0<\d_\t(\g_0,\g)\leq(\g-\g_0)^\t$, $\H^{\g_0,0}_{\d_\t(\g_0,\g)}\left(\Lambda,\phi_0\right)(Q)>\H^{\g_0,0}\left(\Lambda,\phi_0\right)(Q)-(\g-\g_0)^\t$ and $\H^{\g,0}_{\d_\t(\g_0,\g)}\left(\Lambda,\phi_0\right)(Q)>\H^{\g,0}\left(\Lambda,\phi_0\right)(Q)-(\g-\g_0)^\t$. 
Let $(B_m)_{m\leq 0}\in\C^{\g}_{1,\d_\t(\g_0,\g)}(Q)$. Then, by Lemma \ref{hmfd} (i) and Lemma \ref{edl} (i),
\begin{eqnarray*}
&&\H^{\g_0,0}\left(\Lambda,\phi_0\right)(Q)+(\g-\g_0)\Psi^{\g,\g_0}_{1,\d_\t(\g_0,\g)}(Q)+2(\g-\g_0)^\t+3\d_\t(\g_0,\g)\\
&>&\H^{\g,0}\left(\Lambda,\phi_0\right)(Q) + \d_\t(\g_0,\g)\\
&>&\sum\limits_{m\leq 0}\int\limits_{S^mB_m}Z^{\g} d\phi_0\\
&\geq&\sum\limits_{m\leq 0}\int\limits_{S^mB_m}Z^{\g_0} d\phi_0 +(\g-\g_0)\sum\limits_{m\leq 0}\int\limits_{S^mB_m}Z^{\g_0}\log Z d\phi_0\\
&\geq&\sum\limits_{m\leq 0}\int\limits_{S^mB_m}Z^{\g_0} d\phi_0 +(\g-\g_0)\Psi^{\g_0,\g}_{1,\d_\t(\g_0,\g)}(Q).
\end{eqnarray*}
Hence, since, by Lemma \ref{hmfd} (i), $(\g-\g_0)(\Psi^{\g,\g_0}_{1,\d_\t(\g_0,\g)}(Q)-\Psi^{\g_0,\g}_{1,\d_\t(\g_0,\g)}(Q))+2(\g-\g_0)^\t+3\d_\t(\g_0,\g)>\d_\t(\g_0,\g)$,
\begin{equation*}
  (B_m)_{m\leq 0}\in\C^{\g_0}_{1,\e_\t(\g_0,\g)}(Q).
\end{equation*}
That is
\begin{eqnarray}\label{cir}
\C^{\g}_{1,\d_\t(\g_0,\g)}(Q)\subset\C^{\g_0}_{1,\e_\t(\g_0,\g)}(Q).
\end{eqnarray}
Therefore, for every $0<\a\leq 1$,
\begin{equation}\label{saci}
   \Psi^{\a,\g}_{1,\d_\t(\g_0,\g)}(Q)\geq\Psi^{\a,\g_0}_{1,\e_\t(\g_0,\g)}(Q).
\end{equation}
In particular, by setting $\a=\g_0$ and letting $\g\to^+\g_0$, it follows, since\\ $\Psi^{\a,\g}_{1,\d_\t(\g_0,\g)}(Q)\leq \Psi^{\a,\g}_{1}(Q)$, that
\[\Psi^{\g_0,\g_0}_{1}(Q)\leq\liminf\limits_{\g\to^+\g_0}\Psi^{\g_0,\g}_{1}(Q).\]
 Since, by Lemma \ref{hmfd} (i),
 \begin{equation}\label{hmci}
    \Psi^{\g_0,\g}_{1}(Q)\leq\frac{\H^{\g,0}\left(\Lambda,\phi_0\right)(Q)-\H^{\g_0,0}\left(\Lambda,\phi_0\right)(Q)}{\g-\g_0}\leq\Psi^{\g,\g_0}_{1}(Q)
 \end{equation}
 and, by Lemma \ref{eql} (ii), $\lim_{\g\to^+\g_0}\Psi^{\g,\g_0}_{1}(Q)=\Psi^{\g_0,\g_0}_{1}(Q)$, it follows that
\begin{equation}\label{rdhm}
\left.\lim\limits_{\g\to^+\g_0}\Psi^{\g_0,\g}_{1}(Q)=\frac{d_+\H^{\a,0}\left(\Lambda,\phi_0\right)(Q)}{d_+\a}\right|_{\a=\g_0}=\Psi^{\g_0,\g_0}_{1}(Q).
\end{equation}
This proves the right differentiability of $(0,1)\owns\a\longmapsto\H^{\a,0}\left(\Lambda,\phi_0\right)(Q)$. Also, by \eqref{odcl} and \eqref{hmci}, for all $0<\a_0<\a<1$,
\begin{eqnarray*}
	&&\Psi^{\a_0,\a_0}_{1}(Q)+(\a-\a_0)\G_{2}^{\a_0,\a}(Q)\geq\Psi^{\a,\a_0}_{1}(Q)\geq\Psi^{\a_0,\a}_{1}(Q)\\
	&\geq&\Psi^{\a,\a}_{1}(Q)-(\a-\a_0)\G_{2}^{\a_0,\a}(Q)\geq\Psi^{\a_0,\a}_{1}(Q)-(\a-\a_0)\G_{2}^{\a_0,\a}(Q).
\end{eqnarray*}
Thus, by \eqref{rdhm},
\[\lim\limits_{\a\to^+\a_0}\Psi^{\a,\a}_{1}(Q)=\Psi^{\a_0,\a_0}_{1}(Q).\]

Now, let us consider the differentiability from the left.  Let $\e>0$ and $(C_m)_{m\leq 0}\in\dot\C^{\g_0}_{1,\e}(Q)$. By \eqref{hmci}, Lemma \ref{sdpl} (i), Lemma \ref{edl} (i) and Lemma \ref{hfpl}, 
\begin{eqnarray*}
&&\H^{\g,0}\left(\Lambda,\phi_0\right)(Q)+\frac{\g-\g_0}{e\g_0} \Phi(Q)+\e\geq\H^{\g_0,0}\left(\Lambda,\phi_0\right)(Q) + \e\\
&>&\sum\limits_{m\leq 0}\int\limits_{S^mC_m}Z^{\g_0} d\phi_0\geq\sum\limits_{m\leq 0}\int\limits_{S^mC_m}Z^{\g} d\phi_0 +\frac{\g_0-\g}{e(1-\g)}\Lambda\left(X\right),
\end{eqnarray*}
and therefore,
\begin{equation}\label{lsae}
   (C_m)_{m\leq 0}\in\dot\C^{\g}_{1,\frac{\g-\g_0}{e}\left(\frac{\Lambda\left(X\right)}{1-\g}+\frac{\Phi(Q)}{\g_0} \right)+\e}(Q).
\end{equation}
That is
\[\dot\C^{\g_0}_{1,\e}(Q)\subset\dot\C^{\g}_{1,\frac{\g-\g_0}{e}\left(\frac{\Lambda\left(X\right)}{1-\g}+\frac{\Phi(Q)}{\g_0} \right)+\e}(Q).\]
Therefore,  for every $0<\a\leq 1$,
\begin{equation}\label{salc}
   \dot\Psi^{\a,\g_0}_{1,\e}(Q)\geq\dot\Psi^{\a,\g}_{1,\frac{\g-\g_0}{e}\left(\frac{\Lambda\left(X\right)}{1-\g}+\frac{\Phi(Q)}{\g_0} \right)+\e}(Q).
\end{equation}
 Since, by Lemma 10 (ii) (Lemma 6 (ii) in the arXiv version) in \cite{Wer15}, 
 \[\dot\Psi^{\a,\g_0}_{1,\e}(Q)\leq\Psi^{\a,\g_0}_{1,\e}(Q)+\frac{\e}{\a e},\]
 it follows, by \eqref{saci} and \eqref{salc}, that
 \begin{eqnarray}\label{fcsi}
   \Psi^{\a,\g}_{1,\d_\t(\g,\g_0)}(Q)+\frac{\e_\t(\g_0,\g)}{\a e}&\geq&\Psi^{\a,\g_0}_{1,\e_\t(\g_0,\g)}(Q)+\frac{\e_\t(\g_0,\g)}{\a e}\nonumber\\
   &\geq&\Psi^{\a,\g}_{1,\frac{\g-\g_0}{e}\left(\frac{\Lambda\left(X\right)}{1-\g}+\frac{\Phi(Q)}{\g_0} \right)+\e_\t(\g_0,\g)}(Q).
 \end{eqnarray}
 Furthermore,  by \eqref{odeb},
 \begin{equation*}
    \Psi^{\g_0,\g}_{1,\e}(Q)\leq\Psi^{\g,\g}_{1,\e}(Q)\leq\Psi^{\g_0,\g}_{1,\e}(Q)+c(\g_0,\g,\e)(\g-\g_0)+\frac{\e}{\g_0e}
 \end{equation*}
 where
 \[c(\g_0,\g,\e):=\left(\frac{2}{\g_0e}\right)^2\left(\Phi(Q)+\e\right) + \left(\frac{2}{(1-\g)e}\right)^2\Lambda(X).\]
 Therefore, putting $\a=\g_0$ in \eqref{fcsi} implies that
 \begin{equation}\label{pdl}
 \lim\limits_{\g_0\to^-\g}\Psi^{\g_0,\g_0}_{1,\e_\t(\g_0,\g)}(Q)=\Psi^{\g,\g}_{1}(Q).
 \end{equation}
  Also, putting $\a=\g$ in \eqref{fcsi} implies that
\[\lim_{\g_0\to^-\g}\Psi^{\g,\g_0}_{1,\e_\t(\g_0,\g)}(Q)=\Psi^{\g,\g}_{1}(Q).\]

Suppose $\t>1$.  Since, by \eqref{saci} and Lemma \ref{hmfd} (i),
\begin{eqnarray*}
	&&\Psi^{\g_0,\g_0}_{1,\e_\t(\g_0,\g)}(Q)-(\g-\g_0)^{\t-1}-\frac{\d_\t(\g_0,\g)}{\g-\g_0}\\
	&\leq&\Psi^{\g_0,\g}_{1,\d_\t(\g_0,\g)}(Q)-(\g-\g_0)^{\t-1}-\frac{\d_\t(\g_0,\g)}{\g-\g_0}\\
&\leq&\frac{\H^{\g,0}\left(\Lambda,\phi_0\right)(Q)-\H^{\g_0,0}\left(\Lambda,\phi_0\right)(Q)}{\g-\g_0}\\
&\leq&\Psi^{\g,\g_0}_{1,\d_\t(\g_0,\g)}(Q)+(\g-\g_0)^{\t-1}+\frac{\d_\t(\g_0,\g)}{\g-\g_0}
\end{eqnarray*} 
it follows \eqref{ldi}, by \eqref{pdl}.
\hfill$\Box$

\subsubsection{The left derivative of $(0,1)\owns\a\longmapsto\H^{\a,0}\left(\Lambda,\phi_0\right)$}\label{ldsh}

Now, we show that $(0,1)\owns\a\longmapsto\H^{\a,0}\left(\Lambda,\phi_0\right)(Q)$ is also left differentiable for all $Q\in\B$, but its left derivative seems to be, in general, a different function. 

\begin{Definition}\label{hldd}
	Let $0<\a\leq 1$, $0\leq\b\leq 1$, $Q\in\P(X)$ and $\e>0$. Define
	\[\H^{\b,0,1}_{\e}\left(\Lambda,\phi_0\right)(Q):=\inf\limits_{(A_m)_{m\leq 0}\in\C^{0,1}_{\e}(Q)}\sum\limits_{m\leq 0}\int\limits_{S^mA_m}Z^\b d\phi_0\ \ \ \mbox{ and}\]
	\[\H^{\b,0,1}\left(\Lambda,\phi_0\right)(Q):=\lim\limits_{\e\to0}\H^{\b,0,1}_{\e}\left(\Lambda,\phi_0\right)(Q).\]
	As in Lemma \ref{shm} (i), on sees that $\H^{\b,0,1}\left(\Lambda,\phi_0\right)(X)<\infty$, and, by Lemma \ref{icp} (iii), 
	\begin{equation}\label{sche}
	   \H^{\b,0,1}\left(\Lambda,\phi_0\right)(Q)=\H^{\b,0}\left(\Lambda,\phi_0\right)(Q)\ \ \mbox{ for all }Q\in\P(X).
	\end{equation}
	 Now, extending Definition \ref{icfw} one step further, define
	\[\C^{\b,0,1}_{\e}(Q):=\left\{(A_m)_{m\leq 0}\in\C^{0,1}_\e(Q)|\ \sum\limits_{m\leq 0}\int\limits_{S^mA_m}Z^\b d\phi_0<\H^{\b,0,1}\left(\Lambda,\phi_0\right)(Q)+\e\right\},\]
	\[\Xi^{\a,\b}_{1,\e}(Q):=\sup\limits_{(A_m)_{m\leq 0}\in\C^{\b,0,1}_{\e}(Q)}\sum\limits_{m\leq 0}\int\limits_{S^mA_m}Z^\a\log Z d\phi_0\]
	and
	\[\Xi^{\a,\b}_{1}(Q):=\lim\limits_{\e\to 0}\Xi^{\a,\b}_{1,\e}(Q).\]
\end{Definition}

Obviously, by \eqref{hmvi} and \eqref{sche}, for every $Q\in\B$, 
\begin{equation}\label{dmr}
	\Psi^{\a,\b}_{1}(Q)\leq\Xi^{\a,\b}_{1}(Q).
\end{equation} 
However, as the next two lemmas show, the latter shares with the former some of the properties. 

In order to show that it is also a measure, we need the following definition.

\begin{Definition}
	 Let $0\leq\a< 1$, $0\leq\b\leq 1$, $Q\in\P(X)$ and $\e>0$. For $A\in\A_0$, let 
	 \[\omega_\a(A):=\int\limits_{A}\left(-e^{-(1-\a)\log Z}\log Z\right)d\Lambda+\frac{1}{(1-\a) e}\Lambda(A).\]
	   Define
	 \[\Omega^{\a,\b}_\e(Q):=\inf\limits_{(A_m)_{m\leq 0}\in\C^{\b,0,1}_{\e}(Q)}\sum\limits_{m\leq 0}\omega_\a\left(S^mA_m\right)\]
	 and
	 \[\Omega^{\a,\b}(Q):=\lim\limits_{\e\to 0}\Omega^{\a,\b}_\e(Q).\]
\end{Definition}

Let us abbreviate 
\[\G^{\a_0,\a}_{2,\e}(Q):=\G^{\a_0,\a}_2(Q)+\frac{4\e}{e^2}\left(\frac{1}{\a_0^2}+\frac{1}{(1-\a)^2}\right).\]

\begin{lemma}\label{ldp}
	Let $0<\a_0\leq\a< 1$, $0\leq\b\leq 1$ and $Q\in\P(X)$.
	
(i)
\begin{equation*}
	-\frac{\Phi(Q)}{\a e}\leq\Xi^{\a,\b}_{1}(Q)\leq\frac{\Lambda(Q)}{(1-\a)e}.
\end{equation*}

(ii) \[0\leq\Omega^{\a,\b}(Q)=-\Xi^{\a,\b}_{1}(Q)+\frac{1}{(1-\a) e}\Lambda(Q).\]

(iii)  $\Xi^{\a,\b}_{1}$ is a $S$-invariant, signed measure on $\B$.

(iv) For every $\e>0$,
\[0\leq\Xi^{\a,\b}_{1,\e}(Q)-\Xi^{\a_0,\b}_{1,\e}(Q)\leq(\a-\a_0)\G^{\a_0,\a}_{2,\e}(Q).\]
\end{lemma}
{\it Proof.} Let $\e>0$ and $(A_m)_{m\leq 0}\in\C^{\b,0,1}_{\e}(Q)$.

$(i)$  Since $Z^\a\log Z\geq -1/(\a e)$, 
\[\Xi^{\a,\b}_{1,\e}(Q)>-\frac{1}{\a e}\left(\Phi(Q)+\e\right).\]
This implies the first inequality of (i). 

On the other hand, by Lemma \ref{hfpl},
\begin{eqnarray*}
	\sum\limits_{m\leq 0}\int\limits_{S^mA_m}Z^\a\log Z d\phi_0&=&\sum\limits_{m\leq 0}\int\limits_{S^mA_m}e^{-(1-\a)\log Z}\log Z d\Lambda\\
	&\leq&\frac{1}{(1-\a)e}\sum\limits_{m\leq 0}\Lambda(A_m)\leq\frac{1}{(1-\a)e}\left(\Lambda(Q)+\e\right).	
\end{eqnarray*}

Thus, taking the supremum implies the second inequality of (i).

$(ii)$ Observe that, by Lemma \ref{hfpl}, $\omega_\a(A)\geq 0$ for all $A\in\A_0$. Thus, the inequality of (ii) is obvious.

Clearly,
\begin{equation*}
	\Omega^{\a,\b}_\e(Q)<-\sum\limits_{m\leq 0}\int\limits_{S^mA_m}Z^\a\log Z d\phi_0+\frac{1}{(1-\a) e}\left(\Lambda(Q) +\e\right).
\end{equation*}
Hence,
\[\Omega^{\a,\b}_\e(Q)\leq-\Xi^{\a,\b}_{1,\e}(Q)+\frac{1}{(1-\a) e}\left(\Lambda(Q) +\e\right).\]
On the other hand, one readily sees that
\[\sum\limits_{m\leq 0}\omega_\a\left(S^mA_m\right)\geq-\Xi^{\a,\b}_{1,\e}(Q)+\frac{1}{(1-\a) e}\Lambda(Q).\]
Hence,
\[\Omega^{\a,\b}_\e(Q)\geq-\Xi^{\a,\b}_{1,\e}(Q)+\frac{1}{(1-\a) e}\Lambda(Q).\]
Thus, the equality of (ii) follows.

$(iii)$ Since $\omega_\a(A)\geq 0$ for all $A\in\A_0$, it follows, by (i), (ii) and Theorem 16 (ii) (Theorem 4 (ii) in the arXiv version) in \cite{Wer15}, that $\Omega^{\a,\b}$ is a finite, $S$-invariant measure on $\B$, and therefore, $\Xi^{\a,\b}_{1}$ is a $S$-invariant, signed measure on $\B$.

$(iv)$ The first inequality of (iv) is obvious, by the first inequality of Lemma \ref{edl} (ii). 

Now, observe that, by the second inequality of Lemma \ref{edl} (ii), 
\begin{eqnarray*}
&&\sum\limits_{m\leq 0}\int\limits_{S^mA_m}Z^\a\log Z d\phi_0 - \Xi^{\a_0,\b}_{1,\e}(Q)\\
&\leq&\sum\limits_{m\leq 0}\int\limits_{S^mA_m}Z^\a\log Z d\phi_0 - \sum\limits_{m\leq 0}\int\limits_{S^mA_m}Z^{\a_0}\log Z d\phi_0\\
&\leq&(\a-\a_0)\left(\left(\frac{2}{\a_0e}\right)^{2}\sum\limits_{m\leq 0}\phi_0\left(S^mA_m\right)+\left(\frac{2}{(1-\a)e}\right)^{2}\sum\limits_{m\leq 0}\Lambda(A_m)\right)\\
&<&(\a-\a_0)\G^{\a_0,\a}_{2,\e}(Q).\\
\end{eqnarray*}
Thus, taking the supremum and letting $\e\to 0$ implies the second inequality of (iv).
\hfill$\Box$

Also, analogously to Lemma \ref{hmfd} (i), we have the following.

\begin{lemma}\label{ldl}
	Let $0<\a_0<\a\leq 1$, $Q\in\B$ and $\e_0,\e>0$. Let $\d_0,\d>0$ such that $\H^{\a_0,0}_{\d_0}\left(\Lambda,\phi_0\right)(Q) >\H^{\a_0,0}\left(\Lambda,\phi_0\right)(Q)-\e_0$ and $\H^{\a,0}_\d\left(\Lambda,\phi_0\right)(Q)>\H^{\a,0}\left(\Lambda,\phi_0\right)(Q)-\e$. Then
	\begin{eqnarray*}
		(\a-\a_0)\Xi^{\a_0,\a}_{1,\d_0}(Q)-\e_0-\d_0&<&\H^{\a,0}\left(\Lambda,\phi_0\right)(Q)-\H^{\a_0,0}\left(\Lambda,\phi_0\right)(Q)\\
		&<&(\a-\a_0)\Xi^{\a,\a_0}_{1,\d}(Q)+\e+\d.
	\end{eqnarray*}
\end{lemma}	
{\it Proof.}
Let $(A_m)_{m\leq 0}\in\C^{\a,0,1}_{\d_0}(Q)$. Then, by Lemma \ref{edl} (i), \eqref{hmvi} and \eqref{sche},  
\begin{eqnarray*}
&&\left(\a-\a_0\right)\sum\limits_{m\leq 0}\int\limits_{S^mA_m}Z^{\a_0}\log Z d\phi_0\\
&\leq&\sum\limits_{m\leq 0}\int\limits_{S^mA_m}Z^{\a} d\phi_0-\sum\limits_{m\leq 0}\int\limits_{S^mA_m}Z^{\a_0} d\phi_0\\
&<&\H^{\a,0}\left(\Lambda,\phi_0\right)(Q)+\d_0-\H^{\a_0,0}\left(\Lambda,\phi_0\right)(Q)+\e_0.
\end{eqnarray*}
Thus, taking the supremum implies the first inequality.

Now, let $(B_m)_{m\leq 0}\in\C^{\a_0,0,1}_{\d}(Q)$. Then, by \eqref{hmvi}, \eqref{sche} and Lemma \ref{edl} (i),
\begin{eqnarray*}
	&&\H^{\a,0}\left(\Lambda,\phi_0\right)(Q)-\e-\H^{\a_0,0}\left(\Lambda,\phi_0\right)(Q)-\d\\
	&<&\sum\limits_{m\leq 0}\int\limits_{S^mB_m}Z^{\a} d\phi_0-\sum\limits_{m\leq 0}\int\limits_{S^mB_m}Z^{\a_0} d\phi_0\\
	&\leq&\left(\a-\a_0\right)\sum\limits_{m\leq 0}\int\limits_{S^mB_m}Z^{\a}\log Z d\phi_0\\
	&\leq&\left(\a-\a_0\right)\Xi^{\a,\a_0}_{1,\d}(Q).
\end{eqnarray*}
This proves the second inequality.
\hfill$\Box$

Finally, similarly to $\Psi^{\a,\a}_{1}(Q)$, we are only able to show that the introduced set function is a derivative of $(0,1)\owns\a\longmapsto\H^{\a,0}(\Lambda,\phi_0)(Q)$ from one side, but this time the left one. 

In order to clarify the behavior of the left derivative from the right, we need the following definition.
\begin{Definition}
	Let $Q\in\P(X)$, $0<\a_0\leq\a< 1$ and $\t>0$. Define
	\[\e'_\t(\a_0,\a):=(\a-\a_0)\left(\Xi^{\a,\a_0}_{1,\d_\t(\a_0,\a)}(Q)-\Xi^{\a_0,\a}_{1,\d_\t(\a_0,\a)}(Q)\right)+2(\a-\a_0)^\t+3\d_\t(\a_0,\a).\]
\end{Definition}

\begin{theo}\label{ldt} 
	Let $Q\in\B$. Then the function $(0,1)\owns\b\longmapsto\H^{\b,0}(\Lambda,\phi_0)(Q)$ is left differentiable, and
	\begin{eqnarray*}
	  \left.\frac{d_-}{d_-\b}\H^{\b,0}(\Lambda,\phi_0)(Q)\right|_{\b=\a}&=&\Xi^{\a,\a}_{1}(Q)\\
	  &=&\lim\limits_{\b\to^-\a}\Xi^{\a,\b}_{1}(Q)=\lim\limits_{\b\to^-\a}\Xi^{\b,\b}_{1}(Q)\\
	  &=&\lim\limits_{\b\to^-\a}\Psi^{\a,\b}_{1,\d_\t(\b,\a)}(Q)=\lim\limits_{\b\to^-\a}\Psi^{\b,\b}_{1,\d_\t(\b,\a)}(Q)\\
	  &=&\lim\limits_{\b\to^-\a}\Psi^{\a,\b}_{1}(Q)=\lim\limits_{\b\to^-\a}\Psi^{\b,\b}_{1}(Q)
	\end{eqnarray*}
	for all $0<\a<1$ and $\t>1$ where $d_-/d_-\b$ denotes the left derivative. 
	
	From the right, for every $0<\a_0<0$ and $\t>0$,
	 \[\lim\limits_{\a\to^+\a_0}\Xi^{\a,\a}_{1,\e'_\t(\a_0,\a)}(Q)=\lim_{\a\to^+\a_0}\Xi^{\a_0,\a}_{1,\e'_\t(\a_0,\a)}(Q)=\Xi^{\a_0,\a_0}_{1}(Q),\]
	 and, for every $\t>1$,
	 \begin{eqnarray*}
	   &&\lim\limits_{\a\to^+\a_0}\Xi^{\a,\a}_{1,\d_\t(\a_0,\a)}(Q)=\lim_{\a\to^+\a_0}\Xi^{\a_0,\a}_{1,\d_\t(\a_0,\a)}(Q)\\
	   &=&\lim\limits_{\a\to^+\a_0}\Xi^{\a,\a}_{1}(Q)=\lim_{\a\to^+\a_0}\Xi^{\a_0,\a}_{1}(Q)=\Psi^{\a_0,\a_0}_{1}(Q).
	 \end{eqnarray*}
\end{theo}
{\it Proof.}
Let $0<\a_0<\a<1$, $\t>0$ and $(A_m)_{m\leq 0}\in\C^{\a_0,0,1}_{\d_\t(\a_0,\a)}(Q)$. Then, by Lemma \ref{ldl}, \eqref{sche} and Lemma \ref{edl} (i),
\begin{eqnarray*}
 &&\H^{\a,0}\left(\Lambda,\phi_0\right)(Q) - (\a-\a_0)\Xi^{\a_0,\a}_{1, \d_\t(\a_0,\a)}(Q) + 2\d_\t(\a_0,\a)+(\a-\a_0)^\t\\
 &\geq&\H^{\a_0,0}\left(\Lambda,\phi_0\right)(Q)+ \d_\t(\a_0,\a)\\
 &>&\sum\limits_{m\leq 0}\int\limits_{S^mA_m}Z^{\a_0} d\phi_0\\
 &\geq&\sum\limits_{m\leq 0}\int\limits_{S^mA_m}Z^{\a} d\phi_0-\left(\a-\a_0\right)\sum\limits_{m\leq 0}\int\limits_{S^mA_m}Z^{\a}\log Z d\phi_0\\
 &\geq&\sum\limits_{m\leq 0}\int\limits_{S^mA_m}Z^{\a} d\phi_0-\left(\a-\a_0\right)\Xi^{\a,\a_0}_{1,\d_\t(\a_0,\a)}(Q).
\end{eqnarray*}
Therefore, since, by Lemma \ref{ldl}, $(\a-\a_0)\left(\Xi^{\a,\a_0}_{1,\d_\t(\a_0,\a)}(Q)-\Xi^{\a_0,\a}_{1,\d_\t(\a_0,\a)}(Q)\right)+2(\a-\a_0)^\t+3\d_\t(\a_0,\a)>\d_\t(\a_0,\a)$,
\[(A_m)_{m\leq 0}\in\C^{\a,0,1}_{\e'_\t(\a_0,\a)}(Q).\]
That is
\[\C^{\a_0,0,1}_{\d_\t(\a_0,\a)}(Q)\subset\C^{\a,0,1}_{\e'_\t(\a_0,\a)}(Q).\]
Hence, for every $0<\b\leq 1$,
\begin{equation}\label{xsar}
  \Xi^{\b,\a_0}_{1}(Q)\leq\Xi^{\b,\a_0}_{1,\d_\t(\a_0,\a)}(Q)\leq\Xi^{\b,\a}_{1,\e'_\t(\a_0,\a)}(Q).
\end{equation}
Therefore, by Lemma \ref{ldp} (iv) and Lemma \ref{ldl}, in the case $\b=\a$,
\begin{eqnarray*}
\Xi^{\a,\a}_{1}(Q)-(\a-\a_0)\G^{\a_0,\a}_2(Q)&\leq&\Xi^{\a_0,\a}_{1}(Q)\\
&\leq&\frac{\H^{\a,0}\left(\Lambda,\phi_0\right)(Q)-\H^{\a_0,0}\left(\Lambda,\phi_0\right)(Q)}{\a-\a_0}\\
&\leq&\Xi^{\a,\a_0}_{1}(Q)\\
&\leq&\Xi^{\a,\a}_{1,\e'_\t(\a_0,\a)}(Q).
\end{eqnarray*}
Thus 
\[\lim\limits_{\a_0\to^-\a}\frac{\H^{\a,0}\left(\Lambda,\phi_0\right)(Q)-\H^{\a_0,0}\left(\Lambda,\phi_0\right)(Q)}{\a-\a_0}=\Xi^{\a,\a}_{1}(Q),\]
and
\[\lim\limits_{\a_0\to^-\a}\Xi^{\a,\a_0}_{1}(Q)=\Xi^{\a,\a}_{1}(Q).\]
Since, by Lemma \ref{ldp} (iv),
\[\Xi^{\a,\a_0}_{1}(Q)-(\a-\a_0)\G^{\a_0,\a}_2(Q)\leq\Xi^{\a_0,\a_0}_{1}(Q)\leq\Xi^{\a,\a_0}_{1}(Q),\]
it follows also that
\[\lim\limits_{\a_0\to^-\a}\Xi^{\a_0,\a_0}_{1}(Q)=\Xi^{\a,\a}_{1}(Q).\]
Let $\t>1$. Let us abbreviate
\begin{eqnarray*}
 \eta_\t(\a_0,\a)&:=&(\a-\a_0)\left(\left(\frac{2}{\a_0e}\right)^{2}\left(\Phi(Q)+\d_\t(\a_0,\a)\right) +\left(\frac{2}{(1-\a)e}\right)^{2}\Lambda(X)\right)\\
 &&+\d_\t(\a_0,\a)\frac{1}{\a_0 e}.
\end{eqnarray*}
Then, by Lemma \ref{hmfd} (i) and Lemma \ref{eql} (ii),
\begin{eqnarray*}
	&&\frac{\H^{\a,0}\left(\Lambda,\phi_0\right)(Q)-\H^{\a_0,0}\left(\Lambda,\phi_0\right)(Q)}{\a-\a_0}\\
	&\leq&\Psi^{\a,\a_0}_{1,\d_\t(\a_0,\a)}(Q)+(\a-\a_0)^{\t-1}+\frac{\d_\t(\a_0,\a)}{\a-\a_0}\\
	&\leq&\Psi^{\a_0,\a_0}_{1,\d_\t(\a_0,\a)}(Q)+\eta_\t(\a_0,\a) +(\a-\a_0)^{\t-1}+\frac{\d_\t(\a_0,\a)}{\a-\a_0}.
\end{eqnarray*} 
Thus, since $\Psi^{\a_0,\a_0}_{1,\d_\t(\a_0,\a)}(Q)\leq\Psi^{\a_0,\a_0}_{1}(Q)\leq\Psi^{\a,\a_0}_{1}(Q)\leq\Xi^{\a,\a_0}_{1}(Q)$, this implies the remaining equalities from the left.

Now, let us consider the behavior of the function from the right. Let $\t>0$. Putting $\b=\a_0$ in \eqref{xsar} implies that
\begin{equation}\label{xsal}
\Xi^{\a_0,\a_0}_{1}(Q)\leq\Xi^{\a_0,\a_0}_{1,\d_\t(\a_0,\a)}(Q)\leq\Xi^{\a_0,\a}_{1,\e'_\t(\a_0,\a)}(Q).
\end{equation}
Let $\d>0$. Then, similarly to the proof of \eqref{cir}, one verifies, by \eqref{hmvi}, \eqref{sche}, Lemma \ref{sdpl} (i), and Lemma \ref{hfpl}, that
\begin{equation*}
\C^{\a,0,1}_{\d}(Q)\subset\C^{\a_0,0,1}_{1,(\a-\a_0)\left(\frac{\Lambda(Q)}{(1-\a)e}+\frac{\Phi(Q)+\d}{\a_0e}\right)+\d}(Q),
\end{equation*}
and therefore, for every $0<\b\leq 1$,
\begin{equation*}
\Xi^{\b,\a}_{1,\d}(Q)\leq\Xi^{\b,\a_0}_{1,(\a-\a_0)\left(\frac{\Lambda(Q)}{(1-\a)e}+\frac{\Phi(Q)+\d}{\a_0e}\right)+\d}(Q),
\end{equation*}
which combined with \eqref{xsal} implies that
\begin{equation*}
\Xi^{\a_0,\a_0}_{1}(Q)\leq\Xi^{\a_0,\a}_{1,\e'_\t(\a_0,\a)}(Q)\leq\Xi^{\a_0,\a_0}_{1,(\a-\a_0)\left(\frac{\Lambda(Q)}{(1-\a)e}+\frac{\Phi(Q)+\e'_\t(\a_0,\a)}{\a_0e}\right)+\e'_\t(\a_0,\a)}(Q).
\end{equation*}
Thus
\[\lim_{\a\to^+\a_0}\Xi^{\a_0,\a}_{1,\e'_\t(\a_0,\a)}(Q)=\Xi^{\a_0,\a_0}_{1}(Q),\]
and, by Lemma \ref{ldp} (iv), also
\[\lim_{\a\to^+\a_0}\Xi^{\a,\a}_{1,\e'_\t(\a_0,\a)}(Q)=\Xi^{\a_0,\a_0}_{1}(Q).\]

Finally, let $\t>1$. Then, by Lemma \ref{eql} (i) and Lemma \ref{ldl},
\begin{eqnarray*}
 &&\Psi^{\a,\a}_{1}(Q)-(\a-\a_0)\G^{\a_0,\a}_2(Q)\leq\Psi^{\a_0,\a}_{1}(Q)\leq\Xi^{\a_0,\a}_{1}(Q)\leq\Xi^{\a_0,\a}_{1,\d_\t(\a_0,\a)}(Q)\\&\leq&\frac{\H^{\a,0}\left(\Lambda,\phi_0\right)(Q)-\H^{\a_0,0}\left(\Lambda,\phi_0\right)(Q)}{\a-\a_0}+\frac{\d_\t(\a_0,\a)}{\a-\a_0}+(\a-\a_0)^{\t-1}.
\end{eqnarray*}
Thus, by Theorem \ref{dlpl},
\[\lim_{\a\to^+\a_0}\Xi^{\a_0,\a}_{1}(Q)=\lim_{\a\to^+\a_0}\Xi^{\a_0,\a}_{1,\d_\t(\a_0,\a)}(Q)=\Psi^{\a_0,\a_0}_{1}(Q),\]
which, by Lemma \ref{ldp} (iv), implies the final assertion.
\hfill$\Box$

Now, we are able to give a lower bound for $\H^{\a,0}\left(\Lambda,\phi_0\right)(Q)$ in terms of $\Xi^{\a,\a}_{1}(Q)$.

\begin{cor}\label{HlbX}
	Let $0<\a<1$ and $Q\in\B$ such that $\Lambda(Q)>0$ and $\Xi^{\a,\a}_{1}(Q)>0$. Then
	\begin{equation*}
	   \Lambda(Q)e^{W_{-1}\left(\frac{-(1-\a)\Xi^{\a,\a}_{1}(Q)}{\Lambda(Q)}\right)}\leq\H^{\a,0}\left(\Lambda,\phi_0\right)(Q)
	   \leq\Lambda(Q)e^{W\left(\frac{-(1-\a)\Xi^{\a,\a}_{1}(Q)}{\Lambda(Q)}\right)}
	\end{equation*}
	where $W$ and $W_{-1}$ denote the principal and the lower branch of the Lambert function respectively.
\end{cor}
{\it Proof.} 
By Lemma \ref{pchm} (iv), $\H^{\a,0}\left(\Lambda,\phi_0\right)(Q)>0$. By Theorem \ref{ldt} and the second inequality of Lemma \ref{pchm} (v),
\[\Xi^{\a,\a}_{1}(Q)\leq-\frac{1}{1-\a}\H^{\a,0}\left(\Lambda,\phi_0\right)(Q)\log\frac{\H^{\a,0}\left(\Lambda,\phi_0\right)(Q)}{\Lambda(Q)},\]
which is equivalent to
\[\frac{-(1-\a)\Xi^{\a,\a}_{1}(Q)}{\Lambda(Q)}\geq\frac{-(1-\a)\Xi^{\a,\a}_{1}(Q)}{\H^{\a,0}\left(\Lambda,\phi_0\right)(Q)}e^{-\frac{(1-\a)\Xi^{\a,\a}_{1}(Q)}{\H^{\a,0}\left(\Lambda,\phi_0\right)(Q)}}.\]
That is
\[W_{-1}\left(\frac{-(1-\a)\Xi^{\a,\a}_{1}(Q)}{\Lambda(Q)}\right)\leq\frac{-(1-\a)\Xi^{\a,\a}_{1}(Q)}{\H^{\a,0}\left(\Lambda,\phi_0\right)(Q)} \leq W\left(\frac{-(1-\a)\Xi^{\a,\a}_{1}(Q)}{\Lambda(Q)}\right),\]
which is equivalent to 
\[\frac{-(1-\a)\Xi^{\a,\a}_{1}(Q)}{W_{-1}\left(\frac{-(1-\a)\Xi^{\a,\a}_{1}(Q)}{\Lambda(Q)}\right)}\leq\H^{\a,0}\left(\Lambda,\phi_0\right)(Q)\leq\frac{-(1-\a)\Xi^{\a,\a}_{1}(Q)}{W\left(\frac{-(1-\a)\Xi^{\a,\a}_{1}(Q)}{\Lambda(Q)}\right)},\]
which is the assertion, since $x/W(x)=e^{W(x)}$ and $x/W_{-1}(x)=e^{W_{-1}(x)}$.
\hfill$\Box$

It appears that the construction of the (signed) measure $\Xi^{\a,\b}_{1}$ is measure-theoretically new.
We show now that, for $0<\a<1$, it can be also obtained in the standard way of the dynamical measure theory, given by the inductive construction in Subsection 4.1.2 in \cite{Wer15}.
\begin{Definition}\label{ldcc}
	Let $0<\a< 1$, $0\leq\b\leq 1$, $Q\in\P(X)$ and $\e>0$. Define
	\[\C^{\a,\b,0,1}_{\e}(Q):=\left\{(A_m)_{m\leq 0}\in\C^{\b,0,1}_\e(Q)|\ \sum\limits_{m\leq 0}\omega_\a\left(S^mA_m\right)<\Omega^{\a,\b}(Q)+\e\right\},\]
	\[\Upsilon^{\a,\b}_{1,\e}(Q):=\inf\limits_{(A_m)_{m\leq 0}\in\C^{\a,\b,0,1}_{\e}(Q)}\sum\limits_{m\leq 0}\int\limits_{S^mA_m}Z^\a\log Z d\phi_0\]
	and
	\[\Upsilon^{\a,\b}_{1}(Q):=\lim\limits_{\e\to 0}\Upsilon^{\a,\b}_{1,\e}(Q).\]
\end{Definition}

\begin{lemma}\label{coma}
	Let $Q\in\P(X)$, $0<\a< 1$ and $0\leq\b\leq 1$. Then
	\[\Xi^{\a,\b}_{1}(Q)=\Upsilon^{\a,\b}_{1}(Q).\]
\end{lemma}
{\it Proof.} Obviously,
\[\Xi^{\a,\b}_{1}(Q)\geq\Upsilon^{\a,\b}_{1}(Q).\]
Let $\e>0$ and $(A_m)_{m\leq 0}\in\C^{\a,\b,0,1}_{\e}(Q)$. Then
\[\Omega^{\a,\b}(Q)+\e>-\sum\limits_{m\leq 0}\int\limits_{S^mA_m}Z^\a\log Z d\phi_0+\frac{1}{(1-\a) e}\Lambda(Q).\]
Hence, taking the infimum and letting $\e\to 0$ implies that
\[\Omega^{\a,\b}(Q)\geq-\Upsilon^{\a,\b}_{1}(Q)+\frac{1}{(1-\a) e}\Lambda(Q).\]
Thus, the assertion follows by Lemma \ref{ldp} (ii).
\hfill$\Box$

\subsubsection{The differentiability of $(0,1)\owns\a\longmapsto\H^{\a,0}\left(\Lambda,\phi_0\right)$}

 We have seen, by Lemma \ref{pchm} (vii), that the function $(0,1)\owns\a\longmapsto\H^{\a,0}\left(\Lambda,\phi_0\right)$ is Lipschitz on every closed subinterval, and therefore, it is differentiable almost everywhere. Using the well-known Beppo Levi Theorem for both-sided differentiable functions, as in Subsection \ref{hds}, one can conclude from our results much more.

Let us consider the set of exceptional points.
\begin{Definition}
For $Q\in\B$, define
\[\H_Q:=\left\{\a\in(0,1)|\ \Psi^{\a,\a}_{1}(Q) < \Xi^{\a,\a}_{1}(Q)\right\}.\]
\end{Definition}

It has the following properties.
\begin{lemma}\label{nsp}
	(i) $\H_Q=\emptyset$ for all $Q\in\B$ such that there exists $\a\in[0,1]$ with $\H^{\a,0}\left(\Lambda,\phi_0\right)(Q)=0$.
	
	(ii) $\H_A\subset\H_B$ for all $A,B\in\B$ with $A\subset B$.
	
	(iii) $\H_Q=\H_{S^{-1}Q}$ for all $Q\in\B$.
		
	(iv) $\bigcup_{n\in\N}\H_{Q_n}=\H_{\bigcup_{n\in\N}Q_n}$ for all $(Q_n)_{n\in\N}\subset\B$.
	
	(v) $\H_{\bigcup_{n\in\Z}S^nQ}=\H_Q$ for all $Q\in\B$.
\end{lemma}
{\it Proof.} $(i)$ It is obvious, by Theorem \ref{dlpl} and  Theorem \ref{ldt}, since, by Lemma \ref{pchm} (iv), $\H^{\a,0}\left(\Lambda,\phi_0\right)(Q)=0$ for all $\a\in(0,1)$ for such $Q$.

$(ii)$ Let $A,B\in\B$ with $A\subset B$. Let $\a\in\H_A$. Then, since $\Psi^{\a,\a}_{1}$ and $\Xi^{\a,\a}_{1}$ are finite signed measures on $\B$, by \eqref{dmr},
\[\Psi^{\a,\a}_{1}(B) = \Psi^{\a,\a}_{1}(B\setminus A)+\Psi^{\a,\a}_{1}(A)<\Xi^{\a,\a}_{1}(B\setminus A)+\Xi^{\a,\a}_{1}(A)=\Xi^{\a,\a}_{1}(B).\]
Hence,  $\a\in\H_B$.

$(iii)$ It is obvious, since $\Psi^{\a,\a}_{1}$ and $\Xi^{\a,\a}_{1}$ are $S$-invariant.

$(iv)$ Let $(Q_n)_{n\in\N}\subset\B$. By (ii), we only need to show that $\H_{\bigcup_{n\in\N}Q_n}\subset\bigcup_{n\in\N}\H_{Q_n}$. Set $Q'_1:=Q_1$ and $Q'_n:=Q_n\setminus(Q_{n-1}\cup...\cup Q_1)$ for all $n\geq 2$. Let $\a\in\H_{\bigcup_{n\in\N}Q_n}$. Then
\[0<\Xi^{\a,\a}_{1}\left(\bigcup\limits_{n\in\N}Q'_n\right)-\Psi^{\a,\a}_{1}\left(\bigcup\limits_{n\in\N}Q'_n\right)=\sum\limits_{n\in\N}\left(\Xi^{\a,\a}_{1}\left(Q'_n\right)-\Psi^{\a,\a}_{1}\left(Q'_n\right)\right).\]
Hence, by \eqref{dmr}, there exists $n\in\N$ such that $\a\in\H_{Q'_n}\subset\H_{\bigcup_{n\in\N}Q_n}$, by (ii).

$(v)$ It follows immediately by (iii) and (iv).
\hfill$\Box$

\begin{cor}\label{dc}
	The set $\H_X$ is at most countable, and $(0,1)\setminus\H_Q\owns \a\longmapsto\\ \H^{\a,0}\left(\Lambda,\phi_0\right)(Q)$ is continuously differentiable for all $Q\in\B$.
\end{cor}
{\it Proof.}
The assertion follows from  Theorem \ref{dlpl} and  Theorem \ref{ldt} by the Beppo Levi Theorem (e.g. see \cite{To}, p. 143).
\hfill$\Box$

Also, by Lemma \ref{pchm}, the function $[0,1]\owns\a\longmapsto\H^{\a,0}\left(\Lambda,\phi_0\right)$ is almost convex. Since the left derivative of a convex function can not exceed the right, it is necessary to test whether the almost convexity also reverses inequality \eqref{dmr}. It turns out, as the next proposition shows, that it seems only to impose a restriction on the difference of the derivatives. 

Another important conclusion of the next proposition is that, even at the points where the left derivative is greater than the right, the function does not provide the best lower bound for $\Phi(Q)$ by Lemma \ref{shm} (i). 

\begin{prop}\label{ftp}
	Let $Q\in\B$ and $0<\a<1$. Suppose $\H^{\a,0}\left(\Lambda,\phi_0\right)(Q)>0$. Then
(i)
	\[\Xi^{\a,\a}_{1}(Q)-\Psi^{\a,\a}_{1}(Q)\leq-\frac{\H^{\a,0}\left(\Lambda,\phi_0\right)(Q)}{\a(1-\a)}\log\frac{\H^{\a,0}\left(\Lambda,\phi_0\right)(Q)}{\Phi(Q)^{1-\a}\Lambda(Q)^{\a}},\]

(ii)
	\begin{eqnarray*}
		&&e^{W_{-1}\left(\frac{\a(1-\a)\left(\Psi^{\a,\a}_{1}(Q)-\Xi^{\a,\a}_{1}(Q)\right)}{\Phi(Q)^{1-\a}\Lambda(Q)^{\a}}\right)}\Phi(Q)^{1-\a}\Lambda(Q)^{\a}\leq\H^{\a,0}\left(\Lambda,\phi_0\right)(Q)\\
		&\leq& e^{W\left(\frac{\a(1-\a)\left(\Psi^{\a,\a}_{1}(Q)-\Xi^{\a,\a}_{1}(Q)\right)}{\Phi(Q)^{1-\a}\Lambda(Q)^{\a}}\right)}\Phi(Q)^{1-\a}\Lambda(Q)^{\a}
	\end{eqnarray*}
	where $W_{-1}$ and $W$ denote the lower and the principal branch of the Lambert function respectively with $W_{-1}(0):=-\infty$.
\end{prop}
{\it Proof.}
Note that, by Lemma \ref{pchm} (iv), the hypothesis implies that\\ $\H^{x,0}\left(\Lambda,\phi_0\right)(Q)>0$ for all $x\in[0,1]$. 

$(i)$ Let $0\leq\b<\a<\g\leq 1$. Let $\a<y<1$. Then, by Lemma \ref{pchm} (v),
\begin{eqnarray*}
	&&\frac{1}{\a-\b}\H^{\a,0}\left(\Lambda,\phi_0\right)(Q)\log\frac{\H^{\a,0}\left(\Lambda,\phi_0\right)(Q)}{\H^{\b,y,0}\left(\Lambda,\phi_0\right)(Q)}\\
	&\leq&\frac{\H^{y,0}\left(\Lambda,\phi_0\right)(Q)-\H^{\a,0}\left(\Lambda,\phi_0\right)(Q)}{y-\a}.
\end{eqnarray*}
Hence, by Theorem \ref{dlpl},
\[\frac{1}{\a-\b}\H^{\a,0}\left(\Lambda,\phi_0\right)(Q)\log\frac{\H^{\a,0}\left(\Lambda,\phi_0\right)(Q)}{\liminf\limits_{y\to^+\a}\H^{\b,y,0}\left(\Lambda,\phi_0\right)(Q)}\leq\Psi^{\a,\a}_{1}(Q).\]
That is
\[\log\H^{\a,0}\left(\Lambda,\phi_0\right)(Q)\leq\frac{\a-\b}{\H^{\a,0}\left(\Lambda,\phi_0\right)(Q)}\Psi^{\a,\a}_{1}(Q)+\log\liminf\limits_{y\to^+\a}\H^{\b,y,0}\left(\Lambda,\phi_0\right)(Q).\]

Now, let $0<x<\a$. Then, by Lemma \ref{pchm} (v),
\begin{eqnarray*}
	&&\frac{\H^{\a,0}\left(\Lambda,\phi_0\right)(Q)-\H^{x,0}\left(\Lambda,\phi_0\right)(Q)}{\a-x}\\
	&\leq&\frac{1}{\g-x}\H^{\a,0}\left(\Lambda,\phi_0\right)(Q)\log\frac{\H^{\g,x,0}\left(\Lambda,\phi_0\right)(Q)}{\H^{\a,0}\left(\Lambda,\phi_0\right)(Q)}.
\end{eqnarray*}
Therefore, by Theorem \ref{ldt},
\[\Xi^{\a,\a}_{1}(Q)\leq\frac{1}{\g-\a}\H^{\a,0}\left(\Lambda,\phi_0\right)(Q)\log\frac{\liminf\limits_{x\to^-\a}\H^{\g,x,0}\left(\Lambda,\phi_0\right)(Q)}{\H^{\a,0}\left(\Lambda,\phi_0\right)(Q)}.\]
That is 
\[\log\H^{\a,0}\left(\Lambda,\phi_0\right)(Q)\leq\log\liminf\limits_{x\to^-\a}\H^{\g,x,0}\left(\Lambda,\phi_0\right)(Q)-\frac{\g-\a}{\H^{\a,0}\left(\Lambda,\phi_0\right)(Q)}\Xi^{\a,\a}_{1}(Q).\]
Therefore, for $\tau:=(\a-\b)/(\g-\b)$,
\begin{eqnarray*}
	&&\log\H^{\a,0}\left(\Lambda,\phi_0\right)(Q)\\
	&\leq&\frac{(\a-\b)(\g-\a)}{(\g-\b)\H^{\a,0}\left(\Lambda,\phi_0\right)(Q)}\left(\Psi^{\a,\a}_{1}(Q)-\Xi^{\a,\a}_{1}(Q)\right)\\
	&&+\log\left(\liminf\limits_{x\to^-\a}\H^{\g,x,0}\left(\Lambda,\phi_0\right)(Q)^\tau\liminf\limits_{y\to^+\a}\H^{\b,y,0}\left(\Lambda,\phi_0\right)(Q)^{1-\tau}\right).
\end{eqnarray*}
Hence,
\begin{eqnarray*}
	\H^{\a,0}\left(\Lambda,\phi_0\right)(Q)&\leq& e^{\frac{(\a-\b)(\g-\a)}{(\g-\b)\H^{\a,0}\left(\Lambda,\phi_0\right)(Q)}\left(\Psi^{\a,\a}_{1}(Q)-\Xi^{\a,\a}_{1}(Q)\right)}\\
	&&\times\liminf\limits_{x\to^-\a}\H^{\g,x,0}\left(\Lambda,\phi_0\right)(Q)^\tau\liminf\limits_{y\to^+\a}\H^{\b,y,0}\left(\Lambda,\phi_0\right)(Q)^{1-\tau}.
\end{eqnarray*}
Thus, setting $\b=0$ and $\g=1$ gives
\begin{equation*}
	\H^{\a,0}\left(\Lambda,\phi_0\right)(Q)\leq e^{\frac{\a(1-\a)\left(\Psi^{\a,\a}_{1}(Q)-\Xi^{\a,\a}_{1}(Q)\right)}{\H^{\a,0}\left(\Lambda,\phi_0\right)(Q)}}\Phi(Q)^{1-\a}\Lambda(Q)^{\a},
\end{equation*}
which is equivalent to (i).

$(ii)$ Obviously, it only needs to be proved when $\Psi^{\a,\a}_{1}(Q)<\Xi^{\a,\a}_{1}(Q)$, in which case it follows the same way as Corollary \ref{HlbX}.
\hfill$\Box$

\section{Lower bounds for $\Phi$ via the DDMs arising from the Hellinger integral $\J_{\a}\left(\Lambda,\phi_0\right)$}\label{tolb}

Motivated by Proposition \ref{ftp} (ii), we now introduce another DDM arising from the Hellinger integral which naturally suggests itself as the greatest one for the purpose of obtaining a lower bound for $\Phi$ by means of the logic of Lemma \ref{hmlb} (i).

\begin{Definition}
	Let $\a\geq 0$, $Q\in\P(X)$ and $\e>0$. Define
	\[\J_{\a,\e}\left(\Lambda,\phi_0\right)(Q):=\sup\limits_{(A_m)_{m\leq 0}\in\C^{0,1}_{\e}(Q)}\sum\limits_{m\leq 0}\int\limits_{S^mA_m}Z^\a d\phi_0\ \ \ \mbox{ and}\]
	\[\J_{\a}\left(\Lambda,\phi_0\right)(Q):=\lim\limits_{\e\to 0}\J_{\a,\e}\left(\Lambda,\phi_0\right)(Q).\]
\end{Definition}

Obviously,  by \eqref{hmvi}, $\J_{0}\left(\Lambda,\phi_0\right)(Q)=\Phi(Q)$, $\J_{1}\left(\Lambda,\phi_0\right)(Q)=\Lambda(Q)$, and\\ $\H^{\a,0}\left(\Lambda,\phi_0\right)(Q)\leq\J_{\a}\left(\Lambda,\phi_0\right)(Q)$ for all $\a\geq 0$. In order to prove that the latter is also a finite measure for some parameter values, we need the following definition.

\begin{Definition}
	Let $0<\a\leq 1$, $Q\in\P(X)$ and $\e>0$. Define
	\[\mathcal{N}_{\a,\e}(Q):=\inf\limits_{(A_m)_{m\leq 0}\in\C^{0,1}_{\e}(Q)}\sum\limits_{m\leq 0}\int\limits_{S^mA_m}\left(\a Z+1-\a -Z^\a\right) d\phi_0\ \ \ \mbox{ and}\]
	\[\mathcal{N}_{\a}(Q):=\lim\limits_{\e\to 0}\mathcal{N}_{\a,\e}(Q).\]
\end{Definition}

Since
$Z^\a\leq 1 +\a(Z-1)$, it follows, by Theorem 16 (ii) (Theorem 4 (ii) in the arXiv version) in \cite{Wer15}, that $\mathcal{N}_{\a}$ is a $S$-invariant measure on $\B$.

\begin{lemma}\label{tlb}
	(i) For every $0\leq\a\leq 1$,
	\begin{equation*}
	\J_{\a}\left(\Lambda,\phi_0\right)(Q)\leq\Phi(Q)^{1-\a}\Lambda(Q)^{\a}\ \ \ \mbox{ for all }Q\in\B.
	\end{equation*}
	
	(ii) Let $0<\a\leq 1$. Then
	\[\mathcal{N}_{\a}(Q)=\a\Lambda(Q)+(1-\a)\Phi(Q)-\J_{\a}\left(\Lambda,\phi_0\right)(Q)\ \ \ \mbox{ for all }Q\in\B.\]
	
	(iii) $\J_{\a}\left(\Lambda,\phi_0\right)$ is a finite, $S$-invariant measure on $\B$ for all $\a\in[0,1]$.
\end{lemma}
{\it Proof.} Let $Q\in B$, $\e>0$  and $(A_m)_{m\leq 0}\in\C^{0,1}_\e(Q)$.

$(i)$ Observe that, by \eqref{hmvi}, the same way as in Lemma \ref{hmlb} (i),
\begin{eqnarray}\label{Jbi}
   \sum\limits_{m\leq 0}\int\limits_{S^mA_m}Z^{\a} d\phi_0\leq\left(\Phi(Q)+\e\right)^{1-\a}\left(\Lambda(Q)+\e\right)^{\a}.
\end{eqnarray}
Thus, the assertion follows. 

$(ii)$ Now, by \eqref{hmvi},
\[\mathcal{N}_{\a,\e}(Q)\leq\a(\Lambda(Q)+\e)+(1-\a)\left(\Phi(Q)+\e\right)-\sum\limits_{m\leq 0}\int\limits_{S^mA_m}Z^\a d\phi_0.\]
Hence,
\[\mathcal{N}_{\a}(Q)\leq\a\Lambda(Q)+(1-\a)\Phi(Q)-\J_{\a}\left(\Lambda,\phi_0\right)(Q).\]
On the other hand,
\[\sum\limits_{m\leq 0}\int\limits_{S^mA_m}\left(\a Z+1-\a -Z^\a\right) d\phi_0\geq\a\Lambda(Q)+(1-\a)\Phi(Q)-\J_{\a,\e}\left(\Lambda,\phi_0\right)(Q).\]
Thus, (ii) follows.

$(iii)$ It follows immediately from (i) and (ii).
\hfill$\Box$

\begin{Remark}\label{Jccr}
   Observe that, by Lemma \ref{tlb} (ii), for $0\leq\a\leq 1$, $\J_{\a}\left(\Lambda,\phi_0\right)$ can be also obtained as a limit of an outer measure approximation by imposing an additional condition on the set of covers, the same way as in Lemma \ref{coma}.
\end{Remark}

\subsection{The regularity of $\a\longmapsto\J_{\a}\left(\Lambda,\phi_0\right)$}

Having observed an improvement of the regularity of the dependence of a DDM arising from the Hellinger integral on the parameter after the restriction of the set of covers with an additional condition (Lemma \ref{pchm}), one might expect a further improvement of the regularity of $\a\longmapsto\J_{\a}\left(\Lambda,\phi_0\right)$ in view of Remark \ref{Jccr}.

\subsubsection{The log-convexity of $[0,1]\owns\a\longmapsto\J_{\a}\left(\Lambda,\phi_0\right)$}

We show now that in fact, in contrast to $[0,1]\owns\a\longmapsto\H^{\a,0}\left(\Lambda,\phi_0\right)$ (compare with Lemma \ref{pchm} (i)), the new function has a very strong regularity property - it is logarithmically convex. (Recall that a convex function on a closed interval always has its one-sided derivatives in the interior, which are non-decreasing and can disagree only on an at most countable set (which still can be dense though).)

The logarithmic almost convexity of the function $\a\longmapsto\H^{\a,0}\left(\Lambda,\phi_0\right)$ can also be expressed in terms of $\J_{\a}\left(\Lambda,\phi_0\right)$.

\begin{lemma}\label{Jcl}
	 Let $Q\in\P(X)$ and $0\leq\b\leq\a_0\leq\a\leq 1$ such that $\a\neq\b$. \\
	 (i)
	\[\J_{\a_0}(\Lambda,\phi_0)(Q)\leq{\J_{\b}(\Lambda,\phi_0)(Q)}^{1-\frac{\a_0-\b}{\a-\b}}{\J_{\a}\left(\Lambda,\phi_0\right)(Q)}^{\frac{\a_0-\b}{\a-\b}}.\]
 	(ii)
	\[\H^{\a_0,0}(\Lambda,\phi_0)(Q)\leq{\J_{\b}(\Lambda,\phi_0)(Q)}^{1-\frac{\a_0-\b}{\a-\b}}{\H^{\a,0}\left(\Lambda,\phi_0\right)(Q)}^{\frac{\a_0-\b}{\a-\b}}, \mbox{ and}\]
	\[\H^{\a_0,0}(\Lambda,\phi_0)(Q)\leq{\H^{\b,0}(\Lambda,\phi_0)(Q)}^{1-\frac{\a_0-\b}{\a-\b}}{\J_{\a}\left(\Lambda,\phi_0\right)(Q)}^{\frac{\a_0-\b}{\a-\b}}.\]	
\end{lemma}
{\it Proof.}
Let $\e>0$ and $(A_m)_{m\leq 0}\in\C^{0,1}_\e(Q)$. Let $\tau:=(\a_0-\b)/(\a-\b)$. 

$(i)$ Obviously, the inequality is correct for $\a_0=0$. Let $\a_0>0$. Then we also can assume that $\a_0<\a$. In this case, by \eqref{Jbi} and \eqref{hce},
\[\sum\limits_{m\leq 0}\int\limits_{S^mA_m}Z^{\a_0} d\phi_0\leq{\J_{\b,\e}(\Lambda,\phi_0)(Q)}^{1-\tau}{\J_{\a,\e}\left(\Lambda,\phi_0\right)(Q)}^{\tau}.\]
Thus, taking the supremum and letting $\e\to 0$ implies (i).

$(ii)$ It follows the same way as (i) by \eqref{sche}.
\hfill$\Box$

\begin{Remark}\label{lacr}
  Lemma \ref{Jcl}, clearly, suggests the following definition. A function $f:[a,b]\longrightarrow[0,\infty)$ is {\em logarithmically almost convex} iff there exists a logarithmically convex function $g:[a,b]\longrightarrow[0,\infty)$ with $g(a)=f(a)$ and $g(b)=f(b)$ such that 
  \[f(\a_0)\leq\min \left\{{g(\b)}^{1-\frac{\a_0-\b}{\a-\b}}{f(\a)}^{\frac{\a_0-\b}{\a-\b}},{f(\b)}^{1-\frac{\a_0-\b}{\a-\b}}{g(\a)}^{\frac{\a_0-\b}{\a-\b}}\right\}\]
  for all $a\leq\b\leq\a_0\leq\a\leq b$ such that $\a\neq\b$. This raises many questions on properties of such functions and the relation to other notions of almost, approximate and quasi convexity appearing in literature. In particular, the open questions related to this article are the following. Suppose $f$ is logarithmically almost convex with a corresponding logarithmically convex function $g$. Does $f$ always have the one-sided derivatives? Is there a relation of its non-differentiability points to those of $g$? Of course, clarifying them first would have been helpful, but it, probably, would lead us too far aside from our current goal. 
\end{Remark}

\subsubsection{The left derivative of $(0,1)\owns\a\longmapsto\J_{\a}\left(\Lambda,\phi_0\right)$}

Now, we are going to show that the following defines the left derivative of the function (compare with the left derivative of $(0,1)\owns\a\longmapsto\H^{\a,0}\left(\Lambda,\phi_0\right)$, Definition \ref{hldd}).

\begin{Definition}\label{Jld}
	Let $0<\a\leq 1$, $Q\in\P(X)$ and $\e>0$. Define
	\[\F^{\a,0,1}_{\e}(Q):=\left\{(A_m)_{m\leq 0}\in\C^{0,1}_\e(Q)|\ \sum\limits_{m\leq 0}\int\limits_{S^mA_m}Z^\a d\phi_0>\J_{\a}\left(\Lambda,\phi_0\right)(Q)-\e\right\},\]
	\[\T_{\a,\e}(Q):=\inf\limits_{(A_m)_{m\leq 0}\in\F^{\a,0,1}_{\e}(Q)}\sum\limits_{m\leq 0}\int\limits_{S^mA_m}Z^\a\log Z d\phi_0\ \ \ \mbox{ and}\]
	\[\T_\a(Q):=\lim\limits_{\e\to 0}\T_{\a,\e}(Q).\]
\end{Definition}

Despite the fact that the construction of $\T_\a$ is new, we show now that it is still in the realm of the dynamical measure theory developed in \cite{Wer15}. 

\begin{Definition}
	Let $0<\a\leq 1$, $Q\in\P(X)$ and $\e>0$. Let $\F'^{\a,0,1}_{\e}(Q)$ be the set of all $(A_m)_{m\leq 0}\in\C^{0,1}_\e(Q)$ such that 
\[\sum\limits_{m\leq 0}\int\limits_{S^mA_m}\left(\a Z+1-\a-Z^\a\right) d\phi_0<\mathcal{N}_{\a}(Q)+\e.\]
Define
\[\T'_{\a,\e}(Q):=\inf\limits_{(A_m)_{m\leq 0}\in\F'^{\a,0,1}_{\e}(Q)}\sum\limits_{m\leq 0}\int\limits_{S^mA_m}Z^\a\log Z d\phi_0\ \ \ \mbox{ and}\]
\[\T'_\a(Q):=\lim\limits_{\e\to 0}\T'_{\a,\e}(Q).\]
\end{Definition}

Then the construction of $\T'_\a$ is a standard one in the dynamical measure theory (elaborated in Subsection 4.1.2 in \cite{Wer15}). Therefore, the same way as in the proof of Lemma \ref{eml} (iv) and (v), one sees that $\T'_\a$ is a $S$-invariant, signed measure on $B$ for all $0<\a<1$. We show now that it coincides with $\T_\a$ on $\B$.

\begin{lemma}\label{Jldm}
	Let $Q\in B$.
	
	(i) For every $0<\a< 1$,
	\begin{equation*}
	\frac{\J_{\a}\left(\Lambda,\phi_0\right)(Q)-\Phi(Q)}{\a}\leq\T_{\a}(Q)\leq\frac{\Lambda(Q)-\J_{\a}\left(\Lambda,\phi_0\right)(Q)}{1-\a}.
	\end{equation*}
	
	(ii) For every $0<\a\leq 1$,
	  \[\T_{\a}(Q)=\T'_{\a}(Q).\]
	
	(iii) $\T_{\a}$ is a $S$-invariant, signed measure on $\B$ for all $0<\a< 1$.
\end{lemma}
{\it Proof.} Let $0<\a\leq 1$, $\e>0$ and $(A_m)_{m\leq 0}\in\F^{\a,0,1}_{\e}(Q)$.

$(i)$ Let $\a<1$. Since $(Z^\a-1)/\a\leq Z^\a\log Z\leq (Z-Z^\a)/(1-\a)$, 
\begin{eqnarray*}
  \frac{\J_{\a}\left(\Lambda,\phi_0\right)(Q)-\Phi(Q)-2\e}{\a}&<&\sum\limits_{m\leq 0}\int\limits_{S^mA_m}Z^\a\log Z d\phi_0\\
  &<&\frac{\Lambda(Q)-\J_{\a}\left(\Lambda,\phi_0\right)(Q)+2\e}{1-\a}.	
\end{eqnarray*}
This implies the assertion of (i). 

$(ii)$ Let $(A_m)_{m\leq 0}\in\F'^{\a,0,1}_{\e}(Q)$. Then
\[\mathcal{N}_{\a}(Q)+\e>\a\Lambda(Q)+(1-\a)\Phi(Q)-\sum\limits_{m\leq 0}\int\limits_{S^mA_m}Z^\a d\phi_0.\]
Hence, by Lemma \ref{tlb} (ii), $(A_m)_{m\leq 0}\in\F^{\a,0,1}_{\e}(Q)$. That is 
\begin{equation}\label{lddi}
  \F'^{\a,0,1}_{\e}(Q)\subset\F^{\a,0,1}_{\e}(Q).
\end{equation}
 Therefore,
\begin{equation*}
 \T'_\a(Q)\geq\T_\a(Q).
\end{equation*}
Now, let $(B_m)_{m\leq 0}\in\F^{\a,0,1}_{\e}(Q)$. Then, by Lemma \ref{tlb} (ii),
\begin{eqnarray*}
  &&\sum\limits_{m\leq 0}\int\limits_{S^mB_m}\left(\a Z+1-\a-Z^\a\right) d\phi_0\\
  &<&\a(\Lambda(Q)+\e)+(1-\a)\left(\Phi(Q)+\e\right)-\J_{\a}\left(\Lambda,\phi_0\right)(Q)+\e\\
  &=&\mathcal{N}_\a(Q)+2\e.
\end{eqnarray*}
Hence, $(B_m)_{m\leq 0}\in\F'^{\a,0,1}_{2\e}(Q)$, i.e. 
\begin{equation}\label{ldci}
 \F^{\a,0,1}_{\e}(Q)\subset\F'^{\a,0,1}_{2\e}(Q).
\end{equation}
 Therefore,
\[\T_{\a,\e}(Q)\geq\T'_{\a,2\e}(Q).\]
Thus 
\[\T_{\a}(Q)\geq\T'_{\a}(Q),\]
which remained to prove in (ii).

$(iii)$ It follows immediately from (ii).
\hfill$\Box$

The next lemma shows that $\T_{\a}$ is a good candidate for a derivative of $\J_{\a}\left(\Lambda,\phi_0\right)$. 

\begin{lemma}\label{Jldl}
  Let $0<\a_0<\a\leq 1$, $Q\in\B$ and $\e_0,\e>0$. Let $\d_0,\d>0$ such that $\J_{\a_0,\d_0}\left(\Lambda,\phi_0\right)(Q)<\J_{\a_0}\left(\Lambda,\phi_0\right)(Q)+\e_0$ and $\J_{\a,\d}\left(\Lambda,\phi_0\right)(Q)<\J_{\a}\left(\Lambda,\phi_0\right)(Q)+\e$. Then
  \begin{eqnarray*}
  	(\a-\a_0)\T_{\a_0,\d}(Q)-\e-\d&<&\J_{\a}\left(\Lambda,\phi_0\right)(Q)-\J_{\a_0}\left(\Lambda,\phi_0\right)(Q)\\
  	&<&(\a-\a_0)\T_{\a,\d_0}(Q)+\e_0+\d_0.
  \end{eqnarray*}
\end{lemma}
{\it Proof.} 
 Let $(A_m)_{m\leq 0}\in\F^{\a_0,0,1}_{\d}(Q)$. Then
 \begin{eqnarray*}
   (\a-\a_0)\T_{\a_0,\d}(Q)&\leq&(\a-\a_0)\sum\limits_{m\leq 0}\int\limits_{S^mA_m}Z^{\a_0}\log Z d\phi_0\\
   &\leq&\sum\limits_{m\leq 0}\int\limits_{S^mA_m}Z^{\a} d\phi_0-\sum\limits_{m\leq 0}\int\limits_{S^mA_m}Z^{\a_0} d\phi_0\\
   &<&\J_{\a}\left(\Lambda,\phi_0\right)(Q)+\e-\J_{\a_0}\left(\Lambda,\phi_0\right)(Q)+\d.
 \end{eqnarray*}
 This gives the first inequality.
 
 Let $(B_m)_{m\leq 0}\in\F^{\a,0,1}_{\d_0}(Q)$. Then
 \begin{eqnarray*}
   &&\J_{\a}\left(\Lambda,\phi_0\right)(Q)-\d_0-\J_{\a_0}\left(\Lambda,\phi_0\right)(Q)-\e_0\\
   &\leq&\sum\limits_{m\leq 0}\int\limits_{S^mB_m}Z^{\a} d\phi_0-\sum\limits_{m\leq 0}\int\limits_{S^mB_m}Z^{\a_0} d\phi_0\leq(\a-\a_0)\sum\limits_{m\leq 0}\int\limits_{S^mB_m}Z^{\a}\log Z d\phi_0.
 \end{eqnarray*}
 Hence, taking the infimum gives the second inequality.
\hfill$\Box$

Finally, we are only able to show that $\T_{\a}$ is in fact the left derivative of $\J_{\a}\left(\Lambda,\phi_0\right)$. 

\begin{theo}\label{Jldt}
	Let $Q\in\B$. Then 
	\begin{eqnarray*}
		\left.\frac{d_-}{d_-x}\J_x(\Lambda,\phi_0)(Q)\right|_{x=\a}=\T_{\a}(Q)=\lim\limits_{x\to^-\a}\T_{x}(Q)
	\end{eqnarray*}
	for all $0<\a<1$ where $d_-/d_-x$ denotes the left derivative.
\end{theo}
{\it Proof.} Let $0<\a_0<\a<1$, $\e>0$ and $(A_m)_{m\leq 0}\in\F^{\a_0,0,1}_{\e}(Q)$. By Lemma \ref{Jldl}, Lemma \ref{edl} (i) and Lemma \ref{Jldm} (i),
\begin{eqnarray*}
  &&\J_{\a}\left(\Lambda,\phi_0\right)(Q)-\e\\
  &\leq&\J_{\a_0}\left(\Lambda,\phi_0\right)(Q)+(\a-\a_0)\T_{\a}(Q)-\e\\
  &<&\sum\limits_{m\leq 0}\int\limits_{S^mA_m}Z^{\a_0} d\phi_0+(\a-\a_0)\T_{\a}(Q)\\
  &\leq&\sum\limits_{m\leq 0}\int\limits_{S^mA_m}Z^{\a} d\phi_0+(\a-\a_0)\left(\T_{\a}(Q)-\sum\limits_{m\leq 0}\int\limits_{S^mA_m}Z^{\a_0}\log Z d\phi_0\right)\\
  &\leq&\sum\limits_{m\leq 0}\int\limits_{S^mA_m}Z^{\a} d\phi_0+(\a-\a_0)\left(\frac{\Lambda(Q)}{1-\a}+\frac{\Phi(Q)+\e}{\a_0}\right).
\end{eqnarray*}
Hence, $(A_m)_{m\leq 0}\in\F^{\a,0,1}_{(\a-\a_0)\left(\Lambda(Q)/(1-\a)+(\Phi(Q)+\e)/\a_0\right)+\e}(Q)$. That is
\[\F^{\a_0,0,1}_{\e}(Q)\subset\F^{\a,0,1}_{(\a-\a_0)\left(\frac{\Lambda(Q)}{1-\a}+\frac{\Phi(Q)+\e}{\a_0}\right)+\e}(Q).\]
Therefore, by Lemma \ref{edl} (ii) and Lemma \ref{hfpl},
\[ \T_{\a,(\a-\a_0)\left(\frac{\Lambda(Q)}{1-\a}+\frac{\Phi(Q)+\e}{\a_0}\right)+\e}(Q)-\sum\limits_{m\leq 0}\int\limits_{S^mA_m}Z^{\a_0}\log Z d\phi_0\leq(\a-\a_0)\G^{\a_0,\a}_{2,\e}(Q).\]
Hence,
\[ \T_{\a,(\a-\a_0)\left(\frac{\Lambda(Q)}{1-\a}+\frac{\Phi(Q)+\e}{\a_0}\right)+\e}(Q)-(\a-\a_0)\G^{\a_0,\a}_{2,\e}(Q)\leq\T_{\a_0,\e}(Q)\leq\T_{\a_0}(Q).\]
Therefore, by Lemma \ref{Jldl},
\begin{eqnarray*}
&&\T_{\a,(\a-\a_0)\left(\frac{\Lambda(Q)}{1-\a}+\frac{\Phi(Q)+\e}{\a_0}\right)+\e}(Q)-(\a-\a_0)\G^{\a_0,\a}_{2,\e}(Q)\leq\T_{\a_0}(Q)\\
&\leq&\frac{\J_{\a}\left(\Lambda,\phi_0\right)(Q)-\J_{\a_0}\left(\Lambda,\phi_0\right)(Q)}{\a-\a_0}\leq\T_{\a}(Q).
\end{eqnarray*}
Thus, setting $\e:=\a-\a_0$ and letting $\a_0\to\a$ implies the assertion.
\hfill$\Box$

\subsubsection{The right derivative of $(0,1)\owns\a\longmapsto\J_{\a}\left(\Lambda,\phi_0\right)$}

Next, we are going to obtain the right derivative of the function, following the recipe from Subsection \ref{ldsh}.
\begin{Definition}\label{Jrd}
	Let $0<\a\leq 1$, $Q\in\P(X)$ and $\e>0$. Define
	\[\Pi_{\a,\e}(Q):=\sup\limits_{(A_m)_{m\leq 0}\in\F^{\a,0,1}_{\e}(Q)}\sum\limits_{m\leq 0}\int\limits_{S^mA_m}Z^\a\log Z d\phi_0\ \ \ \mbox{ and}\]
	\[\Pi_\a(Q):=\lim\limits_{\e\to 0}\Pi_{\a,\e}(Q).\]
\end{Definition}

We will show that this construction is still covered by the dynamical measure theory \cite{Wer15} for all $0<\a<1$.

\begin{lemma}\label{Jrdp}
	Let $0<\a< 1$ and $Q\in\B$.
	
	(i)
	\begin{equation*}
	\frac{\J_{\a}\left(\Lambda,\phi_0\right)(Q)-\Phi(Q)}{\a}\leq\Pi_{\a}(Q)\leq\frac{\Lambda(Q)-\J_{\a}\left(\Lambda,\phi_0\right)(Q)}{1-\a}.
	\end{equation*}
	
	(ii)  $\Pi_{\a}$ is a $S$-invariant, signed measure on $\B$.
\end{lemma}
{\it Proof.} 
$(i)$  The proof is the same as that of Lemma \ref{Jldm} (i).

$(ii)$ Let $\e>0$. Define
 \[\Omega'_{\a,\e}(Q):=\inf\limits_{(A_m)_{m\leq 0}\in\F'^{\a,0,1}_{\e}(Q)}\sum\limits_{m\leq 0}\omega_\a\left(S^mA_m\right)\]
 and
 \[\Omega'_{\a}(Q):=\lim\limits_{\e\to 0}\Omega'_{\a,\e}(Q).\]
Then, as in the proof of Lemma \ref{ldp} (ii), $\Omega'_{\a}$ is a finite measure on $\B$. 

Now, observe that, by \eqref{lddi},
\begin{eqnarray*}
&&\Omega'_{\a,\e}(Q)\\
&\geq&\inf\limits_{(A_m)_{m\leq 0}\in\F^{\a,0,1}_{\e}(Q)}\left\{\frac{1}{(1-\a) e}\sum\limits_{m\leq 0}\Lambda(A_m)-\sum\limits_{m\leq 0}\int\limits_{S^mA_m}Z^\a\log Z d\phi_0\right\}\\
&\geq&\frac{1}{(1-\a) e}\Lambda(Q)-\Pi_{\a,\e}(Q).
\end{eqnarray*}
Hence,
\[\Omega'_{\a}(Q)\geq\frac{1}{(1-\a) e}\Lambda(Q)-\Pi_{\a}(Q).\]
On the other hand, by \eqref{ldci},
\begin{eqnarray*}
	&&\Omega'_{\a,2\e}(Q)\\
	&\leq&\inf\limits_{(A_m)_{m\leq 0}\in\F^{\a,0,1}_{\e}(Q)}\left\{\frac{1}{(1-\a) e}\sum\limits_{m\leq 0}\Lambda(A_m)-\sum\limits_{m\leq 0}\int\limits_{S^mA_m}Z^\a\log Z d\phi_0\right\}\\
	&\leq&\frac{1}{(1-\a) e}\left(\Lambda(Q)+\e\right)-\Pi_{\a,\e}(Q).
\end{eqnarray*}
This implies the converse inequality, and therefore,
\begin{equation}\label{rdme}
  \Omega'_{\a}(Q)=\frac{1}{(1-\a) e}\Lambda(Q)-\Pi_{\a}(Q),	
\end{equation}
which implies the assertion.
\hfill$\Box$

Observe that, by \eqref{rdme}, one can obtain $\Pi_{\a}$ also via an outer measure approximation for all $0<\a<1$, the same way as in Definition \ref{ldcc}.

Similarly to Lemma \ref{Jldl}, we have the following.
\begin{lemma}\label{Jrdl}
	Let $0<\a_0<\a\leq 1$, $Q\in\B$ and $\e_0,\e>0$. Let $\d_0,\d>0$ such that $\J_{\a_0,\d_0}\left(\Lambda,\phi_0\right)(Q)<\J_{\a_0}\left(\Lambda,\phi_0\right)(Q)+\e_0$ and $\J_{\a,\d}\left(\Lambda,\phi_0\right)(Q)<\J_{\a}\left(\Lambda,\phi_0\right)(Q)+\e$. Then
	\begin{eqnarray*}
		(\a-\a_0)\Pi_{\a_0,\d}(Q)-\e-\d&<&\J_{\a}\left(\Lambda,\phi_0\right)(Q)-\J_{\a_0}\left(\Lambda,\phi_0\right)(Q)\\
		&<&(\a-\a_0)\Pi_{\a,\d_0}(Q)+\e_0+\d_0.
	\end{eqnarray*}
\end{lemma}
{\it Proof.} 
Let $(A_m)_{m\leq 0}\in\F^{\a_0,0,1}_{\d}(Q)$. Then
\begin{eqnarray*}
	(\a-\a_0)\sum\limits_{m\leq 0}\int\limits_{S^mA_m}Z^{\a_0}\log Z d\phi_0&\leq&\sum\limits_{m\leq 0}\int\limits_{S^mA_m}Z^{\a} d\phi_0-\sum\limits_{m\leq 0}\int\limits_{S^mA_m}Z^{\a_0} d\phi_0\\
	&<&\J_{\a}\left(\Lambda,\phi_0\right)(Q)+\e-\J_{\a_0}\left(\Lambda,\phi_0\right)(Q)+\d.
\end{eqnarray*}
Thus, taking the supremum gives the first inequality.

Let $(B_m)_{m\leq 0}\in\F^{\a,0,1}_{\d_0}(Q)$. Then
\begin{eqnarray*}
	&&\J_{\a}\left(\Lambda,\phi_0\right)(Q)-\d_0-\J_{\a_0}\left(\Lambda,\phi_0\right)(Q)-\e_0\\
	&\leq&\sum\limits_{m\leq 0}\int\limits_{S^mB_m}Z^{\a} d\phi_0-\sum\limits_{m\leq 0}\int\limits_{S^mB_m}Z^{\a_0} d\phi_0\leq(\a-\a_0)\sum\limits_{m\leq 0}\int\limits_{S^mB_m}Z^{\a}\log Z d\phi_0\\
	&\leq&(\a-\a_0)\Pi_{\a,\d_0}(Q),
\end{eqnarray*}
which is the second inequality.
\hfill$\Box$

And again, we are only able to show that $\Pi_x$ is the one-sided derivative of $\J_x(\Lambda,\phi_0)$.

\begin{theo}\label{Jrdt}
	Let $Q\in\B$ and $0<\a<1$. Then
	\begin{eqnarray*}
		\left.\frac{d_+}{d_+x}\J_x(\Lambda,\phi_0)(Q)\right|_{x=\a}=\Pi_{\a}(Q)=\lim\limits_{x\to^+\a}\Pi_{x}(Q)=\lim\limits_{x\to^+\a}\T_{x}(Q)
	\end{eqnarray*}
	 where $d_+/d_+x$ denotes the right derivative. Also,
	\[\lim\limits_{x\to^-\a}\Pi_{x}(Q)=\T_{\a}(Q).\]
\end{theo}
{\it Proof.} Let $0<\a_0<\a<1$, $\e>0$ and $(A_m)_{m\leq 0}\in\F^{\a,0,1}_{\e}(Q)$. Then, by Lemma \ref{Jrdl}, Lemma \ref{edl} (i) and Lemma \ref{Jrdp} (i),
\begin{eqnarray*}
  &&\J_{\a_0}\left(\Lambda,\phi_0\right)(Q)-\e\\
  &\leq&\J_{\a}\left(\Lambda,\phi_0\right)(Q)-(\a-\a_0)\Pi_{\a_0}(Q)-\e\\
  &<&\sum\limits_{m\leq 0}\int\limits_{S^mA_m}Z^{\a} d\phi_0-(\a-\a_0)\Pi_{\a_0}(Q)\\
  &\leq&\sum\limits_{m\leq 0}\int\limits_{S^mA_m}Z^{\a_0} d\phi_0+(\a-\a_0)\left(\sum\limits_{m\leq 0}\int\limits_{S^mA_m}Z^{\a}\log Z d\phi_0-\Pi_{\a_0}(Q)\right)\\
  &\leq&\sum\limits_{m\leq 0}\int\limits_{S^mA_m}Z^{\a_0} d\phi_0+(\a-\a_0)\left(\frac{\Lambda(Q)+\e}{1-\a}+\frac{\Phi(Q)}{\a_0}\right).
\end{eqnarray*}
Hence, $(A_m)_{m\leq 0}\in\F^{\a_0,0,1}_{(\a-\a_0)\left((\Lambda(Q)+\e)/(1-\a)+\Phi(Q)/\a_0\right)+\e}(Q)$.
Therefore, by Lemma \ref{edl} (ii) and Lemma \ref{hfpl},
\[ \sum\limits_{m\leq 0}\int\limits_{S^mA_m}Z^{\a}\log Z d\phi_0-\Pi_{\a_0,(\a-\a_0)\left(\frac{\Lambda(Q)+\e}{1-\a}+\frac{\Phi(Q)}{\a_0}\right)+\e}(Q)\leq(\a-\a_0)\G^{\a_0,\a}_{2,\e}(Q).\]
Hence,
\[\Pi_{\a}(Q)\leq\Pi_{\a,\e}(Q)\leq(\a-\a_0)\G^{\a_0,\a}_{2,\e}(Q)+\Pi_{\a_0,(\a-\a_0)\left(\frac{\Lambda(Q)+\e}{1-\a}+\frac{\Phi(Q)}{\a_0}\right)+\e}(Q).\]
Therefore, by Lemma \ref{Jrdl},
\begin{eqnarray*}
	\Pi_{\a_0}(Q)&\leq&\frac{\J_{\a}\left(\Lambda,\phi_0\right)(Q)-\J_{\a_0}\left(\Lambda,\phi_0\right)(Q)}{\a-\a_0}\leq\Pi_{\a}(Q)\\
	&\leq&(\a-\a_0)\G^{\a_0,\a}_{2,\e}(Q)+\Pi_{\a_0,(\a-\a_0)\left(\frac{\Lambda(Q)+\e}{1-\a}+\frac{\Phi(Q)}{\a_0}\right)+\e}(Q).
\end{eqnarray*}
Thus, setting $\e:=\a-\a_0$ and letting $\a\to\a_0$ implies the first two equalities of the assertion.

Now, let us consider the behavior of the right derivative from the left and the left derivative from the right. By the definitions of $\T_{\a_0}(Q)$ and $\Pi_{\a_0}(Q)$, Lemma \ref{Jldl} and Lemma \ref{Jrdl},
\[\T_{\a_0}(Q)\leq\Pi_{\a_0}(Q)\leq\frac{\J_{\a}\left(\Lambda,\phi_0\right)(Q)-\J_{\a_0}\left(\Lambda,\phi_0\right)(Q)}{\a-\a_0}\leq\T_{\a}(Q)\leq\Pi_{\a}(Q).\]
Thus, the remaining two equalities follow by the above and Theorem \ref{Jldt}.
\hfill$\Box$

Similarly to Corollary \ref{HlbX}, here the right derivative can be used to obtain a lower bound for the function.

\begin{cor}\label{rdlb}
	Let $0<\a<1$ and $Q\in\B$ such that $\Lambda(Q)>0$ and $\Pi_{\a}(Q)> 0$. Then
	\begin{equation*}
	\Lambda(Q)e^{W_{-1}\left(\frac{-(1-\a)\Pi_{\a}(Q)}{\Lambda(Q)}\right)}\leq\J_\a\left(\Lambda,\phi_0\right)(Q)
	\leq\Lambda(Q)e^{W\left(\frac{-(1-\a)\Pi_{\a}(Q)}{\Lambda(Q)}\right)}
	\end{equation*}
	where $W$ and $W_{-1}$ denote the principal and the lower branch of the Lambert function respectively.
\end{cor}
{\it Proof.} The proof is the same as that of Corollary \ref{HlbX} (where, instead of Lemma \ref{pchm} and Theorem \ref{ldt}, one should refer to Lemma \ref{Jcl} and Theorem \ref{Jrdt}).
\hfill$\Box$

\subsubsection{The set of non-differentiability points of  $(0,1)\owns\a\longmapsto\J_{\a}\left(\Lambda,\phi_0\right)$}

Now, let us state the properties of the set of non-differentiability points of $(0,1)\owns\a\longmapsto\J_{\a}\left(\Lambda,\phi_0\right)(Q)$.
\begin{Definition}
	For $Q\in\B$, define
	\[\J_Q:=\left\{\a\in(0,1)|\ \T_{\a}(Q) < \Pi_{\a}(Q)\right\}.\]
\end{Definition}

We already know that $\J_Q$ is at most countable, since $(0,1)\owns\a\longmapsto\J_{\a}\left(\Lambda,\phi_0\right)(Q)$ is convex, by Lemma \ref{Jcl}. Analogously to Lemma \ref{nsp}, it has also the following properties.

\begin{lemma}	 
	(i) $\J_Q=\emptyset$ for all $Q\in\B$ such that there exists $\a\in[0,1]$ with $\J_{\a}\left(\Lambda,\phi_0\right)(Q)=0$.
	
	(ii) $\J_A\subset\J_B$ for all $A,B\in\B$ with $A\subset B$.
	
	(iii) $\J_Q=\J_{S^{-1}Q}$ for all $Q\in\B$.
	
	(iv) $\bigcup_{n\in\N}\J_{Q_n}=\J_{\bigcup_{n\in\N}Q_n}$ for all $(Q_n)_{n\in\N}\subset\B$.
	
	(v) $\J_{\bigcup_{n\in\Z}S^nQ}=\J_Q$ for all $Q\in\B$.
\end{lemma}
{\it Proof.} The proof is similar to that of Lemma \ref{nsp}.
\hfill$\Box$

Clearly, as in Remark \ref{lacr}, arises the question on the relation between $\J_Q$ and $\H_Q$, which we leave open here.

\section{The ergodic case for $\H^{\a,0}\left(\Lambda,\phi_0\right)$ and $\J_{\a}\left(\Lambda,\phi_0\right)$}

We continue the analysis of the case of an ergodic $\Lambda$ started in Subsection \ref{satss}, in terms of the absolute continuity relations.

\begin{prop}\label{ecp}
	Suppose $\Lambda$ is an ergodic probability measure. Let $0\leq\a< 1$. Then the following are equivalent.
	
	(i) $\Lambda\ll \H^{\a,0}\left(\Lambda,\phi_0\right)$ on $\B$,
	
	(ii) $\Lambda\ll \J_{\a}\left(\Lambda,\phi_0\right)$ on $\B$, and
	
	(iii) $Z$ is essentially bounded with respect to $\Lambda$.
\end{prop}
{\it Proof.}
The implications $(iii)\Rightarrow (i)\Rightarrow (ii)$ follow by Corollary \ref{fcem}, since\\ $\H_{\a}\left(\Lambda,\phi_0\right)(Q)\leq \H^{\a,0}\left(\Lambda,\phi_0\right)(Q)\leq \J_{\a}\left(\Lambda,\phi_0\right)(Q)$ for all $Q\in\B$.
 
$(ii)\Rightarrow (iii)$: Suppose (iii) is false. Let $B\in\B$  as constructed in the proof of Corollary \ref{fcem}. Then, by Lemma \ref{tlb} (i),
$\J_{\a}\left(\Lambda,\phi_0\right)(B)=0$, since $\Phi(B)=0$, but this contradicts to (ii), since $\Lambda(B)=1$.\hfill$\Box$

Similarly to Corollary \ref{1hme}, we have the following.

\begin{cor}
	Suppose the hypothesis of Corollary \ref{1hme} is satisfied. Let $Y_{\a}\left(\Lambda,\phi_0\right)$ denote $\H^{\a,0}\left(\Lambda,\phi_0\right)$ or $\J_{\a}\left(\Lambda,\phi_0\right)$ for all $0\leq\a\leq 1$. Then the following are equivalent.
	
	(i) $Z$ is essentially bounded with respect to $\Lambda$.
	
	(ii) For every $0\leq\g\leq 1$,  $Y_{\g}\left(\Lambda,\phi_0\right)(X)>0$ and\\ $Y_{\g}\left(\Lambda,\phi_0\right)(Q)/Y_{\g}\left(\Lambda,\phi_0\right)(X) = \Lambda(Q)$ for all $Q\in\B$.
	
	(iii) There exists $0\leq\g<1$ such that  $Y_{\g}\left(\Lambda,\phi_0\right)(X)>0$ and\\ $Y_{\g}\left(\Lambda,\phi_0\right)(Q)/Y_{\g}\left(\Lambda,\phi_0\right)(X) = \Lambda(Q)$ for all $Q\in\B$.
\end{cor}
{\it Proof.} We prove the case $Y_{\a}\left(\Lambda,\phi_0\right)=\H^{\a,0}\left(\Lambda,\phi_0\right)$, the proof in the case  $Y_{\a}\left(\Lambda,\phi_0\right)=\J_{\a}\left(\Lambda,\phi_0\right)$ is the same.

$(i)\Rightarrow (ii)$: Let $0\leq\g<1$. By Corollary \ref{fcem}, $\H^{\g,0}\left(\Lambda,\phi_0\right)(X)>0$. The relation $\H^{\g,0}\left(\Lambda,\phi_0\right)\ll\Lambda$ follows by Lemma \ref{shm} (iii). Hence, (ii) follows the same way as that of Corollary \ref{1hme}.  

$(ii)\Rightarrow (iii)$ is obvious.

$(iii)\Rightarrow (i)$ follows by the implication (i) $\Rightarrow$ (iii) of Proposition \ref{ecp}.
\hfill$\Box$

\section{Explicit computations}

In this section, for the purpose of computing some DDMs explicitly, in terms of $\Lambda$, we will make the following additional assumption on $\Lambda$.

Let $\Lambda(X)=1$. Furthermore, let $\I$ be an at most countable set, $(\Lambda_i)_{i\in \I}$ be a family of distinct ergodic probability measures on $\B$, and $(\lambda_i)_{i\in \I}\subset(0,1]$ such that 
\[\Lambda(Q)=\sum\limits_{i\in \I} \lambda_i\Lambda_i(Q)\ \ \ \mbox{ for all }Q\in\B.\]
For each $i\in \I$, let $Z_i$ be a measurable version of the Radon-Nikodym derivative $d\Lambda_i/d\phi_0$. One easily sees that
\[Z=\sum\limits_{i\in \I} \lambda_iZ_i\ \ \ \phi_0\mbox{-a.e.}\]

We will need the following well-known general lemma, which we give here with a proof for the purpose of completeness.
\begin{lemma}\label{ispl}
	Let $Q\in\B$ such that $\Lambda(Q\Delta S^{-1}Q)=0$. Then there exists $A\in\A_0$ such that $\Lambda(Q\Delta A)=0$.
\end{lemma}
{\it Proof.} 
Let $n\in\N$ and $(A^n_m)_{m\leq 0}\in\C^1_{2^{-n}}(Q)$. Choose $m_n\leq 0$ such that
\[\Lambda\left(\bigcup\limits_{m\leq 0}A^n_m\setminus \bigcup\limits_{m_n\leq m\leq 0}A^n_m\right)<2^{-n}.\]
Set $A_n:=S^{m_n}(\bigcup_{m_n\leq m\leq 0}A^n_m)$. Obviously, $A_n\in\A_0$. Also, by the hypothesis,
\begin{eqnarray*}
	&&\Lambda\left(Q\Delta A_n\right)\\
	&=&\Lambda\left(Q\Delta S^{-m_n}A_n\right) = \Lambda\left(Q\setminus \bigcup\limits_{m_n\leq m\leq 0}A^n_m\right)+\Lambda\left(\bigcup\limits_{m_n\leq m\leq 0}A^n_m\setminus Q\right)\\
	&\leq& \Lambda\left(\bigcup\limits_{m\leq 0}A^n_m\setminus \bigcup\limits_{m_n\leq m\leq 0}A^n_m\right)+\Lambda\left(\bigcup\limits_{m\leq 0}A^n_m\setminus Q\right)< 2^{-n+1}.
\end{eqnarray*}
Now, define $A:=\bigcup_{k\in\N}\bigcap_{n\geq k}A_n$. Then $A\in\A_0$, and
\begin{eqnarray*}
	\Lambda\left(Q\Delta A\right)&=&\Lambda\left(\bigcap\limits_{k\in\N}\bigcup\limits_{n\geq k}Q\setminus A_n\right) + \Lambda\left(\bigcup\limits_{k\in\N}\bigcap\limits_{n\geq k} A_n\setminus Q\right)\\
	&=&\lim\limits_{k\to\infty}\Lambda\left(\bigcup\limits_{n\geq k}Q\setminus A_n\right) +\lim\limits_{k\to\infty}\Lambda\left(\bigcap\limits_{n\geq k} A_n\setminus Q\right)\\
	&\leq&\lim\limits_{k\to\infty}\sum\limits_{n\geq k}2^{-n+1} +\lim\limits_{k\to\infty}2^{-k+1} = 0.
\end{eqnarray*}
\hfill$\Box$

Also, it is well-known that, by Lemma \ref{ispl} and the pairwise singularity of ergodic measures, there exists $(\Omega_i)_{i\in \I}\subset\A_0$ such that $\Lambda_i(\Omega_i)=1$ and $\Lambda_i(\Omega_j)=0$ for all $i\neq j\in \I$. One easily verifies that, for each $i\in\I$,
\[Z1_{\Omega_i}=\lambda_iZ_i\ \phi_0\mbox{-a.e.,  and therefore,}\]
\[Z=\lambda_iZ_i\ \Lambda_i\mbox{-a.e.}\]
This implies a fact, which might be useful to observe in advance. If, for each $i\in\I$, $M_i$ is the $\Lambda_i$-essential supremum of $Z_i$, then
\begin{equation}\label{sis}
    \phi_0(X)\geq\phi_0\{Z>0\}=\int\frac{1}{Z}d\Lambda=\sum\limits_{i\in\I}\lambda_i\int\frac{1}{Z}d\Lambda_i\geq\sum\limits_{i\in\I}\frac{1}{M_i}.
\end{equation}

We will need also the following extension of Definition \ref{icfw}.
\begin{Definition}\label{ddme}
	Let $\xi$ and $\psi$ be non-negative measures on $\A_0$ such that $\Xi(X)<\infty$. Let $Q\in\P(X)$ and $\e>0$.
	Let $\times$ be the placeholder for $\xi$ or the pair ${\psi,\xi}$ if $\Psi^\xi(Q)<\infty$. In the case, when a non-negative measures $\o(f)$ on $\A_0$ is given by $\o(f)(A):=\int_{A}fd\o$ for all $A\in\A_0$ with a non-negative, $\A_0$-measurable function $f$ and a non-negative measures $\o$ on $\A_0$, we will use also the following notation. Define
	\[\underline{\O}^{\times}(f)_{\e}(Q):=\inf\limits_{(A_m)_{m\leq 0}\in\C^\times_\e(Q)}\sum\limits_{m\leq 0}\o(f)\left(S^mA_m\right),\]
	\[\underline{\O}^{\times}(f)(Q):=\lim\limits_{\e\to 0}\underline{\O}^{\times}(f)_{\e}(Q),\]
	\[\overline{\O}^{\times}(f)_{\e}(Q):=\sup\limits_{(A_m)_{m\leq 0}\in\C^{\times}_\e(Q)}\sum\limits_{m\leq 0}\o(f)\left(S^mA_m\right),\mbox{ and}\]
	\[\overline{\O}^{\times}(f)(Q):=\lim\limits_{\e\to 0}\overline{\O}^{\times}(f)_{\e}(Q).\]
	In the case that $f=1_C$ for some $C\in\A_0$, we will also use the abbreviations $\underline{\O}^{\times}_{C,\e}(Q):=\underline{\O}^{\times}(1_C)_{\e}(Q)$, $\underline{\O}^{\times}_{C}(Q):=\underline{\O}^{\times}(1_C)(Q)$, $\overline{\O}^{\times}_{C,\e}(Q):=\overline{\O}^{\times}(1_C)_{\e}(Q)$, and $\overline{\O}^{\times}_{C}(Q):=\overline{\O}^{\times}(1_C)(Q)$.
	Let $\dot\C^\times_\e(Q)$ denote the set of all $(A_m)_{m\leq 0}\in\C^\times_{\e}(Q)$ such that $A_{m}$'s are pairwise disjoint. One easily sees that $\dot\C^\times_{\e}(Q)$ is not empty (it follows from the proof of the next lemma). Define
	$\underline{\dot\O}^{\times}(f)_{\e}(Q)$, $\underline{\dot\O}^{\times}(f)(Q)$, $\overline{\dot\O}^{\times}(f)_{\e}(Q)$, and $\overline{\dot\O}^{\times}(f)(Q)$ the same way as $\underline{\O}^{\times}(f)_{\e}(Q)$, $\underline{\O}^{\times}(f)(Q)$, $\overline{\O}^{\times}(f)_{\e}(Q)$, and $\overline{\O}^{\times}(f)(Q)$ only with $\C^\times_\e(Q)$ replaced by $\dot\C^\times_\e(Q)$.
\end{Definition}

\begin{lemma}\label{del}
	Let $\xi$ and $\psi$ be non-negative measures on $\A_0$ such that $\Xi(X)<\infty$. Let $Q\in\P(X)$ and $\e>0$.
	Let $\times$ be the placeholder for $\xi$ or the pair ${\psi,\xi}$ if $\Psi^\xi(X)<\infty$. Then
	\[\underline{\dot\O}^{\times}(f)_{\e}(Q)=\underline{\O}^{\times}(f)_{\e}(Q).\]
\end{lemma}
{\it Proof.} 
  The proof is straightforward (e.g. see Lemma 9 (Lemma 5 in the arXiv version) in \cite{Wer15}).
\hfill$\Box$

One easily sees that, for every $C\in\A_0$,
\begin{equation}\label{simd}
   \Lambda(Q) =\overline{\Lambda}^{\Lambda}_C(Q)+\underline{\Lambda}^{\Lambda}_{X\setminus C}(Q),\mbox{ and }
\end{equation}
\begin{equation}\label{imsd}
   \Lambda(Q) =\overline{\Lambda}^{\psi,\Lambda}_C(Q)+\underline{\Lambda}^{\psi,\Lambda}_{X\setminus C}(Q)\ \ \mbox{ for all }Q\in\B.
\end{equation}

Also, $\underline{\Lambda}^{\Lambda}(f1_C)$ has the following property.
\begin{lemma}\label{edml}
	Let $\tilde\Lambda$ be an invariant, finite measure on $\A_0$. Let $f$ be a non-negative, $\B$-measurable function such that $\int fd\Lambda<\infty$. For each $i\in\I$, let $C_i\in\A_0$ such that $C_i\subset\Omega_i$, and $\Lambda_i(C_i)<1$. Let $C:=\bigcup_{i\in\I}C_i$. Then
	\[\underline{\Lambda}^{\tilde\Lambda}(f1_C)(Q)=0\ \ \mbox{ for all }Q\in\B.\]
\end{lemma}
{\it Proof.} Let $Q\in\B$. Define
\[\Lambda(f1_C)(Q):=\inf\limits_{(A_m)_{m\leq 0}\in\C(Q)}\sum\limits_{m\leq 0}\int\limits_{C\cap S^mA_m}fd\Lambda.\]
Since $(...,\emptyset,\emptyset, X)\in\C(X)$, $\Lambda(f1_C)(X)\leq\int fd\Lambda<\infty$, and therefore, by Lemma \ref{icp} (i), 
\[\underline{\Lambda}^{\tilde\Lambda}(f1_C)(Q)=\Lambda(f1_C)(Q).\]
Let $k\in\N$. Define $A^k_m:=X\setminus S^{-m}C$ for all $m\leq 0$ such that $m\neq-k$ and $A^k_{-k}:=(\bigcap_{i=0}^kS^iC)\cup(X\setminus S^kC)$. Then $A^k_m\in\A_m$ for all $m\leq 0$ and $\bigcup_{m\leq 0}A^k_m=X$. Hence, $(A^k_m)_{m\leq 0}\in\C(Q)$, and therefore,
\[\underline{\Lambda}^{\tilde\Lambda}(f1_C)(Q)\leq\sum\limits_{m\leq 0}\int\limits_{C\cap S^mA^k_m}fd\Lambda=\int\limits_{S^{-k}\bigcap_{i=0}^kS^iC}fd\Lambda=\int\limits_{\bigcap_{i=0}^kS^{-i}C}fd\Lambda\]
for all $k\in\N$. Hence, since $\int fd\Lambda<\infty$,
\[\underline{\Lambda}^{\tilde\Lambda}(f1_C)(Q)\leq\int\limits_{\bigcap_{i=0}^\infty S^{-i}C}fd\Lambda=\sum\limits_{j\in\I}\lambda_j\int\limits_{\bigcap_{i=0}^\infty S^{-i}C}fd\Lambda_j.\]
Since, for each $j\in\I$,
$\Lambda_j\left(\bigcap_{i=0}^\infty S^{-i}C\right)=\Lambda_j\left(\bigcap_{i=-\infty}^\infty S^{-i}C\right)$, and\\ $\Lambda_j\left(\bigcap_{i=-\infty}^\infty S^{-i}C\right)\leq\Lambda_j(C)=\Lambda_j(C_j)<1$, it follows, by the ergodicity of $\Lambda_j$, that $\Lambda_j\left(\bigcap_{i=0}^\infty S^{-i}C\right)=0$ for all $j\in\I$. Thus, the assertion follows.
\hfill$\Box$

The computation of some DDMs can be deduced from the following theorem.
\begin{theo}\label{ect}
	 Let $f$ be a non-negative, $\A_0$-measurable function such that $\int fd\Lambda<\infty$. For each $i\in\I$, let $l_i$ be the $\Lambda_i$-essential infimum of $f$. Then
	\[\underline{\Lambda}^{\Lambda}(f)(Q)=\sum\limits_{i\in\I}\lambda_i l_i\Lambda_i(Q)\ \ \mbox{ for all }Q\in\B.\]
\end{theo}
{\it Proof.}
Let $Q\in\B$, $\e>0$, and $(A_m)_{m\leq 0}\in\C^1_\e(Q)$. Then
\begin{eqnarray*}
  \sum\limits_{m\leq 0}\int\limits_{S^mA_m}fd\Lambda&=&\sum\limits_{i\in\I}\lambda_i\sum\limits_{m\leq 0}\int\limits_{S^mA_m}fd\Lambda_i\geq\sum\limits_{i\in\I}\lambda_i l_i\sum\limits_{m\leq 0}\Lambda_i\left(A_m\right)\\
  &\geq&\sum\limits_{i\in\I}\lambda_i l_i\Lambda_i\left(Q\right).
\end{eqnarray*}
Hence,
\begin{equation}\label{ectf}
\underline{\Lambda}^{\Lambda}(f)(Q)\geq\sum\limits_{i\in\I}\lambda_i l_i\Lambda_i(Q).
\end{equation}

Now, let $0<c<1$. For each $i\in\I$, define 
\begin{equation*} 
\eta_i(c):=\left\{\begin{array}{cc}
\frac{l_i}{1-c} & \mbox{if }l_i>0\\
c &\mbox{otherwise } 
\end{array}\right.,
\end{equation*}
$C_i:=\{f\geq \eta_i(c)\}\cap\Omega_i$, and $C:=\bigcup_{i\in\I}C_i$.  Let $\d>0$ such that $\int_{B}f d\Lambda<\e$ whenever $\Lambda(B)<\d$ for some $B\in\B$. Let $(B_m)_{m\leq 0}\in\dot\C^1_{\d}(Q)$. Then
\begin{eqnarray*}
  &&\underline{\Lambda}^{\Lambda}(f)_\d(Q)\\
  &\leq&\sum\limits_{i\in\I}\lambda_i\sum\limits_{m\leq 0}\int\limits_{S^mB_m}fd\Lambda_i\\
  &=&\sum\limits_{i\in\I}\lambda_i\sum\limits_{m\leq 0}\left(\int\limits_{C_i\cap S^mB_m}fd\Lambda_i+\int\limits_{\left(S^mB_m\right)\setminus C_i}fd\Lambda_i\right)\\
  &\leq&\sum\limits_{i\in\I}\lambda_i\sum\limits_{m\leq 0}\int\limits_{C\cap S^mB_m}fd\Lambda_i+\sum\limits_{i\in\I}\lambda_i \eta_i(c)\sum\limits_{m\leq 0}\Lambda_i\left(B_m\right)\\
  &=&\sum\limits_{m\leq 0}\int\limits_{C\cap S^mB_m}fd\Lambda+\sum\limits_{i\in\I}\lambda_i \eta_i(c)\Lambda_i\left(Q\right)+\sum\limits_{i\in\I}\lambda_i \eta_i(c)\Lambda_i\left(\bigcup\limits_{m\leq 0}B_m\setminus Q\right).
\end{eqnarray*}
Observe that
\begin{eqnarray*}
 &&\sum\limits_{i\in\I}\lambda_i \eta_i(c)\Lambda_i\left(\bigcup\limits_{m\leq 0}B_m\setminus Q\right)\\
 &=&\frac{1}{1-c}\sum\limits_{i\in\I,\ l_i>0}\lambda_i l_i\Lambda_i\left(\bigcup\limits_{m\leq 0}B_m\setminus Q\right)+c\sum\limits_{i\in\I,\ l_i=0}\lambda_i\Lambda_i\left(\bigcup\limits_{m\leq 0}B_m\setminus Q\right)\\
 &\leq&\frac{1}{1-c}\sum\limits_{i\in\I}\lambda_i \int\limits_{\bigcup_{m\leq 0}B_m\setminus Q}fd\Lambda_i+c\Lambda\left(\bigcup\limits_{m\leq 0}B_m\setminus Q\right)\\
 &<&\frac{\e}{1-c}+c\d.
\end{eqnarray*}
Hence,
\[\underline{\Lambda}^{\Lambda}(f)_\d(Q)\leq\underline{\dot\Lambda}^{\Lambda}(f1_C)_\d(Q)+\sum\limits_{i\in\I}\lambda_i \eta_i(c)\Lambda_i(Q)+\frac{\e}{1-c}+c\d.\]
By Lemma \ref{del} and Lemma \ref{edml}, it follows that
\[\underline{\Lambda}^{\Lambda}(f)(Q)\leq\sum\limits_{i\in\I}\lambda_i \eta_i(c)\Lambda_i(Q)\leq\frac{1}{1-c}\sum\limits_{i\in\I}\lambda_i  l_i\Lambda_i(Q)+c,\]
which implies the converse inequality of \eqref{ectf}.
\hfill$\Box$

\subsection{A computation of $(0,1)\owns\a\longmapsto\H_{\a}\left(\Lambda,\phi_0\right)$}\label{chs1}

In this subsection, we compute the function $(0,1)\owns\a\longmapsto\H_{\a}\left(\Lambda,\phi_0\right)$ explicitly in terms of $Z$. It shows, in particular, that an irregularity of the function can occur only in the case of an uncountable ergodic decomposition of $\Lambda$.

We would like to remind that we are using the definitions $1/\infty:=0$, and $0^0:=1$ (i.e. $0\log 0:=0$).
\begin{theo}\label{aedt} 
	For each $i\in\I$, let $M_i$ be the $\Lambda_i$-essential supremum of $Z_i$. Let $Q\in\B$. Then
	
	(i)\[\Phi(Q)\geq\sum\limits_{i\in \I} \frac{1}{M_i}\Lambda_i(Q)+\overline{\Phi_0}^{\phi_0}_{\{Z=0\}}(Q),\mbox{ and}\]
	 
	(ii) 
	\[\H_{\a}\left(\Lambda,\phi_0\right)(Q)=\sum\limits_{i\in \I} \lambda_i^\a\left(\frac{1}{M_i}\right)^{1-\a}\Lambda_i(Q) \mbox{ for all }0<\a\leq 1.\]
	The equality is also true for $\a=0$ if $\phi_0\ll\Lambda$.
\end{theo}
{\it Proof.}
 $(i)$  Let $\e>0$ and $(A_m)_{m\leq 0}\in\C^0_\e(Q)$. Then, since $Z=\lambda_i Z_i$ $\Lambda_i$-a.e.,
  \begin{eqnarray*}
     &&\Phi(Q)+\e>\sum\limits_{m\leq 0}\phi_0\left(S^mA_m\right)\\
     &=&\sum\limits_{m\leq 0}\int\limits_{S^mA_m}\frac{1}{Z} d\Lambda+\sum\limits_{m\leq 0}\phi_0\left(\{Z=0\}\cap S^mA_m\right)\\
     &=&\sum\limits_{i\in \I} \lambda_i\sum\limits_{m\leq 0}\int\limits_{S^mA_m}\frac{1}{Z}d\Lambda_i+\sum\limits_{m\leq 0}\phi_0\left(\{Z=0\}\cap S^mA_m\right)\\
     &\geq&\sum\limits_{i\in \I}\frac{1}{M_i}\Lambda_i(Q)+\sum\limits_{m\leq 0}\phi_0\left(\{Z=0\}\cap S^mA_m\right).
  \end{eqnarray*}
  This implies the inequality in (i).

  $(ii)$ Obviously, the equality in (ii) is correct for $\a=1$.  
  
  Let $0<\a<1$, or $\a=0$ and $\phi_0\ll\Lambda$. Then
  \[\sum\limits_{m\leq 0}\int\limits_{S^mB_m}Z^\a d\phi_0=\sum\limits_{m\leq 0}\int\limits_{S^mB_m}\left(\frac{1}{Z}\right)^{1-\a} d\Lambda\]
  for all $(B_m)_{m\leq 0}\in\C(Q)$, and therefore, by \eqref{hmvi},
  \[\H_{\a}\left(\Lambda,\phi_0\right)(Q)=\H^{\a,1}\left(\Lambda,\phi_0\right)(Q)=\underline{\Lambda}^{\Lambda}\left(\left(\frac{1}{Z}\right)^{1-\a}\right)(Q).\]
  Observe that, for every $i\in\I$,  the $\Lambda_i$-essential infimum of $(1/Z)^{1-\a}$ is\\ $(1/(\lambda_i M_i))^{1-\a}$. Thus, the equality in (ii) follows by Theorem \ref{ect}.
\hfill$\Box$

By Proposition \ref{ca0}, we know that $\Phi(Q)>0$ if the function $[0,1]\owns\alpha\longmapsto \H_{\a}\left(\Lambda,\phi_0\right)(Q)$ is discontinuous at $0$. Now, we obtain a computation of $\Phi$ in the case of the continuity.

\begin{cor}\label{hecc}
    For each $i\in\I$, let $M_i$ be the $\Lambda_i$-essential supremum of $Z_i$. Let $Q\in\B$.
    
    (i) The function $[0,1]\owns\alpha\longmapsto\H_{\a}\left(\Lambda,\phi_0\right)(Q)$ is continuous at $0$ if and only if
	\[\Phi(Q)=\sum\limits_{i\in \I}\frac{1}{M_i}\Lambda_i(Q).\]
	
	(ii) The function $[0,1]\owns\alpha\longmapsto\H_{\a}\left(\Lambda,\phi_0\right)(Q)$ is continuous at $0$ if $\phi_0\ll\Lambda$.
\end{cor}
{\it Proof.} 
(i) Observe that, for $0<\a<1$,
\begin{eqnarray*}
	&&\sum\limits_{i\in \I} \lambda_i^\a\left(\frac{1}{M_i}\right)^{1-\a}\Lambda_i(Q)\\
	&=&\sum\limits_{i\in\I,\ \lambda_i M_i\leq 1} \lambda_i\left(\frac{1}{\lambda_i M_i}\right)^{1-\a}\Lambda_i(Q)+\sum\limits_{i\in\I,\ \lambda_i M_i>1} \lambda_i\left(\frac{1}{\lambda_i M_i}\right)^{1-\a}\Lambda_i(Q).
\end{eqnarray*}
Hence, by Theorem \ref{aedt} (ii) and the Monotone Convergence Theorem, 
\begin{equation}\label{hzl}
   \lim\limits_{\a\to 0}\H_{\a}\left(\Lambda,\phi_0\right)(Q)=\sum\limits_{i\in \I}\frac{1}{M_i}\Lambda_i(Q).
\end{equation}
This proves (i).

(ii) It follow by Theorem \ref{aedt} (ii) and \eqref{hzl}.
\hfill$\Box$

\begin{cor}\label{hcdc}
		For each $i\in\I$, let $M_i$ be the $\Lambda_i$-essential supremum of $Z_i$. Let $Q\in\B$. Then the function $(0,1)\owns\a\longmapsto\H_{\a}\left(\Lambda,\phi_0\right)(Q)$ is infinitely differentiable, and, for every $\a\in(0,1)$ and $n\in\N$, the $n$-th derivative of it at $\a$,
		\[\left.\frac{d^n}{d^n x}\H_{x}\left(\Lambda,\phi_0\right)(Q)\right|_{x=\a}=\sum\limits_{i\in \I} \lambda_i^\a{M_i}^{\a-1}(\log\left(\lambda_iM_i\right))^n\Lambda_i(Q).\]
\end{cor}
{\it Proof.} 
 The proof follows straightforward from Theorem \ref{aedt} (ii), by induction and the well-known theorem on the differentiation of parameter dependent integrals resulting from the Lebesgue Dominated Conversion Theorem. We give only the induction beginning here. Let $\a\in (0,1)$ and $\d>0$ such that $0<\a-\d$ and $\a+\d<1$, then, for each $i\in\I$,
 \begin{equation*} 
    \left|\left.\frac{d}{d x}\left(\lambda_iM_i\right)^{x-1}\right|_{x=\a}\right|=\left|\left(\lambda_iM_i\right)^{\a-1}\log\left(\lambda_iM_i\right)\right|\leq g_i(\a)
 \end{equation*}
 where
 \[g_i(\a):=\left\{\begin{array}{cc}
 \frac{1}{\d}\left(\lambda_i M_i\right)^{\a-1}\left(\left(\lambda_i M_i\right)^{-\d}-1\right) & \mbox{if }\lambda_i M_i\leq 1\\
 \frac{1}{\d}\left(\lambda_i M_i\right)^{\a-1}\left(\left(\lambda_i M_i\right)^{\d}-1\right) &\mbox{otherwise } 
 \end{array}\right.,\]
 and 
 \begin{eqnarray*}
   &&\sum\limits_{i\in \I} \lambda_i g_i(\a)\Lambda_i(Q)\\
   &\leq&\frac{1}{\d}\left(\sum\limits_{i\in \I} \lambda_i \left(\lambda_i M_i\right)^{\a-\d-1}\Lambda_i(Q)+\sum\limits_{i\in \I} \lambda_i \left(\lambda_i M_i\right)^{\a+\d-1}\Lambda_i(Q)\right)\\
   &=&\frac{1}{\d}\left(\H_{\a-\d}\left(\Lambda,\phi_0\right)(Q)+\H_{\a+\d}\left(\Lambda,\phi_0\right)(Q)\right)<\infty,
 \end{eqnarray*}
 which proves the assertion for $n=1$.
\hfill$\Box$

\begin{cor}\label{pebc}
	 Let $Q\in\B$. Suppose  $\Lambda$ is ergodic, and $[0,1]\owns\alpha\longmapsto \H_{\a}\left(\Lambda,\phi_0\right)(Q)$ is continuous at $0$. Then
	\begin{eqnarray*}
	   \H_\a(\Lambda,\phi_0)(Q)&=&\H^{\a,\a_0}(\Lambda,\phi_0)(Q)=\H^{\a,\b,0}\left(\Lambda,\phi_0\right)(Q)=\J_{\a}\left(\Lambda,\phi_0\right)(Q)\\
	   &=&\Phi(Q)^{1-\a}\Lambda(Q)^\a
	\end{eqnarray*}
	for all $\a,\a_0,\b\in[0,1]$.
\end{cor}
{\it Proof.} 
 For $\a_0=0$, the assertion follows by Theorem \ref{aedt} (ii), Corollary \ref{hecc} (i), Lemma \ref{tlb} (i), and Lemma \ref{pchm} (ii). For $0\leq\a_0\leq\a$, it follows, by Corollary \ref{hecc} (i), Theorem \ref{aedt} (ii), and Lemma \ref{shm} (i), that
\[\H^{\a,\a_0}(\Lambda,\phi_0)(Q)=\H_{\a}(\Lambda,\phi_0)(Q),\]
which implies, by Lemma \ref{hsl}, that also
\[\H^{\a_0,\a}(\Lambda,\phi_0)(Q)=\H_{\a_0}(\Lambda,\phi_0)(Q),\]
which proves the first equality for any other $\a_0$.
\hfill$\Box$

\begin{cor}\label{dhsp}
	Suppose $\Lambda$ is ergodic. Let $M<\infty$ be the $\Lambda$-essential supremum of $Z$. Let $Q\in\B$. Then
	\[\H^{\g,\a_0}(\Lambda,\phi_0)(Q)= M^{\g-\a_0}\H_{\a_0}(\Lambda,\phi_0)(Q)=M^{\g-1}\Lambda(Q).\]
	for all $0<\a_0< 1\leq \g$. If, in addition, $[0,1]\owns\alpha\longmapsto \H_{\a}\left(\Lambda,\phi_0\right)(Q)$ is continuous at $0$, then the equalities are true also for $\a_0=0$.
\end{cor}
{\it Proof.} Clearly, the second equality follows by Theorem \ref{aedt} (ii), as also the first, by Lemma \ref{icp} (ii), for $\g=1$. 

Let $0<\a_0< 1< \g$. Let $\e>0$ and $(A_m)_{m\leq 0}\in\C^{\a_0}_\e(Q)$. Then
\[\H^{\g,\a_0}_\e(\Lambda,\phi_0)(Q)\leq\sum\limits_{m\leq 0}\int\limits_{S^mA_m}Z^{\g-1} d\Lambda\leq M^{\g-\a_0}\sum\limits_{m\leq 0}\int\limits_{S^mA_m}Z^{\a_0-1} d\Lambda.\]
Hence,
\[\H^{\g,\a_0}(\Lambda,\phi_0)(Q)\leq M^{\g-\a_0}\H_{\a_0}(\Lambda,\phi_0)(Q),\]
which implies, in particular, the assertion if $\H_{\a_0}(\Lambda,\phi_0)(Q)=0$.
By Proposition \ref{clbh} (ii) and Theorem \ref{aedt} (ii), it follows that
\begin{eqnarray*}
	\Lambda(Q)&\leq&{\H_{\a_0}(\Lambda,\phi_0)(Q)}^{1-\frac{1-\a_0}{\g-\a_0}}{\H^{\g,\a_0}\left(\Lambda,\phi_0\right)(Q)}^{\frac{1-\a_0}{\g-\a_0}}\leq\H_{\a_0}(\Lambda,\phi_0)(Q)M^{1-\a_0}\\
	&=&\Lambda(Q),
\end{eqnarray*}
which implies the assertion if $\H_{\a_0}(\Lambda,\phi_0)(Q)>0$.
\hfill$\Box$

\begin{cor}
	Suppose $\Lambda$ is an ergodic probability measure. If the $\Lambda$-essential supremum of $Z$ is infinite, then, for every $0<\a<1$,
	\[\K_\a(\Lambda|\phi_0)(Q)=\infty\ \ \mbox{ for all }Q\in\B\mbox{ such that }\Lambda(Q)>0.\]
\end{cor}
{\it Proof.} 
It follows form Theorem \ref{lbre} (ii) and Theorem \ref{aedt} (ii).
\hfill$\Box$

\subsection{A computation of $(1, \g]\owns\a\longmapsto\H_{\a}\left(\Lambda,\phi_0\right)$}\label{chs2}

In this subsection, we make an additional assumption, which is slightly stronger than $\int\log Zd\Lambda<\infty$ and implies the finiteness of $\H_{\a}\left(\Lambda,\phi_0\right)$ for some $\a>1$. This allows us to compute the latter explicitly. The obtained formula reveals a discontinuity of the derivative of the function at $\a=1$.

\begin{theo}\label{ec2} 
	For each $i\in\I$, let $m_i$ be the $\Lambda_i$-essential infimum of $Z_i$. Suppose there exists $\g>1$ such that
	$\int Z^{\g-1} d\Lambda<\infty$. Let $Q\in\B$. Then 
	\[\H_{\a}\left(\Lambda,\phi_0\right)(Q)=\sum\limits_{i\in \I} \lambda_i^\a m_i^{\a-1}\Lambda_i(Q)\mbox{ for all }1\leq\a\leq\g.\]
\end{theo}
{\it Proof.}
 Let $1\leq\a\leq\g$. Then
\begin{equation*}
  \sum\limits_{m\leq 0}\int\limits_{S^mA_m}Z^\a d\phi_0=\sum\limits_{m\leq 0}\int\limits_{S^mA_m}Z^{\a-1} d\Lambda
\end{equation*}
for all $(A_m)_{m\leq 0}\in\C(Q)$. Hence, by \eqref{hmvi},
\[\H_{\a}\left(\Lambda,\phi_0\right)(Q)=\H^{\a,1}\left(\Lambda,\phi_0\right)(Q)=\underline{\Lambda}^{\Lambda}\left(Z^{\a-1}\right)(Q).\]
Since, for every $i\in\I$, the $\Lambda_i$-essential infimum of $Z^{\a-1}$ is $(\lambda_i m_i)^{\a-1}$, the assertion follows by Theorem \ref{ect}.
\hfill$\Box$

\begin{cor}\label{hsdc}
	Suppose $\Lambda$ is ergodic, and there exists $\g>1$ such that $\int Z^{\g-1} d\Lambda<\infty$. Let $m>0$ be the $\Lambda$-essential infimum of $Z$. Let $Q\in\B$. Then
	\[\H^{\b,\a}(\Lambda,\phi_0)(Q)=\left(\frac{1}{m}\right)^{\a-\b}\H_{\a}(\Lambda,\phi_0)(Q)=\left(\frac{1}{m}\right)^{1-\b}\Lambda(Q)\]
	for all $0<\b\leq 1<\a\leq\g$.
\end{cor}
{\it Proof.}  Let $0<\b\leq 1<\a\leq\g$. Obviously, the second equality follows by Theorem \ref{ec2}. Let $\e>0$ and $(A_m)_{m\leq 0}\in\C^{\a}_\e(Q)$. Then, as in Subsection \ref{chm},
\begin{eqnarray*}
    \H^{\b,\a}_\e(\Lambda,\phi_0)(Q)&\leq&\sum\limits_{m\leq 0}\int\limits_{S^mA_m}Z^{\b-1} d\Lambda\leq\left(\frac{1}{m}\right)^{\a-\b}\left(\H_{\a}(\Lambda,\phi_0)(Q)+\e\right),
\end{eqnarray*}
and therefore,
\[\H^{\b,\a}(\Lambda,\phi_0)(Q)\leq\left(\frac{1}{m}\right)^{\a-\b}\H_{\a}(\Lambda,\phi_0)(Q).\]
Hence, the first equality is correct if $\H_{\a}(\Lambda,\phi_0)(Q)=0$. By Proposition \ref{clbh} (i) and Theorem \ref{ec2},
\begin{eqnarray*}
   \Lambda(Q)&\leq&{\H^{\b,\a}\left(\Lambda,\phi_0\right)(Q)}^{1-\frac{1-\b}{\a-\b}}{\H_{\a}(\Lambda,\phi_0)(Q)}^{\frac{1-\b}{\a-\b}}\\
   &\leq&\left(\frac{1}{m}\right)^{(\a-\b)\left(1-\frac{1-\b}{\a-\b}\right)}\H_{\a}(\Lambda,\phi_0)(Q)\\
   &=&\Lambda(Q),
\end{eqnarray*}
which implies the assertion if $\H_{\a}(\Lambda,\phi_0)(Q)>0$.
\hfill$\Box$

\begin{Remark}\label{dmnr}
   A comparison of Theorem \ref{ec2} with Corollary \ref{dhsp} and Theorem \ref{aedt} (ii) with Corollary \ref{hsdc} leaves no doubts that the technique of outer measure approximations (developed in \cite{Wer15}) allows to obtain new measures.
	
   Also, Corollary \ref{dhsp} and Corollary \ref{hsdc} reveal a discontinuity of $\a\longmapsto\H^{\b,\a}(\Lambda,\phi_0)$  at $1$ from the left if $\b> 1$ and from the right if $0<\b< 1$ respectively, as follows. 
   
   Under the assumptions of the first,
   \[\lim\limits_{\a\uparrow 1}\H^{\b,\a}(\Lambda,\phi_0)(Q)=M^{\b-1}\Lambda(Q)\]
   for all $\b\geq 1$, where $M$ is the $\Lambda$-essential supremum of $Z$, whereas by \eqref{hmvi} and Theorem \ref{ec2}, 
   \[\H^{\b,1}(\Lambda,\phi_0)(Q)=m^{\b-1}\Lambda(Q)\]
   for all $\b\geq 1$, where $m$ is the $\Lambda$-essential infimum of $Z$.
   
   Under the assumptions of the second,
   \[\lim\limits_{\a\downarrow 1}\H^{\b,\a}(\Lambda,\phi_0)(Q)=\left(\frac{1}{m}\right)^{1-\b}\Lambda(Q)\]
   for all $0<\b\leq 1$, whereas by \eqref{hmvi} and Theorem \ref{aedt} (ii), 
   \[\H^{\b,1}(\Lambda,\phi_0)(Q)=\left(\frac{1}{M}\right)^{1-\b}\Lambda(Q)\]
   for all $0<\b\leq 1$.
\end{Remark}

\subsection{DDMs from sequences of inconsistent non-additive contents}\label{inc}

In this subsection, we would like also to make some corollaries, which might be useful for a further development of the dynamical measure theory.

\begin{Definition}
	Let $Q\in\P(X)$. For $\a\geq 0$, define
	\[\mathbf{H}_\a(\Lambda,\phi_0)(Q):=\inf\limits_{(A_m)_{m\leq 0}\in\C(Q)}\sum\limits_{m\leq 0,\ \phi_0(S^mA_m)>0}\phi_0\left(S^mA_m\right)\left(\frac{\Lambda(A_m)}{\phi_0\left(S^mA_m\right)}\right)^{\a}.\] 
\end{Definition}

Clearly,  $\mathbf{H}_\a(\Lambda,\phi_0)(Q)\leq\phi_0(X)^{1-\a}$ for all $\a\geq 0$ and $Q\in\P(X)$.
Since $\sum_{m\leq 0,\ \phi_0(S^mA_m)>0}\Lambda\left(A_m\right)=\sum_{m\leq 0}\Lambda\left(A_m\right)$, it is obvious that $\mathbf{H}_0(\Lambda,\phi_0)(Q)=\Phi(Q)$ and $\mathbf{H}_1(\Lambda,\phi_0)(Q)=\Lambda(Q)$.

For other $\a$, in this definition, we deal with a summation of contents which are not only inconsistent to each-other, but also each is non-additive. At the time of writing, the question whether the restriction of this set function on $\B$ is a measure, in general, is beyond of what the dynamical measure theory (developed in \cite{Wer15}) can immediately answer. However, if $\Lambda$ is ergodic, then it is a simple implication from our results that the answer to this question is affirmative in some cases.

\begin{cor}
	Suppose $\Lambda$ is ergodic, and $[0,1]\owns\alpha\longmapsto \H_{\a}\left(\Lambda,\phi_0\right)(Q)$ is continuous at $0$. Let $M$ be the $\Lambda$-essential supremum of $Z$. Let $Q\in\B$. Then
	\begin{eqnarray*}
		\mathbf{H}_\a(\Lambda,\phi_0)(Q)=\left(\frac{1}{M}\right)^{1-\a}\Lambda(Q)
	\end{eqnarray*}
	for all $\a\in[0,1]$. (That is $\mathbf{H}_\a(\Lambda,\phi_0)(Q)=\H_\a(\Lambda,\phi_0)(Q)$ for all $\a\in[0,1]$.)
\end{cor}
{\it Proof.} By Corollary \ref{hecc} (i), we only need to consider $0<\a<1$. For $\e>0$, define
\[\mathbf{H}^{\a,1}_\e(\Lambda,\phi_0)(Q):=\inf\limits_{(A_m)_{m\leq 0}\in\C^1_\e(Q)}\sum\limits_{m\leq 0,\ \phi_0(S^mA_m)>0}\phi_0\left(S^mA_m\right)\left(\frac{\Lambda(A_m)}{\phi_0\left(S^mA_m\right)}\right)^{\a},\]
and $\mathbf{H}^{\a,1}(\Lambda,\phi_0)(Q):=\lim\limits_{\e\to 0}\mathbf{H}^{\a,1}_\e(\Lambda,\phi_0)(Q)$.
Then, by Lemma \ref{lbl},
\[\H_\a(\Lambda,\phi_0)(Q)\leq\mathbf{H}_\a(\Lambda,\phi_0)(Q)\leq\mathbf{H}^{\a,1}(\Lambda,\phi_0)(Q)\leq\Phi(Q)^{1-\a}\Lambda(Q)^{\a}.\]
Thus, the assertion follows by Corollary \ref{pebc}.
\hfill$\Box$

\begin{cor}\label{nhgo}
	Suppose $\Lambda$ is ergodic, and there exists $\g>1$ such that $\int Z^{\g-1} d\Lambda<\infty$. Let $m$ be the $\Lambda$-essential infimum of $Z$. Let $Q\in\B$. If also $\phi_0\ll\Lambda$, then
	\[\mathbf{H}_{\a}(\Lambda,\phi_0)(Q)=m^{\a-1}\Lambda(Q)\]
	for all $1\leq\a\leq\g$. (That is $\mathbf{H}_\a(\Lambda,\phi_0)(Q)=\H_\a(\Lambda,\phi_0)(Q)$ for all $\a\in[1,\g]$.)
\end{cor}
{\it Proof.} Obviously, we only need to give a proof for $1<\a\leq\g$. Let $(A_m)_{m\leq 0}\in\C(Q)$. Since
\[\int\limits_{S^mA_m}Z^{\a}d\phi_0\geq\phi_0\left(S^mA_m\right)\left(\frac{\Lambda\left(A_m\right)}{\phi_0\left(S^mA_m\right)}\right)^{\a}\]
for all $\phi_0\left(S^mA_m\right)>0$, it follows, by Theorem \ref{ec2}, that
\begin{equation*}
	m^{\a-1}\Lambda(Q)=\H_{\a}(\Lambda,\phi_0)(Q)\geq\mathbf{H}_{\a}(\Lambda,\phi_0)(Q).
\end{equation*}
On the other hand, since $Z\geq m$ $\phi_0$-a.e., $\Lambda(A)\geq m\phi(A)$ for all $A\in\A_0$, and therefore,
\[\sum\limits_{m\leq 0,\ \phi_0\left(S^mA_m\right)>0}\phi_0\left(S^mA_m\right)\left(\frac{\Lambda\left(A_m\right)}{\phi_0\left(S^mA_m\right)}\right)^{\a}\geq m^{\a-1}\sum\limits_{m\leq 0}\Lambda\left(A_m\right)\geq m^{\a-1}\Lambda\left(Q\right).\]
Hence,
\[\mathbf{H}_{\a}(\Lambda,\phi_0)(Q)\geq m^{\a-1}\Lambda(Q),\]
what remained to show.
\hfill$\Box$

\subsection{A computation of $(0,\g)\times(0,\g)\owns (\b,\a)\longmapsto\H^{\b,\a}(\Lambda,\phi_0)$}\label{homs}

In this subsection, we complete the computation of $\H^{\b,\a}(\Lambda,\phi_0)$ (from Corollaries \ref{dhsp} and \ref{hsdc}) for some of the remaining parameter values and extend it to the case of a discrete ergodic decomposition of $\Lambda$ using a more systematic method, which also uses the results on the computation of $\H_{\a}(\Lambda,\phi_0)$ from Subsection \ref{chs1} and Subsection \ref{chs2}. In particular, under the additional assumption that $\phi_0\ll\Lambda$, we compute $\H^{\b,0}\left(\Lambda,\phi_0\right)$ for all $\b\in[0,\g)$ where the measure remains finite.

Following Definition \ref{ddme}, for $\e>0$, $\a\geq 0$, $C\in\A_0$, a non-negative measure $\o$ on $\A_0$, and $h_\a(A):=\int_A Z^\a d\phi_0$ for all $A\in\A_0$, let us abbreviate $\underline{\O}_{C,\e}^{\a,1}:=\underline{\O}^{h_\a,\Lambda}_{C,\e}$, $\underline{\O}_{C}^{\a,1}:=\underline{\O}^{h_\a,\Lambda}_C$,
$\overline{\O}_{C,\e}^{\a,1}:=\overline{\O}^{h_\a,\Lambda}_{C,\e}$, and 
$\overline{\O}_{C}^{\a,1}:=\overline{\O}^{h_\a,\Lambda}_C$ (of course, provided that $\H_{\a}(\Lambda,\phi_0)(X)<\infty$). And in order to keep our notation consistent (in particular with Definition \ref{icfw}), we will write
$\H^{\b,\a,1}_\e(\Lambda,\phi_0):=\underline{\Phi_0}^{h_\a,\Lambda}\left(Z^\b\right)_\e$ and $\H^{\b,\a,1}(\Lambda,\phi_0):=\underline{\Phi_0}^{h_\a,\Lambda}\left(Z^\b\right)$.

For the computation of our DDMs arising as outer measure approximations with respect to $\H_{\a}(\Lambda,\phi_0)$, we will need the following two lemmas, which are similar to Lemma \ref{edml} in their functionality.
\begin{lemma}\label{ecal}
	  For each $i\in \I$, let $M_i$ be the $\Lambda_i$-essential supremum and $m_i$ be the $\Lambda_i$-essential infimum of $Z_i$. Let $Q\in\B$.
	 
	 (i) Let $0<c<1$. For each $i\in\I$, define 
	\begin{equation*} 
	\t_i(c):=\left\{\begin{array}{cc}
	M_i(1-c) & \mbox{if }\lambda_i M_i\leq\frac{1}{c}\\
	\frac{1-c}{\lambda_i c} &\mbox{otherwise } 
	\end{array}\right.,
	\end{equation*}
	$C_i:=\{Z_i>\t_i(c)\}\cap\Omega_i$, and $C:=\bigcup_{i\in\I}C_i$.  Then
	\[\underline{\Lambda}_{C}^{\a,1}(Q)=\overline{\Lambda}_{C}^{\a,1}(Q)=\Lambda(Q)\mbox{ for all }0<\a< 1.\] 
	The equalities are also true for $\a=0$ if $[0,1]\owns\alpha\longmapsto \H_{\a}\left(\Lambda,\phi_0\right)(Q)$ is continuous at $0$.
	
	(ii) For each $i\in\I$, let $A_i\in\A_0$ such that $A_i\subset\Omega_i$ and $\Lambda_i(A_i)<1$ and $A:=\bigcup_{i\in\I}A_i$. Then
	\[\underline{\Lambda}_{A}^{1,1}(Q)=0,\mbox{ and } \overline{\Lambda}_{A}^{1,1}(Q)=\Lambda(Q)\mbox{ if }\Lambda_i(A_i)>0\mbox{ for all }i\in\I.\] 
	
	(iii) Suppose $\int Z^{\g-1}d\Lambda<\infty$ for some $\g>1$. Let $0<c^*<1$. For each $i\in\I$, define 
	\begin{equation*} 
	\eta_i(c^*):=\left\{\begin{array}{cc}
	\frac{m_i}{1-c^*} & \mbox{if }\lambda_i m_i\geq c^*\\
	\frac{c^*}{\lambda_i(1-c^*)} &\mbox{otherwise } 
	\end{array}\right.,
	\end{equation*}
	$C^*_i:=\{Z_i<\eta_i(c^*)\}\cap\Omega_i$, and $C^*:=\bigcup_{i\in\I}C^*_i$.  Then
	\[\underline{\Lambda}_{C^*}^{\a,1}(Q)=\overline{\Lambda}_{C^*}^{\a,1}(Q)=\Lambda(Q)\mbox{ for all }1<\a\leq\g.\]
\end{lemma}
{\it Proof.}  
  $(i)$ Let $0<\a<1$, or ($\a=0$, and $[0,1]\owns\alpha\longmapsto \H_{\a}\left(\Lambda,\phi_0\right)(Q)$ is continuous at $0$). Let $\e>0$ and $(A_m)_{m\leq 0}\in\C^{\a,1}_\e(Q)$. Then
  \begin{eqnarray*}
  &&\H_{\a}\left(\Lambda,\phi_0\right)(Q)+\e>\sum\limits_{m\leq 0}\int\limits_{S^mA_m}\left(\frac{1}{Z}\right)^{1-\a}d\Lambda\\
  &=&\sum\limits_{i\in\I}\lambda_i^\a\left[\sum\limits_{m\leq 0}\int\limits_{C_i\cap S^mA_m}\left(\frac{1}{Z_i}\right)^{1-\a}d\Lambda_i+\sum\limits_{m\leq 0}\int\limits_{(X\setminus C_i)\cap S^mA_m}\left(\frac{1}{Z_i}\right)^{1-\a}d\Lambda_i\right]\\
  &\geq&\sum\limits_{i\in\I}\lambda_i^\a\left[\left(\frac{1}{M_i}\right)^{1-\a}\sum\limits_{m\leq 0}\Lambda_i\left(C_i\cap S^mA_m\right)\right.\\
  &&\ \ \ \ \ \ \ \ \ \ \ \ +\left.\left(\frac{1}{\t_i(c)}\right)^{1-\a}\sum\limits_{m\leq 0}\Lambda_i\left((X\setminus C_i)\cap S^mA_m\right)\right]\\
  &\geq&\sum\limits_{i\in\I}\lambda_i^\a\left(\frac{1}{M_i}\right)^{1-\a}\Lambda_i\left(Q\right)\\
  &&+\sum\limits_{i\in\I}\lambda_i\left(\left(\frac{1}{\lambda_i\t_i(c)}\right)^{1-\a}-\left(\frac{1}{\lambda_iM_i}\right)^{1-\a}\right)\sum\limits_{m\leq 0}\Lambda_i\left((X\setminus C_i)\cap S^mA_m\right).
  \end{eqnarray*}
  Let $i\in\I$. Observe that
  \begin{eqnarray*}
    \left(\frac{1}{\lambda_i\t_i(c)}\right)^{1-\a}-\left(\frac{1}{\lambda_iM_i}\right)^{1-\a}&=&\left\{\begin{array}{cc}
    	\frac{1}{\left(\lambda_iM_i(1-c)\right)^{1-\a}}-\frac{1}{\left(\lambda_iM_i\right)^{1-\a}} & \mbox{if }\lambda_i M_i\leq\frac{1}{c}\\
    	\left(\frac{c}{1-c}\right)^{1-\a}-\frac{1}{\left(\lambda_iM_i\right)^{1-\a}} &\mbox{otherwise } 
    \end{array}\right.\\
    &\geq&c^{1-\a}\left(\frac{1}{(1-c)^{1-\a}}-1\right).
  \end{eqnarray*}
Hence, by Theorem \ref{aedt} (ii), or Corollary \ref{hecc} (i) in the case $\a=0$,
  \begin{eqnarray*}
     \e&\geq&c^{1-\a}\left(\frac{1}{(1-c)^{1-\a}}-1\right)\sum\limits_{m\leq 0}\sum\limits_{i\in\I}\lambda_i\Lambda_i\left((X\setminus C_i)\cap S^mA_m\right)\\
       &\geq&\frac{c^{1-\a}\left(1-(1-c)^{1-\a}\right)}{(1-c)^{1-\a}}\sum\limits_{m\leq 0}\Lambda\left((X\setminus C)\cap S^mA_m\right),\\
  \end{eqnarray*}
  which implies that
  \[0=\overline{\Lambda}_{X\setminus C}^{\a,1}(Q)=\underline{\Lambda}_{X\setminus C}^{\a,1}(Q),\]
  which, by \eqref{imsd}, is the assertion of (i).
  
  $(ii)$ It follows immediately by Lemma \ref{edml} that $\underline{\Lambda}_{A}^{1,1}(Q)=0$. Now, suppose $0<\Lambda_i(A_i)<1$ for all $i\in\I$. Then, by Lemma \ref{edml} and \eqref{simd}, $\overline{\Lambda}_{X\setminus\bigcup_{i\in\I}\O_i\setminus A_i}^{1,1}(Q)=\Lambda(Q)$, but, as one sees, $\overline{\Lambda}_{X\setminus\bigcup_{i\in\I}\O_i\setminus A_i}^{1,1}(Q)=\overline{\Lambda}_{A}^{1,1}(Q)$.
  
  $(iii)$ Let $1<\a\leq\g$. Then, for $(A_m)_{m\leq 0}\in\C^{\a,1}_\e(Q)$,
  \begin{eqnarray*}
  	&&\H_{\a}\left(\Lambda,\phi_0\right)(Q)+\e\\
  	&>&\sum\limits_{i\in\I}\lambda_i^\a\left[\sum\limits_{m\leq 0}\int\limits_{C_i^*\cap S^mA_m}Z_i^{\a-1}d\Lambda_i+\sum\limits_{m\leq 0}\int\limits_{(X\setminus C_i^*)\cap S^mA_m}Z_i^{\a-1}d\Lambda_i\right]\\
  	&\geq&\sum\limits_{i\in\I}\lambda_i^\a m_i^{\a-1}\sum\limits_{m\leq 0}\Lambda_i\left(C_i^*\cap S^mA_m\right)\\
  	    &&+\sum\limits_{i\in\I}\lambda_i^\a \eta_i(c^*)^{\a-1}\sum\limits_{m\leq 0}\Lambda_i\left(\left(X\setminus C_i^*\right)\cap S^mA_m\right)\\
  	&\geq&\sum\limits_{i\in\I}\lambda_i^\a m_i^{\a-1}\Lambda_i\left(Q_i\right)\\
  	 	&&+\sum\limits_{i\in\I}\lambda_i\left(\left(\lambda_i\eta_i(c^*)\right)^{\a-1}-\left(\lambda_i m_i\right)^{\a-1}\right)\sum\limits_{m\leq 0}\Lambda_i\left(\left(X\setminus C^*\right)\cap S^mA_m\right).
  \end{eqnarray*}
   Let $i\in\I$. Observe that
  \begin{eqnarray*}
   \left(\lambda_i\eta_i(c^*)\right)^{\a-1}-\left(\lambda_i m_i\right)^{\a-1}&=&\left\{\begin{array}{cc}
    	\left(\frac{\lambda_i m_i}{1-c^*}\right)^{\a-1}-\left(\lambda_i m_i\right)^{\a-1} & \mbox{if }\lambda_i m_i\geq c^*\\
    	\left(\frac{c^*}{1-c^*}\right)^{\a-1}-\left(\lambda_i m_i\right)^{\a-1} &\mbox{otherwise } 
    \end{array}\right.\\
    &\geq&{c^*}^{\a-1}\left(\frac{1}{(1-c^*)^{\a-1}}-1\right).
  \end{eqnarray*}
  Hence, by Theorem \ref{ec2},
  \begin{eqnarray*}
  	\e&>&{c^*}^{\a-1}\left(\frac{1}{(1-c^*)^{\a-1}}-1\right)\sum\limits_{m\leq 0}\Lambda\left((X\setminus C^*)\cap S^mA_m\right),
  \end{eqnarray*}
  which, by \eqref{imsd}, implies the assertion.
\hfill$\Box$

\begin{lemma}\label{ecd}
	For each $i\in\I$, let $C_i\in\A_0$ such that $C_i\subset\Omega_i$ and $C:=\bigcup_{i\in\I}C_i$. Let $Q\in\B$.  Suppose $\underline{\Lambda}_{C}^{\a,1}(Q)=\Lambda(Q)$ for some $\a\geq 0$. Then
	\[\overline{\Lambda_i}_{C_i}^{\a,1}(Q)=\underline{\Lambda_i}_{C_i}^{\a,1}(Q)=\Lambda_i(Q)\mbox{ and }\overline{\Lambda_i}_{X\setminus C_i}^{\a,1}(Q)=0\ \ \mbox{ for all }i\in\I.\]
\end{lemma}
{\it Proof.}
  Let $\e>0$ and $(A_m)_{m\leq 0}\in\C^{\a,1}_\e(Q)$. Let $k\in\I$. Then
\begin{eqnarray*}
	\underline{\Lambda}_{C,\e}^{\a,1}(Q)&\leq&\sum\limits_{i\in\I}\lambda_i\sum\limits_{m\leq 0}\Lambda_i\left(C_i\cap S^mA_m\right)\\
	&=&\sum\limits_{i\in\I}\lambda_i\sum\limits_{m\leq 0}\Lambda_i\left(A_m\right)-\sum\limits_{i\in\I}\lambda_i\sum\limits_{m\leq 0}\Lambda_i\left(\left(X\setminus C_i\right)\cap S^mA_m\right)\\
	&<&\Lambda\left(Q\right)+\e-\lambda_k\sum\limits_{m\leq 0}\Lambda_k\left(\left(X\setminus C_k\right)\cap S^mA_m\right).
\end{eqnarray*} 
Hence,
\[\underline{\Lambda}_{C}^{\a,1}(Q)\leq\Lambda\left(Q\right)-\lambda_k\overline{\Lambda_k}_{X\setminus C_k}^{\a,1}(Q).\]
Thus
\[\overline{\Lambda_k}_{X\setminus C_k}^{\a,1}(Q)=0.\]
Since
\[\Lambda_k(Q)\leq\underline{\Lambda_k}_{C_k}^{\a,1}(Q)+\overline{\Lambda_k}_{X\setminus C_k}^{\a,1}(Q),\]
it follows that
\begin{equation}\label{peka}
	\Lambda_k(Q)\leq\underline{\Lambda_k}_{C_k}^{\a,1}(Q)\mbox{ for all }k\in\I.
\end{equation}
On the other hand,
\[\sum\limits_{k\in\I}\lambda_k\underline{\Lambda_k}_{C_k,\e}^{\a,1}(Q)\leq\sum\limits_{k\in\I}\lambda_k\sum\limits_{m\leq 0}\Lambda_k\left(A_m\right)=\sum\limits_{m\leq 0}\Lambda\left(A_m\right)<\Lambda(Q)+\e.\]
Hence, by the Monotone Convergence Theorem,
\[\sum\limits_{k\in\I}\lambda_k\underline{\Lambda_k}_{C_k}^{\a,1}(Q)\leq\sum\limits_{k\in\I}\lambda_k\sum\limits_{m\leq 0}\Lambda_k\left(Q\right),\]
which together with \eqref{peka} implies that
\begin{equation}\label{reka}
	\underline{\Lambda_i}_{C_i}^{\a,1}(Q)=\Lambda_i(Q)\ \ \mbox{ for all }i\in\I.
\end{equation}
It remains only to show that 
\[\overline{\Lambda_k}_{C_k}^{\a,1}(Q)\leq\Lambda_i(Q).\]
Observe that
\begin{eqnarray*}
	\sum\limits_{i\in\I}\lambda_i\Lambda_i\left(Q\right)+\e&>&\sum\limits_{i\in\I}\lambda_i\sum\limits_{m\leq 0}\Lambda_i\left(A_m\right)\geq\sum\limits_{i\in\I}\lambda_i\sum\limits_{m\leq 0}\Lambda_i\left(C_i\cap S^mA_m\right)\\
	&=&\lambda_k\sum\limits_{m\leq 0}\Lambda_k\left(C_k\cap S^mA_m\right)+\sum\limits_{\I\owns i\neq k}\lambda_i\sum\limits_{m\leq 0}\Lambda_i\left(C_i\cap S^mA_m\right)\\
	&\geq&\lambda_k\sum\limits_{m\leq 0}\Lambda_k\left(C_k\cap S^mA_m\right)+\sum\limits_{\I\owns i\neq k}\lambda_i\underline{\Lambda_i}_{C_i,\e}^{\a,1}(Q).
\end{eqnarray*}
Hence,
\[\sum\limits_{i\in\I}\lambda_i\Lambda_i\left(Q\right)+\e\geq\lambda_k\overline{\Lambda_k}_{C_k,\e}^{\a,1}(Q)+\sum\limits_{\I\owns i\neq k}\lambda_i\underline{\Lambda_i}_{C_i,\e}^{\a,1}(Q),\]
which,  by the Monotone Convergence Theorem and \eqref{reka}, implies the remaining assertion.
\hfill$\Box$

In order to formulate the next theorem, we need the following definitions.

\begin{Definition}
	  For each $i\in\I$, let $m_i$ and $M_i$ be the $\Lambda_i$-essential infimum and the $\Lambda_i$-essential supremum of $Z_i$ respectively. Define
	\[\a(\Lambda|\phi_0):=\inf\left\{0<\b\leq 1|\ \sum\limits_{i\in\I}\lambda_i^{\b}m_i^{\b-1}<\infty\right\},\]
	\[\b(\Lambda|\phi_0):=\sup\left\{\b\geq 1|\ \sum\limits_{i\in\I}\lambda_i^{\b}M_i^{\b-1}<\infty\right\},\mbox{ and}\]
	\[\g(\Lambda|\phi_0):=\sup\left\{\b\geq 1|\ \int Z^{\b-1}d\Lambda<\infty\right\}.\]
\end{Definition}

Obviously, $\a(\Lambda|\phi_0)=0$ if $I$ is finite and all $m_i>0$. Also, $\a(\Lambda|\Lambda)=0$. If some $m_i=0$, then $\a(\Lambda|\phi_0)=1$. If some $M_i=\infty$, then $\b(\Lambda|\phi_0)=1$. If $I$ is finite and $M_i<\infty$ for all $i\in\I$, then $\b(\Lambda|\phi_0)=\infty$. Also, $\b(\Lambda|\Lambda)=\infty$.

\begin{center}
	\unitlength .8cm
	\begin{picture}(12,11)\thinlines
	\put(2,2){\vector(0,1){9}} 
	\put(0,0){\makebox(12,1){Fig. 1: Phase diagram of $\H^{\b,\a}(\Lambda,\phi_0)$}}
	\put(2,2){\vector(1,0){9}} 
	\put(1.7,1.8){\small $0$}
	\put(1.6,10.7){\small $\b$} 
	\put(10.8,1.7){\small $\a$} 
	\put(0.4,7.85){\small $\b(\Lambda|\phi_0)$}
	\put(0.4,9.9){\small $\g(\Lambda|\phi_0)$}
	\put(9.3,1.4){\small $\g(\Lambda|\phi_0)$}
	\put(0.4,3.7){\small $\a(\Lambda|\phi_0)$}
	\put(1.6,4.9){\small $1$} 
	\put(4.9,1.4){\small $1$}
	{ %\thicklines
		\put(1.9,5){\line(1,0){0.2}}
		%\put(2,2){\line(1,1){1.79}}
		\put(1.9,3.8){\line(1,0){0.2}}		
		\put(5,2){\line(0,1){3}}
		\put(4.98,5){\line(0,1){5}}
		\put(2,8){\line(1,0){2.99}}
		\put(10,2){\line(0,1){8.01}}
		\put(2,10){\line(1,0){7.99}}
		\put(2.15,6.03){$\sum\limits_{i\in \I} \lambda_i^\b M_i^{\b-1}\Lambda_i$}
		\put(6.15,6.03){$\sum\limits_{i\in \I} \lambda_i^\b m_i^{\b-1}\Lambda_i$}
		%{\thicklines\put(2,5){\line(1,0){7.99}}}
		%\put(2,2){\framebox(2.97,6){$\sum\limits_{i\in \I} \lambda_i^\b M_i^{\b-1}\Lambda_i$}}
		%\put(5.03,2){\framebox(4.97,6){$\sum\limits_{i\in \I} \lambda_i^\b m_i^{\b-1}\Lambda_i$}}
		%\put(3.79,3.8){\line(1,0){1.21}}
		\put(5.02,3.8){\line(1,0){4.98}}
		\put(3.1,8.85){$(\infty)$}
		\put(7.1,2.78){$(\infty)$}
	}
	\multiput(5.3125,9.755)(0.3125,0){15}{\circle*{0.001}}
	\multiput(5.3125,9.4425)(0.3125,0){15}{\circle*{0.001}}
	\multiput(5.3125,9.13)(0.3125,0){15}{\circle*{0.001}}
	\multiput(5.3125,8.8175)(0.3125,0){15}{\circle*{0.001}}
	\multiput(5.3125,8.505)(0.3125,0){15}{\circle*{0.001}}
	\multiput(5.3125,8.1925)(0.3125,0){15}{\circle*{0.001}}
	\multiput(5.3125,7.88)(0.3125,0){15}{\circle*{0.001}}
    \multiput(5.3125,7.5675)(0.3125,0){15}{\circle*{0.001}}
    \multiput(5.3125,7.255)(0.3125,0){15}{\circle*{0.001}}
    \multiput(5.3125,6.9425)(0.3125,0){15}{\circle*{0.001}}
    \multiput(5.3125,6.63)(0.3125,0){15}{\circle*{0.001}}
    \multiput(5.3125,6.3175)(0.3125,0){15}{\circle*{0.001}}
    \multiput(5.3125,6.005)(0.3125,0){15}{\circle*{0.001}}
    \multiput(5.3125,5.6925)(0.3125,0){15}{\circle*{0.001}}
    \multiput(5.3125,5.38)(0.3125,0){15}{\circle*{0.001}}
    \multiput(5.3125,5.0675)(0.3125,0){15}{\circle*{0.001}}
	\multiput(5.3125,4.755)(0.3125,0){15}{\circle*{0.001}}
	\multiput(5.3125,4.4425)(0.3125,0){15}{\circle*{0.001}}
	\multiput(5.3125,4.1125)(0.3125,0){15}{\circle*{0.001}}
	%\multiput(1.9,5)(0.8,0){4}{\line(1,0){0.2}}
	\multiput(1.9,3.8)(0.8,0){4}{\line(1,0){0.2}}
	%\multiput(2.4,2)(0.4,0){3}{\line(1,1){1.79}}
	%\put(3.6,2){\line(1,1){1.39}}
	%\put(4,2){\line(1,1){0.99}}
	%\put(4.4,2){\line(1,1){0.59}}
	%\put(4.8,2){\line(1,1){0.19}}
	
	\multiput(5.03,3.79)(0.4,0){8}{\line(1,-1){1.79}}
	\put(6.42,2){\line(-1,1){1.38}}
	\put(6.02,2){\line(-1,1){0.98}}
	\put(5.62,2){\line(-1,1){0.58}}
	\put(5.22,2){\line(-1,1){0.18}}
	\put(8.23,3.79){\line(1,-1){1.78}}
	\put(8.63,3.79){\line(1,-1){1.38}}
	\put(9.03,3.79){\line(1,-1){0.98}}
	\put(9.43,3.79){\line(1,-1){0.58}}
	\put(9.83,3.79){\line(1,-1){0.18}}
	
	\multiput(2.01,8)(0.4,0){3}{\line(1,1){1.985}}
	\put(3.21,8){\line(1,1){1.78}}
	\put(3.61,8){\line(1,1){1.37}}
	\put(4.01,8){\line(1,1){0.97}}
	\put(4.41,8){\line(1,1){0.57}}
	\put(4.81,8){\line(1,1){0.17}}
	\put(3.59,9.99){\line(-1,-1){1.59}}
	\put(3.19,9.99){\line(-1,-1){1.19}}
	\put(2.79,9.99){\line(-1,-1){0.79}}
	\put(2.39,9.99){\line(-1,-1){0.39}}
	\end{picture}
\end{center}

The following theorem, in particular, reveals an other kind of discontinuity of the function $(0,\g(\Lambda|\phi_0))\times(0,\g(\Lambda|\phi_0))\owns (\b,\a)\longmapsto\H^{\b,\a}(\Lambda,\phi_0)$ (in addition to those already established in Remark \ref{dmnr}) which can occur at the line $\a=1$ below  $\a(\Lambda|\phi_0)$ and above $\b(\Lambda|\phi_0)$ (see Fig. 1).
\begin{theo}\label{homa}
	For each $i\in\I$, let $m_i$ be the $\Lambda_i$-essential infimum and $M_i$ be the $\Lambda_i$-essential supremum of $Z_i$. Let $Q\in\B$.
	
	(a) Let $0<\a< 1$, or $\a=0$ if $[0,1]\owns\alpha\longmapsto \H_{\a}\left(\Lambda,\phi_0\right)(Q)$ is continuous at $0$.
	
	\hspace{2em}(i) If $\b\in(0,1]\cup\left(1,\b(\Lambda|\phi_0)\right)$, then
	\hspace{2em}\[\H^{\b,\a}(\Lambda,\phi_0)(Q)=\sum\limits_{i\in \I} \lambda_i^\b M_i^{\b-1}\Lambda_i(Q).\]
	\hspace{2em} The equality is also true for $\b=0$ if also $\phi_0\ll\Lambda$. 
	
	\hspace{2em}(ii) If $\b(\Lambda|\phi_0)\leq\b<\g(\Lambda|\phi_0)$ such that
	\hspace{2em}\[\limsup\limits_{L\to\infty}\left(\sum\limits_{i\in\I,\ \lambda_i M_i\leq L}\lambda_i^\b M_i^{\b-1}+L^{\b-1}\sum\limits_{i\in\I,\ \lambda_i M_i>L}\lambda_i\right)=\infty,\mbox{ then}\] 
	\hspace{2em}\[\H^{\b,\a}(\Lambda,\phi_0)(X)=\infty.\]
	
	(b) Let $1<\a<\g(\Lambda|\phi_0)$.
	
	\hspace{2em}(i) If $0<\b\leq\a(\Lambda|\phi_0)$, or ($\b=0$ and $\phi_0\ll\Lambda$) such that 
	\[\limsup\limits_{\d\to 0}\left(\frac{1}{\d}\right)^{1-\b}\sum\limits_{i\in\I,\ \lambda_i m_i<\d}\lambda_i=\infty,\mbox{ then }\H^{\b,\a}(\Lambda,\phi_0)(X)=\infty.\] 
	\hspace{2em}(ii) If $\a(\Lambda|\phi_0)<\b<\g(\Lambda|\phi_0)$, then
	\hspace{2em}\[\H^{\b,\a}(\Lambda,\phi_0)(Q)=\sum\limits_{i\in \I} \lambda_i^\b m_i^{\b-1}\Lambda_i(Q).\]
\end{theo}
{\it Proof.} 
$(a)$ $(i)$  The equality is obviously true for $\b=1$, by Lemma \ref{icp} (ii). Let us first consider the case 
\[0<\b<1,\mbox{ or }\b=0\mbox{ if also }\phi_0\ll\Lambda. \]
Then the part "$\geq$" of the equality obviously follows by Theorem \ref{aedt} (ii). It also proves the equality if $\a=\b$.

Let $\a<\b<1$.  Let $\e>0$ and $(H_m)_{m\leq 0}\in\dot\C^{\a,1}_\e(Q)$. Let $0<c<1$, $\t_i(c)$'s, $C_i$'s, and $C$ be those from Lemma \ref{ecal} (i). Then
\begin{eqnarray*}
   &&\sum\limits_{m\leq 0}\int\limits_{S^mH_m}Z^{\b}d\phi_0=\sum\limits_{i\in\I}\lambda_i^\b\sum\limits_{m\leq 0}\int\limits_{S^mH_m}Z_i^{\b-1}d\Lambda_i\\
   &\leq&\sum\limits_{i\in\I}\lambda_i^\b\left[\sum\limits_{m\leq 0}\int\limits_{(\O_i\setminus C_i)\cap S^mH_m}Z_i^{\b-1}d\Lambda_i+\t_i(c)^{\b-1}\overline{\dot\Lambda_i}^{\a,1}_{C_i,\e}(Q)\right]\\
   &\leq&\sum\limits_{m\leq 0}\int\limits_{(\bigcup_{i\in\I}\O_i\setminus C_i)\cap S^mH_m}Z^{\b-1}d\Lambda+\sum\limits_{i\in\I}\lambda_i^\b\t_i(c)^{\b-1}\overline{\dot\Lambda_i}^{\a,1}_{C_i,\e}(Q).
\end{eqnarray*}
Observe that $\Lambda((\bigcup_{i\in\I}\O_i\setminus C_i)\cap S^mH_m)=\Lambda(H_m)-\Lambda(C\cap S^mH_m)$ for all $m\leq 0$.
Then, as in the proof of Lemma \ref{pchm} (i), the concavity of $[0,\infty)\owns x\longmapsto x^{(\b-\a)/(1-\a)}$ implies that
\begin{eqnarray*}
&&\sum\limits_{m\leq 0}\int\limits_{(\bigcup_{i\in\I}\O_i\setminus C_i)\cap S^mH_m}Z^{\b-1}d\Lambda=\sum\limits_{m\leq 0}\int\limits_{S^mH_m}\left(1_{\bigcup_{i\in\I}\O_i\setminus C_i}Z^{1-\a}\right)^\frac{\b-\a}{1-\a}Z^\a d\phi_0\\
&\leq&\left(\sum\limits_{m\leq 0}\int\limits_{S^mH_m}Z^{\a}d\phi_0\right)^{1-\frac{\b-\a}{1-\a}}\left(\sum\limits_{m\leq 0}\Lambda\left(\bigcup_{i\in\I}\left(\O_i\setminus C_i\right)\cap S^mH_m\right)\right)^\frac{\b-\a}{1-\a}\\
&<&\left(\H_{\a}(\Lambda,\phi_0)(Q)+\e\right)^{1-\frac{\b-\a}{1-\a}}\left(\Lambda(Q)+\e-\sum\limits_{m\leq 0}\Lambda\left(C\cap S^mH_m\right)\right)^\frac{\b-\a}{1-\a}\\
&\leq&\left(\H_{\a}(\Lambda,\phi_0)(Q)+\e\right)^{1-\frac{\b-\a}{1-\a}}\left(\Lambda(Q)+\e-\underline{\Lambda}^{\a,1}_{C,\e}(Q)\right)^\frac{\b-\a}{1-\a}.
\end{eqnarray*}
Hence,
\begin{eqnarray}\label{haub0}
   \sum\limits_{m\leq 0}\int\limits_{S^mH_m}Z^{\b}d\phi_0&<&\left(\H_{\a}(\Lambda,\phi_0)(Q)+\e\right)^{1-\frac{\b-\a}{1-\a}}\left(\Lambda(Q)+\e-\underline{\Lambda}^{\a,1}_{C,\e}(Q)\right)^\frac{\b-\a}{1-\a}\nonumber\\
   &&+\sum\limits_{i\in\I}\lambda_i^\b\t_i(c)^{\b-1}\overline{\dot\Lambda_i}^{\a,1}_{C_i,\e}(Q).
\end{eqnarray}
Since\[\sum\limits_{i\in \I}\lambda_i^\b \t_i(c)^{\b-1}\overline{\dot\Lambda_i}^{\a,1}_{C_i,\e}(Q)\leq \sum\limits_{i\in \I}\lambda_i^\b \t_i(c)^{\b-1}<\infty,\]
 it follows, by Lemma \ref{icp} (i), the Monotone Convergence Theorem, Lemma \ref{ecal} (i), and Lemma \ref{ecd} that
\begin{eqnarray*}
  \H^{\b,\a}(\Lambda,\phi_0)(Q)&\leq&\sum\limits_{i\in\I}\lambda_i^\b\t_i(c)^{\b-1}\Lambda_i(Q)\\
  &\leq&(1-c)^{\b-1}\sum\limits_{i\in\I}\lambda_i^\b M_i^{\b-1}\Lambda_i(Q)+\left(\frac{c}{1-c}\right)^{1-\b}.
\end{eqnarray*}
Thus, letting $c\to 0$ completes the proof of $(a)$ $(i)$ for the case $\a\leq\b<1$.

Now, let $\b<\a<1$. Then, by the above and Theorem \ref{aedt} (ii), $\H^{\a,\b}\left(\Lambda,\phi_0\right)(Q)=\H_{\a}\left(\Lambda,\phi_0\right)(Q)$, but this implies, by Lemma \ref{hsl}, that $\H^{\b,\a}\left(\Lambda,\phi_0\right)(Q)=\H_{\b}\left(\Lambda,\phi_0\right)(Q)$. Thus, by Theorem \ref{aedt} (ii), this completes the proof of $(a)$ $(i)$ for $0<\b<1$, and $\b=0$ if $\phi_0\ll\Lambda$.

Finally, let us consider the case
\begin{equation*}
1<\b<\b(\Lambda|\phi_0).
\end{equation*}
In this case,
\[\sum\limits_{i\in \I}\lambda_i^\b M_i^{\b-1}<\infty.\]
In particular, this implies that $M_i<\infty$ for all $i\in\I$ (and $\int Z^{\b-1} d\Lambda<\infty$). For $B\in\B$, define
\[\overline{\Lambda}^{(\b)}(B):=\sum\limits_{i\in \I}\lambda_i^{\b} M_i^{\b-1}\Lambda_i\left(B\right).\]
Then $\overline{\Lambda}^{(\b)}$ is a finite measure on $\B$ which is absolutely continuous with respect to $\Lambda$. Let $\d>0$ such that $\overline{\Lambda}^{(\b)}(B)<\e$ whenever $\Lambda(B)<\d$ for some $B\in\B$. Let $(B_m)_{m\leq 0}\in\dot\C^{\a,1}_\d(Q)$. Then, since $\Lambda\left(\bigcup_{m\leq 0}B_m\setminus Q\right)<\d$,
\begin{eqnarray}\label{haub2}
\sum\limits_{m\leq 0}\int\limits_{S^mB_m}Z^{\b}d\phi_0&=&\sum\limits_{i\in \I}\lambda_i^\b\sum\limits_{m\leq 0}\int\limits_{S^mB_m}Z_i^{\b-1}d\Lambda_i\nonumber\\
&\leq&\sum\limits_{i\in \I}\lambda_i^\b M_i^{\b-1}\sum\limits_{m\leq 0}\Lambda_i\left(B_m\right)\nonumber\\
&\leq&\sum\limits_{i\in \I}\lambda_i^\b M_i^{\b-1}\Lambda_i\left(Q\right)+\sum\limits_{i\in \I}\lambda_i^{\b} M_i^{\b-1}\Lambda_i\left(\bigcup\limits_{m\leq 0}B_m\setminus Q\right)\nonumber\\
&<&\sum\limits_{i\in \I}\lambda_i^\b M_i^{\b-1}\Lambda_i\left(Q\right)+\e.
\end{eqnarray}
Hence, by Lemma \ref{icp} (i),
\[\H^{\b,\a}(\Lambda,\phi_0)(Q)\leq\sum\limits_{i\in \I}\lambda_i^\b M_i^{\b-1}\Lambda_i\left(Q\right).\]
On the other hand, for every $(C_m)_{m\leq 0}\in\C^{\a,1}_\e(Q)$,
\begin{eqnarray*}
	\sum\limits_{m\leq 0}\int\limits_{S^mC_m}Z^{\b}d\phi_0
	&\geq&\sum\limits_{i\in \I}\lambda_i^\b\t_i(c)^{\b-1}\sum\limits_{m\leq 0}\Lambda_i\left(C_i\cap S^mC_m\right)\\
	&\geq&\sum\limits_{i\in \I}\lambda_i^\b\t_i(c)^{\b-1}\underline{\Lambda_i}_{C_i,\e}^{\a,1}(Q),
\end{eqnarray*}
and therefore,
\begin{equation}\label{halb}
  \H^{\b,\a,1}_\e(\Lambda,\phi_0)(Q)\geq\sum\limits_{i\in \I}\lambda_i^\b\t_i(c)^{\b-1}\underline{\Lambda_i}_{C_i,\e}^{\a,1}(Q).
\end{equation}
Hence, by the Monotone Convergence Theorem, Lemma \ref{ecal} (i), Lemma \ref{ecd}, and Lemma \ref{icp} (iii),
\[\H^{\b,\a}(\Lambda,\phi_0)(Q)\geq\sum\limits_{i\in \I}\lambda_i^\b\t_i(c)^{\b-1}\Lambda_i(Q)\geq(1-c)^{\b-1}\sum\limits_{i\in \I,\  \lambda_i M_i\leq 1/c}\lambda_i^\b M_i^{\b-1}\Lambda_i(Q).\]
Thus, letting $c\to 0$ completes the proof of (a) (i).

$(a)$ $(ii)$ The hypothesis of (a) (ii) implies that $\b>1$. Observe that, the same way as in \eqref{halb},
\[\H^{\b,\a,1}(\Lambda,\phi_0)(Q)\geq\sum\limits_{i\in \I}\lambda_i^\b\t_i(c)^{\b-1}\Lambda_i(Q),\]
 and therefore,
\begin{eqnarray*}
  &&\H^{\b,\a,1}(\Lambda,\phi_0)(X)\\
  &\geq&(1-c)^{\b-1}\left(\sum\limits_{i\in \I,\ \lambda_i M_i\leq 1/c}\lambda_i^\b M_i^{\b-1}+\left(\frac{1}{c}\right)^{\b-1}\sum\limits_{i\in \I,\ \lambda_i M_i> 1/c}\lambda_i\right).
\end{eqnarray*}
Hence, by the hypothesis, $\H^{\b,\a,1}(\Lambda,\phi_0)(X)=\infty$. If $\H^{\b,\a}(\Lambda,\phi_0)(X)<\infty$, then we obtain a contradiction by Lemma \ref{icp} (iii).

$(b)$ $(i)$ The hypothesis of (b) (i) implies that $\b<1$. Let $0<c^*<1$. For each $i\in\I$, let $\eta_i(c^*)$, $C^*_i$, and $C^*$ be defined as in Lemma \ref{ecal} (iii).  Let $\e>0$ and $(G_m)_{m\leq 0}\in\C^{\a,1}_\e(Q)$. Then
\begin{eqnarray*}
	&&\sum\limits_{m\leq 0}\int\limits_{S^mG_m}Z^{\b}d\phi_0=\sum\limits_{i\in\I}\lambda_i^\b\sum\limits_{m\leq 0}\int\limits_{S^mG_m}Z_i^{\b-1}d\Lambda_i\\
	&\geq&\sum\limits_{i\in\I}\lambda_i^\b\eta_i(c^*)^{\b-1}\sum\limits_{m\leq 0}\Lambda_i\left(C^*_i\cap{S^mG_m}\right)\geq\sum\limits_{i\in\I}\lambda_i^\b\eta_i(c^*)^{\b-1}\underline{\Lambda}_{C^*_i,\e}^{\a,1}(Q)\\
	&\geq&\left(\frac{1-c^*}{c^*}\right)^{1-\b}\sum\limits_{i\in\I,\ \lambda_i m_i<c^*}\lambda_i\underline{\Lambda}_{C^*_i,\e}^{\a,1}(Q).
\end{eqnarray*}
Hence, by the Monotone Convergence Theorem, Lemma \ref{ecal} (iii), and Lemma \ref{ecd},
\[\H^{\b,\a,1}(\Lambda,\phi_0)(Q)\geq\left(\frac{1-c^*}{c^*}\right)^{1-\b}\sum\limits_{i\in\I,\ \lambda_i m_i<c^*}\lambda_i\Lambda_i(Q).\]
Suppose $\H^{\b,\a}(\Lambda,\phi_0)(X)<\infty$. Then, by Lemma \ref{icp} (iii),
\[\H^{\b,\a}(\Lambda,\phi_0)(X)\geq\left(\frac{1-c^*}{c^*}\right)^{1-\a}\sum\limits_{i\in\I,\ \lambda_i m_i<c^*}\lambda_i\]
for all $0<c^*<1$, but this contradicts the assumption. Thus, the assertion of (b) (i) is correct.

$(b)$ $(ii)$ Let
\[\a(\Lambda|\phi_0)<\b<\g(\Lambda|\phi_0).\]
Let us first consider the case $\a(\Lambda|\phi_0)<\b<1$.  This implies that $\sum_{i\in\I}\lambda_i^{\b}m_i^{\b-1}<\infty$. Let $\e>0$. Let $\d>0$ such that $\sum_{i\in\I}\lambda_i^{\b}m_i^{\b-1}\Lambda_i(B)<\e$ whenever $\Lambda(B)<\d$ for some $B\in\B$. Let $(D_m)_{m\leq 0}\in\dot\C^{\a,1}_\d(Q)$. Then
\begin{eqnarray*}
	\H^{\b,\a}_\d(\Lambda,\phi_0)(Q)&\leq&\sum\limits_{i\in \I}\lambda_i^\b\sum\limits_{m\leq 0}\int\limits_{S^mD_m}Z_i^{\b-1}d\Lambda_i
	\leq\sum\limits_{i\in \I}\lambda_i^\b m_i^{\b-1}\sum\limits_{m\leq 0}\Lambda_i\left(D_m\right)\\
	&\leq&\sum\limits_{i\in \I}\lambda_i^\b m_i^{\b-1}\Lambda_i(Q)+\sum\limits_{i\in \I}\lambda_i^\b m_i^{\b-1}\Lambda_i\left(\bigcup\limits_{m\leq 0}D_m\setminus Q\right)\\
	&\leq&\sum\limits_{i\in \I}\lambda_i^\b m_i^{\b-1}\Lambda_i(Q)+\e.
\end{eqnarray*}
Hence, 
\[\H^{\b,\a}(\Lambda,\phi_0)(Q)\leq\sum\limits_{i\in \I}\lambda_i^\b m_i^{\b-1}\Lambda_i(Q).\]

Now, let $(E_m)_{m\leq 0}\in\C^{\a,1}_\e(Q)$. Then
\begin{eqnarray*}
	&&\sum\limits_{m\leq 0}\int\limits_{S^mE_m}Z^{\b}d\phi_0=\sum\limits_{i\in \I}\lambda_i^\b\sum\limits_{m\leq 0}\int\limits_{S^mE_m}Z_i^{\b-1}d\Lambda_i\\
	&\geq&\sum\limits_{i\in\I}\lambda_i^\b \eta_i\left(c^*\right)^{\b-1}\sum\limits_{m\leq 0}\Lambda_i\left(C^*_i\cap S^mE_m\right)\\
	&\geq&\left(1-c^*\right)^{1-\b}\sum\limits_{i\in \I,\ \lambda_i m_i\geq c^*}\lambda_i^\b m_i^{\b-1}\underline{\Lambda_i}_{C^*_i,\e}^{\a,1}(Q).
\end{eqnarray*}
Hence, by Lemma \ref{ecal} (iii), Lemma \ref{ecd}, Lemma \ref{icp} (iii), and the Monotone Convergence Theorem,
\[\H^{\b,\a}(\Lambda,\phi_0)(Q)\geq\sum\limits_{i\in \I}\lambda_i^\b m_i^{\b-1}\Lambda_i(Q).\]
This completes the proof of the case $\a(\Lambda|\phi_0)<\b<1$.

The case $\b=1$ follows again by Lemma \ref{icp} (ii).

Finally, let us consider the case \[1<\b<\g(\Lambda|\phi_0).\]
 Since it implies that
$\int Z^{\b-1} d\Lambda<\infty$, it follows by Theorem \ref{ec2} that
\[\H^{\b,\a}(\Lambda,\phi_0)(Q)\geq\sum\limits_{i\in \I}\lambda_i^\b m_i^{\b-1}\Lambda_i(Q).\]
Theorem \ref{ec2} implies also the equality if $\a=\b$.

Now, let $1<\b<\a<\g(\Lambda|\phi_0)$. Let $\d>0$ such that $\int_B Z^{\b-1} d\Lambda<\e$ whenever $\Lambda(B)<\d$ for some $B\in\B$. Let $(F_m)_{m\leq 0}\in\dot\C^{\a,1}_\d(Q)$. Then
\begin{eqnarray*}
    &&\H^{\b,\a}_\d(\Lambda,\phi_0)(Q)\leq\sum\limits_{i\in \I}\lambda_i^\b\sum\limits_{m\leq 0}\int\limits_{S^mF_m}Z_i^{\b-1}d\Lambda_i\\
	&\leq&\sum\limits_{i\in\I}\lambda_i^\b\left[\eta_i(c^*)^{\b-1}\sum\limits_{m\leq 0}\Lambda_i\left(S^mF_m\right)+\sum\limits_{m\leq 0}\int\limits_{\left(\O_i\setminus C^*_i\right)\cap S^mF_m}Z_i^{\b-1}d\Lambda_i\right]\\
	&\leq&\sum\limits_{i\in \I}\lambda_i^\b\eta_i(c^*)^{\b-1}\Lambda_i(Q)+\sum\limits_{i\in \I}\lambda_i^{\b}\eta_i(c^*)^{\b-1}\Lambda_i\left(\bigcup\limits_{m\leq 0}F_m\setminus Q\right)\\
	&&+\sum\limits_{m\leq 0}\int\limits_{\left(\bigcup_{i\in\I}\O_i\setminus C^*_i\right)\cap S^mF_m}Z^{\b-1}d\Lambda.	
\end{eqnarray*}
Observe that
\begin{eqnarray*}
    &&\sum\limits_{i\in\I}\lambda_i^{\b}\eta_i(c^*)^{\b-1}\Lambda_i\left(\bigcup\limits_{m\leq 0}F_m\setminus Q\right)\\
    &=&\left(\frac{1}{1-c^*}\right)^{\b-1}\sum\limits_{i\in\I,\ \lambda_i m_i\geq c^*}\lambda_i^\b m_i^{\b-1}\Lambda_i\left(\bigcup\limits_{m\leq 0}F_m\setminus Q\right)\\
    &&+\left(\frac{c^*}{1-c^*}\right)^{\b-1}\sum\limits_{i\in\I,\ \lambda_i m_i< c^*}\lambda_i\Lambda_i\left(\bigcup\limits_{m\leq 0}F_m\setminus Q\right)\\
    &\leq&\left(\frac{1}{1-c^*}\right)^{\b-1}\sum\limits_{i\in\I}\lambda_i\int\limits_{\bigcup\limits_{m\leq 0}F_m\setminus Q} \left(\lambda_i Z_i\right)^{\b-1}d\Lambda_i+\d\left(\frac{c^*}{1-c^*}\right)^{\b-1}\\
    &<&\e\left(\frac{1}{1-c^*}\right)^{\b-1}+\d\left(\frac{c^*}{1-c^*}\right)^{\b-1}.
\end{eqnarray*}
Also, the same way as in the proof of (a) (i),
\begin{eqnarray*}
  &&\sum\limits_{m\leq 0}\int\limits_{\left(\bigcup_{i\in\I}\O_i\setminus C^*_i\right)\cap S^mF_m}Z^{\b-1}d\Lambda=\sum\limits_{m\leq 0}\int\limits_{S^mF_m}\left(1_{\bigcup_{i\in\I}\O_i\setminus C^*_i}Z^{1-\a}\right)^\frac{\a-\b}{\a-1}Z^\a d\phi_0\\
  &\leq&\left(\sum\limits_{m\leq 0}\int\limits_{S^mF_m}Z^{\a}d\phi_0\right)^{1-\frac{\a-\b}{\a-1}}\left(\sum\limits_{m\leq 0}\Lambda\left(\left(\bigcup_{i\in\I}\O_i\setminus C^*_i\right)\cap S^mF_m\right)\right)^\frac{\a-\b}{\a-1}\\
  &<&\left(\H_{\a}(\Lambda,\phi_0)(Q)+\d\right)^{1-\frac{\a-\b}{\a-1}}\left(\Lambda(Q)+\d-\underline{\Lambda}^{\a,1}_{C^*,\d}(Q)\right)^\frac{\a-\b}{\a-1}.
\end{eqnarray*}
Hence, by Lemma \ref{ecal} (iii),
\begin{eqnarray*}
  &&\H^{\b,\a}(\Lambda,\phi_0)(Q)\leq\sum\limits_{i\in \I}\lambda_i^\b\eta_i(c^*)^{\b-1}\Lambda_i(Q)\\
  &\leq&\left(\frac{1}{1-c^*}\right)^{\b-1}\sum\limits_{i\in \I}\lambda_i^\b m_i^{\b-1}\Lambda_i(Q)+\left(\frac{c^*}{1-c^*}\right)^{\b-1}\sum\limits_{i\in \I,\ \lambda_i m_i< c^*}\lambda_i.
\end{eqnarray*}
Thus, letting $c^*\to 0$ implies
\[\H^{\b,\a}(\Lambda,\phi_0)(Q)\leq\sum\limits_{i\in \I}\lambda_i^\b m_i^{\b-1}\Lambda_i(Q)\]
and completes the proof of (b) (ii) for $1<\b<\a<\g(\Lambda|\phi_0)$.

Finally, the case $1<\a<\b<\g(\Lambda|\phi_0)$ follows then by the above, Theorem \ref{ec2}, and Lemma \ref{hsl}.
\hfill$\Box$

\subsection{A computation of $[0, \g)\owns\a\longmapsto\J_{\a}\left(\Lambda,\phi_0\right)$}\label{homs2}

In this subsection, we compute the function $[0, \b(\Lambda|\phi_0))\owns\a\longmapsto\J_{\a}\left(\Lambda,\phi_0\right)$ under the condition of the equivalence of $\phi_0$ and $\Lambda$.

First, we need to clarify the equivalence of the definition of $\J_{\a}\left(\Lambda,\phi_0\right)$ to that over the disjoint covers, as it is not covered directly by Lemma \ref{del}.
\begin{Definition}
	Let $\b\geq 0$, $Q\in\P(X)$ and $\e>0$. Define
	\[\dot\J_{\b,\e}\left(\Lambda,\phi_0\right)(Q):=\sup\limits_{(A_m)_{m\leq 0}\in\dot\C^{0,1}_{\e}(Q)}\sum\limits_{m\leq 0}\int\limits_{S^mA_m}Z^\b d\phi_0\ \ \ \mbox{ and}\]
	\[\dot\J_{\b}\left(\Lambda,\phi_0\right)(Q):=\lim\limits_{\e\to 0}\dot\J_{\b,\e}\left(\Lambda,\phi_0\right)(Q).\]
\end{Definition}

\begin{lemma}\label{lhmd}
   Let $\b\in[0,1]\cup (1,\b(\Lambda|\phi_0))$ and $Q\in\P(X)$. Then
   \[\dot\J_{\b}\left(\Lambda,\phi_0\right)(Q)=\J_{\b}\left(\Lambda,\phi_0\right)(Q).\]
\end{lemma}
{\it Proof.} 
  The assertion in the case $0\leq\b\leq 1$ follows by Remark \ref{Jccr} and Lemma \ref{del} the same way as below.
  
  Let $1<\b< \b(\Lambda|\phi_0)$. For each $i\in\I$, let $M_i$ be the $\Lambda_i$-essential supremum of $Z_i$. Then
  \[\sum\limits_{i\in \I}\lambda_i^\b M_i^{\b-1}<\infty,\mbox{ and therefore, }\int Z^{\b-1}d\Lambda<\infty.\]
  For $\e>0$, define 
  \[\mathcal{N}_{\b,\e}(Q):=\inf\limits_{(A_m)_{m\leq 0}\in\C^{0,1}_{\e}(Q)}\sum\limits_{m\leq 0}\sum\limits_{i\in\I}\lambda_i^{\b}\int\limits_{S^mA_m}\left(M_i^{\b-1}-Z_i^{\b-1}\right) d\Lambda_i,\mbox{ and }\]
  \[\mathcal{N}_{\b}(Q):=\lim\limits_{\e\to 0}\mathcal{N}_{\b,\e}(Q).\]
  By restricting the infimum on disjoint covers, one immediately sees that $\mathcal{N}_{\b}(Q)<\infty$. Also, by Lemma \ref{del}, $\mathcal{N}_{\b}(Q)=\dot{\mathcal{N}}_{\b}(Q)$. One easily verifies (the same way as Lemma \ref{tlb} (ii)) that
  \[\mathcal{N}_{\b}(Q)=\sum\limits_{i\in \I}\lambda_i^\b M_i^{\b-1}\Lambda_i(Q)-\J_\b\left(\Lambda,\phi_0\right)(Q).\]
  Let $(A_m)_{m\leq 0}\in\C^{0,1}_{\e}(Q)$ such that 
  \[\sum\limits_{m\leq 0}\sum\limits_{i\in\I}\lambda_i^{\b}\int\limits_{S^mA_m}\left(M_i^{\b-1}-Z_i^{\b-1}\right) d\Lambda_i<\mathcal{N}_{\b}(Q) + \e.\]
   Define $B_0:=A_0$ and $B_m:=A_m\setminus(A_{m+1}\cup...\cup A_0)$ for all $m\leq -1$. Then, since $B_m\subset A_m$ for all $m\leq 0$ and $\bigcup_{m\leq 0}B_m=\bigcup_{m\leq 0}A_m$, one easily checks that $(B_m)_{m\leq 0}\in\dot\C^{0,1}_{\e}(Q)$, and therefore, 
  \begin{eqnarray*}
        &&\sum\limits_{i\in \I}\lambda_i^\b M_i^{\b-1}\Lambda_i(Q)-\dot\J_{\b,\e}\left(\Lambda,\phi_0\right)(Q)\leq\sum\limits_{m\leq 0}\int\limits_{S^mB_m}\sum\limits_{i\in\I}\lambda_i^{\b}\left(M_i^{\b-1}-Z_i^{\b-1}\right) d\Lambda_i\\
        &<&\sum\limits_{i\in \I}\lambda_i^\b M_i^{\b-1}\Lambda_i(Q)-\J_\b\left(\Lambda,\phi_0\right)(Q) + \e,
  \end{eqnarray*}
  which implies that 
  \[\dot\J_{\b}\left(\Lambda,\phi_0\right)(Q)\geq\J_{\b}\left(\Lambda,\phi_0\right)(Q).\]
  The converse inequality is obvious.
\hfill$\Box$

\begin{theo}\label{homa2}
	For each $i\in\I$, let $M_i$ be the $\Lambda_i$-essential supremum of $Z_i$. Suppose $\phi_0\ll\Lambda$. Let $Q\in\B$. Then
	\[\J_{\b}(\Lambda,\phi_0)(Q)=\sum\limits_{i\in \I} \lambda_i^\b M_i^{\b-1}\Lambda_i(Q)\mbox{ for all }\b\in[0,1]\cup(1,\b(\Lambda|\phi_0)).\]
\end{theo}
{\it Proof.} 
Since $\J_{\b}\left(\Lambda,\phi_0\right)(Q)\geq\H^{\b,0}\left(\Lambda,\phi_0\right)(Q)$ for all $\b\geq 0$, the part '$\geq$' follows by Theorem \ref{homa} $(a)$ $(i)$. Let $\e>0$ and $(A_m)_{m\leq 0}\in\dot\C^{0,1}_\e(Q)$.
 
First, let
\begin{equation*}
0\leq\b<1.
\end{equation*}
Let $0<c<1$. For each $i\in\I$, let $\t_i(c)$, $C_i$, and $C$ be defined as in Lemma \ref{ecal} (i). Then, by \eqref{haub0}, for $\a=0$,
\begin{eqnarray*}
\sum\limits_{m\leq 0}\int\limits_{S^mA_m}Z^{\b}d\phi_0&<&\left(\H_{\a}(\Lambda,\phi_0)(Q)+\e\right)^{1-\b}\left(\Lambda(Q)+\e-\underline{\Lambda}^{0,1}_{C,\e}(Q)\right)^\b\\
&&+\sum\limits_{i\in\I}\lambda_i^\b\t_i(c)^{\b-1}\overline{\dot\Lambda_i}^{0,1}_{C_i,\e}(Q).
\end{eqnarray*}
Hence, by Lemma \ref{lhmd}, the same way as in the proof of Theorem \ref{homa} $(a)$ $(i)$, it follows that 
\[\J_{\b}\left(\Lambda,\phi_0\right)(Q)\leq\sum\limits_{i\in \I} \lambda_i^\b M_i^{\b-1}\Lambda_i(Q).\]

Obviously, the assertion is correct for $\b=1$. 

Finally, let 
\begin{equation*}
	1<\b<\b(\Lambda|\phi_0).
\end{equation*}
In this case, it follows, by \eqref{haub2} and  Lemma \ref{lhmd}, that again
\[\J_{\b}\left(\Lambda,\phi_0\right)(Q)\leq\sum\limits_{i\in \I} \lambda_i^\b M_i^{\b-1}\Lambda_i(Q).\]
\hfill$\Box$

\subsection{An example}

Now, we will give an example which enables us to learn more about the dynamical measure theory.
\begin{Example}\label{2Mce}
	Consider Example \ref{e1} in the case $N=2$. Let $P$ be irreducible with the invariant probability measure $\mu:=(\mu\{1\}:=\mu_1,\mu\{2\}:=\mu_2)$, and $\Lambda$ be given by
	\[\Lambda\left({_0[i_1,...,i_n]}\right):=\mu\{i_1\}p_{i_1i_2}...p_{i_{n-1}i_n}\]
	for all $_0[i_1,...,i_n]\subset X=\{1,2\}^{\Z}$ and $n\geq0$. Let $\nu\{1\}:=\nu_1>0$ and $\nu\{2\}:=\nu_2>0$ such that $\nu_1\neq\mu_1$. Observe that, in this case,
	\begin{eqnarray*}
		Z(\sigma)=\frac{\mu_{\sigma_0}}{\nu_{\sigma_0}}\ \ \ \mbox{ for }\phi_0\mbox{-a.e. }\sigma\in X.
	\end{eqnarray*} 
	Let $Q\in\B$ and $\a\geq 0$. For $i\in\{1,2\}$, let us abbreviate
	\[\underline{\Lambda}^{\a,1}_{i}:=\underline{\Lambda}^{\a,1}_{_0[i]},\mbox{ and } \overline{\Lambda}^{\a,1}_{i}:=\overline{\Lambda}^{\a,1}_{_0[i]}.\]
	By Lemma \ref{ecal} (i) and \eqref{imsd},
	\begin{equation}\label{mha1}
	\underline{\Lambda}^{\a,1}_{1}(Q)=\overline{\Lambda}^{\a,1}_{1}(Q)=\left\{\begin{array}{cc}
	\Lambda(Q) &\mbox{if } \nu_1<\mu_1\\
	0 & \mbox{if }\nu_1>\mu_1
	\end{array}\right.\mbox{ for all }0\leq\a< 1
	\end{equation}
	 (e.g. if $\nu_1<\mu_1$, choosing $c=(\nu_1/\mu_1)(\mu_1-\nu_1)/(1-\nu_1)$ in Lemma \ref{ecal} (i) gives $C={_0[1]}$).   By Lemma \ref{ecal} (ii), $\underline{\Lambda}^{1,1}_{1}(Q)=0$, and $\overline{\Lambda}^{1,1}_{1}(Q)=\Lambda(Q)$. And by Lemma \ref{ecal} (iii) and \eqref{imsd},
	\begin{equation}\label{mha2}
	\underline{\Lambda}^{\a,1}_{1}(Q)=\overline{\Lambda}^{\a,1}_{1}(Q)=\left\{\begin{array}{cc}
	0 &\mbox{if } \nu_1<\mu_1\\
	\Lambda(Q) & \mbox{if }\nu_1>\mu_1
	\end{array}\right.\mbox{ for all }\a>1
	\end{equation}
	(e.g. if $\nu_1<\mu_1$, choosing $c^*=(\mu_2/\nu_2)(\nu_2-\mu_2)/(1-\mu_2)$ in Lemma \ref{ecal} (iii) gives $C^*={_0[2]}$). 
	
	However, in this example, these measures can be also computed in a simple algebraic way. One easily sees, by \eqref{hmvi}, that
	\[\Lambda(Q)=\underline{\Lambda}^{\a,1}_{1}(Q)+\overline{\Lambda}^{\a,1}_{2}(Q),\]
	\[\Lambda(Q)=\overline{\Lambda}^{\a,1}_{1}(Q)+\underline{\Lambda}^{\a,1}_{2}(Q),\]
	\[\H_\a(\Lambda,\phi_0)(Q)=\left(\frac{\mu_1}{\nu_1}\right)^{\a-1}\underline{\Lambda}^{\a,1}_{1}(Q)+\left(\frac{\mu_2}{\nu_2}\right)^{\a-1}\overline{\Lambda}^{\a,1}_{2}(Q),\mbox{ and}\]
	\[\H_\a(\Lambda,\phi_0)(Q)=\left(\frac{\mu_1}{\nu_1}\right)^{\a-1}\overline{\Lambda}^{\a,1}_{1}(Q)+\left(\frac{\mu_2}{\nu_2}\right)^{\a-1}\underline{\Lambda}^{\a,1}_{2}(Q).\]
	This implies that, for $\a\neq 1$,
	\begin{equation*}
	\underline{\Lambda}^{\a,1}_{1}(Q)=\overline{\Lambda}^{\a,1}_{1}(Q)=\frac{\left(\frac{\mu_2}{\nu_2}\right)^{\a-1}\Lambda(Q)-\H_\a(\Lambda,\phi_0)(Q)}{\left(\frac{\mu_2}{\nu_2}\right)^{\a-1}-\left(\frac{\mu_1}{\nu_1}\right)^{\a-1}}.
	\end{equation*}
	Thus, by Theorem \ref{aedt} (ii) and Theorem \ref{ec2}, it follows \eqref{mha1} and \eqref{mha2}.
	
	Now, we can also compute $\H^{\b,\a}(\Lambda,\phi_0)(Q)$ in this example in a simple way. Since, for all $i\in\{1,2\}$ and $\a\neq 1$, $\underline{\Lambda}^{\a,1}_{i}(Q)=\overline{\Lambda}^{\a,1}_{i}(Q)$, let us denote this number by $\underline{\overline\Lambda}^{\a,1}_{i}(Q)$.  Then, as one easily sees, for every $\b\geq 0$,
    \begin{equation*}
    \H^{\b,\a}(\Lambda,\phi_0)(Q)=\left(\frac{\mu_1}{\nu_1}\right)^{\b-1}\underline{\overline\Lambda}^{\a,1}_{1}(Q)+\left(\frac{\mu_2}{\nu_2}\right)^{\b-1}\underline{\overline\Lambda}^{\a,1}_{2}(Q).
    \end{equation*}
    Hence, for every $\b\geq 0$, by \eqref{mha1} and \eqref{imsd},
    \begin{eqnarray*}
    \H^{\b,\a}(\Lambda,\phi_0)(Q)=\max\left\{\frac{\mu_1}{\nu_1},\frac{\mu_2}{\nu_2}\right\}^{\b-1}\Lambda(Q)\ \ \mbox{ for all }0\leq\a< 1
    \end{eqnarray*}
	 and, by \eqref{mha2} and \eqref{imsd}, 
	\begin{eqnarray*}
	\H^{\b,\a}(\Lambda,\phi_0)(Q)=\min\left\{\frac{\mu_1}{\nu_1},\frac{\mu_2}{\nu_2}\right\}^{\b-1}\Lambda(Q)\ \ \mbox{ for all }\a>1.
	\end{eqnarray*}
	
	We can use this example also to answer the open question whether $\Phi$ remains the same if one takes $\phi_0\circ S^{-1}$ for the initial measure on $\A_0$, instead of $\phi_0$, which is equivalent, in our example, to taking for the construction of $\phi_0$ the initial measure on $\{1,2\}$ corresponding to the next step of the Markov process (open question (1) in \cite{Wer10}, p. 17 (p. 22 in the arXiv version)). Let $\Phi_{(-1)}$ denote $\Phi$ with the initial measure $\nu P$, instead of $\nu$. (By Lemma 4 in \cite{Wer10}, $\Phi_{(-1)}\geq\Phi$.) By Theorem \ref{aedt} (ii),
	\[\Phi_{(-1)}(X)=\min\left\{\frac{\nu_1p_{11}+\nu_2p_{21}}{\mu_1},\frac{\nu_1p_{12}+\nu_2p_{22}}{\mu_2}\right\}.\]
	Suppose $\Phi(X)=\nu_2/\mu_2$, i.e. $\nu_1/\mu_1>\nu_2/\mu_2$.  A simple computation shows that
	\[\frac{\nu_1p_{11}+\nu_2p_{21}}{\mu_1}-\frac{\nu_2}{\mu_2}=p_{11}\left(\frac{\nu_1}{\mu_1}-\frac{\nu_2}{\mu_2}\right),\mbox{ and}\]
	\[ \frac{\nu_1p_{12}+\nu_2p_{22}}{\mu_2}-\frac{\nu_2}{\mu_2}=p_{21}\left(\frac{\nu_1}{\mu_1}-\frac{\nu_2}{\mu_2}\right).\]
	Thus
	\[\Phi_{(-1)}(X)>\Phi(X)\]
	if all entries of $P$ are positive.
\end{Example}

\subsection{On the initial measure of moved covers}

There is something else which we can learn about the dynamical measure theory at this point. One might have the temptation, in the search for a lower bound for $\Phi$, to proceed straightforward through $\sum_{m\leq 0}\phi_0(S^mA_m)\geq\phi_0(\bigcup_{m\leq 0}S^mA_m)$, particularly because of the well-known Chung-Erd\"{o}s inequality. (In fact, it is difficult to find a partition $(A_m)_{m\leq 0}\in\C(\{0,1\}^\Z)$ with pen and paper such that the $\{1/2,1/2\}$-Bernoulli measure of the union at the right-hand side is less than $1/2$ where $S:\{0,1\}^\Z\longmapsto\{0,1\}^\Z$ is the left shift map.) We show now that this would not work even if $\phi_0$ is a Bernoulli measure.

As a consequence of Lemma \ref{edml}, we obtain the following.

\begin{prop}
	Suppose $\Lambda$ is a non-atomic, ergodic probability measure. Then
	\[\inf\limits_{(A_m)_{m\leq 0}\in\C(X)}\Lambda\left(\bigcup\limits_{m\leq 0}S^mA_m\right)=0.\]
\end{prop}
{\it Proof.} Let $\e>0$ and $C\in\A_0$. Define
\[\underrightarrow{\Lambda}_{C,\e}(X):=\inf\limits_{(A_m)_{m\leq 0}\in\C^{1}_{\e}(X)}\Lambda\left(C\cap\bigcup\limits_{m\leq 0} S^mA_m\right),\ \underrightarrow{\Lambda}_{C}(X):=\lim\limits_{\e\to 0}\underrightarrow{\Lambda}_{C,\e}(X),\]
\[\overrightarrow{\Lambda}_{C,\e}(X):=\sup\limits_{(A_m)_{m\leq 0}\in\C^{1}_{\e}(X)}\Lambda\left(C\cap\bigcup\limits_{m\leq 0} S^mA_m\right),\mbox{ and } \overrightarrow{\Lambda}_{C}(X):=\lim\limits_{\e\to 0}\overrightarrow{\Lambda}_{C,\e}(X).\]
Then, obviously,
\[\underrightarrow{\Lambda}_{C}(X)\leq\underline{\Lambda}^\Lambda\left(1_{C}\right)(X),\mbox{ and }\overrightarrow{\Lambda}_{C}(X)\leq\Lambda(C).\]
Also, one readily sees that
\[\underrightarrow{\Lambda}_{X}(X)\leq\underrightarrow{\Lambda}_{X\setminus C}(X)+\overrightarrow{\Lambda}_{C}(X).\]
Then, by Lemma \ref{edml}, for every $C\in\A_0$ such that $0<\Lambda(C)<1$,
\[\underrightarrow{\Lambda}_{X}(X)\leq\Lambda(C).\]
Thus, since $\Lambda$ is non-atomic, $\underrightarrow{\Lambda}_{X}(X)=0$, which implies the assertion.
\hfill$\Box$


\begin{thebibliography}{9}	
	
\bibitem{B}   V. Bogachev, Measure theory Vol. I.,
\emph{Springer} (2007).	
	
\bibitem{Sl} W. Slomczynski, Dynamical entropy, Markov operators, and iterated function
systems, 
\emph{Rozprawy Habilitacyjne Uniwersytetu Jagiello\'{n}skiego Nr
362, Wydawnictwo Uniwersytetu Jagiello\'{n}skiego}, Krak\'{o}w (2003).

\bibitem{To} B. Thomson, Real Functions, Lecture Notes in Mathematics \textbf{1170},
\emph{Springer} (1985).

\bibitem{Wer1} I. Werner, Contractive Markov systems,
\href{https://doi.org/10.1112/S0024610704006088}{
\emph{J. London Math. Soc.}
\textbf{71} (2005) 236--258}.

\bibitem{Wer3} I. Werner, Coding map for a contractive Markov system,
\href{http://dx.doi.org/10.1017/S0305004105009072}{
\emph{Math. Proc. Camb. Phil. Soc.}
\textbf{140} (2) (2006) 333--347}, 
\href{http://arxiv.org/abs/math/0504247}{ arXiv:math/0504247}.

\bibitem{Wer13} I. Werner, Erratum: Coding map for a contractive Markov system,
\href{http://arxiv.org/abs/1410.7545}{ arXiv:1410.7545}.

\bibitem{Wer6} I. Werner, The generalized Markov measure as an equilibrium state,
\emph{Nonlinearity} 
\textbf{18} (2005) 2261--2274,
\href{http://arxiv.org/abs/math/0503644}{ 	arXiv:math/0503644}.

\bibitem{Wer10} I. Werner, Dynamically defined measures and equilibrium states,
\href {http://dx.doi.org/10.1063/1.3666020}{
\emph{J. Math. Phys.} 
\textbf{52}  (2011) 122701},
\href{http://arxiv.org/abs/1101.2623}{ arXiv:1101.2623}.

\bibitem{Wer12} I. Werner, Erratum: Dynamically defined measures and equilibrium states,
\href{http://dx.doi.org/10.1063/1.4736999} {
\emph{J. Math. Phys.} 
\textbf{53}  079902 (2012)},
\href{http://arxiv.org/abs/1101.2623}{ arXiv:1101.2623}.

\bibitem{Wer11}  I. Werner, Equilibrium states and invariant measures  for  random dynamical systems,\href{http://dx.doi.org/10.3934/dcds.2015.35.1285}{
\emph{DCDS-A} 
\textbf{35} (3) (2015) 1285--1326},
\href{http://arxiv.org/abs/1203.6432}{arXiv:1203.6432}.

\bibitem{Wer15} I. Werner, On the Carath\'{e}odory approach to the construction of a measure, \href{http://dx.doi.org/10.14321/realanalexch.42.2.0345}
{\emph{Real Analysis Exchange} 
\textbf{42} (2) (2017) 345--384}, 
\href{http://arxiv.org/abs/1506.04736}{arXiv:1506.04736}.
\end{thebibliography}
\end{document}